\documentclass[a4paper,10pt]{article}
\usepackage[utf8]{inputenc}
\title{Les espaces de Berkovich sont excellents}
\author{\sc Antoine Ducros\\ \small Laboratoire J.-A. Dieudonné \\ \small Université de Nice - Sophia Antipolis, Parc Valrose\\ \small 06108 Nice CEDEX 02~~FRANCE}
\date{}
\usepackage{mathrsfs}

\usepackage[francais]{babel}
\usepackage{amssymb}
\usepackage{hyperref} 
\date{}
\input xy
\xyoption{all} 
\input{xypic} 
\newcommand{\irr}[1]{{\rm Irr}(#1)}
\renewcommand{\leq}{\leqslant}
\newcommand{\sep}{_{\rm sep}}
\renewcommand{\geq}{\geqslant}
\newcommand{\grot}{_{\tiny \mbox{\rm G}}}
\renewcommand{\dot}[1]{#1'}
\renewcommand{\ddot}[1]{#1''}

\newcommand{\hotimes}{\widehat{\otimes}}
\renewcommand{\phi}{\varphi}
\renewcommand{\dim}[1]{\mbox{\rm dim}_{#1}\;}
\newcommand{\codim}[1]{\mbox{\rm codim}_{#1}\;}
\newcommand{\got}[1]{{\mathfrak #1}}
\renewcommand{\Bbb}{\mathbb}    
\renewcommand{\cal}{\mathscr}

\renewcommand{\H}{\mbox{\rm H}}
\newcommand{\RR}{{\Bbb R}}
\newcommand{\KK}{{\Bbb K}}
\newcommand{\ZZ}{{\Bbb Z}}

\newcommand{\FF}{{\Bbb F}}
\newcommand{\DD}{{\Bbb D}}

\newcommand{\Aff}{{\Bbb A}}

\newcommand{\NN}{{\Bbb N}}

\newcommand{\QQ}{{\Bbb Q}}

\renewcommand{\epsilon}{\varepsilon}

\newcommand{\an}{^{\rm an}}
\newcommand{\zero}{^{\mbox{\tiny o}}}
\newcommand{\zeroo}{^{\mbox{\tiny oo}}}

\newcommand{\spec}{\mbox{\rm Spec}\;}

\newcommand{\zar}{_{\tiny\mbox{\rm Zar}}}

\newcommand{\itb}{\medskip \item[$\bullet$]}

\newcommand{\ka}{\widehat{k^a}}
\newcommand{\duit}{_{\rm red}}
\newcommand{\red}{\widetilde}
\newcommand{\hres}{{\cal H}}

\newcommand{\deux}[1]{\refstepcounter{cpt}\label{#1}\noindent {\bf (\thesection.\thecpt)}\hspace{.1cm}}
\newcommand{\trois}[1]{\refstepcounter{cptbis}\label{#1}\noindent {\bf
    (\thesection.\thecpt .\thecptbis)}\hspace{.1cm}}
\newcounter{cpt}
\newcounter{cptbis}

\begin{document}

\maketitle

\noindent
{\bf Abstract.} {\small In this paper, we first study the properties of the local rings of a Berkovich
analytic space from the point of view of the commutative algebra properties ; we
show that those rings are excellent ; we introduce the notion of an
{\em analytically separable} extension of non-archimedean complete fields (it
includes the case of the finite separable extensions, and also that of any
complete extension of a perfect complete non-archimedean field) and show that
the usual commutative algebra properties ($R_m$, $S_m$, Gorenstein, Cohen-Macaulay,
Complete Intersection) are stable under analytically separable ground field
extensions ; we also establish a GAGA principle with respect to those
properties for any finitely generated scheme over an affinoid algebra. The remaining part of the paper deals with more global notions which are closely related to the preceeding ones: the irreducible components of an analytic space, its normalization, and the behaviour of irreducibility and connectedness under base change.}

\tableofcontents

\section*{Introduction}
\addcontentsline{toc}{section}{Introduction}
La première partie de cet article est essentiellement consacrée à l'étude, du point de vue des propriétés usuelles de l'algèbre commutative (régularité, profondeur, excellence...), des anneaux locaux des espaces analytiques {\em au sens de Berkovich} ; dans un second temps, nous nous intéressons à des notions  plus globales comme les composantes irréductibles d'un espace analytique, sa normalisation, ou encore l'irréductibilité et la connexité géométriques.  Décrivons brièvement les principaux résultats obtenus.

\medskip
On fixe un corps ultramétrique complet $k$ ; on note $p$ son exposant caractéristique ; si $p=1$ et si $n\in \NN$, la notation $k^{1/p^n}$ désignera simplement $k$. Rappelons qu'un espace analytique est dit {\em bon} si chacun de ses points a une base de voisinages affinoïdes ; les espaces affinoïdes, les analytifiés de variétés algébriques, les fibres génériques de schémas formels {\em affines ou propres} sont des exemples de bons espaces -- en ce qui concerne les schémas formels propres, ce n'est pas trivial et découle d'un résultat de Temkin (\cite{tmk1}, cor. 4.4) ; la bonté d'un espace garantit que ses anneaux locaux sont noethériens et henséliens (\cite{brk2}, th. 2.1.4 et th. 2.1.5).

\subsection*{Propriétés algébriques des anneaux analytiques}

\begin{itemize} 

\medskip
\item[$\diamond$] {\bf Excellence des anneaux analytiques locaux et globaux.} 

\medskip
\begin{itemize}
\item[]Au paragraphe \ref{EXC}, nous démontrons (th.~\ref{EXC}.\ref{excellence}) que toute algèbre $k$-affinoïde est un anneau excellent, et que les anneaux locaux d'un espace $k$-affinoïde sont excellents.
\end{itemize} 

\medskip
\item[$\diamond$] {\bf Théorèmes de type GAGA.} 

\begin{itemize}
\medskip
\item[] Après avoir établi la régularité géométrique des fibres de certains morphismes entre anneaux locaux algébriques et analytiques (th. \ref{GAGA}.\ref{geomregulschem}) nous en déduisons l'assertion suivante  (th.~\ref{GAGA}.\ref{gagarmsm}) : soit ${\cal A}$ une algèbre $k$-affinoïde, soit $\cal X$ un schéma de type fini sur ${\cal A}$, et soit $V$ un  domaine affinoïde de l'analytifié ${\cal X}\an$ de $\cal X$. Soit $x$ un point de $V$ et soit $\bf x$ son image sur $\cal X$. Soit $\mathsf P$ l'une des propriétés suivantes : être $R_{m}$ (resp. $S_{m}$, resp. régulier, resp. d'intersection complète, resp. de Gorenstein, resp. de Cohen-Macaulay). On a les équivalences suivantes :

$${\cal O}_{V,x}\;\mbox{satisfait} \;\mathsf P\iff {\cal O}_{{\cal X}\an,x}\;\mbox{satisfait} \;\mathsf P\iff {\cal O}_{{\cal X},\bf x}\;\mbox{satisfait} \;\mathsf P\; .$$

\end{itemize} 

\medskip
\begin{itemize}
\item[] Remarquons qu'en vertu des résultats d'excellence mentionnés plus haut, l'ensemble des ${\bf x}\in {\cal X}$ tel que ${\cal O}_{{\cal X},\bf x}$ satisfasse $\mathsf P$ en est un ouvert de Zariski ; l'équivalence ci-dessus assure que l'ensemble des $x\in {\cal X}\an $ tel que ${\cal O}_{{\cal X}\an,x}$ satisfasse $\mathsf P$ en est un ouvert de Zariski ; cette remarque s'applique notamment au cas où ${\cal X}=\spec \cal A$ et où ${\cal X}\an$ est donc égal à ${\cal M}({\cal A})$. 

\medskip
\item[] {\bf À propos du lieu de validité d'une propriété $\mathsf P$ sur un espace analytique. } Soit  $X$ un espace analytique et soit $x\in X$. Soit $\mathsf P$ l'une des propriétés mentionnées ci-dessus. Si $X$ est bon, on dit que $X$ satisfait $\mathsf P$ en $x$ si ${\cal O}_{X,x}$ satisfait $\mathsf P$. Dans le cas général, on dit que  $X$ satisfait $\mathsf P$ en $x$ s'il existe un bon domaine analytique $U$ de $X$ contenant $x$ tel que $U$ satisfasse $\mathsf P$ en $x$, et c'est alors le cas pour {\em tout} bon domaine analytique $U$ de $X$ contenant $x$ (\ref{GAGA}.\ref{defpasbon}). On déduit de ce qui précède que le lieu de validité de $\mathsf P$ sur un espace analytique $X$ est un ouvert de Zariski de $X$. 
\end{itemize}

\medskip
\item [$\diamond$] {\bf Effets de l'extension des scalaires.} 

\medskip
\begin{itemize}
\item[] Soit $X$ un espace $k$-analytique et soit $L$ une extension ultramétrique complète de $k$. Soit $x\in X$, et soit $y$ un point de $X_{L}$ situé au-dessus de $x$, où $X_{L}$ est l'espace déduit de $X$ par extension des scalaires de $k$ à $L$. Soit $\mathsf P$ l'une des propriétés d'algèbre commutative évoquées plus haut. Nous établissons  les faits suivants (le théorème~\ref{GAGA}.\ref{changermsm} constitue leur variante schématique ; ils en découlent modulo le théorème de comparaison ~\ref{GAGA}.\ref{gagarmsm} et les commentaires faits au \ref{GAGA}.\ref{defpasbon}.\ref{encorechangebase}) : 

\medskip
\begin{itemize}
\item[$i)$] si $X_L$ satisfait $\mathsf P$ en $y$ alors $X$ satisfait $\mathsf P$ en $x$ ; 
\item[$ii)$] si $X$ satisfait $\mathsf P$ en $x$ et si $\mathsf P$ est la propriété d'être $S_{m}$ pour un certain $m$ (resp. d'intersection complète, resp. de Gorenstein, resp. de Cohen-Macaulay), alors  $X_L$ satisfait $\mathsf P$ en $y$ ; 
\item[$iii)$] si $X$ satisfait $\mathsf P$ en $x$, si $\mathsf P$ est la propriété d'être $R_{m}$ pour un certain $m$ (resp. régulier) et si $L$ est une extension {\em analytiquement séparable de $k$}, alors  $X_L$ satisfait $\mathsf P$ en $y$.
\end{itemize}

\medskip
\item[]{\bf Commentaires.} La notion d'extension {\em analytiquement séparable} de corps ultramétriques complets est introduite au paragraphe 1 (déf.~\ref{ANSEP}.\ref{fortsep}) ; c'est, aussi bien du point de vue de sa définition que de ses propriétés (telles l'assertion $iii)$ ci-dessus) la variante valuée de la notion classique d'extension séparable.  

\medskip
\item[]{\bf Quelques exemples ({\em cf.} \ref{ANSEP}.\ref{exfortsep}.\ref{parffortsep}, \ref{ANSEP}.\ref{exfortsep}.\ref{sepfortsep}, \ref{ANSEP}.\ref{introunivmult}.\ref{premunivmult} et \ref{ANSEP}.\ref{introunivmult}.\ref{univmult})}. Si $k$ est parfait, toute extension complète de $k$ est analytiquement séparable ; une extension finie de $k$ est analytiquement séparable si et seulement si elle est séparable ; si $\bf r$ est un polyrayon $k$-libre, $k_{\bf r}$ est une extension analytiquement séparable de $k$ (pour la définition d'un polyrayon $k$-libre et le sens de la notation $k_{\bf r}$, voir \ref{RAP}.\ref{polyray}) ; plus généralement, si $\eta$ est l'unique point du bord de Shilov d'un polydisque sur $k$, son corps résiduel complété $\hres(\eta)$ est une extension analytiquement séparable de $k$.

\medskip
\item[]{\bf Régularité géométrique.} Soit $X$ un espace $k$-analytique. Soit $\mathsf P$ l'une des propriétés d'algèbre commutative mentionnées ci-dessus. Si $x\in X$, on dit que $X$ satisfait {\em géométriquement} $\mathsf P$ en $x$ si pour toute extension complète $L$ de $k$, l'espace $X_L$ satisfait $\mathsf P$ en chacun de ses points situés au-dessus de $x$. Au paragraphe \ref{UNIREG} nous démontrons (prop.~\ref{UNIREG}.\ref{proprgeom}) que $X$ est géométriquement régulier en $x$ si et seulement si il est {\em quasi-lisse} en $x$, c'est-à-dire si et seulement si le rang de $\Omega^1_{X/k}$ en $x$ est égal à la dimension de $X$ en $x$ ; nous prouvons ensuite (prop. \ref{UNIREG}.\ref{kunsurp}) que si $P$ est la propriété d'être régulier ou bien d'être $R_m$ pour un certain $m$, alors $X$ satisfait géométriquement $\mathsf P$ en $x$ si et seulement si $X_{k^{1/p}}$ satisfait $\mathsf P$ en son unique point situé au-dessus de $x$. On en déduit (cor. \ref{UNIREG}.\ref{lieugeomreg}) que le lieu de validité géométrique de $\mathsf P$ sur $X$ en est un ouvert de Zariski. 

\medskip
\item[]{\bf Remarque.} En géométrie analytique rigide, les assertions ci-dessus n'auraient pas pu être même {\em énoncées}, puisqu'il n'existe pas dans ce cadre de morphisme de changement de base.

\end{itemize}

\medskip
\item[$\diamond$] {\bf À propos des démonstrations.} 

\medskip
\begin{itemize}
\item[]Les résultats qui précèdent reposent pour l'essentiel sur la proposition~\ref{EXC}.\ref{ouvertsreg} dont l'énoncé, et plus encore la preuve, sont passablement ingrats. Elle assure l'existence, sur le spectre du complété d'un certain anneau local analytique, d'un ouvert de Zariski d'un type particulier qui est non vide et {\em régulier.} Pour l'exhiber, on se fonde sur le lemme de normalisation de Noether (et plus précisément sur ses versions analytique et algébrique, sur la régularité des anneaux locaux de l'espace affine analytique, qui est établie directement au préalable (lemme~\ref{EXC}.\ref{disquereg}), et enfin sur le critère de régularité suivant, dû à Kiehl (\cite{kie}, Folg. 2.3) : {\em si $A$ est un anneau noethérien, si $B$ est une $A$-algèbre de type fini régulière, si $C$ est une $B$-algèbre finie et plate et si $\Omega^1_{B/A}$ et $\Omega^1_{C/A}$ sont tous deux libres de même rang, alors l'anneau $C$ est régulier.} 

\medskip
\item[] {\bf Remarque.} D'après un théorème de Kiehl (\cite{kie}, Satz 1.4), si $F$ est un corps ultramétrique complet de caractéristique $p>0$, et si $E$ est un sous-corps complet de $F$ contenant $F^{p}$ et {\em topologiquement de type dénombrable sur $F^{p}$}, alors $E$ possède une $p$-base topologique sur $F^{p}$. Nous nous servons de ce résultat au cours de la démonstration de la proposition~\ref{EXC}.\ref{ouvertsreg} ; afin de pouvoir l'appliquer, nous sommes amené à plusieurs reprises à remplacer le corps avec lequel nous travaillons par un sous-corps satisfaisant l'hypothèse de dénombrabilité évoquée, puis à utiliser des arguments de limite inductive pour conclure ; ces contorsions techniques que nous n'avons malheureusement pas su éviter alourdissent la rédaction.

\end{itemize}
\end{itemize}
\subsection*{Normalisation et propriétés géométriques globales des espaces analytiques}

\begin{itemize} 

\item[$\diamond$]  {\bf Les composantes irréductibles d'un espace analytique.} 

\medskip
\begin{itemize}
\item[]Si $X$ n'est pas compact, sa topologie de Zariski n'est pas noethérienne, et la théorie classique des composantes irréductibles ne s'applique pas. Nous pallions comme suit cet inconvénient : par des arguments relativement élémentaires reposant {\em in fine} sur la théorie de la dimension nous montrons (th. \ref{COMP}.\ref{existcompirr}) que les parties de $X$ qui peuvent s'écrire comme l'adhérence, relative à la topologie de Zariski de $X$, d'une composante irréductible d'un domaine affinoïde de $X$, sont exactement les fermés irréductibles maximaux de $X$ ; ce sont ces parties que nous appellerons ses composantes irréductibles. Nous établissons les propriétés attendues à leur sujet (du lemme \ref{COMP}.\ref{dimloccasgen} à la fin du paragraphe  \ref{COMP}) ; elles se comportent essentiellement comme les composantes irréductibles usuelles d'un espace topologique noethérien, à condition de remplacer à peu près partout la notion d'ensemble fini de fermés de Zariski par celle d'ensemble {\em G-localement fini} de tels fermés (déf. \ref{COMP}.\ref{defglocfin}).

\end{itemize}

\medskip
\item[$\diamond$] {\bf La normalisation d'un espace analytique.} 

\medskip
\begin{itemize}
\item[]Au début du paragraphe \ref{NOR}, nous introduisons la notion de morphisme {\em quasi-dominant} d'espaces analytiques (définition \ref{NOR}.\ref{pseudom} ; elle repose sur celle de composante irréductible) et nous définissons la normalisation de $X$ comme l'objet final de la catégorie des espaces analytiques normaux munis d'un morphisme quasi-dominant vers $X$.  Nous en montrons l'existence (th. \ref{NOR}.\ref{existnormx}) ; l'approche en termes de propriété universelle a l'avantage de réduire la démonstration au cas affinoïde, en assurant à peu près automatiquement le recollement des constructions locales. Nous montrons que les composantes connexes de la normalisation de $X$ correspondent aux composantes irréductibles de $X$ (th. \ref{NOR}.\ref{desccompnorm}), et que la normalisation est «compatible à la formation des anneaux locaux» dans le cas d'un bon espace (lemme \ref{NOR}.\ref{normannloc}). 

\medskip
Nous montrons l'existence, lorsque $X$ est quasi-compact, d'un entier $n$ tel que la normalisation de $X_{k^{1/p^n}}$ soit géométriquement normale (th. \ref{UNIREG}.\ref{defnormfin}) ; un résultat analogue est établi au préalable concernant l'espace réduit associé à un espace analytique (th. \ref{UNIREG}.\ref{defredfin}). 

\medskip
Nous utilisons enfin le théorème \ref{UNIREG}.\ref{defnormfin} pour prouver, à l'aide du théorème \ref{NOR}.\ref{desccompnorm}, qu'un espace $k$-analytique qui est irréductible le reste après passage à la clôture parfaite de $k$ (prop. \ref{UNIREG}.\ref{perfclos}).

\end{itemize} 
\medskip
\item[$\diamond$]  {\bf Connexité et irréductibilité géométriques.} 

\medskip
\begin{itemize}
\item[] Au paragraphe \ref{GEO}, nous commençons par définir, étant donné un espace $k$-analytique $X$, l'anneau $\got s(X)$ des fonctions analytiques sur $X$ qui annulent localement un polynôme non nul et séparable sur $k$ (\ref{GEO}.\ref{introcx}) ; nous donnons une démonstration du fait (essentiellement dû à Berkovich) que $\got s(X)$ est une extension finie séparable de $k$ lorsque $X$ est connexe et non vide (lemme \ref{GEO}.\ref{cxcorps}). 

\medskip
Le foncteur $X\mapsto \got s(X)$ et le lemme \ref{GEO}.\ref{cxcorps} sont ensuite utilisés pour étudier le comportement de la connexité et de l'irréductibilité par extension du corps de base. Les résultats attendus, à savoir les pendants analytiques des assertions bien connues portant sur les variétés {\em algébriques}, sont établis ; on démontre d'abord ceux qui sont relatifs à la connexité (th. \ref{GEO}.\ref{conngeomfond}) ; suivent leurs analogues concernant l'irréductibilité (th. \ref{GEO}.\ref{irrgeomfond}), qui s'en déduisent à l'aide d'un raisonnement faisant intervenir la notion de normalisation et le théorème  \ref{UNIREG}.\ref{defnormfin} cité ci-dessus. 

\medskip
Nous n'allons pas rappeler ici les énoncé intégraux des théorèmes \ref{GEO}.\ref{conngeomfond} et \ref{GEO}.\ref{irrgeomfond} ; mentionnons simplement qu'ils assurent entre autres que si $k$ est algébriquement clos\footnote{Les assertions correspondantes en géométrie algébrique sont en général énoncées sur un corps {\em séparablement } clos ; mais il se trouve qu'un corps ultramétrique complet dont la valeur absolue n'est pas triviale et qui est séparablement clos est automatiquement algébriquement clos {\em cf.} \cite{bgr}, \S 3.4, prop. 6.} et si $X$ est un espace $k$-analytique connexe (resp. irréductible) alors $X_L$ est connexe (resp. irréductible) pour toute extension complète $L$ de $k$. 

\medskip
Nous terminons ce chapitre en établissant, par un raisonnement fondé sur les théorèmes \ref{GEO}.\ref{conngeomfond} et \ref{GEO}.\ref{irrgeomfond} ainsi que sur une proposition technique (prop. \ref{GEO}.\ref{presquenorm} ; {\em cf.} aussi les commentaires qui la précèdent au \ref{GEO}.\ref{explicanorm}), les faits suivants (th. \ref{GEO}.\ref{corpsdefconn}, $iii)$,  $iv)$ et $vi)$ ; th. \ref{GEO}.\ref{corpsdefirr}, $iii)$, $iv)$ et $vi)$ ): si $X$ est un espace $k$-analytique, si $L$ est une extension complète de $k$, et si $T$ est une composante connexe (resp. irréductible) de $X_L$, alors $T$ possède un \og corps de définition\fg~qui est fini et séparable sur $k$, et l'image de $T$ sur $X$ est une composante connexe (resp. irréductible) de $X$ ; si $X$ n'a qu'un nombre fini de composantes connexes (resp. irréductibles), il existe une extension finie et séparable $K$ de $k$ telle que toutes les composantes connexes (resp. irréductibles) de $X_K$ soient géométriquement connexes (resp. géométriquement irréductibles). 

\medskip
\item[$\diamond$] {\bf Stabilité de certaines propriétés par produit.} 

\medskip
\begin{itemize}
\item[] Au paragraphe \ref{PROD}, nous montrons (th. \ref{PROD}.\ref{prodnorm}) que si $\mathsf P$ est l'une des propriétés d'algèbre commutative mentionnée plus haut et si $X$ et $Y$ sont deux espaces $k$-analytiques la satisfaisant {\em géométriquement}, alors le produit $X\times_k Y$ la satisfait géométriquement\footnote{Nous espérons pouvoir ultérieurement, grâce à un travail en cours sur la platitude, arriver à cette conclusion en nous contentant de la validité {\em simple} de $\mathsf P$ sur l'un des deux espaces en jeu.}.

\medskip
Nous prouvons ensuite (th. \ref{PROD}.\ref{prodconnirr}) que si $X$ est un espace $k$-analytique géométriquement connexe (resp. géométriquement irréductible) et $Y$ un espace $k$-analytique connexe (resp. irréductible) alors $X\times_k Y$ est connexe (resp. irréductible) ; l'assertion relative à la connexité se déduit facilement du théorème \ref{GEO}.\ref{conngeomfond} sur la connexité géométrique, par un argument topologique élémentaire reposant sur la compacité des espaces affinoïdes ; modulo le fait que le produit de deux espaces géométriquement normaux est géométriquement normal (en vertu du théorème  \ref{PROD}.\ref{prodnorm} mentionné ci-dessus), celle relative à s'irréductibilité en découle immédiatement. 
\end{itemize} 
\end{itemize} 
\end{itemize} 

\subsection*{Liens avec divers travaux antérieurs} 

\noindent
{\bf Du cas strictement affinoïde au cas général.} Parmi les différences notables entre la théorie de Berkovich et la géométrie analytique rigide figure la prise en considération, par la première, des polydisques fermés de polyrayon {\em quelconque}, quand la seconde n'autorisait que les polydisques {\em unité} ; la classe des algèbres affinoïdes au sens de Berkovich est ainsi plus large que celle définie par Tate, et les algèbres affinoïdes de ce dernier sont qualifiées de {\em strictement} affinoïdes par Berkovich. 

\medskip
L'intérêt d'accepter n'importe quel polyrayon est particulièrement flagrant lorsqu'on s'intéresse aux espaces analytiques sur un corps muni {\em de la valeur absolue triviale}, espaces dont l'importance a été très récemment mise en lumière dans différents travaux : ceux d'Amaury Thuillier sur la combinatoire du diviseur exceptionnel d'une résolution des singularités (\cite{thui}), ceux de Berkovich et de Nicaise sur les liens entre la géométrie ultramétrique et les variations de structures de Hodge {\em complexes} (\cite{brklimh}, \cite{nic}), ou encore ceux de Jérôme Poineau sur les fondements de la géométrie analytique {\em sur $\ZZ$} et ses applications arithmétiques (\cite{thesejp}). 

\medskip
Il apparaît donc nécessaire de disposer, autant que faire se peut, de résultats algébriques de base au sujet des espaces et algèbres affinoïdes qui soient valables dans le cas général, {\em i.e.} non nécessairement strict. C'est dans cette optique que nous étendons ici aux algèbres affinoïdes quelconques plusieurs assertions déjà connues pour les algèbres strictement affinoïdes ; {\em mentionnons que nous avons choisi, pour la commodité du lecteur, de ne pas utiliser ces assertions déjà connues}, préférant reprendre, lorsque cela s'est avéré nécessaire, certains passages de leurs preuves originales au cours de nos démonstrations. Donnons maintenant quelques détails.

\medskip
\begin{itemize}
\item[$\diamond$] {\bf Propriétés algébriques des anneaux analytiques : ce qui était déjà connu dans le cas strictement affinoïde.} 

\medskip
\begin{itemize}
\itb Le théorème ~\ref{EXC}.\ref{excellence} affirme l'excellence des algèbres affinoïdes ; dans le cas strictement affinoïde, celle-ci a été établie par Kiehl (\cite{kie}). 
\itb Le théorème de comparaison pour les propriétés usuelles de l'algèbre commutative (th.~\ref{GAGA}.\ref{gagarmsm}) énoncé ici met en jeu un schéma de type fini $\cal X$ sur une algèbre affinoïde $\cal A$ et un domaine affinoïde $V$ de ${\cal X}\an$ ; pour certaines de ces propriétés (régularité, normalité, caractère réduit, de Gorenstein, de Cohen-Macaulay, d'intersection complète) il a été démontré par Berkovich lorsque $\cal X={\cal A}$ et $V={\cal X}\an$ (\cite{brk2}, th. 2.2.1) ; {\em lorsque $|k^*|\neq \{1\}$ et lorsque $\cal A$ est strictement $k$-affinoïde}, sa méthode s'applique en fait, même s'il ne le mentionne pas explicitement, à toutes les propriétés d'algèbre commutative considérées ici, à tout $\cal A$-schéma de type fini $\cal X$, et à tout domaine {\em strictement $k$-affinoïde} de ${\cal X}\an$. 

\itb Le théorème \ref{UNIREG}.\ref{kunsurp} assure que si $\mathsf P$ est la propriété d'être régulier ou bien $R_m$ pour un certain $m$, alors la validité {\em géométrique} de $\mathsf P$ en un point d'un espace $k$-analytique $X
$ équivaut à sa validité {\em simple} en l'unique antécédent de ce point sur $X_{k^{1/p}}$ ; dans son article \cite{comprig}, Conrad en démontre un cas particulier (dont l'énoncé nous a inspiré le cas général) pour les espaces strictement $k$-affinoïdes ; plus précisément, le lemme 3.3.1 (resp. le th. 3.3.6) de {\em op. cit.} affirme entre autres que si $\cal A$ est une algèbre strictement $k$-affinoïde telle que $k^{1/p}\hotimes_{k}{\cal A}$ soit réduite (resp. normale), alors pour toute extension complète $F$ de $k$, l'algèbre $F\hotimes_{k}\cal A$ est réduite (resp. normale).
\end{itemize}

\medskip
Concernant les deux premiers résultats évoqués, nos démonstrations sont en grande partie inspirées par celles déjà existantes dans le contexte strictement affinoïde, auquel la partie proprement nouvelle de notre travail consiste justement à se ramener, d'une façon plus ou moins explicite. Comme d'habitude en théorie de Berkovich, on le fait en remplaçant le corps de base $k$ par $k_{\bf r}$, où $\bf r$ est un polyrayon $k$-libre convenable ; il reste à contrôler les effets d'une telle extension des scalaires, ce que permet la proposition~\ref{EXC}.\ref{ouvertsreg}, l'extension $k_{\bf r}$ de $k$ étant analytiquement séparable (exemple \ref{ANSEP}.\ref{introunivmult}.\ref{premunivmult}). 

\medskip
S'agissant du troisième résultat, nous étendons également les scalaires à $k_{\bf r }$ pour un certain polyrayon $k$-libre $\bf r$ afin de se ramener au cas strictement $k$-affinoïde, mais nous traitons ce dernier par une méthode qui diffère de celle utilisée par Conrad, et repose sur le critère de régularité de Kiehl mentionné plus haut. 

\medskip

{\bf Remarque.} En ce qui concerne la normalité, les théorèmes GAGA sur n'importe quelle base affinoïde et la stabilité par l'extension des scalaires à $k_{\bf r}$ ont été prouvés précédemment par l'auteur au moyen de méthodes {\em ad hoc} dans l'appendice de \cite{duc2}.

\medskip
\noindent
\item[$\diamond$] {\bf Propriétés algébriques des anneaux analytiques : ce qui est nouveau même dans le cas strictement affinoïde.} 

\medskip
Il s'agit pour l'essentiel  : 

\medskip
\begin{itemize}
\itb du théorème qui concerne le comportement des propriétés usuelles de l'algèbre commutative par extension des scalaires (th.~\ref{GAGA}.\ref{changermsm}) ; signalons à ce propos que la notion de séparabilité analytique, introduite pour l'occasion, n'avait à notre connaissance pas été considérée jusqu'ici ;
\itb du théorème qui affirme l'excellence des anneaux locaux {\em analytiques} (th.~\ref{EXC}.\ref{excellence}) ;
\itb du théorème de régularité géométrique des fibres analytiques (th. \ref{GAGA}.\ref{geomregulschem}).
\end{itemize}

\medskip
Des résultats partiels à propos des deux premiers points ont été  toutefois établis par Conrad dans {\em op. cit.}

\medskip
Ainsi on déduit de son énoncé rappelé plus haut que si $k$ est parfait (auquel cas $k^{1/p}=k$), toute algèbre strictement $k$-affinoïde réduite (resp. normale) le reste après une extension complète quelconque de $k$ ; ce dernier fait peut se retrouver en corollaire de notre théorème~\ref{GAGA}.\ref{changermsm}, compte-tenu du fait que si $k$ est parfait, toutes ses extensions complètes sont analytiquement séparables (exemple~\ref{ANSEP}.\ref{exfortsep}.\ref{parffortsep}).

\medskip
Conrad établit par ailleurs ({\em op. cit.}, th. 1.1.3) que si $\cal A$ est une algèbre strictement affinoïde sur un corps ultramétrique complet $k$, les anneaux locaux de l'espace ${\cal M}({\cal A})$ en chacun de ses points {\em rigides} sont excellents ; sa démonstration ne s'étend pas au cas de tous les points de ${\cal M}({\cal A})$ : elle utilise en effet de manière essentielle le fait que si $x$ est un point rigide de ${\cal M}({\cal A})$ et si $\got{m}$ désigne son image sur $\spec \cal A$, alors $\widehat{{\cal A}_{\got{m}}}$ et $\widehat{{\cal O}_{{\cal M}({\cal A}),x}}$ sont naturellement isomorphes ; si $x$ n'est pas supposé rigide, cette dernière affirmation est fausse en général. 

\medskip
Mentionnons incidemment que l'excellence des anneaux locaux d'un espace analytique {\em complexe} peut se démontrer au moyen d'un «critère jacobien» ({\em cf.} \cite{mats2}, th. 100, th. 101, et la remarque qui les suit). 

\medskip

\item[$\diamond$] {\bf Normalisation et des propriétés géométriques globales des espaces analytiques : comparaison avec les méthodes et les résultats de Conrad.}

\medskip
\begin{itemize}
\item[] Le contenu de nos paragraphes \ref{COMP}, \ref{NOR}, \ref{UNIREG} et \ref{GEO} constitue pour l'essentiel une extension au cas non nécessairement strict des travaux de Conrad sur les composantes irréductibles et la normalisation d'un espace rigide, qui sont au c\oe ur de l'article \cite{comprig} déjà cité ; toutefois, notre approche diffère de la sienne (y compris dans le cas strictement analytique) sur un point important. 

\begin{itemize}

\itb {\em Différence entre notre approche et celle de Conrad.} Partant d'un espace rigide $X$, ce dernier commence par construire sa normalisation $X'$ (par recollement à partir du cas affinoïde), puis {\em définit} les composantes irréductibles de $X$ comme étant les images des composantes connexes de $X'$. 

\medskip
Nous avons choisi, comme indiqué plus haut, de procéder d'une autre manière, en définissant d'abord les composantes irréductibles d'un espace analytique purement en termes de la topologie de Zariski, en se servant de ces dernières pour introduire la notion de morphisme {\em quasi-dominant}, puis en montrant l'existence d'un objet final dans la catégorie des espaces analytiques normaux munis d'un morphisme quasi-dominant vers un espace $X$ donné ; c'est cet objet final $X'$ que nous appelons normalisation, et nous {\em montrons} qu'il y a correspondance entre les composantes connexes de $X'$ et les composantes irréductibles de $X$. 

\medskip
À notre connaissance cette définition des composantes irréductibles d'un espace analytique et cette caractérisation de sa normalisation par une propriété universelle sont nouvelles, même dans le contexte de la géométrie rigide. 

\itb {\em À propos de la connexité et de l'irréductibilité géométriques.} Pour démontrer nos théorèmes \ref{GEO}.\ref{conngeomfond}, \ref{GEO}.\ref{irrgeomfond}, \ref{GEO}.\ref{corpsdefconn} et \ref{GEO}.\ref{corpsdefirr} nous suivons une stratégie qui est essentiellement la même que celle de Conrad ; modulo les résultats déjà établis à propos de la normalisation d'un espace analytique et certains lemmes intermédiaires tirés de l'article \cite{brk4} de Berkovich (pour la commodité du lecteur, nous les avons redémontrés ici), tout est fondé sur l'assertion suivante : {\em si $k$ est algébriquement clos et si $X$ est un espace $k$-affinoïde connexe alors $X_L$ est connexe pour toute extension complète $L$ de $k$.} 

\medskip
Pour la démontrer, plusieurs méthodes étaient à notre disposition. La plus rapide aurait consisté à invoquer le théorème de Berkovich (\cite{brk2}, th. 7.6.1) sur l'invariance de la cohomologie étale par extension de corps algébriquement clos\footnote{Berkovich, pour établir son théorème, utilise la stabilité de la connexité par extension de corps algébriquement clos... mais uniquement à propos des {\em courbes algébriques} ; nous pourrions donc nous y référer sans risque de raisonnement circulaire.} (en l'appliquant au $\H^0$...). 

\medskip
Ce dernier repose toutefois sur une machinerie extrêmement sophistiquée, et nous avons préféré opter pour une approche plus élémentaire, en nous plaçant tout d'abord dans le cas où $|k^*|\neq\{1\}$ et où $X$ est strictement $k$-affinoïde. Sous ces hypothèses, on peut utiliser la théorie des algèbres strictement $k$-affinoïdes distinguées (et notamment \cite{bgr}, \S 6.4.3, th. 1), la théorie de la réduction et le fait qu'une variété {\em algébrique} connexe sur le corps résiduel $\red k$ (qui est algébriquement clos) reste connexe après une extension quelconque des scalaires, pour conclure ; c'est ce que fait Conrad ({\em op. cit.}, th. 3.2.1) ; nous proposons ici (th. \ref{GEO}.\ref{xkaconnexe}) une démonstration plus directe, qui exploite les propriétés spécifiques des espaces de Berkovich (locale compacité, présence de points à \og gros\fg~ corps résiduels...), qui ne nécessite, en fait d'algèbre commutative normée, que le {\em Nullstellensatz} analytique, et ne fait pas appel à la théorie algébrique de la connexité géométrique. 

\medskip
Le passage au cas non nécessairement strict est fondé sur le corollaire \ref{GEO}.\ref{cxextkr}, qui affirme (sans supposer que $k$ est algébriquement clos) que si $X$ est un espace $k$-analytique connexe et non vide et si $\bf r$ est un polyrayon $k$-libre, alors $\got s(X_{k_{\bf r}})$ s'identifie à $\got s(X)_{\bf r}$. 

\itb {\em Stabilisation de la normalisation à un niveau fini.} Le théorème \ref{UNIREG}.\ref{defnormfin}, qui affirme l'existence, lorsque $X$ est compact, d'un entier $n$ tel que la normalisation de $X_{k^{1/p^n}}$ soit géométriquement normale, et l'assertion analogue concernant l'espace réduit associé à un espace analytique (th. \ref{UNIREG}.\ref{defredfin}), ont été directement inspirés à l'auteur par des résultats de Conrad ({\em op. cit.},  lemme 3.3.1 et prop. 3.3.6), qu'ils étendent au cas non nécessairement strict. Notre preuve est très proche de la sienne ; elle consiste à utiliser la proposition   \ref{UNIREG}.\ref{kunsurp} (Conrad se fondant quant à lui sur le cas particulier de celle-ci qu'il a établi dans {\em op. cit.} et que nous avons évoqué plus haut) et des arguments de noethérianité, appliqués à un module de type fini sur une certaine algèbre affinoïde. 

\itb {\em Stabilité de certaines propriétés par produit.} À notre connaissance, les résultats du paragraphe \ref{PROD} ne figurent pas à ce jour dans la littérature, même dans le cas strictement analytique. 

\end{itemize}

\end{itemize}

\end{itemize} 

\medskip
\noindent
{\bf Remarque.} Cet article reprend une partie d'une prépublication de l'auteur (\cite{ducvert})\footnote{C'est au vu de celle-ci que Berkovich nous a signalé l'article \cite{comprig} de Conrad.} dans laquelle figurent les théorèmes relatifs au comportement des propriétés de l'algèbre commutative vis-à-vis de l'extension des scalaires (la notion de séparabilité analytique y est introduite), les théorèmes de type GAGA, ainsi que ceux qui portent sur la connexité et l'irréductibilité géométriques ; les démonstrations que l'on y trouve ont pour partie été retranscrites ici telles quelles et pour partie profondément remaniées, lorsqu'elles se sont avérées lacunaires, erronées, ou trop compliquées. Sont par contre inédits dans ce qui suit les résultats d'excellence, la définition des composantes irréductibles et la construction de la normalisation, ainsi que la preuve de la stabilité par produit de certaines propriétés.

\subsection*{Remerciements} 

Je voudrais faire part de toute ma gratitude au rapporteur de cet article pour l'attention et la précision extrêmes avec lesquelles il en a relu les différentes versions, ainsi que pour la pertinence de ses remarques et suggestions auxquelles le texte sous sa forme actuelle doit beaucoup. 

\medskip
Je lui sais particulièrement gré de l'intérêt qu'il a montré pour ce travail. Il m'a notamment encouragé à amplement étoffer une première mouture de celui-ci, uniquement consacrée aux propriétés {\em locales} des espaces analytiques, en lui adjoignant les chapitres dévolus à leur géométrie {\em globale} ; qu'il en soit ici vivement remercié. 

\setcounter{section}{-1} 

\section{Rappels et notations}\label{RAP} 

\setcounter{cpt}{0}

\medskip
\deux{siteann} Si $f: (\mathsf T,{\cal O}_{\mathsf T})\to (\mathsf S,{\cal O}_{\mathsf S})$ est un morphisme de sites annelés et si $\cal F$ est un ${\cal O}_{\mathsf S}$-module, on notera $f^*{\cal F}$ le ${\cal O}_{\mathsf T}$-module $f^{-1}{\cal F}\otimes_{f^{-1}{\cal O}_{\mathsf  S}}{\cal O}_{\mathsf T}$.

\medskip
\deux{espacesberk} Dans tout ce texte, nous nous intéresserons aux espaces analytiques  ultramétriques {\em au sens de Berkovich} (\cite{brk1},\cite{brk2}) ; {\em dans cette théorie, la valeur absolue du corps de base peut être triviale}. 

\medskip
\deux{pasdecorps} Deux sortes de questions seront traitées ici : 

\medskip
- celles pour lesquelles le corps de base joue un rôle ; on fixera alors un corps ultramétrique complet $k$, et l'on travaillera dans la catégorie des espaces $k$-analytiques ;

\medskip
- celles pour lesquelles il est secondaire ; on se placera alors dans la catégorie des espaces analytiques (sans référence à un corps), dont les objets sont les couples $(X,k)$ où $k$ est un corps ultramétrique complet et $X$ un espace $k$-analytique ; un morphisme de $(Y,L)$ vers $(X,k)$ est un couple $(\iota, \phi)$ où $\iota$ est une injection isométrique de $k$ dans $L$ et $\phi$ un morphisme d'espaces $L$-analytiques de $Y$ dans $X\times_{\iota}L$ ; bien entendu, on omettra le plus souvent dans ce cadre de mentionner explicitement les corps et les plongements isométriques en jeu. 

\medskip
\noindent
De même, on parlera parfois d'algèbre (ou d'espace) affinoïde, et parfois d'algèbre (ou d'espace) $k$-affinoïde. 

\medskip
\deux{kappax} Si $X$ est un espace analytique et $x$ un point de $X$, on désignera par $\hres(x)$ le {\em corps résiduel complété} de $x$ ; si $V$ est un domaine affinoïde de $X$, on notera ${\cal A}_V$ l'algèbre affinoïde correspondante. Si $\cal X$ est un schéma et $\bf x$ un point de $\cal X$, on désignera par $\kappa({\bf x})$ le corps résiduel de $\bf x$ ; lorsque $\cal X$ est intègre et lorsque $\bf x$ est son point générique, on écrira parfois $\kappa({\cal X})$ au lieu de $\kappa(\bf x)$. 

\medskip
\deux{introespaff} Soit $\cal A$ une algèbre affinoïde, soit $\cal X$ son spectre et soit $X$ l'espace analytique ${\cal M}({\cal A})$. On dispose d'une surjection continue $\rho : X\to \cal X$. La {\em topologie de Zariski} sur $X$ est par définition l'image réciproque par $\rho$ de la topologie de Zariski du schéma (noethérien) $\cal X$ ; autrement dit, un fermé de Zariski de $X$ est une partie que l'on peut décrire comme le lieu des zéros d'un certain idéal $I$ de $\cal A$, c'est-à-dire comme l'ensemble des $x\in X$ tels que $f(x)=0$ pour toute fonction $f$ appartenant à $ I$. 

\medskip
L'adjectif \og irréductible\fg~ appliqué à une partie de $X$ sera toujours relatif à la topologie de Zariski. Par surjectivité de $\rho$, les applications $F\mapsto \rho(F)$ et $G\mapsto \rho^{-1}(G)$ mettent en bijection l'ensemble des fermés de Zariski (resp. des fermés de Zariski irréductibles, resp. des composantes irréductibles) de $X$ et l'ensemble des fermés de Zariski (resp. des fermés de Zariski irréductibles, resp. des composantes irréductibles) de $\cal X$ ; en particulier, $X$ est irréductible si et seulement si $\cal X$ est irréductible. 

\medskip
Soit $Y$ un fermé de Zariski de $X$, et soit $I$ un idéal dont $Y$ est le lieu des zéros. La surjection ${\cal A}\to {\cal A}/I$ induit un homéomorphisme ${\cal M}({\cal A}/I)\simeq Y$, ce qui permet de munir $Y$ d'une structure d'espace affinoïde {\em qui dépend de l'idéal $I$ choisi}. L'ensemble des idéaux dont le lieu des zéros coïncide avec $Y$ a un plus grand élément, à savoir l'idéal formé de {\em toutes} des fonctions qui s'annulent en tout point de $Y$. Par surjectivité de $\rho$ et par le fait correspondant en théorie des schémas, cet idéal n'est autre que $\sqrt I$.

\medskip
\deux{faisceaucoh} Si $X$ est un espace analytique, il est muni d'une topologie de Grothendieck ensembliste, plus fine que sa topologie usuelle, qui est appelée la G-{\em topologie} (\cite{brk2}, \S 1.3). Le site correspondant est noté $X\grot$, il est naturellement annelé ; l'on désignera par ${\cal O}_{X\grot}$ son faisceau structural ; si $V$ est un domaine affinoïde de $X$, alors ${\cal O}_{X\grot}(V)={\cal A}_V$. Lorsque $X$ est {\em bon} (c'est-à-dire lorsque chacun de ses points a une base de voisinages affinoïdes), on désigne par ${\cal O}_{X}$ la restriction de ${\cal O}_{X\grot}$ à la catégorie des ouverts de $X$, et ${\cal O}_X$ fait alors de $X$ un espace localement annelé, dont les anneaux locaux sont noethériens et henséliens (\cite{brk2}, th. 2.1.4 et th. 2.1.5). Le lemme suivant est certainement bien connu mais, à la connaissance de l'auteur, ne figure pas explicitement dans la littérature. 

\medskip
\deux{oxoxgcoh} {\bf Lemme.} {\em Soit $X$ un espace analytique. Le faisceau ${\cal O}_{X\grot}$ est cohérent. Si $X$ est bon, le faisceau ${\cal O}_X$ est cohérent.} 

\medskip
{\em Démonstration.} On traite séparément chacune des deux assertions.

\setcounter{cptbis}{0}
\medskip
\trois{oxgcoh} {\em Preuve de la cohérence de ${\cal O}_{X\grot}$}. Soit $V$ un domaine analytique de $X$, soit $n\in \NN$ et soit $\phi : {\cal O}_{V\grot}^n\to {\cal O}_{V\grot}$ une surjection ; notons $\cal K$ son noyau. La section $1$ de ${\cal O}_{V\grot}(V)$ appartenant G-localement à l'image de $\phi$, il existe un G-recouvrement $(V_i)$ de $V$ par des domaines affinoïdes tels que ${\cal A}_{V_i}^n\to{\cal A}_{V_i}$ soit surjective pour tout $i$. Fixons $i$. L'anneau ${\cal A}_{V_i}$ étant noethérien, le noyau de la flèche ${\cal A}_{V_i}^n\to{\cal A}_{V_i}$, qui n'est autre que ${\cal K}(V_i)$, est un ${\cal A}_{V_i}$-module de type fini. Choisissons une famille génératrice $(e_1,\ldots,e_m)$ de ce module. Soit $W$ un domaine affinoïde de $V_i$. Par platitude de la ${\cal A}_{V_i}$-algèbre ${\cal A}_W$ (\cite{brk1}, prop. 2.2.4), la suite $$0\to{\cal K}(V_i)\otimes_{{\cal A}_{V_i}}{\cal A}_W\to{\cal A}_W^n\to{\cal A}_W\to 0$$ est exacte ; dès lors, ${\cal K}(W)$ s'identifie à ${\cal K}(V_i)\otimes_{{\cal A}_{V_i}}{\cal A}_W$. Par conséquent, le morphisme ${\cal O}_{V_{i,\grot}}^m\to{\cal K}_{|V_i}$ défini par les $e_j$ induit une surjection ${\cal A}_W^m\to{\cal K}(W)$ ; ceci valant pour tout domaine affinoïde $W$ de $V_i$, le faisceau ${\cal K}_{|V_i}$ est de type fini ; ceci valant pour chacun des $V_i$, le faisceau ${\cal O}_{X\grot}$ est cohérent.

\medskip
\trois{oxcoh} {\em Preuve de la cohérence de ${\cal O}_X$ lorsque $X$ est bon.} Supposons $X$ bon, soit $U$ un ouvert de $X$, soit $n\in \NN$, et soit $\phi : {\cal O}_U^n\to {\cal O}_U$ une surjection ; notons $\cal K$ son noyau. Il existe une famile finie $(g_1,\ldots,g_n)$ de fonctions analytiques sur $U$ telles que pour tout ouvert $\Omega$ de $U$ et tout $n$-uplet $(f_1,\ldots,f_n)$ de fonctions analytiques sur $\Omega$, l'on ait $\phi(f_1,\ldots,f_n)=\sum f_ig_{i|\Omega}$. Soit $x$ un point de $U$, soit $U'$ un voisinage ouvert de $x$ dans $U$ tel que la section $1$ de ${\cal O}_{U'}(U')$ appartienne à l'image de ${\cal O}_U'^n\to {\cal O}_U'$, et soit $V$ un voisinage affinoïde de $x$ dans $U'$ ; on note $\Omega$ l'intérieur topologique de $V$ dans $U$. Le morphisme ${\cal A}^n_{V}\to{\cal A}_{V}$ qui envoie $(f_1,\ldots,f_n)$ sur $\sum f_ig_{i|V}$ est surjectif en vertu de l'inclusion $V\subset U'$. Par noethérianité de ${\cal A}_{V}$, le noyau de ce morphisme est un ${\cal A}_{V}$-module de type fini $M$. Soit $(e_1,\ldots,e_m)$ une famille génératrice de ce dernier. Pour tout $j$, écrivons $e_j=(e_{1,j},\ldots,e_{n,j})$. Comme $\sum_i e_{i,j}g_{i|V}=0$ quel que soit $j$, la famille $(e_{1|\Omega},\ldots,e_{m|\Omega})$ appartient à ${\cal K}(\Omega)$. Elle définit donc un morphisme ${\cal O}_{\Omega}^{m}\to {\cal K}_{|\Omega}$. Soit $\omega\in \Omega$ et soit $W$ un voisinage affinoïde de $\omega$ dans $\Omega$. La platitude de la ${\cal A}_{V}$-algèbre ${\cal A}_{W}$  (\cite{brk1}, prop. 2.2.4) assure que $$0\to M\otimes_{{\cal A}_V}{\cal A}_W\to {\cal A}_{W}^n\to {\cal A}_W\to 0$$ est exacte. Ceci valant pour tout voisinage affinoïde $W$ de $\omega$, on obtient par passage à la limite inductive l'exactitude de $$0\to M\otimes_{{\cal A}_V}{\cal O}_{\Omega,\omega}\to {\cal O}_{\Omega,\omega}^n\to {\cal O}_{\Omega,\omega}\to 0.$$ Autrement dit, le noyau de ${\cal O}_{\Omega,\omega}^n\to {\cal O}_{\Omega,\omega}$ est engendré par les germes des $e_j$ ; mais cela signifie exactement que ${\cal O}_{\Omega,\omega}^m\to {\cal K}_{\omega}$ est surjective. Ceci ayant été établi pour tout $\omega \in \Omega$, le morphisme ${\cal O}_{\Omega}^{m}\to {\cal K}_{|\Omega}$ est surjectif, ce qui achève la démonstration.~$\Box$ 

\medskip
\deux{remdefox} {\bf Remarque.} La définition de ${\cal O}_{X}$ donnée ci-dessus garde un sens pour $X$ quelconque ({\em i.e.} non nécessairement bon), et $(X,{\cal O}_{X})$ est toujours un espace localement annelé ; toutefois, on ne sait à peu près rien de ${\cal O}_X$ dans ce cas. À la connaissance de l'auteur, on ignore s'il est cohérent, si ses anneaux locaux sont noethériens, ou s'il peut arriver (en dimension strictement positive) qu'ils soient réduits au corps de base... 

\medskip
\deux{faisccoh} Soit $X$ un espace analytique (resp. un bon espace analytique). Le lemme \ref{RAP}.\ref{oxoxgcoh} garantit qu'un ${\cal O}_{X\grot}$-module (resp. un ${\cal O}_X$-module) $\cal F$ est cohérent si et seulement si il existe un G-recouvrement (resp. un recouvrement ouvert) $(X_i)$ de $X$ tel que pour tout $i$, l'on puisse écrire ${\cal F}_{|X_i}$ comme le conoyau d'une flèche de la forme ${\cal O}_{{X_{i}}\grot}^m\to {\cal O}_{{X_{i}}\grot}^n$ (resp. ${\cal O}_{X_i}^m\to {\cal O}_{X_i}^n$), où $m$ et $n$ appartiennent à $\NN$. 

\medskip
\setcounter{cptbis}{0}
\trois{theokiehl} Soit $\cal A$ une algèbre affinoïde ; posons $X={\cal M}({\cal A})$. D'après un théorème de Kiehl ({\em cf.} \cite{brk2}, \S 1.2), le foncteur ${\cal F}\mapsto {\cal F}(X)$ établit une équivalence entre la catégorie des ${\cal O}_{X\grot}$-modules cohérents et celle des ${\cal A}$-modules de type fini. On dispose d'un quasi-inverse explicite $M\mapsto \red M$, où $\red M$ est uniquement déterminé par le fait que si $V$ est un domaine affinoïde de $X$, alors $\red M(V)=M\otimes_{\cal A}{\cal A}_V$. 

\medskip
\trois{oxequivoxg} Soit $X$ un bon espace analytique, soit $\pi : X\grot\to X$ le morphisme naturel de sites annelés. Berkovich a démontré les résultats suivants (\cite{brk2}, prop. 1.3.4 et 1.3.6) : les foncteurs $\pi_*$ et $\pi^* $ établissent une équivalence entre la catégorie des ${\cal O}_{X\grot}$-modules cohérents et celle des ${\cal O}_X$-modules cohérents ; ils envoient faisceaux localement libres sur faisceaux localement libres ; ils préservent les groupes de cohomologie.

\medskip
\deux{convfasccoh} {\bf Convention.} Si $X$ est un espace analytique, un {\em faisceau cohérent sur $X$} sera un ${\cal O}_{X\grot}$-module cohérent.

\setcounter{cptbis}{0}
\medskip
\trois{rangxbon} Soit $X$ un bon espace analytique, soit $\cal F$ un faisceau cohérent sur $X$ et soit $x\in X$. On notera ${\cal F}\otimes{\cal O}_{X,x}$ la fibre en $x$ du faisceau $\cal F$ restreint à la catégorie des ouverts de $X$ ; elle coïncide avec ${\cal F}(V)\otimes_{{\cal A}_V}{\cal O}_{X,x}$ pour tout voisinage affinoïde $V$ de $x$. 

\medskip
\trois{rangxgen} Soit $X$ un espace analytique et soit $\cal F$ un faisceau cohérent sur $X$. Soit $x\in X$ ; si $V$ est un domaine affinoïde de $X$ contenant $x$, le $\hres(x)$-espace vectoriel ${\cal F}(V)\otimes_{{\cal A}_V}\hres(x)$ ne dépend pas de $V$ et sera noté ${\cal F}\otimes\hres(x)$ ; si $Y$ est un domaine analytique de $X$ contenant $x$, alors ${\cal F}_{|Y}\otimes \hres(x)$ coïncide avec ${\cal F}\otimes \hres(x)$ ; si $X$ est bon, ${\cal F}\otimes \hres(x)\simeq( {\cal F}\otimes{\cal O}_{X,x})\otimes_{{\cal O}_{X,x}}\hres(x)$. 

\medskip
\trois{fyx} Soit $f :Y\to X$ un morphisme d'espaces analytiques. Si $\cal F$ est un faisceau cohérent sur $X$, alors $f^*{\cal F}$ est un faisceau cohérent sur $Y$. Si $y\in Y$ (resp. si $Y$ est bon et si $y\in Y$), on écrira ${\cal F}\otimes\hres(y)$ (resp. ${\cal F}\otimes{\cal O}_{Y,y}$) au lieu de $f^*{\cal F}\otimes \hres(y)$ (resp. $f^*{\cal F}\otimes{\cal O}_{Y,y}$).

\medskip
\deux{zarcasgen} Soit $X$ un espace analytique et soit $\cal I$ un faisceau cohérent d'idéaux de ${\cal O}_{X\grot}$. Le {\em lieu des zéros de $\cal I$} est l'ensemble des points $x$ de $X$ tels que pour tout domaine analytique $V$ de $X$ contenant $x$ et toute fonction $f$ appartenant à ${\cal I}(V)$ l'on ait $f(x)=0$. 

\medskip
\setcounter{cptbis}{0}
\trois{zarcasgenaff} Soit $X$ un espace affinoïde et soit $I$ un idéal de son algèbre de fonctions ; le lieu des zéros de $I$ coïncide avec le lieu des zéros du faisceau cohérent $\red I$. 

\medskip
\trois{zarcasgendom}  Soit $f: Y\to X$ un morphisme entre espaces analytiques. Soit $\cal I$ un faisceau cohérent d'idéaux sur $X$ ; si $Z$ désigne le lieu des zéros de $\cal I$, alors $f^{-1}(Z)$ est le lieu des zéros de $f^*\cal I$. En particulier, si $Y$ est un domaine analytique de $X$, alors $Z\cap Y$ est le lieu des zéros de ${\cal I}_{|Y}$.

\medskip
\trois{topzarcasgen} Soit $X$ un espace analytique. La somme d'une famille quelconque (resp. le produit d'une famille finie) de faisceaux cohérents d'idéaux sur $X$ est un faisceau cohérent d'idéaux : on le teste sur un G-recouvrement de $X$ par des domaines affinoïdes, et cela se déduit alors de la noethérianité des algèbres affinoïdes.  Les parties qui peuvent se décrire comme le lieu des zéros d'un faisceau cohérent d'idéaux  sur $X$ sont donc les fermés d'une topologie appelée {\em topologie de Zariski} sur $X$ ; en vertu du \ref{RAP}.\ref{zarcasgen}.\ref{zarcasgenaff} et du  \ref{RAP}.\ref{faisccoh}.\ref{theokiehl}, cette définition est compatible avec la précédente lorsque $X$ est affinoïde ; un fermé de Zariski est donc fermé pour la topologie usuelle de $X$ (on vérifie cette propriété G-localement, ce qui permet de supposer $X$ affinoïde) ; {\em lorsqu'on parlera d'un espace analytique irréductible, ce sera toujours en référence à sa topologie de Zariski.} 

\medskip
\deux{supportfaisccoh} {\bf Exemple.} Soit $X$ un espace analytique et soit  $\cal F$ un faisceau cohérent sur $X$. Le {\em support} de $\cal F$ est l'ensemble des points $x$ de $X$ tels que ${\cal F}\otimes\hres(x)\neq 0$. C'est un fermé de Zariski : en effet, montrons que si $\cal I$ désigne l'annulateur de $\cal F$, c'est-à-dire le noyau de la flèche naturelle ${\cal O}_{X\grot}\to {\cal End}\;\cal F$, alors le support de $\cal F$ est le lieu des zéros de $\cal I$. On peut le vérifier G-localement et donc supposer que $X$ est bon. Un point $x\in X$ appartient au lieu des zéros de $\cal I$ si et seulement si ${\cal I}\otimes{\cal O}_{X,x}$ est inclus dans l'idéal maximal de ${\cal O}_{X,x}$, donc si et seulement si ${\cal F}\otimes{\cal O}_{X,x}\neq 0$, donc si et seulement si (par le lemme de Nakayama) ${\cal F}\otimes\hres(x)\neq 0$. 

\medskip
\deux{imfinizar} {\bf Exemple.} Soit $f:Y\to X$  un morphisme fini d'espaces analytiques. Le ${\cal O}_{X\grot}$-module $f_*{\cal O}_{Y\grot}$ est cohérent. Le sous-ensemble $f(Y)$ de $X$ coïncide avec le support de $f_*{\cal O}_{Y\grot}$ en vertu du lemme 2.1.6 de \cite{brk2} (ainsi que du lemme de Nakayama). En conséquence, $f(Y)$ est, d'après l'exemple \ref{RAP}.\ref{supportfaisccoh} ci-dessus, un fermé de Zariski de $X$.

\medskip
\deux{immferm} Soit $X$ un espace analytique, soit $\cal I$ un faisceau cohérent d'idéaux sur $X$ et soit $Y$ son lieu des zéros. La catégorie  des flèches $f: Z\to X$ telles que $f^*{\cal I}=0$ admet un objet final $\iota :{\mathsf V}({\cal I})\to X$ ; le $X$-espace {\em topologique} sous-jacent à $\mathsf V({\cal I})$ s'identifie à $Y$. Pour tout domaine affinoïde $V$ de $X$, l'intersection $V\cap Y$, vue comme sous-ensemble de ${\mathsf V}({\cal I})$, en est un domaine affinoïde ; la ${\cal A}_V$-algèbre correspondante est ${\cal A}_V/{\cal I}(V)$ ; en conséquence, la ${\cal O}_{X\grot}$-algèbre $\iota_*{\cal O}_{\mathsf V({\cal I})\grot}$ s'identifie à ${\cal O}_{X\grot}/{\cal I}$.  On dira de ${\mathsf V}({\cal I})$ qu'il est le {\em sous-espace analytique fermé de $X$ défini par $\cal I$} ; les foncteurs $\iota_*$ et $\iota^*$ établissent une équivalence entre la catégorie des faisceaux cohérents sur $\mathsf V({\cal I})$ et celle des faisceaux cohérents sur $X$ dont l'idéal annulateur contient $\cal I$. Si $g:T\to X$ est un morphisme d'espaces analytiques, alors $\mathsf V({\cal I})\times_XT\simeq \mathsf V(g^{-1}{\cal I}.{\cal O}_{T\grot})$. 

\medskip
Par abus,  il arrivera que l'on écrive $Y$ au lieu de $\mathsf V({\cal I})$, mais on aura pris soin au préalable de préciser que  l'on munit $Y$ de la structure définie par $\cal I$. Dans ce contexte, si $V$ est un domaine analytique de $X$, l'ensemble $V\cap Y$ sera {\em a priori} considéré comme muni de sa structure de domaine analytique de $Y$ (qui l'identifie à $Y\times_XV$) ; cette structure est aussi celle du sous-espace analytique fermé de $V$ défini par ${\cal I}_{|V}$, comme on le vérifie à l'aide des propriété universelles respectives de ces deux espaces. 

\medskip
\deux{defimmeferm} Un morphisme $f: Z\to X$ entre espaces analytiques est appelé une {\em immersion fermée} s'il existe un faisceau cohérent d'idéaux $\cal I$ sur $X$ et un $X$-isomorphisme entre $Z$ et ${\mathsf V}({\cal I})$ ; l'idéal $\cal I$ est dans ce cas uniquement déterminé, puisque il s'identifie au noyau de ${\cal O}_{X\grot}\to f_*{\cal O}_{Z\grot}$. La composée de deux immersions fermées est une immersion fermée, et le fait d'être une immersion fermée est stable par schangement de base. 

\medskip
\deux{zarind} Soit $X$ un espace analytique et soit $Y$ un fermé de Zariski de $X$. Choisissons un faisceau cohérent d'idéaux $\cal I$ dont $Y$ est le lieu des zéros, et munissons $Y$ de la structure correspondante. Soit $Z$ une partie de $Y$ ; alors $Z$ est un fermé de Zariski de $X$ si et seulement si $Z$ est un fermé de Zariski de $Y$. En effet, si $Z$ est un fermé de Zariski de $X$, c'est aussi un fermé de Zariski de $Y$ en vertu du \ref{RAP}.\ref{convfasccoh}.\ref{fyx} ; et si $Z$ est un fermé de Zariski de $Y$, il s'écrit comme le lieu des zéros d'un faisceau cohérent d'idéaux $\cal J$ de ${\cal O}_{Y\grot}$ ; si $\cal H$ désigne l'image réciproque de $\iota_*{\cal J}$ par ${\cal O}_{X\grot}\to \iota_*{\cal O}_{Y\grot}$, où $\iota$ est le morphisme naturel $Y\to X$, il est immédiat que $Z$ est le lieu des zéros du faisceau cohérent d'idéaux $\cal H$ de ${\cal O}_{X\grot}$. La topologie de Zariski de $Y$ est donc, indépendamment du faisceau cohérent d'idéaux utilisé pour faire de $Y$ un sous-espace analytique fermé de $X$, induite par la topologie de Zariski de $X$. 

\medskip
\deux{defomega} Soit $k$ un corps ultramétrique complet et soit $X$ un espace $k$-analytique. La diagonale $X\to X\times_{k}X$ est G-localement sur $X$ une immersion fermée (\cite{brk2}, \S 1.4) ; son {\em faisceau conormal} (\cite{brk2}, \S 1.3) est un faisceau cohérent sur $X$ qui sera noté $\Omega^1_{X/k}$ ; si $V$ est un domaine affinoïde de $X$ alors $\Omega^1_{X/k}(V)$ s'identifie au module des $k$-différentielles {\em bornées} de ${\cal A}_V$ (\cite{brk2},\S 3.3). Signalons que par souci de simplicité, nous avons choisi de noter $\Omega^1_{X/k}$ ce que Berkovich désigne par $\Omega^1_{X\grot/k}$.

\medskip
\deux{extcompl} Si $k$ est un corps ultramétrique complet, on notera $\red k$ son corps résiduel. On appellera {\em extension complète} de $k$ la donnée d'un corps ultramétrique complet $L$ et d'une injection isométrique $k\hookrightarrow L$. Si $L$ est une extension complète de $k$, on notera $\mathsf{d}(L/k)$ l'élément $$\mbox{deg. tr.}\; (\red L/\red k)+\dim \QQ \QQ\otimes_\ZZ ( |L^*|/|k^*|)\in \NN\cup \{+\infty\}.$$ Notons que $\mathsf{d}(./.)$ satsifait de manière évidente une règle de transitivité : si $L$ est une extension complète de $k$ et si $F$ est une extension complète de $L$ alors $\mathsf{d}(F/k)=\mathsf{d}(F/L) +\mathsf{d}(L/k)$. 

\medskip
\deux{composext} Soit $k$ un corps ultramétrique complet et soient $F$ et $K$ deux extensions complètes de $k$. Une {\em extension complète de $k$ composée de $F$ et $K$} est la donnée d'une extension complète $L$ de $k$ et de deux injections isométriques $F\hookrightarrow L$ et $K\hookrightarrow L$ qui coïncident sur $k$ et dont la réunion des images engendre un sous-corps dense de $L$. Une telle extension composée existe toujours : en effet, $F\otimes_kK$ s'injecte dans $F\hotimes_kK$ (\cite{gru}, th. 1, $4^{\mbox{o}}$), qui est donc un anneau de Banach non nul. Dès lors, ${\cal M}(F\hotimes_kK)$ est non vide (\cite{brk1}, th. 1.2.1) ; choisissons $y\in {\cal M}(F\hotimes_kK)$. Le corps $\hres(y)$ est alors une extension complète composée de $F$ et de $K$. 

\medskip
\deux{sensxl} Soit $k$ un corps ultramétrique complet et soit $L$ une extension complète de $k$. 

\setcounter{cptbis}{0}
\medskip
\trois{defal} Si $\cal A$ est une algèbre $k$-affinoïde, on notera ${\cal A}_{L}$ l'algèbre $L$-affinoïde $L\hotimes_{k}{\cal A}$ ; c'est une ${\cal A}$-algèbre fidèlement plate (\cite{brk2}, lemme 2.1.2). Si $\cal X$ est un $\cal A$-schéma de type fini, on posera ${\cal X}_{L}={\cal X}\otimes_{\cal A}{\cal A}_{L}$. 

\medskip
\trois{defxl} Si $X$ est un espace $k$-analytique, on notera $X_{L}$ l'espace $L$-analytique déduit de $X$ par extension des scalaires de $k$ à $L$. On dispose d'une application continue $X_{L}\to X$ qui est surjective ({\em cf.} \cite{duc3}, 0.5) ; lorsque $X$ est bon, $X_{L}\to X$ est plat en tant que morphisme d'espaces localement annelés (\cite{brk2}, cor. 2.1.3). 

\medskip
\trois{densehomeolimproj} Si $X$ est un espace $k$-analytique {\em compact} et s'il existe une famille fitrante $(L_i)$ de sous-corps complets de $L$ contenant $k$ telle que le corps $\bigcup L_i$ soit dense dans $L$, alors l'application continue canonique $X_L\to \lim\limits_{\leftarrow} X_{L_i}$ est un homéomorphisme. En effet comme $X_L$ est compact et $\lim\limits_{\leftarrow} X_{L_i}$ séparé (puisque $X$, étant compact, est {\em par définition} topologiquement séparé), il suffit de vérifier que cette application est bijective ; on peut, pour ce faire, supposer que $X$ est $k$-affinoïde ; notons $\cal A$ l'algèbre des fonctions analytiques sur $X$. Si $\phi$ et $\psi$ sont deux semi-normes multiplicatives bornées sur ${\cal A}_L$ dont les restrictions à ${\cal A}_{L_i}$ coïncident pour tout $i$, elles sont égales par densité de la réunion des images des ${\cal A}_{L_i}$ dans $\cal A$, ce qui établit l'injectivité de $X_L\to \lim\limits_{\leftarrow} X_{L_i}$. L'image $\mathsf I$ de $X_L$ par cette flèche est compacte, et donc fermée (rappelons que le but est séparé, {\em cf.} {\em supra}). Soit $\mathsf U$ son ouvert complémentaire ; s'il était non vide, il existerait un indice $j$ et un ouvert non vide $V$ de $X_{L_j}$ tel que $\mathsf U$ contienne l'image réciproque $\mathsf V$ de $V$ dans $\lim\limits_{\leftarrow} X_{L_i}$ ; mais $\mathsf V$ contient l'image de $V_L$, laquelle est un sous-ensemble non vide de $\mathsf I$ ; en conséquence, $\mathsf U\cap \mathsf I\neq \emptyset$, ce qui est absurde ; dès lors $\mathsf U$ est vide et $X_L\to \lim\limits_{\leftarrow} X_{L_i}$ est surjective, et finalement bijective, ce qu'il fallait démontrer. 

\medskip
\trois{raddensehomeo} Signalons un cas particulier important (que l'on aurait pu établir directement) de la situation générale traitée au \ref{RAP}.\ref{sensxl}.\ref{densehomeolimproj} ci-dessus : celui où chacun des $L_i$ est {\em une extension finie purement inséparable} de $k$ ; l'application naturelle $X_{L_i}\to X$ induit alors pour tout $i$ un homéomorphisme entre les espaces topologiques sous-jacents : en effet, comme $X_{L_i}$ est compact et $X$ séparé, il suffit de vérifier qu'elle est ensemblistement bijective, ce que l'on voit en considérant ses fibres. Puisque $X_L\simeq \lim\limits_{\leftarrow} X_{L_i}$ par le \ref{RAP}.\ref{sensxl}.\ref{densehomeolimproj},  $X_L\to X$ induit un homéomorphisme entre les espaces topologiques sous-jacents ; {\em cette assertion s'étend immédiatement, par réduction au cas affinoïde, à un espace $k$-analytique $X$ quelconque, {\em i.e.} non nécessairement compact.} 

\medskip
\deux{domaffplat} Soit ${\cal A}$ une algèbre affinoïde, et soit ${\cal B}$ l'anneau des fonctions d'un domaine affinoïde de ${\cal M}({\cal A})$. La ${\cal A}$-algèbre $\cal B$ est plate (\cite{brk1}, prop. 2.2.4). Par un argument de limite inductive, on en déduit le fait suivant : si $Y$ est un bon domaine analytique d'un bon espace analytique $X$, la flèche $Y\hookrightarrow X$ est un morphisme plat d'espaces localement annelés. 

\medskip
\deux{aletcalxl} Soit $k$ un corps ultramétrique complet, soit $\cal A$ une algèbre $k$-affinoïde et soit $\cal X$ un $\cal A$-schéma de type fini. On désignera par ${\cal X}\an$ l'analytifié de $\cal X$ (\cite{brk2}, prop. 2.6.1) ; c'est un bon espace $k$-analytique, muni d'une flèche surjective ${\cal X}\an \to {\cal X}$, qui est plate en tant que morphisme d'espaces localement annelés (\cite{brk2}, prop. 2.6.2) ; si ${\cal X}=\spec \cal A$ alors ${\cal X}\an={\cal M}({\cal A})$. Si $L$ est une extension complète de $k$, il résulte de la propriété universelle qui définit ${\cal X}\an$ que $({\cal X}_{L})\an$ et $({\cal X}\an)_{L}$ sont canoniquement isomorphes ; on pourra de ce fait écrire ${\cal X}\an_{L}$ sans qu'il y ait ambiguïté. 

\medskip
Soit $\cal F$ un faisceau cohérent sur $\cal X$ ; si $\bf x$ est un point de $\cal X$, on notera ${\cal F}\otimes{\cal O}_{{\cal X},\bf x}$ la fibre de $\cal F$ en $\bf x$. Si $\rho$ désigne le morphisme de sites annelés ${\cal X}\an \grot\to {\cal X}$, alors $\rho^*{\cal F}$ est un faisceau cohérent sur ${\cal X}\an$, que l'on notera ${\cal F}\an$ ; si $x$ est un point de ${\cal X}\an$ d'image $\bf x$ sur $\cal X$, on dispose d'un isomorphisme naturel ${\cal F}\an\otimes{\cal O}_{{\cal X}\an,x}\simeq ({\cal F}\otimes{\cal O}_{{\cal X},\bf x})\otimes_{ {\cal O}_{{\cal X},\bf x}}{\cal O}_{{\cal X}\an,x}.$ Pour cette raison, on écrira ${\cal F}\otimes{\cal O}_{{\cal X}\an,x}$ au lieu de ${\cal F}\an\otimes{\cal O}_{{\cal X}\an,x}$. 

\medskip
\deux{polyray} Un {\em polyrayon} sera une famille finie de réels strictement positifs. Si $k$ est un corps ultramétrique complet et si $\bf r$ est un polyrayon, on notera $k_{\bf r}$ l'algèbre des fonctions analytiques du domaine affinoïde de l'espace affine défini par la condition $|{\bf T}|=\bf r$ ; on dira que $\bf r$ est {\em $k$-libre} s'il constitue une famille libre de $\QQ\otimes_{\ZZ}(\RR^{*}_{+}/|k^{*}|)$. Si $\bf r$ est $k$-libre, $k_{\bf r}$ est un corps, et plus précisément une extension complète de $k$ ({\em cf.} \cite{duc3}, 1.2.2) ; dans ce cas, on mettra simplement $\bf r$ en indice au lieu de $k_{\bf r}$ pour noter l'extension des scalaires de $k$ à $k_{\bf r}$. Si $\cal A$ (resp. $X$) est une algèbre $k$-affinoïde (resp. un espace $k$-analytique) on dira qu'un polyrayon {\em déploie} $\cal A$ (resp. $X$) si $\bf r$ est $k$-libre, si $|k_{\bf r}^{*}|\neq\{1\}$ et si ${\cal A}_{\bf r}$ (resp. $X_{\bf r}$) est strictement $k_{\bf r}$-affinoïde (resp. strictement $k_{\bf r}$-analytique). 

\medskip
\deux{remanalkr} Soit $k$ un corps ultramétrique complet et soit $\bf r$ un polyrayon $k$-libre. Soit $\cal A$ une algèbre $k_{\bf r}$-affinoïde et soit $\cal X$ un $\cal A$-schéma de type fini. L'algèbre $\cal A$, vue simplement comme une $k$-algèbre de Banach, est $k$-affinoïde. Il résulte de la propriété universelle de l'analytification que l'espace ${\cal X}\an$ ne dépend pas du fait que l'on considère $\cal A$ comme $k$-affinoïde ou comme $k_{\bf r}$-affinoïde.

\medskip
\deux{secshilov} Soit $k$ un corps ultramétrique complet, soit $\bf r$ un polyrayon $k$-libre et soit $X$ un espace $k$-analytique. Soit $x\in X$. La fibre de $X_{\bf r}\to X$ en $X$ s'identifie à l'espace $\hres(x)$-affinoïde ${\cal M}(\hres(x)\hotimes_{k}k_{\bf r})$ ; on notera $\sigma(x)$ l'unique point du bord de Shilov de cet espace ; il correspond à la semi-norme $\sum a_I{\bf T}^i\mapsto  \max |a_I|{\bf r}^ I$ (où $a_I\in \hres(x)$ pour tout $I$). L'application $x\mapsto \sigma(x)$ est une section {\em continue} de $X_{\bf r}\to X$ (\cite{brk1}, lemma 3.2.2 $(i)$), que l'on appellera la {\em section de Shilov}. 

\medskip
Soit $\cal X$ un schéma de type fini sur une algèbre $k$-affinoïde $\cal A$ et soit $\cal U$ un ouvert de Zariski de ${\cal X}_{\bf r}$. Soit $\sigma :{\cal X}\an\to{\cal X}_{\bf r}\an$ la section de Shilov. Il existe alors un ouvert de Zariski $\cal V$ de $\cal X$ tel que $\sigma^{-1}({\cal U}\an)={\cal V}\an$ (remarquons qu'un tel $\cal V$ est nécessairement unique, puisqu'il coïncide alors forcément avec l'image de $\sigma^{-1}({\cal U}\an)$ sur $\cal X$). Pour le voir on se ramène, la question étant locale sur $\cal X$, au cas où ce dernier est affine ; soit $\cal B$ son anneau des fonctions. Il existe une famille $(f_j)$ d'éléments de ${\cal B}\otimes_{\cal A}{\cal A}_{\bf r}$ telle que $\cal U$ soit réunion des lieux d'inversibilité des $f_j$. Pour tout $j$, on peut écrire $f_j$ sous la forme $\sum a_{I,j}{\bf T}^{I}$, où les $a_{I,j}$ appartiennent à $\cal B$. Par définition de $\sigma$, l'ouvert $\sigma^{-1}({\cal U}\an)$ est l'ensemble des $x\in {\cal X}\an$ tels que $\max |a_ {I,j}(x)|{\bf r}^I\neq 0$ pour au moins un indice $j$, soit encore l'ensemble des $x\in {\cal X}\an$ tels que $a_{I,j}(x)\neq 0$ pour au moins un couple $(I,j)$ ; il s'identifie donc à ${\cal V}\an$, où $\cal V$ est l'ouvert de Zariski de $\cal X$ défini comme la réunion des lieux d'inversibilité des $a_{I,j}$.

\medskip
\deux{dimanal} Il existe une bonne théorie de la dimension dans le contexte des espaces de Berkovich, qui se comporte comme on s'y attend ; pour ses définitions (dues à Berkovich) et les démonstrations de ses principales propriétés, on pourra se reporter au paragraphe 1 de \cite{duc3}. Soit $X$ un espace analytique sur un corps ultramétrique complet $k$.

\medskip
On sait définir la {\em dimension $k$-analytique} de $X$, notée $\dim k X$, ainsi que la dimension $k$-analytique de $X$ en l'un de ses points $x$, notée $\dim {k,x}X$. Remarquons que $\dim k X$ et $\dim {k,x}X$ dépendent en général effectivement de $k$, et pas seulement de $X$ et $x$; ainsi si $r$ est un réel $k$-libre, $\dim k {\cal M}(k_r)=1$ et $\dim {k_r}{\cal M}(k_r)=0$. Rappelons maintenant quelques faits ; l'ordre dans lequel nous les présentons ici diffère en partie de celui dans lequel ils sont démontrés dans la littérature. 

\medskip
\setcounter{cptbis}{0} 

\trois{defdimstr} Si $X$ est strictement $k$-affinoïde alors $\dim k X$ coïncide avec la dimension de Krull de $X$ pour la topologie de Zariski.

\medskip
\trois{defcasgen} Dans le cas général, $\dim k X=\sup\limits_V \dim k V$, où $V$ parcourt l'ensemble des domaines affinoïdes de $X$.

\medskip
\trois{dimetd} L'on démontre que $\dim k X=\sup\limits_{x\in X}\mathsf{d}(\hres(x)/k)$. Il en découle que si $X\to {\cal M}(k)$ se factorise par ${\cal M}(L)$ pour une certaine extension {\em finie} $L$ de $k$, auquel cas $X$ peut naturellement être vu comme un espace $L$-analytique, alors $\dim L X=\dim k X$ : cela provient simplement du fait que $\mathsf{d}(L/k)=0$, et de la transitivité de $\mathsf{d}(./.)$.

\medskip
\trois{dimxl} Si $L$ est une extension complète de $k$ alors $\dim k X=\dim L X_L$.

\medskip
\trois{dimyv} Soit $Y$ un fermé de Zariski de $X$ et soit $V$ un domaine analytique de $X$ ; munissons $Y$ d'une structure de sous-espace analytique fermé de $X$. On a $\dim k V\leq \dim k X$ et $\dim k Y\leq \dim k X$. Comme $V\cap Y$ peut être vu aussi bien comme un sous-espace analytique fermé de $V$ que comme un domaine analytique que $Y$, on en déduit que $\dim k (V\cap Y)\leq \dim k V$ et $\dim k (V\cap Y)\leq \dim k Y$. 

\medskip
{\em On vérifie, par exemple à l'aide de la formule donnée au {\rm \ref{RAP}.\ref{dimanal}.\ref{dimetd}}, que $\dim k Y $ et $\dim k (V\cap Y)$ ne dépendent pas de la structure de sous-espace analytique fermé dont on a muni $Y$ ; pour cette raison, nous nous permettrons le plus souvent de les évoquer sans prendre la peine de fixer une telle structure.}

\medskip
\trois{defdimloc} Si $x\in X$ on pose $\dim{k,x}X=\inf\limits_V \dim k V$, où $V$ parcourt l'ensemble des domaines analytiques de $X$ qui contiennent $x$ ; on voit immédiatement que $\dim k X=\sup\limits_{x\in X}\dim {k,x}X$, et que si $V$ est un domaine analytique de $X$ et si $x$ est un point de $V$ alors $\dim{k,x}V=\dim{k,x}X$.

\medskip
\trois{dimlocaff} Si $X$ est $k$-affinoïde alors pour tout $x\in X$ l'entier $\dim {k,x}X$ coïncide avec le maximum des dimensions $k$-analytiques des composantes irréductibles de $X$ qui contiennent $x$, et $\dim k X$ est donc égal au maximum des dimensions $k$-analytiques des composantes irréductibles de $X$ 

\medskip 
\trois{dimlocxl} Si $L$ est une extension complète de $k$, si $x$ est un point de $X$ et si $y$ est un point de $X_L$ situé au-dessus de $x$ alors $\dim{k,x}X=\dim {L,y}X_L$. Notons qu'au vu du \ref{RAP}.\ref{dimanal}.\ref{dimlocaff} ci-dessus, ceci entraîne que si $X$ est affinoïde et irréductible de dimension $k$-analytique $d$, alors toutes les composantes irréductibles de $X_L$ sont de dimension $k$-analytique $d$ ; en réalité, on commence par démontrer ce fait pour en déduire l'invariance de la dimension locale par extension des scalaires.

\medskip
\trois{dimcodim} Si $X$ est $k$-affinoïde irréductible et si $Y$ est un fermé de Zariski irréductible de $X$ alors $\dim k Y +\codim{\mathsf {Krull}}(Y,X)= \dim k X$ ; dans le cas strictement $k$-affinoïde remarquons que cette égalité s'écrit tout simplement $\dim {\mathsf {Krull}} Y +\codim{\mathsf {Krull}}(Y,X)= \dim  {\mathsf {Krull}} X$, et c'est d'ailleurs ce cas particulier que l'on prouve en premier lieu avant d'en déduire la validité de la formule en général. 

\medskip
Notons deux conséquences de ce qui précède. La première est immédiate : si $X$ est $k$-affinoïde irréductible et si $Y$ est un fermé de Zariski strict de $X$, alors $\dim k Y<\dim k X$ ; la seconde se fonde sur le fait qu'un espace affinoïde irréductible reste équidimensionnel après extension des scalaires ({\em cf.} \ref{RAP}.\ref{dimanal}.\ref{dimlocxl}) : si $X$ est un espace $k$-affinoïde, si $Y$ est un fermé de Zariski de $X$ et si $L$ est une extension complète de $k$, alors  $\codim{\mathsf {Krull}}(Y_L,X_L)=\codim{\mathsf {Krull}}(Y,X)$.

\medskip
\trois{oxxkrull} Supposons que $|k^*|\neq \{1\}$ et que $X$ est strictement $k$-affinoïde, et soit $x$ un point {\em rigide} de $X$. Notons $\cal X$ le spectre de l'anneau des fonctions analytiques sur $X$ et désignons par $\bf x$ l'image de $x$ sur $\cal X$. On a alors $\dim{\mathsf {Krull}}{\cal O}_{X,x}=\dim{\mathsf {Krull}}{\cal O}_{{\cal X},\bf x}=\dim x X$. Cette égalité est établie au cours de la démonstration du lemme 1.23 de \cite{duc3}, avec une référence au lemme 1.12 de {\em op. cit.} Signalons que dans la preuve de ce dernier l'égalité du \ref{RAP}.\ref{dimanal}.\ref{dimcodim} est utilisée {\em implicitement} (dans le cas strictement $k$-affinoïde, où elle ne met en jeu que des dimensions de Krull).

\medskip
\deux{convdimk} Si le corps de base est déterminé sans ambiguïté (ou éventuellement à une extension finie près, {\em cf.} \ref{RAP}.\ref{dimanal}.\ref{dimetd}) on parlera simplement de dimension, et non de dimension $k$-analytique, et l'on utilisera des notations du type $\dim{}X$ ou $\dim x X$ de préférence à $\dim k X$ et $\dim {k,x}X$. Nous commettrons notamment,  sauf mention expresse du contraire, les abus suivants :

\medskip
\begin{itemize}
 \itb si l'on travaille avec un espace $X$ dont on a précisé qu'il était $k$-analytique pour un certain corps ultramétrique complet $k$, la notion de dimension, lorsqu'elle concernera $X$, l'un de ses domaines analytiques, l'un de ses fermés de Zariski, etc. , sera toujours à prendre au sens de \og dimension $k$-analytique\fg~; si $L$ est une extension complète de $k$ la notion de dimension, lorsqu'elle concernera $X_L$, l'un de ses domaines analytiques, l'un de ses fermés de Zariski, etc. , sera toujours à prendre au sens de \og dimension $L$-analytique\fg~;

\itb si l'on travaille avec un espace analytique $X$ {\em sans mention d'un corps de base}, la notion de dimension, lorsqu'elle concernera $X$, l'un de ses domaines analytiques, l'un de ses fermés de Zariski, etc. sera toujours à prendre au sens de \og dimension $k$-analytique\fg~, où $k$ est le corps ultramétrique complet qui fait {\em implicitement} partie de la donnée de $X$.  
\end{itemize}

\medskip
\deux{dimimfini} Soit $f: Y\to X$ un morphisme fini d'espaces analytiques. On a vu (exemple \ref{RAP}.\ref{imfinizar}) que $f(Y)$ est un fermé de Zariski de $X$. {\em Sa dimension est égale à celle de $Y$} : en effet si $y\in Y$, alors $\hres(y)$ est une extension finie de $\hres(f(y))$, ce qui entraîne l'égalité $\mathsf{d}(\hres(y)/\hres(f(y))=0$ ; on conclut à l'aide de la transitivité de $\mathsf{d}(./.)$. 

\medskip
\deux{remzarintvide} {\bf Remarque.} Soit $X$ un espace affinoïde irréductible, soit $d$ sa dimension et soit $Z$ un fermé de Zariski strict de $X$. Soit $V$ un domaine analytique non vide de $X$. Comme $X$ est irréductible, il est purement de dimension $d$ et $V$ est donc de dimension $d$ ; des inégalités $\dim {} V\cap Z\leq \dim {} Z < d$, l'on déduit que $V$ n'est pas contenu dans $Z$ ; ceci implique notamment que l'intérieur {\em topologique} de $Z$ dans $X$ est vide. 

\subsection*{Algèbre commutative}

\bigskip
\deux{listepropri} {\bf Les propriétés usuelles de l'algèbre commutative.} On définit trois ensembles de propriétés : 

\begin{itemize}
\itb l'ensemble ${\cal S}$, qui comprend les propriétés $S_{m}$ pour $m$ dans $\NN$, ainsi que celle d'être de Cohen-Macaulay ;  ce sont des propriétés des modules de type fini sur un anneau local noethérien ;
\itb l'ensemble $\cal Q$, qui comprend la propriété d'être de Gorenstein, ainsi que celle d'être une intersection complète ; ce sont des propriétés des anneaux locaux noethériens ; 
\itb l'ensemble $\cal R$, qui comprend les propriétés $R_{m}$ pour $m$ dans $\NN$, ainsi que celle d'être régulier ; ce sont des propriétés des anneaux locaux noethériens. 
\end{itemize}

\medskip
\deux{platdirect} Soit $f:{\cal Y}\to {\cal X}$ un morphisme {\em plat} de schémas localement noethériens, soit ${\bf y}\in \cal Y$ et soit $\bf x$ son image sur $\cal X$. Soit $\cal F$ un faisceau cohérent sur $\cal X$. Soit ${\mathsf P}\in {\cal S}\cup {\cal Q}\cup {\cal R}$ ; si ${\mathsf P}\in{\cal Q}\cup {\cal R}$, on suppose que ${\cal F}={\cal O}_{\cal X}$. Si ${\cal F}\otimes {\cal O}_{{\cal Y},\bf y}$ satisfait $\mathsf P$, alors ${\cal F}\otimes {\cal O}_{{\cal X},\bf x}$ satisfait $\mathsf P$ ; si ${\cal F}\otimes {\cal O}_{{\cal X},\bf x}$ et $ {\cal F}\otimes {\cal O}_{f^{-1}({\bf x}),\bf y}$ satisfont $\mathsf P$, alors ${\cal F}\otimes {\cal O}_{{\cal Y},\bf y}$ satisfait $\mathsf P$.

\medskip
Ce sont des assertions classiques ; pour leurs  justifications, on pourra par exemple se référer à  \cite{ega42}, prop. 6.4.1 $i)$ et $ii)$ concernant les propriétés appartenant à $\cal S$, à {\em cf.} \cite{ega42}, prop. 6.5.3 $i)$ et $ii)$ concernant celles appartenant à $\cal R$, à  {\em cf.} \cite{mats}, th. 23.4 concernant le caractère de Gorenstein et enfin à  \cite{avr} pour le caractère d'intersection complète. 

\medskip
\deux{defgeomregul} Soit $f:{\cal Y}\to {\cal X}$ un morphisme de schémas localement noethériens. On dira que $f$ est {\em géométriquement régulier} s'il est plat et si pour tout ${\bf x}\in \cal X$ et toute extension finie purement inséparable $F$ de $\kappa({\bf x})$ le schéma localement noethérien $f^{-1}({\bf x})\otimes_{\kappa({\bf x})}F$ est régulier. 

\medskip
\deux{defexc} La définition de l'excellence que nous donnons ci-dessous est celle qui figure dans EGA IV (\cite{ega42}, \S 7.8 et notamment déf. 7.8.2) ; les quelques propriétés que nous énonçons ensuite sont démontrées dans {\em loc. cit.}

\medskip
Un anneau noethérien $\mathsf A$ est dit {\em excellent} s'il satisfait les propriétés suivantes :

\medskip
\begin{itemize}

\item[$i)$]  $\mathsf A$ est universellement caténaire ;

\medskip
\item [$ii)$] pour tout idéal premier $\got p$ de $A$, le morphisme $\spec \widehat {\mathsf A_{\got p}}\to \spec \mathsf A_{\got p}$ est géométriquement régulier\footnote{Sa platitude est automatique ; c'est donc uniquement sur la régularité géométrique des fibres que porte cette condition.} ;

\medskip
\item[$iii)$] pour tout idéal premier $\got p$ de $\mathsf A$ et pour toute extension finie radicielle $F$ de $\kappa(\got p)$ il existe une sous $\mathsf A/\got p$-algèbre finie $\mathsf B$ de $F$, de corps des fractions égal à $F$, et telle que le lieu régulier de $\spec \mathsf B$ contienne un ouvert non vide. 

\end{itemize}

\medskip
Si $\mathsf A$ est local, il suffit, pour qu'il soit excellent, qu'il satisfasse $i)$, et que $ii)$ soit vraie lorsque $\got p$ est son idéal {\em maximal}.  Si $\mathsf A$ est excellent, alors tout localisé d'une $\mathsf A$-algèbre de type fini est encore excellent ; tout anneau excellent est universellement japonais ; si un schéma localement noethérien $\cal X$ admet un recouvrement par des ouverts affines dont les anneaux de fonctions sont excellents alors ${\cal O}_{\cal X}({\cal U})$ est excellent pour {\em tout} ouvert affine $\cal U$ de $\cal X$ ; on dit dans ce cas que $\cal X$ est excellent. 

\medskip
\deux{exexcell} {\bf Exemples}. Il résulte immédiatement de la définition que tout corps est excellent et que $\ZZ$ est excellent ; on peut montrer que les anneaux locaux d'un espace analytique {\em complexe} sont excellents (on le voit au moyen d'un «critère jacobien» ; {\em cf.} \cite{mats2}, th. 100, th. 101, et la remarque qui les suit). 

\medskip
\deux{lieuxouverts}Soit $\cal X$ un schéma excellent, soit $\mathsf P\in{\cal S}\cup {\cal Q}\cup {\cal R}$ et soit $\cal F$ un faisceau cohérent sur $\cal X$ ; si ${\mathsf P}\in{\cal Q}\cup {\cal R}$, on suppose que ${\cal F}={\cal O}_{\cal X}$. L'ensemble $\cal U$ des points $\bf x$ de $\cal X$ tels que ${\cal F}\otimes {\cal O}_{{\cal X},\bf x}$ satisfasse $\mathsf P$ est alors un ouvert de Zariski de $\cal X$. En effet, hormis pour les propriétés d'être de Gorenstein ou d'intersection complète, cela résulte de \cite{ega42}, scholie 7.8.3, $iv)$ ; en ce qui concerne le fait d'être de Gorenstein, cela découle d'un résultat de Greco et Marinari (\cite{grecmar}, {\em cf.} \cite{mats}, th. 24.6) combiné à la validité de notre assertion pour la régularité, au fait que celle-ci est vraie en le point générique d'un schéma intègre, et à celui qu'elle entraîne le caractère de Gorenstein ; la justification est analogue en ce qui concerne le caractère d'intersection complète (le résultat de Greco et Marinari à utiliser dans ce cas figure aussi dans \cite{grecmar}, et est mentionné sans preuve dans \cite{mats} juste après le th. 24.6). 

\medskip
\deux{reminterscompl} {\bf Remarque.} Si $\cal X$ est un schéma localement isomorphe à un sous-schéma fermé d'un schéma régulier alors l'ensemble des points $\bf x$ de $\cal X$ tels que ${\cal O}_{{\cal X},\bf x}$ soit d'intersection complète est un ouvert de $\cal X$ et ce, sans hypothèse d'excellence (\cite{ega44}, cor. 19.3.3). 

\subsection*{Quelques lemmes techniques}

\medskip
\deux{tenseurkrintegre} {\bf Lemme.} {\em Soit $k$ un corps ultramétrique complet et soit $\bf r$ un polyrayon $k$-libre. Si $\cal X$ est un schéma de type fini sur une algèbre $k$-affinoïde alors $\cal X$ est intègre si et seulement si ${\cal X}_{\bf r}$ est intègre.}

\medskip
{\em Démonstration.} On se ramène aussitôt au cas où $\cal X$ est affine ; l'on note $\cal B$ l'anneau des fonctions du schéma $\cal X$. De la description de ${\cal B}_{\bf r}$ comme un anneau de séries formelles à coefficients dans $\cal B$ l'on déduit immédiatement que $\cal B$ s'injecte dans ${\cal B}_{\bf r}$ ; il en découle que si ${\cal B}_{\bf r}$ est intègre alors $\cal B$ est intègre. 

\medskip
Supposons que $\cal B$ est intègre. Puisque ${\cal B}\hookrightarrow {\cal B}_{\bf r}$, l'anneau  ${\cal B}_{\bf r}$ est non nul. Soient $f$ et $g$ deux éléments de ${\cal B}_{\bf r}$ tels que $fg=0$ ; écrivons $f=\sum f_I{\bf T}^I$ et $g=\sum g_I{\bf T}^I$ où les $f_I$ et $g_I$ appartiennent à $\cal B$. Soit $x\in {\cal X}\an$. La fibre en $x$ de ${\cal X}_{\bf r}\an \to {\cal X}\an$ s'identifie à ${\cal M}(\hres(x)_{\bf r})$ et l'image de $f$ (resp. $g$) dans $\hres(x)_{\bf r}$ modulo cette identification n'est autre que $\sum f_I(x){\bf T}^I$ (resp. $\sum g_I(x){\bf T}^I$) ; l' on a donc $(\sum f_I(x){\bf T}^I)(\sum g_I(x){\bf T}^I)=0$ dans $\hres(x)_{\bf r}$ ; comme ce dernier est intègre (\cite{duc3}, 1.2.1), l'on a ou bien $f_I(x)=0$ pour tout $I$, ou bien $g_I(x)=0$ pour tout $I$. Soit ${\cal I}_1$ (resp. ${\cal I}_2$) l'idéal de $\cal B$ engendré par les $f_I$ (resp. les $g_I$) et soit ${\cal Z}_1$ (resp. ${\cal Z}_2)$ le fermé de Zariski correspondant de $\cal X$. 

\medskip
Par ce qui précède, ${\cal X}\an={\cal Z}_1\an\cup {\cal Z}_2\an$ ; comme ${\cal X}\an\to {\cal X}$ est surjective, ${\cal X}={\cal Z}_1\cup {\cal Z}_2$. Le schéma $\cal X$ étant intègre, on a ${\cal I}_1=0$ ou ${\cal I}_2=0$ ; si l'on est dans le premier (resp. dans le second) cas alors tout les $f_I$ (resp. tous les $g_I$) sont nuls, et $f=0$ (resp. $g=0$) ; l'anneau  ${\cal B}_{\bf r}$ est donc intègre.~$\Box$

\medskip
\deux{prodalgplat} {\bf Lemme (Berkovich).} {\em Soit $k$ un corps ultramétrique complet, soit $\cal A$ une $k$-algèbre de Banach et soit $\cal E$ un $k$-espace de Banach. Le $\cal A$-module ${\cal A}\hotimes_k {\cal E}$ est plat et ${\cal A}\otimes_k {\cal E}\to{\cal A}\hotimes_k {\cal E}$ est injective.}

\medskip
{\em Démonstration.} C'est le lemme 2.1.2 de \cite{brk2} ; il est énoncé lorsque $\cal A$ est $k$-affinoïde et lorsque $\cal E$ est l'espace de Banach sous-jacent à une extension complète $K$ de $k$, mais sa démonstration n'utilise pas ces hypothèses et établit en réalité notre assertion.~$\Box$ 

\medskip
\deux{dimproduit} {\bf Lemme.} {\em Soit $k$ un corps ultramétrique complet et soient $X$ et $Y$ deux espaces $k$-analytiques. Soit $z$ un point de $X\times_k Y $ et soient $x$ et $y$ ses images respectives sur $X$ et $Y$. On a $\dim z X\times_kY=\dim x X+\dim y Y$.} 

\medskip
{\em Démonstration.} On se ramène aussitôt au cas où $X$ et $Y$ sont $k$-affinoïdes. On peut prouver le résultat après extension du corps de base, ce qui autorise à supposer que $|k^*|\neq \{1\}$, que $X$ et $Y$ sont strictement $k$-affinoïdes, et que $z$ est un point $k$-rationnel. On pose $Z=X\times_kY$ ; l'on désigne par $\cal X$ le spectre de l'anneau des fonctions analytiques sur $X$ et par $\bf x$ l'image de $x$ sur $\cal X$ ; l'on définit de manière analogue ${\cal Y},{\cal Z},{\bf y}$ et $\bf z$ ; on appelle $\got m$ l'idéal maximal de ${\cal O}_{{\cal X},\bf x}$. Le schéma ${\cal Z}$ est plat    d'après le lemme \ref{RAP}.\ref{prodalgplat} ci-dessus, et l'on a donc $$\dim{\mathsf{Krull}}{\cal O}_{{\cal Z},\bf z}=\dim{\mathsf{Krull}}{\cal O}_{{\cal X},\bf x}+\dim{\mathsf{Krull}}{\cal O}_{{\cal Z},\bf z}/\got m.$$ Comme $z$ est rigide, l'anneau local  ${\cal O}_{{\cal Z},\bf z}/\got m$ s'identifie à ${\cal O}_{{\cal Y},\bf y}$. L'égalité ci-dessus peut alors, en vertu du \ref{RAP}.\ref{dimanal}.\ref{oxxkrull}, se réécrire $\dim z Z=\dim x X+\dim y Y$, ce qu'il fallait démontrer.~$\Box$

\medskip
\deux{localintegre} {\bf Lemme.} {\em Soit $X$ un espace affinoïde et soit $x$ un point de $X$ tel que ${\cal O}_{X,x}$ soit intègre. Le point $x$ n'est situé que sur une seule composante irréductible de $X$.} 

\medskip
{\em Démonstration.} Soit $Y$ une composante irréductible de $X$ contenant $x$, et soit $Z$ la réunion des composantes irréductibles de $X$ distinctes de $Y$. Soit $f$ une fonction analytique sur $X$ qui est nulle en tout point de $Y$ et dont le lieu des zéros ne contient aucune autre composante, et soit $g$ une fonction analytique sur $X$ qui est nulle en tout point de $Z$ et dont le lieu des zéros ne contient pas $Y$. Le produit $fg$ étant nilpotent et ${\cal O}_{X,x}$ étant intègre, l'une au moins des deux fonctions $f$ et $g$ est nulle au voisinage de $x$. Comme $g_{|Y}$ n'est pas identiquement nulle et comme $Y$ est irréductible, le lieu des zéros de $g$ ne contient aucun ouvert non vide de $Y$ (remarque \ref{RAP}.\ref{remzarintvide}) ; on en conclut que $f$ est nulle au voisinage de $x$. Si $T$ est une composante irréductible de $X$ distincte de $Y$, alors $f_{|T}$ n'est pas identiquement nulle ; par le même argument que précédemment, le lieu des zéros de $f$ ne contient aucun ouvert non vide de $T$. Il en découle que $x\notin T$.~$\Box$

\medskip
\deux{remjpoineau} {\bf Remarque.} Sous les hypothèses du lemme ci-dessus, Jérôme Poineau a montré (\cite{jpoi1}, cor. 4.7) que le point $x$ possède une base de voisinages affinoïdes irréductibles ; nous n'aurons pas besoin de cette assertion, dont la preuve est très délicate. 

\medskip
\deux{quotoxx} {\bf Remarque.} Conservons les hypothèses et notations du lemme ci-dessus. Soit $\cal A$ l'algèbre des fonctions analytiques sur $X$  et soit $\got p$ l'idéal premier de $\cal A$ qui correspond à l'unique composante irréductible de $X$ contenant $x$. Chaque élément de $\got p$ s'annule alors en tout point d'un voisinage de $x$, donc est nilpotent au voisinage de $x$, donc d'image nulle dans ${\cal O}_{X,x}$ puisque ce dernier est intègre. En conséquence, ${\cal O}_{{\cal M}({\cal A}/\got p),x}={\cal O}_{X,x}$. 

\medskip
\deux{idminoxx} {\bf Lemme.} {\em Soit $X$ un espace affinoïde, soit $x\in X$ et soient $\got p_1,\ldots,\got p_r$ les idéaux premiers minimaux de ${\cal O}_{X,x}$. Il existe un voisinage affinoïde $V$ de $x$ dans $X$ et des idéaux $\got q_1,\ldots,\got q_r$ de ${\cal A}_V$ tels que les propriétés suivantes soient vérifiées, en désignant pour tout $i$ par $F_i$ le lieu des zéros de $\got q_i$ dans $V$ :

\medskip
$i)$ $\got q_i {\cal O}_{X,x}=\got p_i$ pour tout $i$ ;

\medskip
$ii)$ pour tout voisinage affinoïde $W$ de $x$ dans $V$ et pour tout indice $i$ le point $x$ n'appartient qu'à une composante irréductible $G_i$ de $F_i\cap W$ ; les $G_i$ sont deux à deux distinctes et sont exactement les composantes irréductibles de $W$ contenant $x$.} 

\medskip
{\em Démonstration.} Par noethérianité de ${\cal O}_{X,x}$ il existe un voisinage affinoïde $V$ de $x$ et $r$ idéaux  $\got q_1,\ldots,\got q_r$ de ${\cal A}_V$ tels que $\got q_i {\cal O}_{X,x}=\got p_i$ pour tout $i$ et tels que $\prod \got q_i$ soit nilpotent. Pour tout $i$, désignons par $F_i$ le lieu des zéros de $\got q_i$ dans $V$.

\medskip
Soit $W$ un voisinage affinoïde de $x$ dans $V$. Munissons $F_i\cap W$ de sa structure définie par l'idéal $\got q_i{\cal A}_W$. L'anneau local ${\cal O}_{F_i\cap W,x}$ s'identifie à ${\cal O}_{X,x}/\got p_i$ et est donc intègre. Le lemme \ref{RAP}.\ref{localintegre} ci-dessus assure que $x$ n'appartient qu'à une composante irréductible $G_i$ de $F_i\cap W$. Comme $\prod \got q_i$ est nilpotent, $V=\bigcup F_i$ et $W=\bigcup W\cap F_i$. Si ${\cal E}_i$ désigne pour tout $i$ l'ensemble des composantes irréductibles de $W\cap F_i$, alors l'ensemble des composantes irréductibles de $W$ est l'ensemble des éléments maximaux de $\bigcup {\cal E}_i$. Il suffit donc pour conclure de s'assurer que si $G_i\subset G_j$ alors $i=j$. 

\medskip
Soient donc $i$ et $j$ deux indices tels que $G_i \subset G_j$ et soit $f\in \got q_j$. Il existe un voisinage affinoïde $W'$ de $x$ dans $V$ tel que $W'\cap F_i\subset G_i$. La fonction $f$ est nulle sur $G_j$, donc sur $G_i$, donc sur $W'\cap F_i$. L'image de $f$ dans ${\cal A}_W'$ appartient donc à $\sqrt{\got q_i {\cal A}_W'}$. On en déduit que $\got p_j\subset \sqrt{\got p_i}=\got p_i$ ; par conséquent, $i=j$.~$\Box$ 

\medskip
Le lemme ci-dessous reprend un résultat de Berkovich (\cite{brk2}, lemme 2.1.6) sur les anneaux locaux analytiques et en énonce une conséquence immédiate concernant leurs complétés. 

\medskip
\deux{finicomplet} {\bf Lemme ( Berkovich).}  {\em Soit $\cal A$ une algèbre affinoïde et soit ${\cal Y}\to {\cal X}$ un morphisme fini entre deux $\cal A$-schémas de type fini. Soit $x$ un point de ${\cal X}\an$ et soit $\{y_{1},\ldots,y_{n}\}$ l'ensemble de ses antécédents sur ${\cal Y}\an$.

\medskip
$i)$ On dispose d'un isomorphisme naturel $${\cal Y}\times_{\cal X} \spec {\cal O}_{{\cal X}\an,x}\simeq \coprod \spec {\cal O}_{{\cal Y}\an,y_{i}}.$$

\medskip
$ii)$ On dispose d'un isomorphisme naturel $${\cal Y}\times_{\cal X} \spec \widehat{{\cal O}_{{\cal X}\an,x}}\simeq \coprod \spec \widehat{{\cal O}_{{\cal Y}\an,y_{i}}}.$$}

{\em Démonstration.} La question étant locale sur ${\cal X}\an$, on se ramène au cas où ${\cal X}=\spec \cal A$ ; 
l'assertion $i)$ est exactement le lemme 2.1.6 de \cite{brk2}, l'assertion $ii)$ en découle trivialement.~$\Box$

\medskip
\deux{remdenom} {\bf Remarque à propos des objets de type dénombrable.} Soit $k$ un corps ultramétrique complet et soit $k_{0}$ un sous-corps fermé de $k$. Nous aurons à plusieurs reprises recours au principe général suivant : {\em si $\mathsf O$ est un objet défini sur $k$ dont la description ne nécessite qu'une quantité (au plus) dénombrable de paramètres, alors $\mathsf O$ est déjà défini sur un sous-corps complet de $k$ contenant $k_{0}$ et topologiquement de type dénombrable sur ce dernier} ; la plupart du temps, nous nous contenterons d'énoncer lorsque nous en aurons besoin une déclinaison précise de cet énoncé un peu vague, en laissant au lecteur le soin de la justifier par les arguments standard (et laborieux); seul un cas, concernant le complété d'un anneau local analytique, nous a paru requérir des explications détaillées (lemme~\ref{EXC}.\ref{ouvertsreg}.\ref{vllimvf}). 

\section{Extensions analytiquement séparables}\label{ANSEP}

\setcounter{cpt}{0}

\bigskip
\deux{famscal} On dira que deux familles $(a_{i})_{i\in I}$ et $(b_{i})_{i\in I}$ de réels strictement positifs sont \textit{équivalentes} s'il existe un couple $(A,B)$ de $(\RR^{*}_{+})^{2}$ tel que $Aa_{i}\leq b_{i}\leq Ba_{i}$ pour tout $i$. Soit $k$ un corps ultramétrique complet ; pour tout réel $a>0$, soit ${\cal L}_{k,a}$ le $k$-espace de Banach dont l'espace vectoriel sous-jacent est $k$ et dont la norme envoie $1$ sur $a$ ; on notera $1_{k,a}$ le scalaire $1$ {\em vu comme appartenant à ${\cal L}_{k,a}$.} 

\medskip 
Soit $M$ un $k$-espace de Banach et soit ${\bf m}=(m_{i})$ une famille d'éléments de $M$. Si ${\bf a}=(a_{i})$ est équivalente à $||{\bf m}||:=(||{\bf m}_i||)$, il existe une unique application linéaire bornée $\Phi_{\bf a}$ de $\widehat{\bigoplus}{\cal L}_{k,a_{i}}$ vers $M$ qui envoie $1_{k,a_{i}}$ sur $m_{i}$ pour tout $i$.

\medskip
Rappelons qu'une application linéaire bornée $\phi : M\to N$ entre deux $k$-espaces de Banach est dite  {\em admissible} si la bijection naturelle $\mbox{Im}\;\phi\to M/\mbox{Ker}\;\phi$ est bornée, l'espace de gauche (resp. de droite) étant muni de la norme induite (resp. quotient).

\bigskip
\deux{injadm} \textbf{Lemme}. {\em Les propositions suivantes sont équivalentes. 

\begin{itemize}

\medskip
\item[$i)$] Pour toute famille $\bf a $ équivalente à $||\bf m||$, l'application $\Phi_{\bf a}$ est injective (resp. bijective) et admissible. 

\medskip
\item[$ii)$] L'application $\Phi_{||\bf m ||}$ est injective (resp. bijective) et admissible.

\medskip
\item[$iii)$] Il existe une famille $(b_{i})$ et une injection (resp. bijection) linéaire admissible $\Phi$ de $\widehat{\bigoplus}{\cal L}_{k,b_{i}}$ vers $M$ qui envoie $1_{k,b_{i}}$ sur $m_{i}$ pour tout $i$. 

\medskip

\end{itemize}}

\textit{Démonstration}. Il est trivial que $i)\Rightarrow ii) \Rightarrow iii)$. Si $iii)$ est vraie, alors $(b_{i})$ est équivalente à $||\bf m||$. Si ${\bf a}=(a_{i})$ est équivalente à $||\bf m||$, elle est équivalente à $(b_{i})$, et il existe donc un isomorphisme admissible $$\Psi : \widehat{\bigoplus}{\cal L}_{k, a_{i}}\simeq
\widehat{\bigoplus}{\cal L}_{k,b_{i}}$$ qui envoie $1_{k,a_{i}}$ sur $1_{k,b_{i}}$ pour tout $i$ ; par construction, $\Phi_{\bf a}=\Phi\circ \Psi$, et $\Phi_{\bf a}$ est donc injective (resp. bijective) et admissible.~$\Box$

\bigskip
\deux{introklibre} {\bf Définition.} On dira que $\bf m$ est {\em une famille topologiquement $k$-libre} (resp. est une {\em $k$-base topologique}) de $M$ si $\Phi_{||\bf m ||}$ est injective (resp. bijective) et admissible. 

\setcounter{cptbis}{0}

\bigskip
\trois{changenorme} Toute injection (resp. bijection) admissible transforme une famille topologiquement $k$-libre (resp. une $k$-base topologique) en une famille topologiquement $k$-libre (resp. en une $k$-base topologique).

\bigskip
\trois{changeklibre} Soit $L$ une extension complète de $k$. Le foncteur $L\hotimes_{k}$ transforme suites exactes courtes admissibles en suites exactes courtes admissibles (\cite{gru}, \S 2, th. 1 ; \cite{brk2}, preuve du lemme 2.1.2) ; comme ${\cal L}_{k,a}\hotimes_{k}L\simeq {\cal L}_{L,a}$ pour tout $a$, on en déduit que si $\bf m$ est une famille topologiquement $k$-libre  (resp. est une $k$-base topologique) de $M$, alors $1 \hotimes{\bf m}$ est une famille topologiquement $L$-libre (resp. est une $L$-base topologique) de $L\hotimes_{k}M$.

\bigskip
\trois{expoklibre} Soit $r$ un réel strictement positif. Soit $\Phi : N\to M$ une application linéaire entre deux $k$-espaces de Banach. Si $\Phi$ est une injection (resp. une bijection) admissible, c'est encore une injection (resp. une bijection) admissible lorsqu'on la considère comme une application $(k,|.|^{r})$-linéaire de $(N,||.||^{r})$ vers $(M,||.||^{r})$. Par ailleurs, pour tout $a>0$, le $(k,|.|^{r})$-espace de Banach $({\cal L}_{k,a},||.||^{r})$ s'identifie à ${\cal L}_{(k,|.|^{r}),a^{r}}$.

On déduit de ces remarques que si $\bf m $ est une famille topologiquement $k$-libre (resp. est une $k$-base topologique) du $k$-espace de Banach $M$, alors $\bf m$ est une famille topologiquement $(k,|.|^{r})$-libre (resp. est une $(k,|.|^{r})$-base topologique) du $(k,|.|^{r}) $-espace de Banach $(M,||.||^{r})$.

\bigskip
\deux{defplibre} Soit $k$ un corps ultramétrique complet de caractéristique $p>0$. Soit $\cal B$ une $k$-algèbre de Banach, et soit ${\bf a}=(a_{i})_{i \in I}$ une famille d'éléments de $\cal B$, tels que $a_{i}^{p}$ appartienne pour tout $i$ à l'image de $k$. Soit $\cal N$ l'ensemble des familles, indexées par $I$ et à support fini, d'entiers compris entre $0$ et $p-1$. On dira que $\bf a$ est \textit{une famille topologiquement $p$-libre} (resp. est \textit{une $p$-base topologique})  de $\cal B$ sur $k$ si $({\bf a}^{\bf n})_{{\bf n} \in {\cal N}}
$ est une famille topologiquement $k$-libre (resp. est une $k$-base topologique) de $\cal B$. 

\bigskip
\deux{kiehlpbase} {\bf Exemple}. Soit $F$ une sous-$k$-extension complète de $k^{1/p}$, topologiquement de type dénombrable sur $k$. Alors $F$ possède une $p$-base topologique sur $k$ : si $|k^{*}|=\{1\}$, il n'y a qu'à prendre une $p$-base algébrique ; sinon, c'est un théorème de Kiehl (\cite{kie} Satz 1.4). 

\medskip
\deux{pbasestable} {\bf Lemme.} {\em Soit $k$ un corps ultramétrique complet de caractéristique $p>0$. Soit $L$ une extension complète de $k$, et soit $F$ une sous-$k$-extension complète de $k^{1/p}$ ; on considère $k^{1/p}$ et $F$ comme plongés dans $L^{1/p}$. Supposons que $F$ possède une $p$-base topologique $\bf a=(a_{i})$ sur $k$ qui  est topologiquement $p$-libre sur $L$. La semi-norme $$\sum l_{i}\otimes x_{i}\mapsto |\underbrace{\sum l_{i}x_{i}}_{\in L^{1/p}}|$$ de $L\hotimes_{k}F$ est alors une norme, équivalente à la norme tensorielle.} 

\bigskip
{\em Démonstration.} Soit $\cal N$ l'ensemble des familles, indexées par $I$ et à support fini, d'entiers compris entre $0$ et $p-1$. Il existe un isomorphisme admissible $\Phi: \widehat{\bigoplus\limits_{{\bf n}\in {\cal N}}} {\cal L}_{k,|{\bf a}^{\bf n}|}\simeq F$, qui envoie $1_{k,|{\bf a}^{\bf n}|}$ sur ${\bf a}^{\bf n}$ pour tout $\bf n$ ; notons $\Phi_{L}$ la flèche déduite de $\Phi$ après tensorisation par $L$ ; c'est encore un isomorphisme admissible. 

\medskip
Posons $\mu=\sum l_{i}\otimes x_{i}\mapsto \sum l_{i}x_{i}$. La composée $\mu \circ \Phi_{L}$ est une application $L$-linéaire bornée de $\widehat{\bigoplus}{\cal L}_{L,|{\bf a}^{\bf n}|}$ vers $L^{1/p}$ qui envoie $1_{L,|{\bf a}^{\bf n}|}$ sur ${\bf a}^{\bf n}$ pour tout $\bf n$. Or ${\bf a}$ est une famille topologiquement $p$-libre de la $L$-algèbre $L^{1/p}$ ; en conséquence, $\mu\circ \Phi_{L}$ est une injection admissible, d'où il résulte que $\mu$ est elle aussi injective et admissible ; ceci achève la démonstration.~$\Box$ 

\bigskip
\deux{normequiv} {\bf Proposition.} {\em Soit $k$ un corps ultramétrique complet de caractéristique $p>0$. Soit $L$ une extension complète de $k$ ; on considère $k^{1/p}$ comme plongé dans $L^{1/p}$. Les propriétés suivantes sont équivalentes : 

\bigskip
\begin{itemize}

\item[$i)$] la semi-norme $\sum l_{i}\otimes{x_{i}}\mapsto |\underbrace{\sum l_{i}x_{i}}_{\in L^{1/p}}|$ de l'algèbre $L\hotimes_{k}k^{1/p}$ est une norme équivalente à la norme tensorielle ; 

\bigskip

\item[$ii)$] toute famille d'éléments de $k^{1/p}$ qui est topologiquement $p$-libre sur $k$ est topologiquement $p$-libre sur $L$ ; 

\bigskip
\item[$iii)$] toute famille dénombrable d'éléments de $k^{1/p}$ qui est topologiquement $p$-libre sur $k$ est topologiquement $p$-libre sur $L$. 

\end{itemize}
\medskip
De plus, pour que ces propriétés équivalentes soient satisfaites, il est suffisant, et nécessaire si $k^{1/p}$ est topologiquement de type dénombrable sur $k$, qu'il existe une $p$-base topologique de la $k$-algèbre $k^{1/p}$ qui soit topologiquement $p$-libre sur $L$. 

}

\setcounter{cptbis}{0}
\bigskip

{\em Démonstration.} Soit $\mu$ l'application $\sum l_{i}\otimes{x_{i}}\mapsto \sum l_{i}x_{i}$.

\bigskip
\trois{1impli2} \textbf{Prouvons que $i)$ implique $ii)$}. On suppose que $i)$ est vérifiée. Comme $|.|\circ \mu$ est d'après l'hypothèse $i)$ une norme équivalente à la norme tensorielle, $\mu$ est une injection $L$-linéaire admissible de $L\hotimes_{k}k^{1/p}$ dans $L^{1/p}$. 

\medskip
Soit $(a_{i})$ une famille d'éléments de $k^{1/p}$ topologiquement $p$-libre sur $k$. Il découle du~\ref{ANSEP}.\ref{introklibre}.\ref{changeklibre} que $(1\otimes a_{i})$ est une famille d'éléments de $L\hotimes_{k}k^{1/p}$ qui est topologiquement $p$-libre sur $L$ ; en vertu du~\ref{ANSEP}.\ref{introklibre}.\ref{changenorme}, $(a_{i})=(\mu(1\otimes a_{i}))$ constitue une famille d'éléments de $L^{1/p}$ qui est topologiquement $p$-libre sur $L$.

\bigskip
\trois{2impli3} L'implication $ii)\Rightarrow iii)$ est évidente. 

\bigskip
\trois{3impli1} \textbf{Prouvons que $iii)\Rightarrow i)$}. On 
fait l'hypothèse que $iii)$ est vérifiée. 

\medskip
{\em Remarque préliminaire.} Soit $F$ une sous-$k$-extension complète de $k^{1/p}$, topologiquement de type dénombrable sur $k$. Elle possède ({\em cf.} exemple~\ref{ANSEP}.\ref{kiehlpbase}) une $p$-base topologique dénombrable $(a_{i})$ sur $k$. Par l'hypothèse $iii)$, $(a_{i})$ est une famille topologiquement $p$-libre sur $L$. Le lemme~\ref{ANSEP}.\ref{pbasestable} assure alors que la semi-norme $|.|\circ \mu$ de $L\hotimes_{k}F$ est une norme équivalente à la norme tensorielle.

\medskip
{\em Suite de la preuve.} Soit $(y_{n})$ une suite d'éléments de $L\hotimes_{k}k^{1/p}$ telle que $\mu(y_{n})$ tende vers zéro quand $n$ tend vers l'infini. Il existe une sous-$k$-extension $F$ de $k^{1/p}$, topologiquement de type dénombrable et telle que $y_{n}$ appartienne à $L\hotimes_{k}F$ pour tout $n$. D'après la remarque ci-dessus, $|.|\circ \mu$ est une norme de $L\hotimes_{k}F$, équivalente à la norme tensorielle. On en déduit que $y_{n}$ tend vers zéro au sens de la norme tensorielle sur $L\hotimes_{k}F$ ; elle tend {\em a fortiori} vers zéro au sens de la norme tensorielle sur $L\hotimes_{k}k^{1/p}$, et $i)$ est démontré.

\bigskip
\trois{condnecsuf}En ce qui concerne la dernière assertion, elle découle immédiatement du lemme~\ref{ANSEP}.\ref{pbasestable} et de l'exemple~\ref{ANSEP}.\ref{kiehlpbase}.~$\Box$

\bigskip
\deux{fortsep} {\bf Définition.} On dira qu'une extension de corps ultramétriques complets est \textit{analytiquement séparable} s'ils sont de caractéristique nulle ou si les trois propriétés équivalentes de la proposition~\ref{ANSEP}.\ref{normequiv} sont satisfaites.

\bigskip
\deux{exfortsep} {\bf Exemples.} Soit $k$ un corps ultramétrique complet de caractéristique $p>0$. 

\setcounter{cptbis}{0}

\bigskip
\trois{parffortsep} Si $k$ est parfait, toute extension de $k$ est analytiquement séparable, puisque $k^{1/p}$ s'identifie alors à $k$. 

\bigskip
\trois{sousansep} Si $L$ est une extension analytiquement séparable de $k$, et si $F$ est une sous-$k$-extension complète de $L$, alors $F$ est une extension analytiquement séparable de $k$ ; on le voit à l'aide de la propriété $ii)$ de la proposition ~\ref{ANSEP}.\ref{normequiv}.

\bigskip
\trois{compfortsep} Si $L$ est une extension analytiquement séparable de $k$, et si $M$ est une extension analytiquement séparable de $L$, alors $M$ est une extension analytiquement séparable de $k$ ; on le déduit de la propriété $ii)$ de la proposition~\ref{ANSEP}.\ref{normequiv}.

\bigskip
\trois{sepfortsep} Soit $L$ une extension {\em finie} de $k$. La flèche naturelle $L\otimes_{k}k^{1/p}\to L\hotimes_{k}k^{1/p}$ est alors un isomorphisme, et le corps $L$ est par ailleurs séparable sur $k
$ si et seulement si $L\otimes_{k}k^{1/p}$ est un corps ({\em cf.} \cite{mats}, th. 26.2). On va montrer, en se fondant sur ces deux remarques, que $L$ est séparable sur $k$ si et seulement si $L$ est analytiquement séparable sur $k$. 

\begin{itemize}

\itb Si $L$ est séparable sur $k$, alors $L\hotimes_{k}k^{1/p}=L\otimes_{k}k^{1/p}$ est un corps de dimension finie sur $k^ {1/p}$. La semi-norme multiplicative $\sum l_{i}\otimes x_{i}\mapsto |\sum l_{i}x_{i}|^{1/p}$ est de ce fait injective, et est donc une norme. Elle est équivalente à la norme tensorielle : on le vérifie directement si $|k^{*}|=\{1\}$ (auquel cas les deux sont triviales), et sinon cela résulte de l'équivalence des normes en dimension finie. En conséquence, $L$ est analytiquement séparable sur $k$.

\itb Si $L$ est analytiquement séparable sur $k$, la semi-norme multiplicative $\sum l_{i}\otimes x_{i}\mapsto |\sum l_{i}x_{i}|^{1/p}$ est une norme, et la $L$-algèbre $L\hotimes_{k}k^{1/p}=L\otimes_{k}k^{1/p}$ est donc intègre. Comme elle est par ailleurs entière, c'est un corps, et $L$ est dès lors séparable sur $k$. 

\end{itemize}
\bigskip
\deux{introunivmult} {\bf Universelle multiplicativité et séparabilité analytique.} Soit $L$ une extension complète de $k$. Supposons que la norme tensorielle de $L\hotimes_{k}k^{1/p}$ soit multiplicative. La $L$-algèbre $L\hotimes_{k}k^{1/p}$ est dès lors intègre, et par ailleurs entière ; c'est donc un corps. La norme tensorielle est une valeur absolue de ce corps, l'application $\sum l_{i}\otimes x_{i}\mapsto |\sum l_{i}x_{i}|$ en est une autre, majorée par la précédente ; on en déduit que ces deux valeurs absolues coïncident, et $L$ est de ce fait analytiquement séparable sur $k$. 

\medskip
On est notamment dans cette situation si la valeur absolue de $L$ est {\em universellement multiplicative} ({\em peaked} dans la terminologie de Berkovich), c'est-à-dire si la norme tensorielle de $L\hotimes_{k}F$ est multiplicative pour toute extension complète $F$ de $k$.

\setcounter{cptbis}{0}
\medskip
\trois{premunivmult} {\bf Exemple.} Si $\bf r$ est un polyrayon $k$-libre, alors la valeur absolue de $k_{\bf r}$ est, par un calcul explicite immédiat, universellement multiplicative ; par conséquent, $k_{\bf r}$ est une extension analytiquement séparable de $k$.

\medskip
\trois{peakedalg} {\bf Lemme (Berkovich, \cite{brk1}, lemma 5.5.2.)} {\em Soit $\cal A$ une algèbre $k$-affinoïde dont la norme est universellement multiplicative, et soit $x$ l'unique point du bord de Shilov de ${\cal M}({\cal A})$. La valeur absolue de l'extension $\hres(x)$ de $k$ est alors universellement multiplicative.} 

\medskip
{\em Démonstration\footnote{Notre lemme est moins général que celui de Berkovich, mais la preuve de ce dernier reposant sur l'affirmation non justifiée qu'une certaine flèche est une injection isométrique, nous l'avons modifiée pour éviter le recours à cet argument ; et ce, au prix d'un affaiblissement de l'énoncé.}.} Soit $F$ une extension complète de $k$. Par hypothèse, la norme de ${\cal A}\hotimes_{k}F$ est multiplicative. Le bord de Shilov de l'espace affinoïde ${\cal M}({\cal A}\hotimes_kF)$ est donc un singleton $\{y\}$; le point $y$ étant situé au-dessus de $x$, il peut être vu comme appartenant à ${\cal M}(\hres(x)\hotimes_kF)$. Soit $\phi$ le morphisme canonique ${\cal A}\hotimes_kF\to \hres(x)\hotimes_k F$ et soit $f\in {\cal A}\hotimes_{k}F$. On a $$|f(y)|=||f||\geq ||\phi(f)||\geq |\phi(f)(y)|=|f(y)|\; ;$$on en déduit que $||\phi(f)||=||f||$. 

\medskip
Soit $\mathsf A$ la sous-algèbre de $\hres(x)\hotimes_k F$ constituée des éléments de la forme $a(x)^{-1}\phi(f)$, où $a$ est une fonction appartenant à $\cal A$ et inversible en $x$ ({\em i.e.} non nulle) et où $f\in {\cal A}\hotimes_kF$. Donnons-nous un élément $a(x)^{-1}\phi(f)$ de $\mathsf A$. On a $$||a(x)^{-1}\phi(f)||=|a(x)|^{-1}.||\phi(f)||=|a(x)|^{-1}.||f||\;\;\;\;\;\;\;(*)\;,$$ la première égalité provenant du fait que la norme tensorielle sur $\hres(x)\hotimes_kF$ est une norme de $\hres(x)$-algèbre, et la seconde découlant de ce qui précède. 

\medskip
Soient $\mathsf a=a(x)^{-1}\phi(f)$ et $\mathsf b=b(x)^{-1}\phi(g)$ deux éléments de $\mathsf A$. On peut écrire $$\begin{array}{ccc}||\mathsf a\mathsf b||&=&||a(x)^{-1}b(x)^{-1}\phi(f)\phi(g)||\\ &=&||a(x)^{-1}b(x)^{-1}\phi(fg)||\\&=&|a(x)|^{-1}.|b(x)|^{-1}.||fg||\;,\end{array}$$ la dernière égalité se fondant sur $(*)$. Comme la norme $||.||$ est multiplicative, $|a(x)|^{-1}.|b(x)|^{-1}.||fg||=|a(x)|^{-1}.|b(x)|^{-1}.||f||.||g||$ ; en utilisant à nouveau $(*)$, il vient $$|a(x)|^{-1}.|b(x)|^{-1}.||f||.||g||=||a(x)^{-1}\phi(f)||.||b(x)^{-1}\phi(g)||=||\mathsf a||.||\mathsf b||.$$ En conclusion, la restriction de la norme de $\hres(x)\hotimes_kF$ à $\mathsf A$ est multiplicative ; comme $\mathsf A$ est dense dans $\hres(x)\hotimes_kF$, la norme de $\hres(x)\hotimes_kF$ est multiplicative, ce qu'il fallait démontrer.~$\Box$ 

\medskip
\trois{univmult} {\bf Exemples.} \nopagebreak

\begin{itemize}
\itb Soit $\bf r$ un polyrayon, soit $X$ le polydisque correspondant, et soit $\eta$ l'unique point du bord de Shilov de $X$. La norme spectrale de l'algèbre des fonctions analytiques sur $X$ est, par un calcul direct, universellement multiplicative. En vertu du lemme ci-dessus, la valeur absolue de $\hres(\eta)$ est universellement multiplicative, et $\hres(\eta)$ est donc une extension analytiquement séparable de $k$. Notons que si $\bf r$ est $k$-libre, $\hres(\eta)=k_{\bf r}$ ; on retrouve ainsi l'exemple \ref{ANSEP}.\ref{introunivmult}.\ref{premunivmult}. 

\itb Supposons que la valeur absolue de $k$ n'est pas triviale et soit $\cal A$ une algèbre strictement $k$-affinoïde {\em distinguée} (\cite{bgr}, 6.4.3, def. 2). Supposons que la $\red k$-algèbre ${\cal A}\zero/{\cal A}\zeroo$ est géométriquement intègre, et soit $\eta$ l'unique point du bord de Shilov de ${\cal M}({\cal A})$. D'après la proposition 5.2.5 de \cite{brk1}, la norme spectrale de $\cal A$ est universellement multiplicative. En vertu du lemme ci-dessus, la valeur absolue de $\hres(\eta)$ est universellement multiplicative, et $\hres(\eta)$ est donc une extension analytiquement séparable de $k$. 
\end{itemize}

\section{Excellence des algèbres et espaces affinoïdes}\label{EXC}
\setcounter{cpt}{0}

\subsection*{Exhibition d'un ouvert régulier}

\medskip
\deux{disquereg}  {\bf Lemme.} {\em Soit $n$ un entier, soit $V$ un domaine affinoïde de $\Aff^{n,an}_{k}$ et soit $\cal A$ l'algèbre associée à $V$. L'anneau $\cal A$ et l'espace localement annelé $V$ sont réguliers.}

\medskip
{\em Démonstration.} En vertu de la fidèle platitude de $V\to \spec \cal A$, il suffit de montrer l'assertion relative à l'espace localement annelé $V$. Soit $v\in V$ ; on va démontrer que ${\cal O}_{V,v}$ est régulier. Si $L$ est une extension complète de $k$, le morphisme d'espaces localement annelés $V_L\to V$ est fidèlement plat ; il suffit donc d'établir la régularité de ${\cal O}_{V_L,w}$ pour n'importe quel antécédent $w$ de $v$ sur $V_L$ ; ceci permet de se ramener au cas où $v$ est un $k$-point de $V$. Il existe un polydisque fermé $\DD$ qui contient $V$ ; le point rigide $v$ appartient à l'intérieur topologique de $V$ dans $\DD$ et l'on a donc ${\cal O}_{V,v}={\cal O}_{\DD,v}$. Comme $v$ est un $k$-point, il possède une base de voisinages dans $\DD$ formée de polydisques, et ${\cal O}_{\DD,v}$ s'identifie de ce fait à l'anneau des séries convergentes en $n$ variables $\tau_1,\ldots,\tau_n$ dont il est immédiat que $(\tau_1,\ldots,\tau_n)$ constitue un système régulier de paramètres.~$\Box$

\bigskip
\deux{ouvertsreg} {\bf Proposition.}  {\em Soit $\cal X$ un schéma de type fini intègre sur une algèbre $k$-affinoïde $\cal A$ et soit $L$ une extension complète (resp. une extension complète {\em analytiquement séparable}) de $k$. 

\medskip
\begin{itemize}

\item[$a)$] Il existe un ouvert de Zariski non vide ${\cal U}$ de ${\cal X}$ tel que les anneaux locaux de ${\cal U}_{L}\an$ et ${\cal U}_{L}$ soient d'intersection complète (resp. réguliers). 

\medskip
\item[$b)$] Si ${\cal A}$ est strictement $k_{\bf r}$-affinoïde pour un certain polyrayon $k$-libre $\bf r$, il existe deux ouverts de Zariski ${\cal V}$ et $\cal W$ de $\cal X$ qui possèdent les propriétés suivantes : 

\medskip
\begin{itemize}

\item[$b_{0})$] $\emptyset\neq {\cal V}\subset {\cal W}$ ;

\medskip
\item[$b_{1})$] si ${\cal X}=\spec {\cal A}$, alors ${\cal W}={\cal X}$ ;  

\item[$b_{2})$] pour tout $x\in {\cal W}_{L}\an$, l'image réciproque de ${\cal V}_{L}$ sur $\spec\widehat{{\cal O}_{{\cal W}_{L}\an,x}}$ est un schéma qui est d'intersection complète en chacun de ses points (resp. un schéma régulier) ; 

\medskip
\item[$b_{3})$] pour tout ${\bf x}\in {\cal W}_{L}$, l'image réciproque de ${\cal V}_{L}$ sur $\spec\widehat{{\cal O}_{{\cal W}_{L},\bf x}}$ est un schéma qui est d'intersection complète en chacun de ses points (resp. un schéma régulier).

\end{itemize}
\end{itemize} }

\bigskip
{\em Démonstration.} Commençons par quelques remarques.

\medskip
\setcounter{cptbis}{0}
\trois{b2impliqueb3} Soit ${\cal W}$ un ouvert de Zariski de ${\cal X}$, soit $x\in{\cal W}_{L}\an$ et soit $\bf x$ l'image de $x$ sur ${\cal W}_{L}$. Comme ${\cal O}_{{\cal W}_{L}\an,x}$ est une ${\cal O}_{{\cal W}_{L},\bf x}$-algèbre plate, $\widehat{{\cal O}_{{\cal W}_{L}\an,x}}$ est une $\widehat{{\cal O}_{{\cal W}_{L},\bf x}}$-algèbre plate ; en conséquence, $b_{2})\Rightarrow b_{3})$ (\ref{RAP}.\ref{platdirect}). 

\medskip
\trois{bimpliquea} Supposons que $b)$ a été prouvée. Soit $\bf r$ un polyrayon déployant $\cal A$ ; l'assertion $b)$, appliquée au schéma ${\cal X}_{\bf r}$ qui est intègre (lemme \ref{RAP}.\ref{tenseurkrintegre}) et de type fini sur l'algèbre $k$-affinoïde ${\cal A}_{\bf r}$, laquelle est strictement $k_{\bf r}$-affinoïde, fournit en particulier un ouvert non vide ${\cal V}$ de ${\cal X}_{\bf r}$ tel que pour tout $x\in {\cal V}_{L}\an$, l'anneau $\widehat{{\cal O}_{{\cal V}_{L}\an,x}}$ soit d'intersection complète (resp. régulier). Désignons par ${\cal U}$ l'ouvert de ${\cal X}$ tel que $\sigma^{-1}({\cal V}\an)={\cal U}\an$, où $\sigma$ est la section de Shilov de ${\cal X}_{\bf r}\an\to{\cal X}\an$ ({\em cf.} \ref{RAP}.\ref{secshilov}). Soit $y$ appartenant à ${\cal U}_{L}\an$ ; notons $\bf y$ son image sur ${\cal U}_{L}$. Comme $\sigma(y)$ appartient à ${\cal V}_{L}\an$, l'anneau $\widehat{{\cal O}_{{\cal V}_{L}\an,\sigma(y)}}$ est d'intersection complète (resp. régulier). La platitude de ${\cal O}_{{\cal U}_{L},\bf y}\to {\cal O}_{{\cal U}_{L}\an,y}$ et de ${\cal O}_{{\cal U}_{L}\an,y}\to\widehat{{\cal O}_{{\cal V}_{L}\an,\sigma(y)}}$ implique (\ref{RAP}.\ref{platdirect}) que les anneaux locaux de ${\cal U}_{L}\an$ (resp. ${\cal U}_{L})$ en $y$ (resp. $\bf y$) sont d'intersection complète (resp. réguliers), et $a)$ est démontrée. 

\medskip
\trois{onprouveb1b2} D'après le~\ref{EXC}.\ref{ouvertsreg}.\ref{bimpliquea} ci-dessus il suffit, pour prouver la proposition, d'établir l'assertion $b)$ ; l'on peut même, en vertu du~\ref{EXC}.\ref{ouvertsreg}.\ref{b2impliqueb3}, se contenter d'exhiber un couple $({\cal V},{\cal W})$ d'ouverts de ${\cal X}$ qui satisfont $b_{0}),\;b_{1})$ et $b_{2})$. 

\medskip
{\em À partir de maintenant, on suppose donc que ${\cal A}$ est strictement $k_{\bf r}$-affinoïde pour un certain polyrayon $k$-libre $\bf r$.}

\bigskip
 \trois 
{cmregnormal} Soit $\cal Z$ l'adhérence de Zariski de l'image de $\cal X$ sur $\spec \cal A$. La version {\em analytique} du lemme de normalisation de Noether assure l'existence d'une flèche finie et surjective $\cal Z\to \spec k_{\bf r}\{T_{1},\ldots,T_{l}\}$ pour un certain entier $l$. On pose ${\cal D}=\spec k_{\bf r}\{T_{1},\ldots,T_{l}\}$. En appliquant la version {\em algébrique} de ce lemme sur le corps des fonctions de ${\cal D}$, on obtient l'existence d'un entier $n$ et d'un ouvert non vide $\cal W$ de ${\cal X}$ muni d'un morphisme fini et dominant sur un ouvert $\dot{\cal W}$ de $\Aff^{n}_{\cal D}$ ;  si ${\cal X}=\spec \cal A$ l'on a $n=0$ et l'on peut prendre ${\cal W}={\cal X}$ et ${\cal W}'={\cal D}$. La fibre générique de ${\cal W}\to {\cal W}'$ est le spectre d'un corps, et est en conséquence d'intersection complète ; en vertu du corollaire 19.3.8 de \cite{ega44}, il existe un ouvert non vide $\dot{\cal V}$ de $\Aff^{n}_{\cal D}$, inclus dans $\dot{\cal W}$, et tel que ${\cal W}\times_{\Aff^{n}_{\cal D}}\dot{\cal V}\to\dot {\cal V}$ soit un morphisme d'intersection complète. On pose ${\cal V}={\cal W}\times_{\Aff^{n}_{\cal D}}\dot{\cal V}={\cal W}\times_{\dot{\cal W}}\dot{\cal V}$. 

\bigskip
Soit $x$ appartenant à ${\cal W}_{L}\an$ et soit $\dot{x}$ son image sur ${\dot{\cal W}_{L}}\an\subset \Aff^{n,an}_{{\cal D}_{L}}$. On note $\mathsf S_{x}$ (resp. $\mathsf S_{\dot{x}}$) le spectre de l'anneau local de ${\cal W}_{L}\an$ (resp. ${\dot{\cal W}_{L}}\an$) en $x$ (resp. $\dot{x}$) ; on désigne par $\widehat{\mathsf S_{x}}$ et $\widehat{\mathsf S_{\dot{x}}}$ les spectres des complétés correspondants, et par $\mathsf V_{L}$ (resp. $\dot{\mathsf V}_{L}$) l'image réciproque de ${\cal V}_{L}$ sur $\widehat{\mathsf S_{x}}$ (resp. $\widehat{\mathsf S_{\dot{x}}}$). 

On dispose d'un diagramme commutatif $$\xymatrix{&\mathsf V_{L}\ar[rr]\ar[d]&&\widehat{\mathsf S_{x}}\ar[d]\\ &\widehat {\mathsf S_{\dot{x}}}\times_{\dot{\cal W}_{L}}{\cal V}_{L}\ar[rr]\ar[dl]\ar@{.>}[dd]&&\widehat{\mathsf S_{\dot{x}}}\times_{\dot{\cal W}_{L}}{\cal W}_{L}\ar[dl]\ar[dd]\\{\cal V}_{L}\ar[dd]\ar[rr]&&{\cal W}_{L}\ar[dd]&\\ & \dot{\mathsf V}_{L}\ar@{.>}[rr]\ar@{.>}[ld]&& \widehat{\mathsf S_{\dot{x}}}\ar[ld]\\ \dot{\cal V}_{L}\ar[rr]&&\dot{\cal W}_{L}}$$ dans lequel les carrés sont cartésiens, et dont les flèches horizontales sont des immersions ouvertes. Quant aux flèches verticales, celles du parallélépipède inférieur sont par construction des morphismes d'intersection complète, et les deux du haut sont des immersions à la fois ouvertes et fermées : il suffit de le vérifier pour celle de droite, pour laquelle cela découle du lemme~\ref{RAP}.\ref{finicomplet}. 

\medskip
Le schéma $\mathsf S_{\dot{x}}$ étant régulier d'après le lemme~\ref{EXC}.\ref{disquereg}, $\widehat{\mathsf S_{\dot{x}}}$ l'est aussi, et l'ouvert $\dot{\mathsf V}_{L}$ de ce dernier également. On en déduit que $\mathsf V_{L}$ est d'intersection complète en chacun de ses points.

\bigskip
{\em À partir de maintenant, on fait l'hypothèse que $L$ est analytiquement séparable sur $k$ ; le but de ce qui suit est de prouver que l'on peut restreindre les différents ouverts de Zariski construits ci-dessus de sorte que $\mathsf V_L$ soit régulier, ce qui achèvera la démonstration.}

\bigskip
\trois{dembacrzero} {\bf Supposons $k$ de caractéristique nulle.} On peut restreindre $\dot{\cal V}$ de sorte que  ${\cal V}\to \dot {\cal V}$ soit étale. Dans ce cas, $\mathsf V_{L}$ est étale sur le schéma régulier $\dot{\mathsf V}_{L}$, et est de ce fait régulier. 

\bigskip
\trois{fortsepcarp} {\bf Supposons $k$ de caractéristique $p>0$.} Il existe ({\em cf.} \ref{RAP}.\ref{remdenom}) :

\begin{itemize}
 
\itb une sous-$k^{p}$-extension complète $k_\aleph$ de $k$ qui est topologiquement de type dénombrable sur $k^{p}$ ;
\itb deux ouverts ${\cal V}'_\aleph\subset {\cal W}'_\aleph$ de $\Aff^n_{{\cal D}_\aleph}$, où ${\cal D}_\aleph=\spec k_{\aleph,\bf r}\{T_1,\ldots,T_l\}$ ;
\itb un morphisme fini et plat ${\cal W}_\aleph\to{\cal W}'_\aleph$ ; on pose ${\cal V}_\aleph={\cal W}_\aleph\times_{{\cal W}'_\aleph}{\cal V}'_\aleph$ ;
\itb un $k_{\aleph,\bf r}\{T_1,\ldots,T_l\}$-schéma intègre de type fini ${\cal X}_\aleph$, et une immersion ouverte ${\cal W}_\aleph\hookrightarrow {\cal X}_\aleph$ ;
\itb des isomorphismes ${\cal X}_{{\aleph},k}\simeq {\cal X}$, etc. , modulo lesquels ${\cal V}'_{\aleph}\hookrightarrow {\cal W}'_{\aleph}\hookrightarrow \Aff^n_{{\cal D}_{\aleph}}$, ${\cal W}_{\aleph}\to {\cal W}'_{\aleph}$ et ${\cal V}_{\aleph}\hookrightarrow {\cal W}_{\aleph}\hookrightarrow {\cal X}_{\aleph}$ induisent ${\cal V}'\hookrightarrow {\cal W}'\hookrightarrow \Aff^n_{{\cal D}}$, ${\cal W}\to {\cal W}'$ et ${\cal V}\hookrightarrow {\cal W}\hookrightarrow {\cal X}$.
\end{itemize}

\medskip
On notera $S_1,\ldots,S_n$ les fonctions coordonnées sur $\Aff^n_{{\cal D}_{\aleph}}$. 

\medskip
D'après l'exemple \ref{ANSEP}.\ref{kiehlpbase}, $k_{\aleph}$ possède une $p$-base topologique $\bf a$ sur $k^p$. Pour toute partie ${\bf b}$ de $\bf a$, notons $k^{p}[[{\bf b}]]$ le complété du sous-corps de $k_\aleph$ engendré par $k^{p}$ et ${\bf b}$ (remarquons que ${\bf b}$ constitue une $p$-base topologique de $k^{p}[[{\bf b}]]$ sur $k^p$, et que $k_\aleph=k^p[[{\bf a}]]$), et $\Delta_{\bf b}$ le schéma $$\spec k^{p}[[{\bf b}]]_{{\bf r}^{p}}\{T_{1}^{p},\ldots,T_{l}^{p}\}[S_{1}^{p},\ldots,S_{n}^{p}].$$ Par construction, $k_\aleph$ est une extension complète de $k^{p}[[{\bf b}]]$ dont ${\bf a}-{\bf b}$ est une $p$-base topologique, et elle est donc finie si et seulement $\bf b$ est cofinie, c'est-à-dire si ${\bf a}-{\bf b}$ est finie ; dans ce cas $\Aff^{n}_{{\cal D}_{\aleph}}$ est un $\Delta_{\bf b}$-schéma fini, radiciel et plat.  

\medskip
\begin{itemize}
\itb {\bf Une première remarque.} L'intersection $\bigcap\limits_{{\bf b}\;\tiny \mbox{cofinie}} \kappa(\Delta_{\bf b})$ est le corps des fractions de $k^{p}_{{\bf r}^{p}}\{T_{1}^{p},\ldots,T_{l}^{p}\}[S_{1}^{p},\ldots,S_{n}^{p}]$, qui coïncide avec $\kappa(\Aff^{n}_{\cal D})^{p}$. 

\medskip
En effet, soit $f\in \bigcap\limits_{{\bf b}\;\tiny \mbox{cofinie}} \kappa(\Delta_{\bf b})$. Écrivons-la sous la forme $g/h$, où $g$ et $h$ appartiennent tous deux à $k_{\aleph,{\bf r}^{p}}\{T_{1}^{p},\ldots,T_{l}^{p}\}[S_{1}^{p},\ldots,S_{n}^{p}]$ (c'est-à-dire à l'anneau des fonctions de $\Delta_{\bf a}$), et où $h\neq 0$. Quitte à remplacer $h$ par $h^p$ et $g$ par $gh^{p-1}$ on peut supposer que $h\in k^{p}_{{\bf r}^{p}}\{T_{1}^{p},\ldots,T_{l}^{p}\}[S_{1}^{p},\ldots,S_{n}^{p}]$.  On va montrer que $g$ appartient aussi à $k^{p}_{{\bf r}^{p}}\{T_{1}^{p},\ldots,T_{l}^{p}\}[S_{1}^{p},\ldots,S_{n}^{p}]$, ce qui permettra de conclure. 

\medskip
Soit $\bf b$ une partie cofinie de $\bf a$. 
La fonction $g=hf$ appartient à $\kappa(\Delta_{\bf b})$ ; on va prouver qu'elle est définie sur $\Delta_{\bf b}$ tout entier. Soit $\xi$ un point de codimension $1$ du schéma noethérien $\Delta_{\bf b}$ et soit $\varpi$ une uniformisante de l'anneau ${\cal O}_{\Delta_{\bf b},\xi}$, qui est de valuation discrète puisque $\Delta_{\bf b}$ est normal en vertu du lemme \ref{EXC}.\ref{disquereg}. Si $g$ est nulle, elle est évidemment définie au voisinage de $\xi$ ; sinon, elle s'écrit $u\varpi^m$ pour un certain entier relatif $m$ et une certaine fonction $u$ définie et inversible au voisinage de $\xi$. Supposons $m<0$ et soit $\omega$ un point de $\Delta_{\bf a}$ situé au-dessus de $\xi$. L'égalité $g\varpi^{-m}=u$ met en jeu des fonctions définies sur $\Delta_{\bf a}$ tout entier. Évaluée en $\omega$, elle conduit à la contradiction $0=u(\xi)$. Par conséquent, $m\geq 0$ et $g\in {\cal O}_{\Delta_{\bf b},\xi}$. Ceci vaut pour tout point $\xi$ de codimension $1$ du schéma noethérien normal $\Delta_{\bf b}$ ; dès lors, $g\in {\cal O}_{\Delta_{\bf b}}(\Delta_{\bf b})$.

\medskip
En conclusion, $$g\in\bigcap\limits_{{\bf b}\;\tiny \mbox{cofinie}}{\cal O}_{\Delta_{\bf b}}(\Delta_{\bf b})=\bigcap\limits_{{\bf b}\;\tiny \mbox{cofinie}}k^{p}[[{\bf b}]]_{{\bf r}^{p}}\{T_{1}^{p},\ldots,T_{l}^{p}\}[S_{1}^{p},\ldots,S_{n}^{p}]$$ $$\;\;\;\;\;\;\;=k^{p}_{{\bf r}^{p}}\{T_{1}^{p},\ldots,T_{l}^{p}\}[S_{1}^{p},\ldots,S_{n}^{p}],$$ ce qu'il fallait démontrer. 

\itb {\bf Une seconde remarque.} Si $(\alpha_{i})$ est une base de $\kappa({\cal X}_{\aleph})$ sur $\kappa(\Aff^{n}_{{\cal D}_{\aleph}})$, c'est encore une base de $\kappa({\cal X})$ sur $\kappa(\Aff^{n}_{\cal D})$, puisque $\kappa({\cal X})$ s'identifie à $$\kappa({\cal X}_{\aleph})\otimes_{\kappa(\Aff^{n}_{{\cal D}_{\aleph}})}\kappa(\Aff^{n}_{\cal D}).$$ En conséquence, $(\alpha_{i}^{p})$ est une famille  $\kappa(\Aff^{n}_{\cal D})^{p}$-libre.
\end{itemize}

\medskip
\noindent
On est dans les conditions d'application d'un théorème de Kiehl (\cite{kie}, Satz 2.4). Il en découle l'existence d'une partie cofinie $\bf b$ de $\bf a$ telle que les rangs (calculés respectivement sur $\kappa({\cal X}_{\aleph})$ et $\kappa(\Aff^{n}_{{\cal D}_{\aleph}}$)) des espaces vectoriels de dimension finie $$\Omega^1(\kappa({\cal X}_{\aleph})/ \kappa(\Delta_{\bf b}))\;\mbox{et}\;\Omega^1(\kappa(\Aff^{n}_{{\cal D}_{\aleph}})/ \kappa(\Delta_{\bf b}))$$ coïncident (cette assertion est en tout point analogue au corollaire 3.8 de \cite{kie}, dont nous avons décalqué la preuve {\em mutatis mutandis.}). 

\bigskip
On peut donc restreindre $\dot{\cal V}_{\aleph}$ (ainsi bien entendu que ${\cal V}_{\aleph}$, $\dot{\cal V}$ et ${\cal V}$) de sorte que la condition suivante soit vérifiée, en notant $\ddot{\cal V}_{\bf b}$ l'ouvert de $\Delta_{\bf b}$ égal à l'image de $\dot{\cal V}_{\aleph}$ : {\em les faisceaux cohérents $\Omega^1_{\dot{\cal V}_{\aleph}/\ddot{\cal V}_{\bf b}}$ et $\Omega^1_{{\cal V}_{\aleph}/\ddot{\cal V}_{\bf b}}$ sont libres, et tous deux de même rang.} On désigne  par $\ddot{\cal W}_{\bf b}$ l'ouvert de $\Delta_{\bf b}$ égal à l'image de $\dot{\cal W}_{\aleph}$. 

\medskip
\trois{rembfa} {\bf Remarque à propos des notations.} En ce qui concerne les symboles choisis pour désigner les schémas, la présence d'un $\aleph$ (resp. d'un $\bf b$) en indice vise à rappeler que le schéma concerné est de type fini sur l'algèbre $k_\aleph$-affinoïde $k_{\aleph,\bf r}\{T_1,\ldots,T_l\}$ (resp. sur l'algèbre $k^p[[{\bf b}]]$-affinoïde $k^{p}[[{\bf b}]]_{{\bf r}^p}\{T_1^p,\ldots,T_l^p\}$).

\medskip
\trois{introcorpsf} Afin de pouvoir utiliser l'exemple~\ref{ANSEP}.\ref{kiehlpbase}, nous allons travailler avec la famille des sous-corps complets de $L$ qui contiennent $L^{p}$ et $k_\aleph$ et qui sont topologiquement de type dénombrable sur $L^{p}$ ; on qualifiera d'{\em admissible} un tel sous-corps.

\medskip
{\bf C'est dans la preuve du lemme ci-dessous, et uniquement là, qu'est utilisée l'hypothèse de séparabilité analytique de l'extension $L/k$.} 

\medskip
\trois{fvientdek} {\bf Lemme.} {\em Soit $F$ un sous-corps admissible de $L$. Il existe une extension complète $F^{\sharp}$ de $k^{p}[[{\bf b}]]$ telle que $F$ s'identifie à $k_\aleph\otimes_{k^{p}[[{\bf b}]]}F^{\sharp}$.}

\medskip
{\em Démonstration.} Comme $L$ est une extension analytiquement séparable de $k$, la famille $\bf a$ d'éléments de $L$ est topologiquement $p$-libre sur le corps $L^{p}$ (la définition de la séparabilité analytique est formulée en termes des extensions $k\hookrightarrow k^{1/p}$ et $L\hookrightarrow L^{1/p}$, et non $k^p\hookrightarrow k$ et $L^p\hookrightarrow L$, mais le~\ref{ANSEP}.\ref{introklibre}.\ref{expoklibre} assure l'équivalence des deux points de vue) ; en conséquence, $\bf a$ est une $p$-base topologique sur $L^p$ du sous-corps complet $\widehat{k_\aleph L^p}$ de $L$ engendré par $L^p$ et $k_\aleph$. Puisque $F$ est topologiquement de type dénombrable sur $L^{p}$, il possède, toujours grâce à l'exemple~\ref{ANSEP}.\ref{kiehlpbase}, une $p$-base topologique $\bf c$ sur $\widehat{k_\aleph L^p}$ ; par construction, ${\bf a}\coprod {\bf c}$ est une $p$-base topologique de $F$ sur le corps $L^{p}$. Soit $F^{\sharp}$ le complété du sous-corps de $L$ engendré par $L^{p}$ et ${\bf b}\coprod {\bf c}$. Alors ${\bf a}-{\bf b}$ est une $p$-base finie de $F$ sur $F^{\sharp}$, et $F$ s'identifie en conséquence à $k_\aleph\otimes_{k^{p}[[{\bf b}]]}F^{\sharp}$.~$\Box$

\medskip
\noindent
{\em On fixe un sous-corps admissible $F$ de $L$, et on choisit $F^{\sharp}$ comme dans le lemme ci-dessus.} 

\medskip
\trois{commentalephb} {\bf Encore quelques remarques à propos des notations.} Si $\Xi_\aleph$ (resp. $\Xi_{\bf b}$) est un schéma de type fini sur l'algèbre $k_\aleph$-affinoïde $k_{\aleph,\bf r}\{T_1,\ldots,T_l\}$ (resp. sur l'algèbre $k^p[[{\bf b}]]$-affinoïde $k^{p}[[{\bf b}]]_{{\bf r}^p}\{T_1^p,\ldots,T_l^p\}$), alors conformément à nos conventions générales $\Xi_{\aleph,F}$ (resp. $\Xi_{{\bf b},F^\sharp}$) désignera $$\Xi_{\aleph}\otimes_{k_{\aleph,\bf r}\{T_1,\ldots,T_l\}}F_{\bf r}\{T_1,\ldots,T_l\}\;\mbox{(resp. }\Xi_{\bf b}\otimes_{k^{p}[[{\bf b}]]_{{\bf r}^p}\{T_1^p,\ldots,T_l^p\}}F^\sharp_{{\bf r}^p}\{T_1^p,\ldots,T_l^p\}\;).$$

\medskip
\trois{dfdeltab} L'on peut écrire : $$ \hskip-20pt k_{\aleph,\bf r}\{T_1,\ldots,T_l\}[S_1,\ldots,S_n]\otimes_{k^{p}[[{\bf b}]]_{{\bf r}^{p}}\{T_{1}^{p},\ldots,T_{l}^{p}\}[S_{1}^{p},\ldots,S_{n}^{p}] } F^\sharp_{{\bf r}^{p}}\{T_{1}^{p},\ldots,T_{l}^{p}\}[S_{1}^{p},\ldots,S_{n}^{p}]$$ $$\simeq k_{\aleph,\bf r}\{T_1,\ldots,T_l\}[S_1,\ldots,S_n]\otimes_{k^{p}[[{\bf b}]]_{{\bf r}^{p}}\{T_{1}^{p},\ldots,T_{l}^{p}\}}F^\sharp_{{\bf r}^{p}}\{T_{1}^{p},\ldots,T_{l}^{p}\}$$ $$\simeq\left(k_{\aleph,\bf r}\{T_1,\ldots,T_l\}\otimes_{k^{p}[[{\bf b}]]_{{\bf r}^{p}}\{T_{1}^{p},\ldots,T_{l}^{p}\}}F^\sharp_{{\bf r}^{p}}\{T_{1}^{p},\ldots,T_{l}^{p}\}\right)[S_1,\ldots,S_n]$$ $$\simeq\left(k_{\aleph,\bf r}\{T_1,\ldots,T_l\}\hotimes_{k^{p}[[{\bf b}]]_{{\bf r}^{p}}\{T_{1}^{p},\ldots,T_{l}^{p}\}}F^\sharp_{{\bf r}^{p}}\{T_{1}^{p},\ldots,T_{l}^{p}\}\right)[S_1,\ldots,S_n]$$ (puisque $k_{\aleph,\bf r}\{T_1,\ldots,T_l\}$ est une $k^{p}[[{\bf b}]]_{{\bf r}^{p}}\{T_{1}^{p},\ldots,T_{l}^{p}\}$-algèbre de Banach {\em finie}) $$\simeq\left(k_{\aleph,\bf r}\{T_1,\ldots,T_l\}\hotimes_{k^{p}[[{\bf b}]]}F^\sharp\right)[S_1,\ldots,S_n]$$ $$\simeq F_{\bf r}\{T_1,\ldots,T_l\}[S_1,\ldots,S_n]$$ en vertu de la définition de $F^\sharp$. La flèche naturelle $\Aff^n_{{\cal D}_{{\aleph},F}}\to \Aff^n_{{\cal D}_{\aleph}}\times_{\Delta_{\bf b}}{\Delta_{{\bf b},F^\sharp}}$ est donc un isomorphisme. 

\medskip
\trois{notprekiehl} {\bf Quelques notations.} 

\begin{itemize}

\itb on appelle $x_{{\aleph},F}$ (resp. $\dot{x}_{{\aleph},F}$, resp. $\ddot{x}_{{\bf b},F^{\sharp}}$) l'image de $x$ sur ${\cal W}_{{\aleph},F}\an$ (resp. ${\dot{\cal W}_{{\aleph},F}}\an$, resp. ${\ddot{\cal W}_{{\bf b},F^{\sharp}}}\an$), on note $\mathsf S_{x_{{\aleph}, F}}$, $\mathsf S_{\dot{x}_{{\aleph},F}}$ et $\mathsf S_{\ddot{x}_{{\bf b},F^{\sharp}}}$ les spectres des anneaux locaux correspondants, et $\widehat{\mathsf S_{x_{{\aleph},F}}}$, $\widehat{\mathsf S_{\dot{x}_{{\aleph},F}}}$ et $\widehat{\mathsf S_{\ddot{x}_{{\bf b},F^{\sharp}}} }$ ceux de leurs complétés ;
\itb on désigne par $\mathsf V_{{\aleph},F}$ (resp. $ \dot{\mathsf V}_{{\bf b},F}$, resp. $\ddot{\mathsf V}_{{\bf b},F^{\sharp}}$)  l'image réciproque de ${\cal V}_{{\aleph},F}$ (resp. $\dot{\cal V}_{{\aleph},F}$, resp. $\ddot{\cal V}_{{\bf b},F^{\sharp}}$) sur $\widehat{\mathsf S_{x_{{\aleph},F}}}$ (resp. $\widehat{\mathsf S_{\dot{x}_{{\aleph},F}}}$, resp. $\widehat{\mathsf S_{\ddot{x}_{{\bf b},F^{\sharp}}} }$). 

\end{itemize}

\medskip
\trois{criterekiehl} Dans le diagramme commutatif $$\xymatrix{{\cal V}_{{\aleph},F}\ar[r]\ar[d]&{\cal V}_{\aleph}\ar[d]\\ \dot{\cal V}_{{\aleph},F}\ar[r]\ar[d]&\dot{\cal V}_{\aleph}\ar[d]\\ \ddot{\cal V}_{{\bf b},F^{\sharp}}\ar[r]&\ddot{\cal V}_{\bf b}}$$ les deux carrés sont cartésiens : c'est évident pour celui du haut, et pour celui du bas, cela résulte du fait que $\Aff^n_{{\cal D}_{{\aleph},F}}\to \Aff^n_{{\cal D}_{{\aleph}}}\times_{\Delta_{\bf b}}{\Delta_{{\bf b},F^\sharp}}$ est un isomorphisme d'après le \ref{EXC}.\ref{ouvertsreg}.\ref{dfdeltab}.

\medskip
Comme $\Omega^1_{\dot{\cal V}_{\aleph}/\ddot{\cal V}_{\bf b}}$ et $\Omega^1_{{\cal V}_{\aleph}/\ddot{\cal V}_{\bf b}}$ sont tous deux libres de même rang, $\Omega^1_{{\cal V}_{{\aleph},F}/\ddot{\cal V}_{{\bf b},F^\sharp}}$ et $\Omega^1_{\dot{\cal V}_{{\aleph},F}/\ddot{\cal V}_{{\bf b},F^{\sharp}}}$ sont tous deux libres de même rang. 

\medskip
Dans le diagramme commutatif $$\xymatrix{ &\mathsf V_{{\aleph},F}\ar[rr]\ar[d]&&\widehat{\mathsf S_{x_{{\aleph},F}}}\ar[d]\\ &\widehat {\mathsf S_{\dot{x}_{{\aleph},F}}}\times_{\dot{\cal W}_{{\aleph},F}}{\cal V}_{{\aleph},F}\ar[rr]\ar[dl]\ar@{.>}[dd]&&\widehat{\mathsf S_{\dot{x_{{\aleph},F}}}}\times_{\dot{\cal W}_{{\aleph},F}}{\cal W}_{{\aleph},F}\ar[dl]\ar[dd]\\{\cal V}_{{\aleph},F}\ar[dd]\ar[rr]&&{\cal W}_{{\aleph},F}\ar[dd]&\\ & \dot{\mathsf V}_{{\aleph},F}\ar@{.>}[rr]\ar@{.>}[ld]\ar@{.>}[dd]& &\widehat{\mathsf S_{\dot{x}_{{\aleph},F}}}\ar[ld]\ar[dd]\\ \dot{\cal V}_{{\aleph},F}\ar[rr]\ar[dd]&&\dot{\cal W}_{{\aleph},F}\ar[dd]&\\ &\ddot{\mathsf V}_{{\bf b},F^{\sharp}}\ar@{.>}[rr]\ar@{.>}[dl]&&\widehat{\mathsf S_{\ddot{x}_{{\bf b},F^{\sharp}}}} \ar[dl]\\ \ddot{\cal V}_{{\bf b},F^{\sharp}}\ar[rr]&&\ddot{\cal W}_{{\bf b},F^{\sharp}} }$$ les carrés sont cartésiens : en ce qui concerne les faces antérieures des deux parallélépipèdes, cela provient du faut que $\Aff^n_{{\cal D}_{{\aleph},F}}\to \Aff^n_{{\cal D}_{\aleph}}\times_{\Delta_{\bf b}}{\Delta_{{\bf b},F^\sharp}}$ est un isomorphisme d'après le \ref{EXC}.\ref{ouvertsreg}.\ref{dfdeltab} ; en ce qui concerne la face de droite du parallélépipède inférieur, cela résulte du lemme~\ref{RAP}.\ref{finicomplet} et du fait que $\dot{x}_{{\aleph},F}$ est l'unique antécédent de $\ddot{x}_{{\bf b},F^{\sharp}}$ par la flèche radicielle ${\dot{\cal W}_{{\aleph},F}}\an\to {\ddot{\cal W}_{{\bf b},F^{\sharp}}}\an$ ; c'est évident pour les autres carrés. Par ailleurs, les flèches allant de gauche à droite sont des immersions ouvertes ; quant aux flèches verticales, elles sont finies et plates, les deux flèches verticales du haut étant plus précisément des immersions à la fois ouvertes et fermées en vertu, là encore, du lemme~\ref{RAP}.\ref{finicomplet}.  

\medskip
On en déduit que $\Omega^1_{\mathsf V_{{\aleph},F}/\ddot{\mathsf V}_{{\bf b},F^{\sharp}}}$ et $\Omega^1_{\dot{\mathsf V}_{{\aleph},F}/\ddot{\mathsf V}_{{\bf b},F^{\sharp}}}$ sont tous deux libres de même rang. Par ailleurs, l'espace localement annelé ${\dot{\cal W}_{{\aleph},F}}\an$ est régulier d'après le lemme~\ref{EXC}.\ref{disquereg} ; en conséquence, le schéma $\mathsf S_{\dot{x}_{{\aleph},F}}$, son complété $\widehat{\mathsf S_{\dot{x}_{{\aleph},F}}}$ et l'ouvert $\dot{\mathsf V}_{{\aleph},F}$ de ce dernier sont réguliers.

\medskip
Récapitulons : $\mathsf V_{{\aleph},F}$ est fini et plat sur le schéma régulier $\dot{\mathsf V}_{{\aleph},F}$, ce dernier est également fini et plat, et {\em a fortiori} de type fini, sur $\ddot{\mathsf V}_{{\bf b},F^{\sharp}}$, et $\Omega^1_{{\mathsf V}_{{\aleph},F}/\ddot{\mathsf V}_{{\bf b},F^{\sharp}}}$ et $\Omega^1_{\dot{\mathsf V}_{{\aleph},F}/\ddot{\mathsf V}_{{\bf b},F^{\sharp}}}$ sont tous deux libres de même rang. 

\medskip
On en déduit, à l'aide d'un critère de Kiehl (\cite{kie}, Satz 2.2 ; pour une preuve, {\em cf.} \cite{comprig}, th. 1.1.1) que $\mathsf V_{{\aleph},F}$ est régulier. {\em Ceci vaut pour tout sous-corps admissible $F$ de $L$.} Les morphismes de transition du système projectif des $\mathsf V_{{\aleph},F}$ sont plats ; en vertu de la proposition 5.13.7 de \cite{ega42}, la régularité de $\mathsf V_{L}$ découle du lemme ci-dessous.~$\Box$

\medskip
\trois{vllimvf} {\bf Lemme.} {\em Le schéma $\mathsf V_{L}$ s'identifie à la limite projective des $\mathsf V_{{\aleph},F}$, où $F$ parcourt la famille des sous-corps admissibles de $L$.} 

\medskip
{\em Démonstration.} Comme $\mathsf V_{L}$ est pour tout sous-corps admissible $F$ de $L$ l'image réciproque de $\mathsf V_{{\aleph},F}$ par la flèche $\widehat{\mathsf S_{x}}\to \widehat{\mathsf S_{x_{{\aleph},F}}}$, il suffit de montrer que $\widehat{\mathsf S_{x}}$ s'identifie à la limite projective des $\widehat{\mathsf S_{x_{{\aleph},F}}}$ ou, ce qui revient au même, que $\widehat{{\cal O}_{{\cal V}_{L}\an,x}}$ s'identifie à la limite inductive des $\widehat{{\cal O}_{{\cal V}_{{\aleph},F}\an,x_{{\aleph},F}}}$. Pour tout $F$, la flèche $\widehat{{\cal O}_{{\cal V}_{{\aleph},F}\an,x_{{\aleph},F}}}\to \widehat{{\cal O}_{{\cal V}_{L}\an,x}}$ est fidèlement plate, et par conséquent injective ; il y a donc simplement à s'assurer que tout élément de $\widehat{{\cal O}_{{\cal V}_{L}\an,x}}$ provient de $\widehat{{\cal O}_{{\cal V}_{{\aleph},F}\an,x_{{\aleph},F}}}$ pour un certain $F$. La question est locale sur le schéma ${\cal V}_{\aleph}$, on peut donc le supposer affine. 

\medskip
Il existe un voisinage affinoïde $\Xi$ de $x$ dans ${\cal V}_{L}\an$ et une famille finie $(f_{1},\ldots,f_{d})$ de fonctions sur $\Xi$ qui engendrent l'idéal maximal de l'anneau local noethérien ${\cal O}_{{\cal V}_{L}\an,x}$. Par la {\em construction explicite} de l'analytifié d'un schéma affine donnée par exemple au paragraphe 1.4 de {\em cf.} \cite{ducbourb}, on peut supposer que $\Xi$ est décrit par un nombre fini d'inégalités larges portant sur des normes de fonctions du {\em schéma} ${\cal V}_{L}$, et que les inégalités {\em strictes} correspondantes sont satisfaites en $x$. Comme ${\cal V}_{L}$ est la limite projective des ${\cal V}_{{\aleph},F}$, il existe un sous-corps admissible $F_{0}$ de $L$ tel que les fonctions rationnelles évoquées proviennent de ${\cal O}_{{\cal V}_{{\aleph},F_{0}}}({\cal V}_{{\aleph},F_{0}})$ ; dès lors, $\Xi$ est l'image réciproque d'un voisinage affinoïde $\Xi_0$ de $x_{{\aleph},F_{0}}$ dans ${\cal V}_{{\aleph},F_{0}}\an$. Quitte à agrandir $F_{0}$, on peut faire l'hypothèse que $f_{i}$ provient de ${\cal O}_{\Xi_0}(\Xi_0)$ pour tout $i$ ; en particulier, $f_{i}$ provient de ${\cal O}_{{\cal V}_{{\aleph},F_{0}}\an,x_{{\aleph},F_{0}}}$ pour tout $i$. 

\medskip
Donnons-nous un élément $g$ de $\widehat{{\cal O}_{{\cal V}_{L}\an,x}}$. Soit $n\in \NN^{*}$. Choisissons dans ${\cal O}_{{\cal V}_{L}\an,x}$ un élément $g_{n}$ congru à $g$ modulo $(f_{1},\ldots,f_{d})^{n}$. Par le même raisonnement que précédemment, il existe un sous-corps admissible $F_{n}$ de $L$ tel que $g_{n}$ provienne de ${\cal O}_{{\cal V}_{{\aleph},F_{n}}\an,x_{{\aleph},F_{n}}}$. 

\medskip
Soit $F$ le complété du sous-corps de $L$ engendré par les $F_{n}$, où $n$ parcourt $\NN$ ; c'est un sous-corps admissible de $L$. Pour tout $n$ non nul, $g_{n}$ provient de ${\cal O}_{{\cal V}_{{\aleph},F}\an,x_{{\aleph},F}}$ ; par ailleurs, le corps $F$ contient $F_ {0}$, d'où il découle que $f_{i}$ provient de ${\cal O}_{{\cal V}_{{\aleph},F}\an,x_{{\aleph},F}}$ quel que soit $i$. 

\medskip
L'anneau ${\cal O}_{{\cal V}_{L}\an,x}$ est fidèlement plat sur ${\cal O}_{{\cal V}_{{\aleph},F}\an,x_{{\aleph},F}}$. Pour $n$ et $i$ convenables, définissons les éléments $g'_{n}$ et $f'_{i}$ de ${\cal O}_{{\cal V}_{{\aleph},F}\an,x_{{\aleph},F}}$ comme les antécédents de $g_n$ et $f_i$. La flèche ${\cal O}_{{\cal V}_{{\aleph},F}\an,x_{{\aleph},F}}/(f'_{1},\ldots,f'_{d})\to {\cal O}_{{\cal V}_{L}\an,x}/(f_{1},\ldots,f_{d})$ est un morphisme d'anneaux fidèlement plat dont le but est un corps ; sa source est donc également un corps ; autrement dit, les $f'_{i}$ engendrent l'idéal maximal de ${\cal O}_{{\cal V}_{{\aleph},F}\an,x_{{\aleph},F}}$. 

\medskip
Pour tout entier $n>0$, les images de $g_{n}$ et $g_{n+1}$ dans ${\cal O}_{{\cal V}_{L}\an,x}/(f_{1},\ldots,f_{d})^{n}$ coïncident. Par fidèle platitude, ${\cal O}_{{\cal V}_{{\aleph},F}\an,x_{{\aleph},F}}/(f'_{1},\ldots,f'_{d})^{n}\to{\cal O}_{{\cal V}_{L}\an,x}/(f_{1},\ldots,f_{d})^{n}$ est injective ; de ce fait, les images de $g'_{n}$ et $g'_{n+1}$ dans ${\cal O}_{{\cal V}_{{\aleph},F}\an,x_{{\aleph},F}}/(f'_{1},\ldots,f'_{d})^{n}$ coïncident.

\medskip
La famille $(g'_{n}\;\mbox{mod}\;(f'_{1},\ldots,f'_{d})^{n}\;)_{n}$ définit ainsi un élément de $\widehat{{\cal O}_{{\cal V}_{{\aleph},F}\an,x_{{\aleph},F}}}$ dont l'image dans  $\widehat{{\cal O}_{{\cal V}_{L}\an,x}}$ est par construction égale à $g$. Ceci achève la preuve du lemme.~$\Box$ 

\subsection*{Vérification des propriétés constitutives de l'excellence} 

\medskip
\deux{geomregul} {\bf Proposition.} {\em Soit $\cal A$ une algèbre $k$-affinoïde ; on pose ${\cal X}=\spec \cal A$ et $X={\cal M}({\cal A})$. Soit $x\in X$ et soit $\bf x$ son image sur $\cal X$. Les morphismes $$\spec\widehat{{\cal O}_{X,x}}\to \spec {\cal O}_{X,x},\;\spec\widehat{{\cal O}_{{\cal X},\bf x}}\to \spec {\cal O}_{{\cal X},\bf x}\;\mbox{et}\;\spec {\cal O}_{X,x}\to \spec {\cal O}_{{\cal X},\bf x}$$ sont géométriquement réguliers.}

\medskip
{\em Démonstration.} On procède en plusieurs temps ; tout au long de la preuve, nous utiliserons sans les rappeler explicitement à chaque occurence les résultats du \ref{RAP}.\ref{platdirect} relatifs au comportement de la régularité vis-à-vis d'un morphisme plat. 

\medskip
Les flèches étudiées étant plates, il reste à montrer les assertions suivantes :

\medskip
$i)$ si $\got{p}$ est un idéal premier de ${\cal O}_{X,x}$ et si $F$ est une extension finie radicielle de $\kappa(\got{p})$, l'anneau $\widehat{{\cal O}_{X,x}}\otimes_{{\cal O}_{X,x}}F$ est régulier ;

\medskip
$ii)$ si $\got{p}$ est un idéal premier de ${\cal O}_{{\cal X},\bf x}$ et si $F$ est une extension finie radicielle de $\kappa(\got{p})$, l'anneau $\widehat{{\cal O}_{{\cal X},\bf x}}\otimes_{{\cal O}_{{\cal X},\bf x}}F$ est régulier ;

\medskip
$iii)$ si $\got{p}$ est un idéal premier de ${\cal O}_{{\cal X},\bf x}$ et si $F$ est une extension finie radicielle de $\kappa(\got{p})$, l'anneau ${\cal O}_{X,x}\otimes_{{\cal O}_{{\cal X},\bf x}}F$ est régulier.

\medskip
\setcounter{cptbis}{0}
\trois{geomregulgen} {\bf Simplification des assertions à établir.} On va tout d'abord expliquer pourquoi l'on peut, pour chacune d'elles, supposer que $\cal A$ est intègre, que $\got{p}=0$ et que $F=\kappa(\got{p})$.  
\medskip
\begin{itemize}

\itb {\em L'assertion $i)$.} On peut supposer, quitte à restreindre $X$, que $\got{p}$ provient d'un idéal $\got{q}$ de $\cal A$ puis, en quotientant $\cal A$ par $\got{q}$, que $\got{p}=0$. Soit $\mathsf B$ une sous ${\cal O}_{X,x}$-algèbre finie et radicielle de $F$ de corps des fractions égal à $F$. On peut à nouveau restreindre $X$ de sorte que $\mathsf B$ soit de la forme ${\cal B}\otimes_{\cal A}{\cal O}_{X,x}$, où $\cal B$ est une $\cal A$-algèbre finie et radicielle. Soit $y$ l'unique antécédent de $x$ sur $Y={\cal M}({\cal B})$ ; l'anneau local ${\cal O}_{Y,y}$ s'identifie à $\mathsf B$ ({\em cf.} lemme \ref{RAP}.\ref{finicomplet}) et l'on dispose d'un isomorphisme $\widehat{{\cal O}_{Y,y}}\otimes_{{\cal O}_{Y,y}}\mbox{Frac}\;{\cal O}_{Y,y}\simeq \widehat{{\cal O}_{X,x}}\otimes_{{\cal O}_{X,x}}F.$ En remplaçant ${\cal A}$ par ${\cal B}$ et $x$ par $y$, on se ramène au cas où $F=\mbox{Frac}\;{\cal O}_{X,x}$. Comme ${\cal O}_{X,x}$ est intègre, le lemme~\ref{RAP}.\ref{localintegre} assure que $x$ n'est situé que sur une composante irréductible de ${\cal M}({\cal A})$. En substituant à ${\cal A}$ son quotient par l'idéal premier correspondant à la composante en question ({\em cf.}  \ref{RAP}.\ref{quotoxx}), on se ramène finalement au cas où $\cal A$ est intègre, où $\got{p}=0$, et où $F=\mbox{Frac}\;{\cal O}_{X,x}$.

 \itb {\em Les assertions $ii)$ et $iii)$.} Quitte à quotienter $\cal A$ par son idéal premier qui correspond à $\got{p}$, on peut supposer que $\cal A$ est intègre et que $\got{p}=0$. Le corps $F$ est maintenant une extension finie radicielle du corps des fractions de $\cal A$. Il existe une sous-$\cal A$-algèbre finie et radicielle $\cal B$ de $F$ dont $F$ est le corps des fractions. Posons $Y={\cal M}({\cal B})$ et ${\cal Y}=\spec \cal B$. Soit $\bf y$ (resp. $y$) l'unique antécédent de $\bf x$ (resp. $x$) sur ${\cal Y}$ (resp. sur $Y$). On dispose d'isomorphismes canoniques $$\widehat{{\cal O}_{{\cal Y},\bf y}}\simeq \widehat{{\cal O}_{{\cal X},\bf x}}\otimes_{{\cal O}_{{\cal X},\bf x}}{{\cal O}_{{\cal Y},\bf y}}\;\mbox{et}\;{\cal O}_{Y,y}\simeq{\cal O}_{X,x}\otimes_{{\cal O}_{{\cal X},\bf x}}{{\cal O}_{{\cal Y},\bf y}}$$ (pour le second, {\em cf.} lemme \ref{RAP}.\ref{finicomplet}). On peut donc écrire $$\widehat{{\cal O}_{{\cal X},\bf x}}\otimes_{{\cal O}_{{\cal X},\bf x}}F\simeq \widehat{{\cal O}_{{\cal Y},\bf y}}\otimes_{{\cal O}_{{\cal Y},\bf y}}\mbox{Frac}\:{\cal O}_{{\cal Y},\bf y}$$ $$\mbox{et}\;{\cal O}_{X,x}\otimes_{{\cal O}_{{\cal X},\bf x}}F\simeq {\cal O}_{Y,y}\otimes_{{\cal O}_{{\cal Y},\bf y}}\mbox{Frac}\;{\cal O}_{{\cal Y},\bf y}\; ;$$ en remplaçant ${\cal A}$ par $\cal B$ et $x$ par $y$, on se ramène finalement au cas où $\cal A$ est intègre, où $\got{p}=0$, et où $F=\mbox{Frac}\;{\cal O}_{{\cal X},\bf x}$.

\end{itemize}

\medskip
Remarquons que si $|k^{*}|\neq\{1\}$ et si $\cal A$ est strictement $k$-affinoïde, on peut procéder aux réductions ci-dessus en préservant le caractère strictement $k$-affinoïde de $\cal A$. 

\medskip
\trois{assersimplestric} {\bf Preuve des assertions simplifiées, et donc de la proposition, dans le cas où $|k^{*}|\neq\{1\}$ et où $\cal A$ est strictement $k$-affinoïde.} On applique à $\cal X$ la proposition~\ref{EXC}.\ref{ouvertsreg}, $b)$  avec $L=k$ et ${\bf r}=\emptyset$ ; elle fournit un ouvert non vide $\cal V$ de $\cal X$ satisfaisant $b_{2})$ et $b_{3})$ avec ${\cal W}={\cal X}$. 

\medskip
{\em L'assertion $i)$.} Soit $\got{q} \in\spec \widehat{{\cal O}_{X,x}}$ situé au-dessus du point générique de $\spec {\cal O_{X,x}}$. Par platitude de $\spec {\cal O_{X,x}}\to \cal X$, le point $\got{q}$ s'envoie sur le point générique de $\cal X$, qui appartient à $\cal V$. Il résulte alors de $b_{2})$ que $\got{q}$ appartient au lieu régulier de $\spec \widehat{{\cal O}_{X,x}}$.

\medskip
{\em L'assertion $ii)$.} Soit $\got{q}\in\spec \widehat{{\cal O}_{{\cal X},\bf x}}$ situé au-dessus du point générique de $\spec {\cal O}_{{\cal X},\bf x}$ ; le point  $\got{q}$ s'envoie sur le point générique de $\cal X$, qui appartient à $\cal V$. Il résulte alors de $b_{3})$ que $\got{q}$ appartient au lieu régulier de $\spec \widehat{{\cal O}_{{\cal X},\bf x}}$.

\medskip
{\em L'assertion $iii)$.} Soit $\got{q}
\in\spec {\cal O}_{X,x}$ situé au-dessus du point générique de $\spec {\cal O_{{\cal X},\bf x}}$ et soit $\got{q}'$ un point de $\spec \widehat{{\cal O}_{X,x}}$ situé au-dessus de $\got{q}$ ; le point $\got{q}'$ s'envoie sur le point générique de $\cal X$, qui appartient à $\cal V$. Il résulte alors de $b_{2})$ que $\got{q}'$ appartient au lieu régulier de $\spec \widehat{{\cal O}_{X,x}}$ ; par platitude de la complétion, $\got{q}$ appartient au lieu régulier de $\spec {\cal O}_{X,x}$.

\medskip
\trois{assersimplgen} {\bf Preuve des assertions simplifiées, et donc de la proposition, dans le cas général.} Soit $\bf r$ un polyrayon déployant $\cal A$ ; notons que comme $\cal A$ est intègre, il en va de même de ${\cal A}_{\bf r}$ (\ref{RAP}.\ref{tenseurkrintegre}). On applique à $\cal X$ la proposition~\ref{EXC}.\ref{ouvertsreg}, $a)$  avec $L=k_{\bf r}$ ; elle fournit un ouvert non vide ${\cal U}$ de ${\cal X}$ tel que les espaces localement annelés ${\cal U}_{\bf r}$ et ${\cal U}_{\bf r}\an$ soient réguliers. Soit $z$ un antécédent de $x$ sur $X_{\bf r}$, et soit $\bf z$ l'image de $z$ sur ${\cal X}_{\bf r}$.  

\medskip
{\em L'assertion $i)$.} Soit $\got{q}\in \spec \widehat{{\cal O}_{X,x}}$ situé au-dessus du point générique de $\spec {\cal O_{X,x}}$ ; par platitude de $\spec {\cal O}_{X,x}\to \cal X$, le point $\got{q}$ s'envoie sur le point générique de $\cal X$, qui appartient à $\cal U$. Soit $\got{q}'\in\spec \widehat{{\cal O}_{X_{\bf r},z}}$ situé au-dessus de $\got{q}$ ; le point $\got{q}'$ est situé au-dessus de ${\cal U}_{\bf r}$, donc au-dessus du lieu régulier de $\spec {\cal O}_{{\cal X}_{\bf r},\bf z}$. Les flèches $\spec \widehat{{\cal O}_{X_{\bf r},z}}\to \spec {\cal O}_{X_{\bf r},z}\to \spec {\cal O}_{{\cal X}_{\bf r},\bf z}$ sont plates, et à fibres régulières en vertu du cas strictement affinoïde déjà traité. On en déduit que $\got{q}'$ appartient au lieu régulier de $\spec \widehat{{\cal O}_{X_{\bf r},z}}$ ; par platitude, $\got{q}$ appartient au lieu régulier de $\spec \widehat{{\cal O}_{X,x}}$.

\medskip
{\em L'assertion $ii)$.} Soit $\got{q}\in\spec \widehat{{\cal O}_{{\cal X},\bf x}}$ situé au-dessus du point générique de $\spec {\cal O}_{{\cal X},\bf x}$ ; le point  $\got{q}$ s'envoie sur le point générique de $\cal X$, qui appartient à $\cal U$. Soit $\got{q}'\in\spec \widehat{{\cal O}_{{\cal X}_{\bf r},\bf z}}$ situé au-dessus de $\got{q}$ ; le point $\got{q}'$ est situé au-dessus de ${\cal U}_{\bf r}$, donc au-dessus du lieu régulier de $\spec {\cal O}_{{\cal X}_{\bf r},\bf z}$. Comme la flèche $\spec \widehat{{\cal O}_{{\cal X}_{\bf r},\bf z}}\to \spec {\cal O}_{{\cal X}_{\bf r},\bf z}$ est plate, et à fibres régulières en vertu du cas strictement affinoïde déjà traité, $\got{q}'$ appartient au lieu régulier de $\spec \widehat{{\cal O}_{{\cal X}_{\bf r},\bf z}}$ ; par platitude, $\got{q}$ appartient au lieu régulier de $\spec \widehat{{\cal O}_{{\cal X},\bf x}}$.

\medskip
{\em L'assertion $iii)$.} Soit $\got{q}
\in\spec {\cal O}_{X,x}$ situé au-dessus du point générique de $\spec {\cal O_{{\cal X},\bf x}}$ ; le point  $\got{q}$ s'envoie sur le point générique de $\cal X$, qui appartient à $\cal U$. Soit $\got{q}'\in\spec {\cal O}_{X_{\bf r}, z}$ situé au-dessus de $\got{q}$ ; le point $\got{q}'$ est situé au-dessus de ${\cal U}_{\bf r}$, donc au-dessus du lieu régulier de $\spec {\cal O}_{{\cal X}_{\bf r},\bf z}$. Comme la flèche $\spec {\cal O}_{X_{\bf r},z}\to \spec {\cal O}_{{\cal X}_{\bf r},\bf z}$ est plate, et à fibres régulières en vertu du cas strictement affinoïde déjà traité, $\got{q}'$ appartient au lieu régulier de $\spec {\cal O}_{X_{\bf r},z}$ ; par platitude, $\got{q}$ appartient au lieu régulier de $\spec {\cal O}_{X,x}$.~$\Box$ 

\medskip
\deux{rempropiii} {\bf Remarque.} La régularité géométrique de $\spec {\cal O}_{X,x}\to \spec {\cal O}_{{\cal X},\bf x}$ n'est pas {\em a priori} requise pour établir les résultats d'excellence que nous avons en vue ; mais nous l'avons insérée dans l'énoncé car nous l'utilisons lors de la preuve de la régularité géométrique de $ \spec\widehat{{\cal O}_{{\cal X},\bf x}}\to \spec {\cal O}_{{\cal X},\bf x}$ dans le cas général ; nous l'étendrons un peu plus bas au cas d'un schéma de type fini sur une algèbre affinoïde (th. \ref{GAGA}.\ref{geomregulschem}). 

\medskip
\deux{lieureg} {\bf Proposition.} {\em Soit $\cal A$ une algèbre affinoïde, soit $\got{p}\in \spec \cal A$ et soit $F$ une extension finie et radicielle de $\kappa(\got{p})$. Il existe une sous-${\cal A}/\got{p}$-algèbre finie et radicielle $\cal B$ de $F$, de corps des fractions égal à $F$, tel que le lieu régulier de $\spec \cal B$ contienne un ouvert non vide.} 

\medskip
{\em Démonstration.}  Quitte à quotienter $\cal A$ par $\got{p}$, on peut supposer que $\cal A$ est intègre et que $\got{p}=0$. Le corps $F$ est alors une extension finie et radicielle de $\mbox{Frac}\;\cal A$. Il existe une sous-$\cal A$-algèbre finie et radicielle ${\cal B}$ de $F$ dont le corps des fractions est égal à $F$ ; la proposition~\ref{EXC}.\ref{ouvertsreg} assure qu'il existe un ouvert non vide de $\spec \cal B$ qui est régulier.~$\Box$

\medskip
\deux{excellence} {\bf Théorème.} {\em Soit $\cal A$ une algèbre $k$-affinoïde. L'algèbre $\cal A$ est excellente, et les anneaux locaux de ${\cal M}({\cal A})$ sont excellents.} 

\medskip
{\em Démonstration.} Compte-tenu des propositions~\ref{EXC}.\ref{geomregul} et ~\ref{EXC}.\ref{lieureg}, il reste simplement à s'assurer de l'universelle caténarité de $\cal A$ et de chacun des anneaux locaux de ${\cal M}({\cal A})$. Or $\cal A$ peut s'écrire comme un quotient de $k\{{\bf r}^{-1}{\bf T}\}$ pour un certain polyrayon $\bf r$. Le lemme~\ref{EXC}.\ref{disquereg} assure que $k\{{\bf r}^{-1}{\bf T}\}$ est régulier, et que l'espace localement annelé ${\cal M}(k\{{\bf r}^{-1}{\bf T}\})$ l'est également, ce qui permet de conclure ({\em cf.} \cite{mats}, th. 17.9).~$\Box$ 

\medskip
\deux{refconradkiehl} {\bf Remarque.} La première assertion du théorème ci-dessus a été démontrée par Kiehl dans le cas d'une algèbre strictement affinoïde (\cite{kie}, th. 3.3) ; la seconde l'a été par Conrad (\cite{comprig}, th. 1.1.3) dans le cas des anneaux locaux en les points {\em rigides} d'un espace strictement affinoïde.
 
\section{Changement de corps de base et théorèmes de comparaison}\label{GAGA}

\setcounter{cpt}{0}

{\em On utilisera dans ce qui suit les ensembles de propriétés $\cal P, \cal Q$ et $\cal R$ définis au} \ref{RAP}.\ref{listepropri}.

\subsection*{Le changement de base} 

\bigskip
\deux{changermsm} {\bf Théorème.} {\em Soit $\cal X$ un schéma de type fini sur une algèbre $k$-affinoïde $\cal A$. Soit $L$ une extension complète de $k$, soit $\bf x$ un point de $\cal X$ et soit $\bf y$ un point de ${\cal X}_{L}$  situé au-dessus de $\bf x$.  Soit $\cal F$ un faisceau cohérent sur $\cal X$. Soit ${\mathsf P}\in {\cal S}\cup {\cal Q}\cup {\cal R}$ ; si ${\mathsf P}\in{\cal Q}\cup {\cal R}$, on suppose que ${\cal F}={\cal O}_{\cal X}$. 

\bigskip
\begin{itemize}

\item[$a)$] Si ${\cal F}\otimes{\cal O}_{{\cal X}_ {L},{\bf y}}$ vérifie $\mathsf P$, alors ${\cal F}\otimes{\cal O}_{{\cal X},{\bf x}}$ vérifie $\mathsf P$.

\medskip
\item[$b)$] Si ${\mathsf P}\in {\cal S}\cup \cal Q$ et si ${\cal F}\otimes{\cal O}_{{\cal X},{\bf x}}$ vérifie $\mathsf P$, alors ${\cal F}\otimes{\cal O}_{{\cal X}_ {L},{\bf y}}$ vérifie $\mathsf P$. 

\medskip
\item[$c)$] Si ${\mathsf P}\in {\cal R}$, si ${\cal F}\otimes{\cal O}_{{\cal X},{\bf x}}$ vérifie $\mathsf P$, et si $L$ est une extension analytiquement séparable de $k$, alors ${\cal F}\otimes{\cal O}_{{\cal X}_ {L},{\bf y}}$ vérifie $\mathsf P$.

\end{itemize}

}

\bigskip
{\em Démonstration.} La platitude de ${\cal X}_{L} \to {\cal X}$ entraîne immédiatement l'assertion $a)$ (\ref{RAP}.\ref{platdirect}). 

\medskip
Soit ${\bf t}\in \cal X$, soit ${\bf s}\in {\cal X}_{L}$ situé au-dessus de $\bf t$, et soit $\cal Y$ l'adhérence réduite de $\bf t$ dans $\cal X$. La proposition~\ref{EXC}.\ref{ouvertsreg} assure l'existence d'un ouvert non vide $\cal U$ de $\cal Y$ tel que ${\cal U}_{L}$ soit d'intersection complète en chacun de ses points, et régulier si $L$ est analytiquement séparable sur $k$. On en déduit que ${\cal O}_{{\cal X}_{L},\bf s}\otimes_{{\cal O}_{{\cal X},\bf t}}\kappa({\bf t})$ est d'intersection complète en chacun de ses points, et régulier si $L$ est analytiquement séparable sur $k$. En vertu de la platitude de ${\cal X}_{L}\to {\cal X}$, les assertions $b)$ et $c)$ s'en déduisent aussitôt (\ref{RAP}.\ref{platdirect}).~$\Box$

\subsection*{Régularité géométrique des fibres analytiques} 

\bigskip
\deux{gagaregul} {\bf Lemme.} {\em Soit $\cal X$ un schéma de type fini sur une algèbre $k$-affinoïde $\cal A$, soit $V$ un domaine affinoïde de ${\cal X}\an$, et soit ${\cal V}$ le spectre de l'anneau des fonctions analytiques sur $V$. Soit $\bf v \in {\cal V}$ et soit $\bf x$ son image sur $\cal X$. Si $\cal X$ est régulier en $\bf x$ alors $\cal V$ est régulier en $\bf v$.}

\medskip
{\em Démonstration.} Supposons donc $\cal X$ régulier en $\bf x$. Soit $v$ un point de $V$ situé au-dessus de $\bf v$. Puisque $\cal A$ est un anneau excellent d'après le théorème \ref{EXC}.\ref{excellence}, le lieu de régularité $\cal U$ de $\cal X$ en est un ouvert de Zariski ; on a ${\bf x}\in \cal U$ et $v\in {\cal U}\an$. Soit $\bf r$ un polyrayon déployant $\cal A$ et $V$. Comme $k_{\bf r}$ est une extension analytiquement séparable de $k$, le théorème précédent assure que le schéma ${\cal U}_{\bf r}$ est régulier. 

\medskip
Soit $w$ un antécédent de $v$ sur ${\cal U}_{\bf r}\an$ et soit $W$ un voisinage strictement $k_{\bf r}$-affinoïde de $w$ dans le bon espace strictement $k_{\bf r}$-analytique $V_{\bf r}\cap {\cal U}_{\bf r}\an$. Soit ${\cal W}$ le spectre de l'algèbre des fonctions analytiques sur $W$. Le morphisme ${\cal W}\to {\cal V}$ est plat et son image contient $\bf v$ ; il suffit donc (\ref{RAP}.\ref{platdirect}) de démontrer que $\cal W$ est régulier. 

\medskip
Soit $\bf t$ un point fermé de $\cal W$. Il correspond par le {\em Nullstellensatz} à un point $k_{\bf r}$-rigide $t$ de $W$, dont on note $\bf u$ l'image sur ${\cal U}_{\bf r}$. L'anneau local complété $\widehat{{\cal O}_{{\cal U}_{\bf r},\bf u}}$ est régulier ; l'on dispose par ailleurs d'isomorphismes naturels $$\widehat{{\cal O}_{{\cal U}_{\bf r},\bf u}}\simeq \widehat{{\cal O}_{{\cal U}_{\bf r}\an, t}}\simeq \widehat{{\cal O}_{W, t}}\simeq \widehat{{\cal O}_{{\cal W},\bf t}} \;;$$ le premier et le dernier proviennent du lemme 2.6.3 de \cite{brk2}, et celui du milieu du fait que le point $k_{\bf r}$-rigide $t$ appartient à l'intérieur topologique de $W$ dans ${\cal U}_{\bf r}\an$ ; on en déduit que $\cal W$ est régulier en $\bf t$. Ceci valant pour tout point fermé $\bf t$ de $\cal W$, le schéma $\cal W$ est régulier. ~$\Box$

\medskip
La proposition ci-dessous constitue entre autres une généralisation de l'assertion $iii)$ de la proposition \ref{EXC}.\ref{geomregul} à un schéma de type fini sur une algèbre affinoïde.

\medskip
\deux{geomregulschem} {\bf Théorème.} {\em Soit $\cal X$ un schéma de type fini sur une algèbre affinoïde. Soit $V$ un domaine affinoïde de ${\cal X}\an$ et soit $\cal V$ le spectre de son algèbre des fonctions. Soit $x\in V$ et soient $\bf v$ et $\bf x$ ses images respectives sur $\cal V$ et $\cal X$. 

\medskip
Dans le diagramme commutatif $$\xymatrix { \spec  {\cal O}_{V,x} \ar[r]\ar[d]\ar[rd] &\spec {\cal O}_{{\cal V},\bf v}\ar[d]\\ \spec{\cal O}_{{\cal X}\an,x}\ar[r]&\spec {\cal O}_{{\cal X},\bf x} } $$ toutes les flèches sont géométriquement régulières.}

\medskip
{\em Démonstration.} On sait déjà que les deux flèches horizontales et la flèche verticale de gauche sont fidèlement plates ; la fidèle platitude de la flèche diagonale et de la flèche verticale de droite s'ensuivent immédiatement. 

\medskip
Venons-en à la régularité géométrique des fibres ; tout au long de la preuve nous utiliserons sans les rappeler explicitement à chaque occurence les résultats du \ref{RAP}.\ref{platdirect} relatifs au comportement de la régularité vis-à-vis d'un morphisme plat. 

\medskip

\setcounter{cptbis}{0}
\medskip
\trois{geomregulzero} {\em  Régularité géométrique de $\spec  {\cal O}_{V,x} \to\spec {\cal O}_{{\cal V},\bf v}$.} C'est l'une des assertions de la proposition \ref{EXC}.\ref{geomregul}.

\medskip
\trois{geomregulun} {\em Régularité géométrique de $\spec  {\cal O}_{V,x} \to \spec {\cal O}_{{\cal X},\bf x}$.} Soit $\got p$ un idéal premier de ${\cal O}_{{\cal X},\bf x}$ et soit $F$ une extension finie purement inséparable de $\kappa(\got p)$ ; nous allons démontrer que ${\cal O}_{V,x}\otimes_{{\cal O}_{{\cal X},\bf x}}F$ est régulier. Soit $\cal Z$ le fermé de Zariski de $\cal X$ égal à l'adhérence du point qui correspond à $\got p$, muni de sa structure réduite. Quitte à remplacer $\cal X$ par $\cal Z$, on peut supposer que $\cal X$ est intègre et que $\got p=0$. Il existe une sous-${\cal O}_{{\cal X},\bf x}$-algèbre finie et radicielle $\mathsf B$ de $F$ dont $F$ est le corps des fractions ; on peut donc restreindre $\cal X$ de sorte qu'il existe  un schéma intègre $\cal Y$ et un morphisme ${\cal Y}\to \cal X$ fini, dominant et radiciel de fibre générique isomorphe à $\spec F$. Soit $\bf y$ (resp. $y$) l'unique point de $\cal Y$ (resp. ${\cal Y}\an$) situé au-dessus de $\bf x$ (resp. $x$) et soit $W$ l'image réciproque de $V$ sur ${\cal Y}\an$. 

\medskip
L'on dispose ({\em cf.} lemme \ref{RAP}.\ref{finicomplet}) d'un isomorphisme ${\cal O}_{W,y}\simeq{\cal O}_{V,x}\otimes_{{\cal O}_{{\cal X},\bf x}}{\cal O}_{{\cal Y},\bf y}$, et l'on a donc ${\cal O}_{W,y}\otimes_{{\cal O}_{{\cal Y},\bf y}}\mbox{Frac}\;{\cal O}_{{\cal Y},\bf y}\simeq {\cal O}_{V,x}\otimes_{{\cal O}_{{\cal X},\bf x}}F.$ Ceci  permet, en remplaçant $\cal X$ par $\cal Y$, $x$ par $y$ et $V$ par $W$, de se ramener finalement au cas où $\cal X$ est intègre, où $\got p=0$ et où $F=\kappa(\got p)$. 

\medskip
Soit $\got q$ un idéal premier de ${\cal O}_{V,x}$ situé au-dessus du point générique de $\cal X$. Notons $\bf z$ son image sur $\cal V$ ; par le lemme \ref{GAGA}.\ref{gagaregul} ci-dessus, ${\cal O}_{{\cal V},\bf z}$ est régulier. En vertu de la régularité géométrique de $\spec {\cal O}_{V,x}\to \spec O_{{\cal V},\bf v}$ (\ref{GAGA}.\ref{geomregulschem}.\ref{geomregulzero}), $\got q$ appartient au lieu régulier de $\spec {\cal O}_{V,x}$, ce qu’il fallait démontrer.

\medskip
\trois{geomreguldeux} {\em Régularité géométrique de $\spec  {\cal O}_{V,x} \to \spec {\cal O}_{{\cal X}\an, x}$.} Quitte à remplacer ${\cal X}\an$ par un voisinage affinoïde de $x$, et $V$ par son intersection avec celui-ci, on peut supposer que $\cal X$ est le spectre d'une algèbre affinoïde $\cal A$ ; l'on a alors ${\cal X}\an={\cal M}({\cal A})$ ; nous allons également désigner cet espace affinoïde par $X$. Soit $\got p$ un idéal premier de ${\cal O}_{X,x}$ et soit $F$ une extension finie radicielle de $\kappa(\got p)$ ; on cherche à montrer que ${\cal O}_{V,x}\otimes_{{\cal O}_{X,x}}F$ est régulier.

\medskip
On peut supposer, quitte à restreindre $X$, que $\got{p}$ provient d'un idéal de $\cal A$ puis, en quotientant $\cal A$ par ce dernier, que $\got{p}=0$. Soit $\mathsf B$ une sous ${\cal O}_{X,x}$-algèbre finie et radicielle de $F$ de corps des fractions égal à $F$. On peut à nouveau restreindre $X$ de sorte que $\mathsf B$ soit de la forme ${\cal B}\otimes_{\cal A}{\cal O}_{X,x}$, où $\cal B$ est une $\cal A$-algèbre finie et radicielle. Soit $y$ l'unique antécédent de $x$ sur $Y={\cal M}({\cal B})$, et soit $W$ l'image réciproque de $V$ sur $Y$.  Comme ${\cal O}_{Y,y}$ s'identifie à $\mathsf B$, et ${\cal O}_{W,y}$  à ${\cal O}_{V,x}\otimes_{{\cal O}_{X,x}} \mathsf B$ ({\em cf.} lemme \ref{RAP}.\ref{finicomplet}), on a $${\cal O}_{V,x}\otimes_{{\cal O}_{X,x}}F\simeq {\cal O}_{W,y}\otimes_{{\cal O}_{Y,y}}\mbox{Frac}\; {\cal O}_{Y,y}.$$ En remplaçant $X$ par $Y$, $x$ par $y$ et $V$ par $W$, on se ramène au cas où $\got p=0$ et où l'on a de surcroît $F=\mbox{Frac}\;{\cal O}_{X,x}$. Comme ${\cal O}_{X,x}$ est intègre, le lemme~\ref{RAP}.\ref{localintegre} assure que $x$ n'est situé que sur une composante irréductible de $X$. En substituant à $\cal A$ son quotient par l'idéal premier correspondant à la composante en question ({\em cf.} remarque \ref{RAP}.\ref{quotoxx}), on se ramène finalement au cas où $\cal A$ est intègre, où $\got{p}=0$, et où $F=\mbox{Frac}\;{\cal O}_{X,x}$. 

\medskip
Notons $\cal X$ le spectre de $\cal A$. Soit $\got q$ un idéal premier de ${\cal O}_{V,x}$ situé au-dessus du point générique de $\spec {\cal O}_{X,x}$. Par platitude de $\spec {\cal O}_{X,x}\to \cal X$, l'image de $\got q$ sur le schéma intègre $\cal X$ est le point générique de celui-ci ; en vertu de la régularité géométrique de $\spec  {\cal O}_{V,x} \to \spec {\cal O}_{{\cal X},\bf x}$, établie au \ref{GAGA}.\ref{geomregulschem}.\ref{geomregulun} ci-dessus, $\got q$ appartient au lieu régulier de $\spec {\cal O}_{V,x}$.

\medskip
\trois{geomregultrois} {\em Régularité géométrique de $\spec {\cal O}_{{\cal V},\bf v}\to \spec {\cal O}_{{\cal X},\bf x}$.} Elle découle immédiatement de la platitude de $\spec  {\cal O}_{V,x}\to\spec {\cal O}_{{\cal V},\bf v}$ et de la régularité géométrique de $\spec  {\cal O}_{V,x} \to \spec {\cal O}_{{\cal X},\bf x}$ établie au \ref{GAGA}.\ref{geomregulschem}.\ref{geomregulun} ci-dessus.

\medskip
\trois{geomregulquatre} {\em Régularité géométrique de $\spec {\cal O}_{{\cal X}\an, x}\to \spec{\cal O}_{{\cal X},\bf x}$.} C'est un cas particulier de la régularité géométrique de $\spec  {\cal O}_{V,x} \to \spec {\cal O}_{{\cal X}\an, x}$, établie au \ref{GAGA}.\ref{geomregulschem}.\ref{geomreguldeux} ci-dessus : celui où $V$ est un  {\em voisinage} affinoïde de $x$ dans ${\cal X}\an$.~$\Box$ 

\subsection*{Le théorème de comparaison} 

{\em On rappelle que les ensembles de propriétés $\cal P, \cal Q$ et $\cal R$ ont été définis au} \ref{RAP}.\ref{listepropri}.

\bigskip
\deux{gagarmsm} {\bf Théorème.} {\em Soit $\cal X$ un schéma de type fini sur une algèbre affinoïde $\cal A$. Soit $V$ un domaine affinoïde de ${\cal X}\an$, soit $x$ un point de $V$ et soit $\bf x$ son image sur $\cal X$. Soit $\cal F$ un faisceau cohérent sur $\cal X$. Soit ${\mathsf P}\in {\cal S}\cup {\cal Q}\cup {\cal R}$ ; si ${\mathsf P}\in {\cal Q}\cup {\cal R}$, on suppose que ${\cal F}={\cal O}_{\cal X}$. Sous ces hypothèses : 

\bigskip
\begin{itemize}

\medskip
\item[$A)$]  l'ensemble $\cal U$ des points $\bf z$ de $\cal X$ tels que ${\cal F}\otimes{\cal O}_{{\cal X},{\bf z}}$ vérifie $\mathsf P$ est un ouvert de Zariski ; 
\medskip
\item[$B)$] les propositions suivantes sont équivalentes : 

\medskip
\begin{itemize}

\item[$B')$] ${\cal F}\otimes{\cal O}_{V,x}$ vérifie $\mathsf P$ ;
\item[$B'')$] ${\cal F}\otimes{\cal O}_{{\cal X}\an,x}$ vérifie $\mathsf P$ ;
\item[$B''')$] ${\cal F}\otimes{\cal O}_{{\cal X},\bf x}$ vérifie $\mathsf P$.

\end{itemize}
\end{itemize}
}

\bigskip

{\em Démonstration.} L'assertion $A)$ provient du fait que $\cal A$ est un anneau excellent par la proposition \ref{EXC}.\ref{excellence}, et des résultats rappelés au \ref{RAP}.\ref{lieuxouverts} ; notons que comme $\cal A$ est un quotient d'un anneau régulier par le lemme \ref{EXC}.\ref{disquereg}, la situation est justiciable de la remarque \ref{RAP}.\ref{reminterscompl} lorsque $\mathsf P$ est la propriété d'être d'intersection complète.

\medskip
Les morphismes $$\spec {\cal O}_{V,x}\to \spec {\cal O}_{{\cal X}\an,x}\; \mbox{et}\;  \spec {\cal O}_{{\cal X}\an,x}\to \spec {\cal O}_{{\cal X},\bf x}$$ sont géométriquement réguliers d'après le théorème \ref{GAGA}.\ref{geomregulschem} ci-dessus ; en vertu du \ref{RAP}.\ref{platdirect}, les équivalences $B')\iff B'')\iff B''')$ s'ensuivent immédiatement.~$\Box$ 

\medskip
\deux{gagaberkorappel} {\bf Remarque.} Lorsque ${\cal X}={\cal A}$, l'équivalence $B'')\iff B''')$ dans le théorème ci-dessus a été prouvée par Berkovich (\cite{brk2}, th. 2.2.1) ; on se convainc facilement que {\em lorsque $|k^*|\neq \{1\}$ et lorsque $\cal A$ est strictement $k$-affinoïde} sa démonstration s'étend {\em mutatis mutandis} au cas où $\cal X$ est un $\cal A$-schéma de type fini quelconque ; de même, celle du corollaire 2.2.8 de \cite{brk2} permet en fait d'établir $B'\iff B'')$ dès que $|k^*|\neq \{1\}$ et que $\cal A$ et $V$ sont strictement $k$-affinoïdes ; nous nous sommes en partie inspiré ici de ces différentes preuves.

\subsection*{Les propriétés usuelles de l'algèbre commutative dans le cas non nécessairement bon} 

{\em On travaille toujours avec  les ensembles de propriétés $\cal P, \cal Q$ et $\cal R$ définis au} \ref{RAP}.\ref{listepropri}.

\medskip
\deux{defpasbon} Soit $X$ un espace analytique, soit $\cal F$ un faisceau cohérent sur $X$, et soit $\mathsf P\in{\cal S}\cup {\cal Q}\cup {\cal R}$ ; si ${\mathsf P}\in {\cal Q}\cup {\cal R}$, on suppose que ${\cal F}={\cal O}_{X\grot}$. Soit $x$ un point de $X$ et soient $V$ et $W$ deux domaines affinoïdes de $X$ contenant $x$ et tels que $W\subset V$. Il résulte du théorème~\ref{GAGA}.\ref{gagarmsm} que ${\cal F}\otimes {\cal O}_{V,x}$ satisfait $\mathsf P$ si et seulement si ${\cal F}\otimes {\cal O}_{W,x}$ satisfait $\mathsf P$. On en déduit l'équivalence des deux propositions suivantes : 

\medskip
\begin{itemize}
\item[$i)$]  il existe un bon domaine analytique $U$ de $X$ contenant $x$ tel que ${\cal F}\otimes{\cal O}_{U,x}$ satisfasse $\mathsf P$ ;
\item[$ii)$]  ${\cal F}\otimes{\cal O}_{U,x}$ satisfait $\mathsf P$ pour {\em tout} bon domaine analytique $U$ de $X$ contenant $x$.
\end{itemize}
\medskip

Lorsque ces deux propriétés sont vérifiées, on dit que {\em $\cal F$ satisfait $\mathsf P$ en $x$}.

\setcounter{cptbis}{0}

\medskip
\trois{remsixbon} Si $X$ est bon, il résulte de la définition que $\cal F$ satisfait $\mathsf P$ en $x$ si et seulement si ${\cal F}\otimes{\cal O}_{X,x}$ satisfait $\mathsf P$. 

\medskip
\trois{ppropglocale} Si $U$ est un domaine analytique de $X$ contenant $x$, il est immédiat que ${\cal F}$ satisfait $\mathsf P$ en $x$ si et seulement si ${\cal F}_{|U}$ satisfait $\mathsf P$ en $x$.

\medskip
\trois{espreduitnorm} On dira que ${\cal F}$ satisfait $\mathsf P$ (et même, si ${\cal F}={\cal O}_{X\grot}$, que $X$ satisfait $\mathsf P$) si $\cal F$ satisfait $\mathsf P$ en tout point de $X$. Ceci vaut bien entendu également pour la conjonction de plusieurs propriétés ; on parlera ainsi d'espace analytique normal ({\em i.e.} satisfaisant $R_1$ et $S_2$, d'après le critère de Serre), réduit  ({\em i.e.} satisfaisant $R_0$ et $S_1$, d'après le critère de Serre), etc.

\medskip
\trois{encorechangebase} Certains résultats des théorèmes~\ref{GAGA}.\ref{changermsm} et \ref{GAGA}.\ref{gagarmsm} peuvent être reformulés dans ce nouveau contexte. Les énoncés obtenus se résument comme suit : les propriétés appartenant à $\cal S \cup \cal Q$ descendent et montent par extension quelconque des scalaires ; celles qui appartiennent à $\cal R$ descendent par extension quelconque et montent par extension {\em analytiquement séparable} des scalaires.

\medskip
{\bf Les théorèmes \ref{GAGA}.\ref{changermsm} et \ref{GAGA}.\ref{gagarmsm}, ainsi que leurs déclinaisons dans le cas non nécessairement bon (\ref{GAGA}.\ref{defpasbon} et \ref{GAGA}.\ref{defpasbon}.\ref{ppropglocale}-\ref{GAGA}.\ref{defpasbon}\ref{encorechangebase}) seront désormais utilisés librement, sans rappel des références.} 

\section{Les composantes irréductibles en géométrie analytique}\label{COMP}

\subsection*{Caractère G-local de la topologie de Zariski et structures réduites} 
\setcounter{cpt}{0}

\medskip
\deux{fcoher} {\bf Lemme.} {\em Soit $\cal A$ une algèbre affinoïde et soit $X$ l'espace ${\cal M}({\cal A})$. Soit $ I$ un idéal de $\cal A$, et soit $Y$ le lieu des zéros de $I$ sur $X$. Soit $\cal F$ le faisceau d'idéaux de ${\cal O}_{X\grot}$ qui associe à un domaine analytique $V$ de $X$ l'ensemble des $f\in {\cal O}_{X\grot}(V)$ telles que $f(x)=0$ pour tout $x\in V\cap Y$. Le faisceau $\cal F$ coïncide avec le faisceau cohérent $\red{\sqrt{I}}$.} 

\medskip
{\em Démonstration.} L'inclusion $\red{\sqrt{I}}\subset {\cal F}$ est triviale. Pour établir la réciproque, il suffit de vérifier que ${\cal F}(V)\subset\red{\sqrt{I}}(V)=\sqrt{I}.{\cal A}_V$ pour tout domaine affinoïde $V$ de $X$ ; donnons-nous un tel $V$. 

\medskip
On peut décrire $Y$ comme le lieu des zéros de  $\red{\sqrt{I}}$ ; munissons-le de la structure de sous-espace analytique fermé correspondante. Il s'identifie alors à ${\cal M}({\cal A}/\sqrt{I})$, et $V\cap Y$ en est un domaine affinoïde, dont l'algèbre des fonctions analytiques est ${\cal A}_V/(\sqrt{I}.{\cal A}_V)$ ; l'anneau ${\cal A}/\sqrt{ I}$ étant réduit, ${\cal A}_V/(\sqrt{ I}.{\cal A}_V)$ l'est aussi. Soit $f$ appartenant à ${\cal F}(V)$. Par hypothèse, $f$ est nulle en tout point de $V\cap Y$. Son image dans ${\cal A}_V/(\sqrt{I}.{\cal A}_V)$ est donc nilpotente, et partant  nulle puisque cet anneau est réduit. Cela signifie que $f\in \sqrt{ I}.{\cal A}_V$.~$\Box$ 

\medskip
\deux{zargloc} {\bf Proposition.} {\em Soit $X$ un espace analytique et soit $(X_i)$ un G-recouvrement de $X$ par des domaines analytiques. Soit $Y$ une partie de $X$ telle que $Y\cap X_i$ soit un fermé de Zariski de $X_i$ pour tout $i$. Soit $\cal F$ le faisceau d'idéaux de ${\cal O}_{X\grot}$ qui associe à un domaine analytique $V$ de $X$ l'ensemble des $f\in {\cal O}_{X\grot}(V)$ telles que $f(x)=0$ pour tout $x\in V\cap Y$. 

\medskip
$i)$ Le faisceau $\cal F$ est cohérent ; le sous-ensemble $Y$ de $X$ s'identifie à son lieu des zéros et est  donc un fermé de Zariski ;

\medskip
$ii)$ le faisceau $\cal F$ est le plus grand faisceau cohérent d'idéaux dont $Y$ est le lieu des zéros ;

\medskip
$iii)$ le sous-espace analytique fermé $Y\duit$ de $X$ défini par $\cal F$ est un objet final de la catégorie des espaces analytiques réduits munis d'un morphisme vers $X$ dont l'image est incluse dans $Y$.} 

\medskip
{\em Démonstration.} Quitte à raffiner le recouvrement $(X_i)$, on peut supposer que chacun des $X_i$ est affinoïde. Le lemme \ref{COMP}.\ref{fcoher} ci-dessus assure que ${\cal F}_{|X_i}$ est cohérent pour tout $i$ ; par conséquent, ${\cal F}$ est cohérent. 

\medskip
Fixons $i$. Comme $Y\cap X_i$ est un fermé de Zariski de l'espace affinoïde $X_i$, il est tautologique que $Y\cap X_i$ est le lieu des zéros de ${\cal F}_{|X_i}$ ; on en déduit que $Y$ est le lieu des zéros de $\cal F$ (\ref{RAP}.\ref{zarcasgen}.\ref{zarcasgendom}), ce qui termine la preuve de $i)$.

\medskip
L'assertion  $ii)$ résulte immédiatement de la définition de $\cal F$. L'assertion $iii)$ se vérifie G-localement sur $X$, ce qui autorise à le supposer affinoïde ; soit $\cal A$ l'algèbre des fonctions analytiques sur $X$ et soit $I={\cal F}(X)$ ; on a $Y\duit={\cal M}({\cal A}/I)$ ; comme $I=\sqrt I$ par définition du faisceau $\cal F$, l'anneau ${\cal A}/I$ est réduit, et l'espace $Y\duit$ est de ce fait réduit ; l'image de $Y\duit\to X$ est évidemment incluse dans (et même égale à) $Y$.

\medskip
 Il suffit, en raisonnant G-localement sur la source des morphismes, de vérifier que $Y\duit$ est un objet final de la catégorie des espaces {\em affinoïdes} réduits munis d'un morphisme vers $X$ dont l'image est incluse dans $Y$. Soit $\cal B$ une algèbre affinoïde réduite et soit ${\cal A}\to {\cal B}$ un morphisme borné tel que l'image de ${\cal M}({\cal B})$ sur ${\cal M}({\cal A})$ soit incluse dans $Y$. Si $f$ est un élément de $I$, il est nul en tout point de $Y$ ; son image dans $\cal B$ est donc nulle en tout point de ${\cal M}({\cal B})$, et est par conséquent nilpotente ; comme $\cal B$ est réduite, elle est nulle. On en déduit que ${\cal B}\to {\cal A}$ se factorise de manière unique par un morphisme, nécessairement borné, de ${\cal A}/I$ dans $\cal B$.~$\Box$ 

\setcounter{cptbis}{0} 

\medskip
\trois{yredabus} La structure définie par $\cal F$ est appelée la {\em structure réduite} du fermé de Zariski $Y$ ; si l'on a préalablement pris la peine de préciser que l'on munit $Y$ de sa structure réduite, on s'autorisera à écrire $Y$ au lieu de $Y\duit$. Si $Y=X$, le faisceau $\cal F$ est celui des fonctions G-localement nilpotentes (qui se trouve donc être cohérent), et $X\duit$ est appelé {\em l'espace analytique réduit associé à $X$}. Notons que $X=X\duit$ si et seulement si $X$ est réduit. 

\medskip
\deux{conszargloc} {\bf Remarque.} La proposition ci-dessus assure en particulier que la propriété d'être ouvert (ou fermé) de Zariski est de nature G-locale. Ce fait permet de globaliser certains résultats obtenu dans le cas affinoïde ; donnons-en une illustration. 

\medskip
\deux{lieupasbonfermzar} {\bf Théorème.} {\em Soit $X$ un espace analytique, soit $\cal F$ un faisceau cohérent sur $X$ et soit $\mathsf P$ une propriété   appartenant à ${\cal S} \cup {\cal Q}\cup {\cal R}$ (les ensembles ${\cal S}, {\cal Q}$ et $ {\cal R}$ ont été définis au {\rm \ref{RAP}.\ref{listepropri}} ) ; si $\mathsf P$ appartient à ${\cal Q}\cup {\cal R}$, on suppose que ${\cal F}={\cal O}_{X\grot}$. L'ensemble $Z$ des points de $X$ en lesquels $\cal F$ satisfait $\mathsf P$ est un ouvert de Zariski de $X$.} 

\medskip
{\em Démonstration.} Si $V$ est un domaine affinoïde de $X$,  l'intersection $V\cap Z$ est l'ensemble des points de $V$ en lesquels ${\cal F}_{|V}$ satisfait $\mathsf P$ ; on déduit du théorème \ref{GAGA}.\ref{gagarmsm} (appliqué à ${\cal X}=\spec {\cal A}_V$) que $V\cap Z$ est un ouvert de Zariski de $V$. Ceci valant pour tout domaine affinoïde $V$ de $X$, la proposition \ref{COMP}.\ref{zargloc} garantit que $Z$ est un ouvert de Zariski de $X$.~$\Box$ 

\medskip
\deux{remformu} {\bf Remarque.} Lorsque $\mathsf P\in {\cal S}\cup {\cal Q}$, il est immédiat que la formation de l'ouvert de Zariski évoqué ci-dessus commute à tout changement du corps de base.

\medskip
{\bf Le théorème ci-dessus sera utilisé librement dans la suite, sans rappel ni référence.}

\medskip
\deux{remzred} {\bf Remarque.} Soit $X$ un espace analytique, soit $\cal I$ un faisceau cohérent d'idéaux sur $X$ et soit $Y$ son lieu des zéros ; munissons-le de la structure associée à $\cal I$. Soit $Z$ un fermé de Zariski de $Y$. La notation $Z\duit$ peut {\em a priori} faire référence à $Z$ ou bien en tant que fermé de Zariski de l'espace analytique $Y$, ou bien en tant que fermé de Zariski de $X$ ; cette ambiguïté n'en est en réalité pas une, car ces deux espaces notés $Z\duit$ sont canoniquement isomorphes, comme on le vérifie à l'aide de leurs propriétés universelles.

\medskip
\deux{redetdom} {\bf Remarque.} Soit $X$ un espace analytique, soit $Y$ un fermé de Zariski de $X$ et soit $U$ un domaine analytique de $X$ ; il résulte des propriétés universelles des espaces en jeu que les sous-ensembles analytiques fermés $U\times_X Y\duit$ et $(U\cap Y)\duit$ de $U$ sont canoniquement isomorphes.

\medskip
\deux{redannloc} {\bf Proposition.} {\em Soit $X$ un bon espace analytique et soit $x\in X$. Le morphisme naturel $({\cal O}_{X,x})\duit\to{\cal O}_{X\duit,x}$ est un isomorphisme.}

\medskip
{\em Démonstration.} On peut supposer que $X$ est affinoïde ; soit $\cal A$ l'algèbre affinoïde correspondante, et soit $\cal I$ l'idéal des éléments nilpotents de $\cal A$. L'anneau ${\cal O}_{X\duit,x}$ est réduit et s'identifie au quotient de ${\cal O}_{X,x}$ par son idéal nilpotent ${\cal I}{\cal O}_{X,x}$ ; il est donc naturellement isomorphe à $({\cal O}_{X,x})\duit$.~$\Box$ 

\medskip
\deux{xredl} {\bf Lemme.} {\em Soit $X$ un espace $k$-analytique et soit $L$ une extension complète de $k$. Il existe un morphisme naturel $(X_L)\duit\to (X\duit)_L$ qui est une immersion fermée, et un isomorphisme si et seulement si $(X\duit)_L$ est réduit.}

\medskip
{\em Démonstration.} Comme l'espace $(X_L)\duit$ est réduit, le morphisme composé $(X_L)\duit\to X_L\to X$ se factorise par $X\duit$ et donc par $(X\duit)_L$ ;  les flèches $(X_L)\duit\to X_L$ et $(X\duit)_L\to X_L$ sont des immersions fermées ; il en résulte que $(X_L)\duit\to (X\duit)_L$ est une immersion fermée.

\medskip
Soit $Z$ un espace analytique réduit, et soit $Z\to X_L$ un morphisme. Comme $Z$ est réduit, la composée $Z\to X_L\to X$ se factorise d'une unique manière par $X\duit$. L'ensemble des $k$-morphismes de $Z$ vers $X\duit$ étant en bijection naturelle avec celui des $L$-morphismes de $Z$ vers $(X\duit)_L$, la flèche $Z\to X_L$ se factorise de manière unique par $(X\duit)_L$. Si celui-ci est réduit, il s'identifie donc à $(X_L)\duit$ ; la réciproque est triviale.~$\Box$

\subsection*{Autour des composantes irréductibles d'un espace affinoïde}

\medskip
\deux{ztequiv} {\bf Proposition.} {\em Soit $V$ un espace affinoïde et soit $W$ un domaine affinoïde de $V$ ; soit $Z$ un fermé de Zariski de $V$ et soit $T$ un fermé de Zariski de $W$. Les assertions suivantes sont équivalentes :

\medskip
$\alpha)$ $Z$ est une composante irréductible de $V$ et $T$ une composante irréductible de $W$ incluse dans $Z$ ; 

\medskip
$\beta)$ $Z$ est une composante irréductible de $V$ et $T$ une composante irréductible de $Z\cap W$ ; 

\medskip
$\gamma)$ $Z$ est une composante irréductible de $V$ et $T$ est inclus dans $Z$, irréductible et de même dimension que $Z$ ; 

\medskip
$\delta)$ $Z$ est égal à l'adhérence $\overline{T}^{V\zar}$ de $T$ dans $V$ pour la topologie de Zariski de ce dernier, et $T$ est une composante irréductible de $W$.}

\medskip
{\em Démonstration.} Supposons $\alpha)$ vraie et montrons $\beta)$. Comme $T$ est une composante irréductible de $W$ incluse dans $Z$, c'est un fermé de Zariski irréductible maximal de $W$ inclus dans $Z$ ; c'est {\em  a fortiori} un fermé irréductible maximal de $Z\cap W$, et donc une composante irréductible de ce dernier. 

\medskip
Supposons $\beta)$ vraie et montrons $\gamma)$. Soit $d$ la dimension de $Z$. Comme $Z$ est irréductible, il est purement de dimension $d$, et $W\cap Z$ est donc purement de dimension $d$ ; ses composantes irréductibles sont donc toutes de dimension $d$ ; c'est en particulier le cas de $T$. 

\medskip
Supposons $\gamma)$ vraie et montrons $\delta)$. Soit $Y$ un fermé de Zariski de $V$ contenant $T$. Le fermé de Zariski $Y\cap Z$ de $V$ contient $T$, donc est au moins de dimension $d$. Comme c'est un fermé de $Z$, lequel est lui-même irréductible de dimension $d$, il coïncide avec $Z$. L'on a donc bien $\overline{T}^{V\zar}=Z$. Si $S$ est un fermé de Zariski irréductible de $W$ contenant $T$, alors $\overline{S}^{V\zar}$ est un fermé irréductible de $V$ contenant $Z$ ; ce dernier étant une composante irréductible de $V$, l'on a $\overline{S}^{V\zar}=Z$. La dimension de $S$ est donc majorée par $d$, ce qui implique que $T=S$ ; ainsi, $T$ est une composante irréductible de $W$.  

\medskip
Supposons $\delta)$ vraie et montrons $\alpha)$. L'irréductibilité de $T$ entraîne celle de $Z$. Soit $Y$ une composante irréductible de $V$ contenant $Z$. Comme $T$ est une composante irréductible de $W$ contenue dans $Y$, l'implication $\alpha)\Rightarrow \gamma)$ déjà établie assure que $\dim{}T=\dim{}Y$ ; d'autre part $T\subset Z$ et $Z\subset Y$, d'où l'encadrement $\dim{}T \leq \dim{}Z\leq \dim{}Y$ ; en conséquence $\dim{}Z=\dim {}Y$. Puisque $Y$ est irréductible, $Y=Z$.~$\Box$

\medskip
\deux{corunion} {\bf Corollaire.} {\em Soit $V$ un espace affinoïde et soit $W$ un domaine affinoïde de $V$. Si $Z$ est une composante irréductible de $V$, alors $Z\cap W$ est réunion de composantes irréductibles de $W$ ; si $T$ est une composante irréductible de $W$, elle est contenue dans une seule composante irréductible de $V$ qui est égale à $\overline T^{V\zar}$.}

\medskip
{\em Démonstration.} La première assertion résulte de l'implication $\beta)\Rightarrow \alpha)$ dans la proposition ci-dessus, et la seconde de $\alpha)\Rightarrow \delta)$. ~$\Box$ 

\medskip
\deux{commentvideun} {\bf Commentaire.} Notons que dans le corollaire ci-dessus, $Z\cap W$ peut très bien être égal à la {\em réunion vide} de composantes irréductibles de $W$.

\medskip
\deux{propvd} {\bf Proposition.} {\em Soit $X$ un espace analytique, soit $Y$ un fermé de Zariski de $X$ et soit $d$ un entier. Il existe deux fermés de Zariski $Y_d^+$ et $Y_d^-$ de $X$ caractérisés par le fait que pour tout domaine affinoïde $V$ de $X$, l'intersection $Y_d^+\cap V$ (resp. $Y_d^-\cap V$) est la réunion des composantes irréductibles de $V$ contenues dans $Y$ et de dimension $d$ (resp. non contenues dans $Y$ ou de dimension différente de $d$).} 

\medskip
{\em Démonstration.} Pour tout domaine affinoïde $V$ de $X$, notons ${\cal C}(V)$ l'ensemble des composantes irréductibles de $V$, et ${\cal C}^+(V)$ l'ensemble des composantes irréductibles de $V$ qui sont contenues dans $Y$ et de dimension $d$ ; posons ${\cal E}^+(V)=\bigcup\limits_{T\in {\cal C}^+(V)}T$ et  ${\cal E}^-(V)=\bigcup\limits_{T\in {\cal C}(V)-{\cal C}^+(V)}T$. 

\medskip
Soit $V$ un domaine affinoïde de $X$, et soit $W$ un domaine affinoïde de $V$. Soit $Z$ une composante irréductible de $V$ et soit $T$ une composante irréductible de $W$ contenue dans $Z$. Le couple $(Z,T)$ satisfaisant la condition $\alpha)$ de la proposition \ref{COMP}.\ref{ztequiv}, il satisfait également $\gamma)$ et $\delta)$ ; on en déduit que $T\in{\cal C}^+(W)$ si et seulement si $Z\in{\cal C}^+(V)$. Il découle alors du corollaire \ref{COMP}.\ref{corunion} que ${\cal E}^+(W)={\cal E}^+(V)\cap W$ et ${\cal E}^-(W)={\cal E}^-(V)\cap W$. 

\medskip
Soit $Y_d^+$ (resp. $Y_d^-$) la réunion des ${\cal E}^+(V)$ (resp. des ${\cal E}^-(V)$) où $V$ parcourt l'ensemble des domaines affinoïdes de $X$. En vertu de ce qui précède, l'on a pour tout domaine affinoïde $V$ de $X$ les égalités $Y_d^+\cap V={\cal E}^+(V)$ et $Y_d^-\cap V={\cal E}^-(V)$ ; puisque ${\cal E}^+(V)$ et ${\cal E}^-(V)$ sont des fermés de Zariski de $V$ pour tout domaine affinoïde $V$ de de $X$, la proposition \ref{COMP}.\ref{zargloc} assure que $Y_d^+$ et $Y_d^-$ sont des fermés de Zariski de $X$.~$\Box$ 

\medskip
\deux{propadhtzar} {\bf Proposition.} {\em Soit $X$ un espace analytique, soit $V$ un domaine affinoïde non vide de $X$ et soit $T$ une composante irréductible de $V$ dont on note $d$ la dimension.

\medskip
$i)$  Pour tout domaine affinoïde $W$ de $X$, l'intersection $\overline{T}^{X\zar}\cap W$ est réunion de composantes irréductibles de $W$ de dimension $d$ ; 

\medskip
$ii)$ $\overline{T}^{X\zar}$ est purement de dimension $d$ ;

\medskip
$iii)$ $\overline{T}^{X\zar}$ est l'unique fermé de Zariski irréductible de $X$ qui contient $T$ ; si $X$ est irréductible, on a donc $\overline{T}^{X\zar}=X$.}

\medskip
{\em Démonstration.} Posons $Y=\overline{T}^{X\zar}$. Il est clair que $Y$ est un fermé de Zariski irréductible de $X$ ; on utilise les notations $Y_d^+$ et $Y_d^-$ de la proposition \ref{COMP}.\ref{propvd} ci-dessus ; on a clairement $X=Y_d^+\cup Y_d^-$.

\medskip
Soit $Z$ un fermé de Zariski irréductible de $X$ contenant $Y$. L'irréductibilité de $Z$ et l'égalité $Z=(Z\cap Y_d^+)\cup(Z\cap Y_d^-)$ impliquent que $Z\subset Y_d^+$ ou $Z\subset Y_d^-$. Comme $T\subset Z$, l'inclusion $Z\subset Y_d^-$ est impossible ; dès lors, $Z\subset Y_d^+$.

\medskip
Appliquons en particulier ceci lorsque $Z=Y$. On a donc $Y\subset Y_d^+$ ; comme il est clair que $Y_d^+\subset Y$, on a $Y=Y_d^+$. Cela signifie exactement que pour tout domaine affinoïde $W$ de $X$, l'intersection $Y\cap W$ est réunion de composantes irréductibles de $W$ de dimension $d$, ce qui prouve $i)$ ; l'assertion $ii)$ en est une conséquence triviale.

\medskip
Revenons à un fermé de Zariski irréductible quelconque $Z$ qui contient $Y$. L'on a $Z\subset Y_d^+$ et $Y_d^+\subset Y$ ; de ce fait, $Z=Y$, ce qui prouve $iii)$.~$\Box$ 

\medskip
\deux{commentvidedeux} {\bf Commentaire.} Notons que dans l'énoncé $i)$ de la proposition ci-dessus, $\overline{T}^{X\zar}\cap W$ peut très bien être la {\em réunion vide} de composantes irréductibles de $W$. 

\medskip
\deux{dimespirr} {\bf Corollaire.} {\em Soit $X$ un espace analytique irréductible. Il existe $d\in \NN$ tel que $X$ soit purement de dimension $d$ ; si $Z$ est un fermé de Zariski strict de $X$, alors $\dim{} Z<d$ et $Z$ ne contient aucun domaine analytique non vide de $X$ ; il est en particulier d'intérieur vide.} 

\medskip
{\em Démonstration.} Faisons une remarque préliminaire : si $T$ est une composante irréductible d'un domaine a	ffinoïde $V$ de $X$ si l'on note $d$ la dimension de $T$, alors $\overline{T}^{X\zar}$ est en vertu de la proposition ci-dessus purement de dimension $d$ et coïncide avec $X$ puisque celui-ci est irréductible. 

\medskip
La première assertion provient alors du fait que $X$ étant irréductible, il est non vide ; il possède donc un domaine affinoïde non vide $V$, lequel a au moins une composante irréductible. 

\medskip
Établissons la seconde assertion. Soit $Z$ un fermé de Zariski strict de $X$. Soit $V$ un domaine affinoïde de $X$. Il résulte de la remarque préliminaire que $Z$ ne peut contenir aucune composante irréductible de $V$. Comme $X$ est purement de dimension $d$, son domaine affinoïde $V$ est purement de dimension $d$ et $Z\cap V$ est donc de dimension strictement inférieure à $d$. Ceci valant pour tout domaine affinoïde $V$ de $X$, la dimension de $Z$ est strictement inférieure à $d$.

\medskip
Soit $V$ un domaine analytique non vide de $X$. Comme $X$ est purement de dimension $k$-analytique $d$, la dimension de $V$ vaut $d$ ; 
comme $\dim {}Z<d$, on a $\dim {}Z\cap V<d$ ; par conséquent, $V$ n'est pas inclus dans $Z$.~$\Box$ 

\subsection*{Ensembles G-localement finis de fermés de Zariski et composantes irréductibles}

\medskip
\deux{defglocfin} On dira qu'un ensemble $\cal E$ de fermés de Zariski d'un espace analytique $X$ est {\em G-localement fini} si  tout domaine affinoïde $V$ de $X$ ne rencontre qu'un nombre fini d'éléments de $\cal E$. Notons que sous cette hypothèse $(\bigcup\limits_{E\in\cal E}E)\cap V$ est pour tout domaine affinoïde $V$ de $X$ une réunion {\em finie} de fermés de Zariski de $V$, et est donc fermé de Zariski de $V$ ; la proposition  \ref{COMP}.\ref{zargloc} assure alors que $\bigcup\limits_{E\in\cal E}E$ est un fermé de Zariski de $X$. 

\medskip
\deux{testglocfin} Soit $f: Y\to X$ un morphisme d'espaces analytiques et soit $\cal E$ un ensemble G-localement fini de fermés de Zariski de $X$. Si $V$ est un domaine affinoïde de $Y$ alors $f(V)$ est une partie {\em quasi}-compacte de $X$ (qui n'est pas supposé topologiquement séparé), et est en conséquence contenue dans une réunion finie de domaines affinoïdes de $X$ ; il en résulte que $f(V)$ ne rencontre qu'un nombre fini d'éléments de $\cal E$, et $\{f^{-1}(E)\}_{E\in \cal E}$ constitue par conséquent un ensemble G-localement fini de fermés de Zariski de $Y$. Notons deux cas particuliers importants : si $Y$ est un domaine analytique ou un sous-espace analytique fermé de $X$, alors $\{E\cap Y\}_{E\in \cal E}$ est un ensemble G-localement fini de fermés de Zariski de $Y$. 

\medskip
\deux{irrlocfin} {\bf Lemme.} {\em Soit $X$ un espace analytique irréductible et soit $\cal E$ un ensemble G-localement fini de fermés de Zariski de $X$ tels que $X=\bigcup\limits_{E\in \cal E} E$. On a alors $X\in {\cal E}$. }

\medskip
{\em Démonstration.} Comme $X$ est irréductible il est non vide, et possède donc un domaine affinoïde non vide $V$. L'ensemble $\cal E$ est G-localement fini ; dès lors, l'ensemble $\{E\cap V\}_{E\in \cal E}$ constitue un recouvrement {\em fini} de $V$ par des fermés de Zariski ; ceci entraîne, $V$ étant non vide, l'existence d'une composante irréductible $T$ de $V$ et de $E\in \cal E$ tels que $T\subset E$. Comme $X$ est irréductible, la proposition \ref{COMP}.\ref{propadhtzar} assure que $\overline{T}^{X\zar}=X$, d'où il découle que $E=X$.~$\Box$

 \medskip
 \deux{uniccompirr} {\bf Corollaire.} {\em Soit $X$ un espace analytique, et soit ${\cal E}$ un ensemble G-localement fini de fermés de Zariski irréductibles de $X$ deux à deux non comparables pour l'inclusion ; supposons que $X=\bigcup\limits_{E\in{\cal E}}E$. Tout fermé de Zariski irréductible de $X$ est alors contenu dans un élément $E$ de $\cal E$, et l'ensemble $\cal E$ est exactement l'ensemble des  fermés de Zariski irréductibles maximaux de $X$.}

 \medskip
{\em Démonstration.} Soit $Z$ un fermé de Zariski irréductible de $X$ ; munissons-le (par exemple) de sa structure réduite. Il résulte du \ref{COMP}.\ref{testglocfin} que $\{E\cap Z\}_{E\in \cal E}$ est un ensemble G-localement fini de fermés de Zariski de $Z$, et il est immédiat que $Z=\bigcup\limits_{E\in \cal E} E\cap Z$ ; d'après le lemme \ref{COMP}.\ref{irrlocfin}, il existe alors $E\in{\cal E}$ tel que $Z=E\cap Z$, c'est-à-dire tel que $Z\subset E$ ; on en déduit aussitôt que tout fermé de Zariski irréductible maximal de $X$ appartient à $\cal E$.

\medskip
Soit $T\in{\cal E}$ et soit $Z$ un fermé de Zariski irréductible de $X$ contenant $T$ ; d'après ce qui précède, il existe $E\in {\cal E}$ tel que $Z\subset E$ ; on a donc $T\subset E$, et comme les éléments de $\cal E$ sont deux à deux non comparables pour l'inclusion, $T=E$ ; par conséquent, $Z=T$ et $T$ est bien un fermé de Zariski irréductible maximal de $Y$.~$\Box$

 \medskip
 \deux{existcompirr} {\bf Théorème (description des composantes irréductibles).} {\em Soit $X$ un espace analytique. Soit $\irr X$ l'ensemble des parties de $X$ qui sont de la forme $\overline{T}^{X\zar}$, où $T$ est une composante irréductible d'un domaine affinoïde de $X$ ; l'ensemble $\irr X$ est constitué de fermés de Zariski irréductibles deux à deux non comparables pour l'inclusion  et est G-localement fini ; l'on a $\bigcup\limits_{E\in\irr X}E=X$. Tout fermé de Zariski de $X$ est contenu dans un élément de $\irr X$, et $\irr X$ est exactement l'ensemble des fermés de Zariski irréductibles maximaux de $X$.} 
 
 \medskip
 {\em Démonstration.} Il résulte de la définition de $\irr X$ que ses éléments sont des fermés irréductibles de $X$ dont la réunion est égale à $X$.
 
\medskip
 Soit $E\in \irr X$ et soit $V$ un domaine affinoïde de $X$. La proposition \ref{COMP}.\ref{propadhtzar} $iii)$ assure que $E\cap V$ est réunion de composantes irréductibles de $V$ et garantit que si $T$ est l'une de celles-ci, alors $E=\overline{T}^{X\zar}$ ; on en déduit que $V$ ne rencontre qu'un nombre fini d'éléments de $\irr X$ ; ce dernier est donc bien G-localement fini.
 
 \medskip
 Soient $E$ et $E'$ deux éléments de $\irr X$ tels que $E\subset E'$. Choisissons un domaine affinoïde $V$ de $X$ et une composante irréductible $T$ de $V$ tels que $E=\overline{T}^{X\zar}$ ; comme $E'$ est un fermé de Zariski irréductible contenant $E$ et donc $T$, il coïncide avec $\overline{T}^{X\zar}$ d'après la proposition \ref{COMP}.\ref{propadhtzar} ; autrement dit $E=E'$, et les éléments de $\irr X$ sont donc deux à deux non comparables pour l'inclusion. Les deux dernières assertions découlent du corollaire  \ref{COMP}.\ref{uniccompirr}. ~$\Box$ 
 
 \medskip
 \deux{defcompirr} Les éléments de l'ensemble $\irr X$ défini dans l'énoncé du théorème ci-dessus seront appelés les {\em composantes irréductibles} de $X$ ; notons que si $X$ est affinoïde, on retrouve bien les composantes irréductibles usuelles. Si $Y$ est une composante irréductible de $X$, la proposition \ref{COMP}.\ref{propadhtzar} assure que pour tout domaine affinoïde $V$ de $X$, l'intersection $V\cap Y$ est réunion de composantes irréductibles de $V$. 

\medskip
Il résulte de la caractérisation topologique des composantes irréductibles (comme fermés de Zariski irréductibles maximaux) que si $Y$ est un fermé de Zariski de $X$, alors l'ensemble des composantes irréductibles de $Y$ ne dépend pas de la structure de sous-espace analytique fermé dont on le munit ; on se permettra donc d'évoquer ces composantes sans fixer une telle structure. Par exemple, soit ${\cal E}$ un sous-ensemble de $\irr X$ et soit $Y=\bigcup\limits_{E\in{\cal E}}E$ ; on sait d'après le \ref{COMP}.\ref{defglocfin} que $Y$ est un fermé de Zariski de $X$. Il découle alors du \ref{COMP}.\ref{testglocfin}, du corollaire \ref{COMP}.\ref{uniccompirr} et de la définition même des composantes irréductibles que ${\cal E}$ est l'ensemble des composantes irréductibles de $Y$.

\medskip
\deux{dimloccasgen} {\bf Lemme.} {\em Soit $X$ un espace analytique. Si $x$ est un point de $X$, alors $\dim x X$ est le maximum des dimensions des composantes irréductibles de $X$ qui contiennent $x$.}

\medskip
{\em Démonstration.} Soient $X_1,\ldots,X_n$ les composantes irréductibles de $X$ contenant $x$, soient $d_1,\ldots,d_n$ leurs dimensions  respectives et soit $V$ un domaine affinoïde de $X$ contenant $x$.  Pour tout $i$, la proposition \ref{COMP}.\ref{propadhtzar} assure que $X_i\cap V$ est réunion de composantes irréductibles de $V$ de dimension $d_i$ ; par ailleurs, toute composante irréductible $T$ de $V$ contenant $x$ est contenue dans $\overline{T}^{X\zar}$, qui est une composante irréductible de $X$ et est donc l'une des $X_i$. On en déduit que l'ensemble des dimensions des composantes irréductibles de $V$ qui contiennent $x$ est égal à $\{d_1,\ldots,d_n\}$ ; par conséquent, $\dim x X=\dim xV=\max_i d_i$.~$\Box$ 

\medskip
Grâce à la proposition \ref{COMP}.\ref{dimespirr} et au lemme \ref{COMP}.\ref{dimloccasgen}, on peut généraliser la proposition \ref{COMP}.\ref{ztequiv} et le corollaire \ref{COMP}.\ref{corunion} en conservant {\em mutatis mutandis} les mêmes démonstrations.

\medskip
\deux{ztequivgen} {\bf Proposition.} {\em Soit $V$ un espace analytique et soit $W$ un domaine analytique de $V$ ; soit $Z$ un fermé de Zariski de $V$ et soit $T$ un fermé de Zariski de $W$. Les assertions suivantes sont équivalentes :

\medskip
$\alpha)$ $Z$ est une composante irréductible de $V$ et $T$ une composante irréductible de $W$ incluse dans $Z$ ; 
 
\medskip
$\beta)$ $Z$ est une composante irréductible de $V$ et $T$ une composante irréductible de $Z\cap W$ ; 

\medskip
$\gamma)$ $Z$ est une composante irréductible de $V$ et $T$ est inclus dans $Z$, irréductible et de même dimension que $Z$ ; 

\medskip
$\delta)$ $Z$ est égal à l'adhérence $\overline{T}^{V\zar}$ de $T$ dans $V$ pour la topologie de Zariski de ce dernier, et $T$ est une composante irréductible de $W$.~$\Box$}

\medskip
\deux{coruniongen} {\bf Corollaire.} {\em Soit $V$ un espace analytique et soit $W$ un domaine analytique de $V$. Si $Z$ est une composante irréductible de $V$, alors $Z\cap W$ est réunion de composantes irréductibles de $W$ ; si $T$ est une composante irréductible de $W$, elle est contenue dans une seule composante irréductible de $V$ qui est égale à $\overline T^{V\zar}$.}~$\Box$

\medskip
\deux{commentvidetrois} {\bf Commentaire.} Notons que dans le corollaire ci-dessus, $Z\cap W$ peut parfaitement être égal à la {\em réunion vide} de composantes irréductibles de $W$ ; il peut tout aussi bien être une réunion {\em infinie} de telles composantes. 

\medskip
\deux{compouvdense} {\bf Lemme.} {\em Soit $X$ un espace analytique et soit $U$ un ouvert de Zariski de $X$. Les propriétés suivantes sont équivalentes :

\medskip
$i)$ $U$ est dense pour la topologie de Zariski de $X$ ;

\medskip
$ii)$ $U$ rencontre chacune des composantes irréductibles de $X$ ;

\medskip
$iii)$ $U$ est dense pour la topologie usuelle de $X$.

}

\medskip
{\em Démonstration.} Supposons $i)$  vraie, et soit $T$ une composante irréductible de $X$. Soit $Z$ la réunion des autres composantes irréductibles de $X$ ; par construction, $X-Z$ est un ouvert de Zariski non vide de $X$, qui rencontre donc $U$ ; comme $(X-Z)\subset T$, l'intersection de $U$ et $T$ est non vide, ce qui montre $ii)$. 

\medskip
Supposons $ii)$ vraie, et soit $\Omega$ un ouvert de $X$ contenu dans le complémentaire $F$ de $U$. Supposons que $\Omega$ est non vide, et soit $Y$ une composante irréductible de $X$ rencontrant $\Omega$ ; soit $d$ la dimension de $Y$. Comme $\Omega\cap Y$ est un ouvert non vide de $Y$, laquelle est purement de dimension $d$ (cor. \ref{COMP}.\ref{dimespirr}), il est lui-même de dimension $d$. En conséquence, $\dim {} F\cap Y \geq d$, ce qui implique que $F\cap Y=Y$ (cor. \ref{COMP}.\ref{dimespirr}) et contredit le fait que $U$ rencontre $Y$ ; dès lors $\Omega$ est vide et $iii)$ est établie.

\medskip
Supposons $iii)$ vraie et soit $V$ un ouvert de Zariski non vide de $X$ ; il existe une composante irréductible $T$ de de $X$ rencontrant $V$. Soit $Z$ la réunion des autres composantes irréductibles de $X$. L'ouvert non vide $X-Z$ de $X$ est contenu dans $T$ et rencontre $U$ puisque celui-ci est dense. Par conséquent, $U\cap T$ et $V\cap T$ sont deux ouverts de Zariski non vide du fermé irréductible $T$ ; on en déduit que $U\cap V\neq \emptyset$, ce qui prouve $i)$.~$\Box$ 

\medskip
\deux{irrcompextscal} {\bf Lemme.} {\em Soit $X$ un espace $k$-analytique et soit $(X_i)$ la famille de ses composantes irréductibles. Soit $L$ une extension complète de $k$. Pour tout $i$, notons $(X_{i,j})$ la famille des composantes irréductibles de $X_{i,L}$. 

\medskip
$i)$ Pour tout $i$ et tout $j$, l'on a $\dim {} X_{i,j}=\dim {} X_i$ ;

\medskip
$ii)$ les $X_{i,j}$ sont deux à deux non comparables pour l'inclusion et sont les composantes irréductibles de $X_L$.}

\medskip
{\em Démonstration.} Le corollaire \ref{COMP}.\ref{dimespirr} assure que $X_i$ est purement de dimension $d_i$ pour un certain entier $d_i$. L'espace $X_{i,L}$ est donc purement de dimension $d_i$. En vertu du lemme \ref{COMP}.\ref{dimloccasgen}, les composantes irréductibles de $X_{i,L}$ sont toutes de dimension $d_i$ ; l'assertion $i)$ est ainsi démontrée.

\medskip
Soient $(i,j)$ et $(i',j')$ deux couples d'indices, et supposons que $X_{i,j}\subset X_{i',j'}$. La composante $X_{i,j}$ est alors incluse dans $(X_i\cap X_{i'})_L$. Si $i\neq i'$ alors $(X_i\cap X_{i'})$ est un fermé de Zariski strict de $X_i$, donc est de dimension strictement inférieure à $d_i$ par le corollaire  \ref{COMP}.\ref{dimespirr} ; en conséquence $\dim {} (X_i\cap X_{i'})_L<d_i$, ce qui contredit l'égalité $\dim {}X_{i,j}=d_i$. Dès lors $i=i'$, et $j=j'$ puisque les $X_{i,j}$ à $i$ fixé sont {\em par définition} deux à deux non comparables pour l'inclusion. 

\medskip
Les $X_{i,j}$ forment un ensemble de fermés de Zariski de $X_L$ qui le recouvrent, sont irréductibles et sont deux à deux non comparables pour l'inclusion. Cet ensemble est par ailleurs G-localement fini puisque celui des $X_{i,L}$ l'est en vertu du \ref{COMP}.\ref{testglocfin}. Les $X_{i,j}$ sont donc les composantes irréductibles de $X_L$.~$\Box$ 

\medskip
\deux{ouvdenseextscal} {\bf Corollaire.} {\em Soit $X$ un espace $k$-analytique et soit $U$ un ouvert de Zariski de $X$ qui est dense dans $X$. Pour tout extension complète $L$ de $k$, l'ouvert $U_L$ est dense dans $X_L$.} 

\medskip
{\em Démonstration.} D'après le lemme \ref{COMP}.\ref{compouvdense}, il suffit de démontrer que pour toute composante irréductible $T$ de $X_L$ l'ouvert $U_L$ rencontre $T$. Soit donc $T$ une telle composante ; notons $F$ le complémentaire de $U$ dans $X$. Il existe en vertu du lemme \ref{COMP}.\ref{irrcompextscal} ci-dessus une composante irréductible $S$ de $X$ tel que $T$ soit une composante irréductible de $S_L$. Si $d$ désigne la dimension de $S$, alors $\dim{}  T=d$. Comme $U\cap S\neq\emptyset$ par l'hypothèse faite sur $U$ et le lemme \ref{COMP}.\ref{compouvdense}, on a l'inégalité $\dim {} F\cap S< d$. Dès lors, $\dim {} F_L\cap S_L< d$, et $F_L\cap S_L$ ne peut donc contenir $T$ ; autrement dit, $U_L$ rencontre $T$, ce qui achève la démonstration.~$\Box$ 
  
\medskip
{\bf Tous les résultats que l'on a établis dans cette section seront désormais utilisés librement, sans références précises.} 

\section{La normalisation d'un espace analytique}\label{NOR}

\subsection*{Les morphismes quasi-dominants} 

\setcounter{cpt}{0}
\medskip
\deux{pseudom} {\bf Définition.} Soit $f: Y\to X$ un morphisme d'espaces analytiques. On dira que $f$ est {\em quasi-dominant} si pour tout couple $(V,U)$ formé d'un domaine analytique $V$ de $Y$ et d'un domaine analytique $U$ de $X$ contenant $f(V)$, et pour toute composante irréductible $T$ de $V$, $\overline{f(T)}^{U\zar}$ est une composante irréductible de $U$.

\medskip
Mentionnons tout de suite quelques propriétés élémentaires des morphismes quasi-dominants.

\setcounter{cptbis}{0}
\medskip
\trois{pseudomtauto} Si $Y$ est un domaine analytique d'un espace analytique $X$, la flèche naturelle $Y\hookrightarrow X$ est quasi-dominante ; en particulier $\mathsf{Id}_X$ et $\emptyset\hookrightarrow X$ sont quasi-dominants ; la composée de deux morphismes quasi-dominants est un morphisme quasi-dominant. 

\medskip
\trois{irrfonct} Si $\mathsf C$ désigne la catégorie dont les objets sont les espaces analytiques et dont les flèches les morphismes quasi-dominants alors $X\mapsto \irr X$ est de manière évidente un foncteur $\mathsf C\to\mathsf {Ens}$ : si $f: Y\to X$ est une flèche de $\mathsf C$, l'application $\irr Y\to \irr X$ correspondante associe à une composante irréductible $T$ de $Y$ l'adhérence $\overline{f(T)}^{X\zar}$ ou, ce qui revient au même, l'unique composante irréductible de $X$ contenant $f(T)$. 

\medskip
\trois{domanalqd} Soit $f: Y\to X$ un morphisme d'espaces analytiques, et soit $U$ un domaine analytique de $X$ tel que $f(Y)\subset U$. Soit $g: Y\to U$ le morphisme induit par $f$. Si $f$ est quasi-dominant il est immédiat que $g$ est quasi-dominant ; réciproquement si $g$ est quasi-dominant alors $f$ est quasi-dominant, puisque c'est la composée de $g$ et $U\hookrightarrow X$.

\medskip
\trois{qdomprofib} Soit $f: Y\to X$ un morphisme d'espaces analytiques, et soit $U$ un domaine analytique de $X$. Si $f$ est quasi-dominant, alors $Y\times_X U\to U$ est quasi-dominant : en effet, $Y\times_X U$ est un domaine analytique de $Y$ ; par conséquent, $Y\times_X U\hookrightarrow Y$ est quasi-dominant ; le morphisme composé $Y\times_X U\hookrightarrow Y\to X$ l'est également. Il se factorise par $U$, et il découle alors du \ref{NOR}.\ref{pseudom}.\ref{domanalqd} ci-dessus que $Y\times_X U\to U$ est quasi-dominant. 

\medskip
\deux{qdomestloc} {\bf Vérification locale de la quasi-dominance.} En pratique, on peut parfois se contenter de vérifier la condition qui définit la quasi-dominance en se limitant à {\em certains} domaines analytiques bien choisis de la source ou du but ; donnons quelques exemples. 

\setcounter{cptbis}{0}
\medskip
\trois{qdomfamille} Soit $f:Y\to X$ un morphisme d'espaces analytiques, et soit $\cal V$ (resp. $\cal U$) un ensemble de domaines analytiques de $Y$ (resp. $X$) possédant la propriété suivante : {\em tout domaine analytique de $Y$ (resp. $X$) est réunion d'éléments de $\cal V$ (resp. $\cal U$).} Le morphisme $f$ est alors quasi-dominant si et seulement si pour tout $V\in \cal V$, pour tout $U\in \cal U$ contenant $f(V)$, et pour toute composante irréductible $T$ de $V$, $\overline{f(T)}^{U\zar}$ est une composante irréductible de $U$. 

\medskip
La condition est en effet clairement nécessaire. Vérifions qu'elle est suffisante ; on la suppose donc satisfaite. Soit $V$ un domaine analytique de $Y$, soit $T$ une composante irréductible de $V$, et soit $U$ un domaine analytique de $X$ tel que $f(V)\subset U$. Par définition de $\cal U$, il existe un élément $U'$ de $\cal U$ inclus dans $U$ et tel que $f^{-1}(U')$ rencontre $T$ ; par définition de $\cal V$, il existe un élément $V'$ de $\cal V$ contenu dans $f^{-1}(U')$ et rencontrant $T$. Soit $S$ une composante irréductible de $T\cap V'$ ; on a $T=\overline{S}^{V\zar}$. On en déduit que $\overline{f(T)}^{U\zar}=\overline{f(S)}^{U\zar}=\overline{\overline{f(S)}^{U'\zar}}^{U\zar}$. Par hypothèse, $\overline{f(S)}^{U'\zar}$ est une composante irréductible de $U'$ ; en conséquence, $\overline{\overline{f(S)}^{U'\zar}}^{U\zar}$ est une composante irréductible de $U$ et $f$ est quasi-dominant. 

\medskip
\trois{qdomglocs} Soit $f: Y\to X$ un morphisme d'espaces analytiques, et soit $\cal V$ un ensemble de domaines analytiques de $Y$ tel que tout domaine analytique de $Y$ soit réunion d'éléments de $\cal V$. Le morphisme $f$ est alors quasi-dominant si et seulement si $V\to X$ est quasi-dominant pour tout $V\in \cal V$ : cette condition est en effet trivialement nécessaire, et est suffisante d'après le  \ref{NOR}.\ref{qdomestloc}.\ref{qdomfamille} ci-dessus.

\medskip
\trois{qdomglocb} Soit $f: Y\to X$ un morphisme d'espaces analytiques et soit ${\cal U}$  un ensemble de domaines analytiques de $X$ tel que $X$ soit réunion d'éléments de $\cal U$. Le morphisme $f$ est alors quasi-dominant si et seulement si $Y\times_X U\to U$ est quasi-dominant pour tout $U\in \cal U$. La condition est en effet nécessaire d'après le \ref{NOR}.\ref{pseudom}.\ref{qdomprofib}. Pour voir qu'elle est suffisante, supposons-la satisfaite. Soit $\cal V$ la famille des domaines analytiques $V$ de $Y$ tels qu'il existe $U\in \cal U$ par lequel $V\to X$ se factorise. Comme $X$ est réunion d'éléments de $\cal U$, tout domaine analytique de $Y$ est réunion d'éléments de $\cal V$ (ses intersections avec les images réciproques des différents éléments de $\cal U$, par exemple) ; de plus, si $V\in \cal V$ et si $U$ est un élément de $\cal U$ tel que $V\to Y$ se factorise par $U$, alors $V\to U$ est quasi-dominant, puisque $Y\times_X U\to U$ l'est par hypothèse ; on en déduit que $V \to X$ est quasi-dominant. Le \ref{NOR}.\ref{qdomestloc}.\ref{qdomglocs} ci-dessus permet de conclure que $f$ est quasi-dominant. 

\medskip
\deux{qdometfermzar} {\bf Quasi-dominance et immersions fermées.} Nous allons appliquer ce qui précède à l'étude des immersions fermées dont l'image est réunion de composantes irréductibles du but. Pour ce faire, nous aurons besoin du résultat qui suit sur les domaines analytiques d'un sous-espace analytique fermé d'un espace analytique. 

\setcounter{cptbis}{0}
\medskip
\trois{gerrgrauert} Soit $X$ un espace analytique, soit $Z$ un sous-espace analytique fermé de $X$ et soit $V$ un domaine analytique de $Z$. Il existe alors un G-recouvrement de $V$ par des domaines analytiques de la forme $W\cap Z$, où $W$ est un domaine analytique de $X$. La question est en effet G-locale sur $V$, ce qui permet de supposer que $V$ est affinoïde et qu'il est contenu dans un domaine affinoïde $U$ de $X$ ; mais l'assertion cherchée est alors une conséquence directe de la variante en théorie de Berkovich du théorème de Gerritzen-Grauert, variante due à Temkin (\cite{tmk}, th. 3.1).

\medskip
\trois{fermzarqdom} Soient $Y$ et $X$ deux espaces analytiques. Soit $\cal I$ un faisceau cohérent d'idéaux sur $X$ dont le lieu des zéros $Z$ est réunion (éventuellement vide) de composantes irréductibles de $X$ ; on munit $Z$ de la structure de sous-espace analytique fermé définie par $\cal I$. Soit $f$ un morphisme de $Y$ vers $Z$. Alors $f$ est quasi-dominant si et seulement si le morphisme composé $g:Y\to Z\hookrightarrow X$ est quasi-dominant (ceci entraîne notamment, comme on le voit en prenant $f$ égal à l'identité de $Z$, que $Z\hookrightarrow X$ est quasi-dominant). 

\medskip
En effet, supposons $f$ quasi-dominant, soit $V$ un domaine analytique de $Y$, soit $U$ un domaine analytique de $X$ tel que $g(V)\subset U$, et soit $T$ une composante irréductible de $V$.  Le morphisme $V\to U$ se factorise canoniquement par la flèche $V\to U\cap Z$ induite par $f$. Comme $f$ est quasi-dominant, $\overline{f(T)}^{(U\cap Z)\zar}$ est une composante irréductible de $U\cap Z$. Comme $Z$ est réunion de composantes irréductibles de $X$, toute composante irréductible de $U\cap Z$ est une composante irréductible de $U$ ; dès lors, $\overline{g(T)}^{U\zar}=\overline{f(T)}^{(U\cap Z)\zar}$ est une composante irréductible de $U$, et $g$ est quasi-dominant.

\medskip
Réciproquement, supposons $g$ quasi-dominant. Soit $V$ un domaine analytique de $Y$, et soit $U$ un domaine analytique de $X$ tel que $f(V)\subset U\cap Z$.  Soit $T$ une composante irréductible de $V$. L'on a $\overline{f(T)}^{(U\cap Z)\zar}=\overline{g(T)}^{U\zar}$ ; comme $g$ est quasi-dominant, ce dernier fermé de Zariski est une composante irréductible de $U$, donc un fermé de Zariski irréductible maximal de $U$ ; c'est {\em a fortiori} un fermé de Zariski irréductible maximal de $U\cap Z$, c'est-à-dire une composante irréductible de ce dernier. 

\medskip
Comme tout domaine analytique de $Z$ est, en vertu du \ref{NOR}.\ref{qdometfermzar}.\ref{gerrgrauert}, réunion de parties de la forme $U\cap Z$ où $U$ est un domaine analytique de $X$, on déduit de ce qui précède, en se fondant sur le \ref{NOR}.\ref{qdomestloc}.\ref{qdomfamille}, que $f$ est quasi-dominant.

\medskip
\trois{qdomtop} Soit $f: Y\to X$ un morphisme d'espaces analytiques. Le morphisme $f$ est quasi-dominant  si et seulement si le morphisme $f\duit: Y\duit\to X\duit$ induit par $f$ est quasi-dominant. 

\medskip
En effet, supposons $f\duit$ quasi-dominant. Soit $V$ un domaine analytique de $Y$, soit $U$ un domaine analytique de $X$ tel que $f(V)\subset U$, et soit $T$ une composante irréductible de $V$ ; l'hypothèse de quasi-dominance de $f\duit$ appliquée à $V\duit$, $U\duit$ et à la composante irréductible $T$ de $V\duit$ assure immédiatement que $\overline{f(T)}^{U\zar}$ est une composante irréductible de $U$. 

\medskip
Supposons maintenant $f$ quasi-dominant. L'immersion fermée $Y\duit\hookrightarrow Y$ est quasi-dominante en vertu du \ref{NOR}.\ref{qdometfermzar}.\ref{fermzarqdom} ; la composée $Y\duit\to X$ est dès lors quasi-dominante. Elle se factorise par le sous-espace analytique fermé $X\duit$ de $X$ ; en utilisant à nouveau le \ref{NOR}.\ref{qdometfermzar}.\ref{fermzarqdom}, on voit que la flèche induite $Y\duit\to X\duit$ est quasi-dominante. 

\medskip
\deux{qdometplat} {\bf Quasi-dominance et platitude.} Nous allons maintenant expliquer pourquoi les morphismes plats, dans un sens que nous précisons ci-dessous, sont quasi-dominants, et mentionner un cas particulier important. 

\setcounter{cptbis}{0} 
\medskip
\trois{qdomflat} Soit $f:Y\to X$ un morphisme d'espaces analytiques possédant la propriété suivante : {\em pour tout domaine affinoïde $V$ de $Y$ et tout domaine affinoïde $U$ de $X$ tel que $f(V)\subset U $, la ${\cal A}_U$-algèbre ${\cal A}_V$ est plate.} Le morphisme $f$ est alors quasi-dominant. En effet, soit $V$ un domaine affinoïde de $Y$, soit $U$ un domaine affinoïde de $X$ tel que $f(V)\subset U$, et soit $T$ une composante irréductible de $V$. Comme ${\cal A}_V$ est plat sur ${\cal A}_U$, le morphisme $\spec{\cal A}_V\to \spec {\cal A}_U$ envoie le point générique de toute composante irréductible de $\spec {\cal A}_V$ sur le point générique d'une composante irréductible de $\spec {\cal A}_U$. En conséquence, $\overline{f(T)}^{U\zar}$ est une composante irréductible de $U$ ; le \ref{NOR}.\ref{qdomestloc}.\ref{qdomfamille} permet de conclure que $f$ est quasi-dominant. 

\medskip
\trois{qdomext} Soit $X$ un espace $k$-analytique et soit $L$ une extension complète de $k$. On déduit du \ref{NOR}.\ref{qdometplat}.\ref{qdomflat} ci-dessus que $X_L\to X$ est quasi-dominant.

\medskip
\deux{qdomdim} {\bf Un critère dimensionnel de quasi-dominance}. Soit $f: Y\to X$ un morphisme d'espaces $k$-analytiques possédant les deux propriétés suivantes :

\begin{itemize}
\itb $f$ est de dimension relative égale à zéro ;
\itb pour toute composante irréductible $T$ de $Y$, il existe une composante irréductible $S$ de $X$ telle que $\dim {} S=\dim {} T$ et telle que $f(T)\subset S$.
\end{itemize}

\medskip
Alors pour toute extension complète $L$ de $k$, le morphisme $Y_L\to X_L$ induit par $f$ est quasi-dominant. En effet, soit $L$ une telle extension. La dimension relative de $Y_L\to X_L$ est encore égale à zéro. Soit $Z$ une composante irréductible de $Y_L$ ; en vertu du lemme \ref{COMP}.\ref{irrcompextscal}, $Z$ est une composante irréductible de $T_L$ pour une certaine composante irréductible $T$ de $Y$ ; si $d$ désigne la dimension de $T$, alors $Z$ est de dimension $d$, toujours d'après le lemme \ref{COMP}.\ref{irrcompextscal}. Par hypothèse, il existe une composante irréductible $S$ de $X$ de dimension $d$ telle que $f(T)\subset S$. L'image de $Z$ est contenue dans une composante irréductible de $S_L$, laquelle est, une fois de plus d'après le lemme \ref{COMP}.\ref{irrcompextscal}, une composante irréductible de $X_L$ de dimension $d$. 

\medskip
Il en résulte que $f_L$ satisfait les mêmes hypothèses que $f$ ; on peut donc, pour établir notre assertion, supposer que $L=k$. Soit $V$ un domaine analytique de $Y$, soit $U$ un domaine analytique de $X$ tel que $f(V)\subset U$ et soit $T$ une composante irréductible de $V$. L'adhérence $\overline{T}^{Y\zar}$ de $T$ est une composante irréductible de $Y$ ; soit $d$ sa dimension ; elle coïncide avec celle de $T$. L'image de $ \overline{T}^{Y\zar}$ est contenue dans une composante irréductible $S$ de $X$ de dimension $d$. Comme $T$ est de dimension $d$, il existe $t\in T$ tel que $\mathsf{d}(\hres(t)/k)=d$ ; puisque $f$ est de dimension relative nulle, $\mathsf{d}(\hres(f(t))/k)=d$. L'adhérence $\overline{f(T)}^{U\zar}$ est donc de dimension au moins égale à $d$ ; comme elle est incluse dans $U\cap S$ qui est purement de dimension $d$, c'est une composante irréductible de $U\cap S$, ce qu'il fallait démontrer. 	

\subsection*{Définition de la normalisation} 

\deux{defnomr} {\bf Définition.} Soit $X$ un espace analytique. On appellera {\em normalisation} de $X$ tout objet final de la catégorie des espaces analytiques normaux munis d'une flèche pseudo-dominante vers $X$ ; il découle de cette définition que si $X'$ et $X''$ sont deux normalisations de $X$ alors il existe un {\em} unique $X$-isomorphisme entre $X'$ et $X''$ ; il est tautologique que $X$ est une normalisation de lui-même si et seulement si il est normal. 

\medskip
\deux{normdom} {\bf Lemme.} {\em Soit $X$ un espace analytique et soit $U$ un domaine analytique de $X$. Supposons que $X$ possède une normalisation $X'$. Sous cette hypothèse, $U\times_X X'$ est une normalisation de $U$.}

\medskip
{\em Démonstration.} Puisque $U\times_X X'$ est un domaine analytique de $X'$, il est normal, et $U\times_X X'\to U$ est quasi-dominant d'après le \ref{NOR}.\ref{pseudom}.\ref{qdomprofib}. 

\medskip
Il reste à montrer que tout morphisme quasi-dominant d'un espace normal $Z$ vers $U$ se factorise de manière unique par $U\times_X X'\to U$. Ceci est une conséquence formelle des faits suivants : 

\medskip
\begin{itemize}
\itb si $Z\to U$ est un morphisme, alors il est quasi-dominant si et seulement si le morphisme $Z\to X$ induit est quasi-dominant ;

\itb il existe une bijection naturelle entre l'ensemble des morphismes de $Z$ dans $U$ et l'ensemble des morphismes de $Z$ dans $X$ dont l'image est contenue dans $U$ ; il existe une bijection naturelle entre l'ensemble des morphismes de $Z$ dans $U\times_XX'$ et l'ensemble des morphismes de $Z$ dans $X'$ dont l'image est contenue dans $U\times_XX'$.~$\Box$
\end{itemize}

\subsection*{Existence de la normalisation : le cas affinoïde} 

Rappelons qu'une algèbre affinoïde est un anneau excellent par le théorème \ref{EXC}.\ref{excellence} ; c'est en particulier un anneau japonais\footnote{Ce dernier fait a été montré plus simplement par Berkovich (\cite{brk1},prop. 2.1.14,; $i)$) à partir du caractère japonais des algèbres strictement affinoïdes (\cite{bgr}, 6.1.2/4).}. 

\medskip
\deux{normaliaff} {\bf Proposition.} {\em Soit $\cal A$ une algèbre $k$-affinoïde et soient $\got p_1,\ldots,\got p_r$ ses idéaux premiers minimaux. Pour tout $i$, soit ${\cal B}_i$ la fermeture intégrale de ${\cal A}/\got p_i$ dans son corps des fractions. On munit ${\cal B}:=\prod {\cal B}_i$  de sa structure de $\cal A$-algèbre de Banach finie. L'on pose $X={\cal M}({\cal A})$ et $Y={\cal M}({\cal B})$, et l'on note $f$ le morphisme naturel $Y\to X$.

\medskip
$i)$ Le morphisme $f$ est fini et surjectif ; pour toute extension complète $L$ de $k$, la flèche induite $Y_L\to X_L$ est quasi-dominante ; si $X$ est purement de dimension $d$ pour un certain $d$, alors $Y$ est purement de dimension $d$ ; 

\medskip
$ii)$ $Y$ est une normalisation de $X$.}

\medskip
{\em Démonstration.} Établissons tout d'abord $i)$. Il est immédiat que $f$ est fini. Pour tout $i$, notons $X_i$ (resp. $Y_i$) l'espace ${\cal M}({\cal A}/\got p_i)$ (resp. ${\cal M}({\cal B}_i)$). Les $X_i$ sont les composantes irréductibles de $X$, munies de leur structure réduite ; les $Y_i$ sont les composantes connexes de $Y$ et sont toutes irréductibles par normalité de $\spec \cal B$. Pour tout $i$, le morphisme fini naturel $Y_i\to X_i$ est surjectif puisque ${\cal A}/\got p_i$ s'injecte dans ${\cal B}_i$ ({\em cf.} \cite{duc3}, 0.8) ; il en résulte que $\dim {}X_i=\dim{} Y_i$.  On en déduit que $f$ est fini et surjectif et, en vertu du \ref{NOR}.\ref{qdomdim}, que $f_L$ est  quasi-dominant pour toute extension complète $L$
 de $k$. Si les $X_i$ sont tous de dimension $d$ pour un certain $d\in \NN$, alors par ce qui précède  les $Y_i$ sont tous de dimension $d$ ; l'assertion $i)$ est ainsi démontrée 

\medskip
Établissons maintenant $ii)$. 

\setcounter{cptbis}{0}
\medskip
\trois{normaffqdom} Le $i)$ assure que $f$ est quasi-dominant, et $Y$ est normal par construction.

\medskip
\trois{objinitaff} Soit $Z$ un espace analytique normal et soit $Z\to X$ un morphisme quasi-dominant ; montrons que $Z\to X$ se factorise d'une unique manière par un morphisme $Z\to Y$. En raisonnant G-localement sur $Z$, on se ramène au cas où celui-ci est affinoïde ; notons $\cal C$ l'algèbre des fonctions analytiques sur $Z$ et ${\cal Z}_1,\ldots,{\cal Z}_n$ les composantes connexes de $\spec \cal C$, qui sont intègres puisque $\cal C$ est normal. Soit $j$ un entier compris entre $1$ et $n$. Comme $Z\to X$ est quasi-dominant, le point générique de ${\cal Z}_j$ s'envoie sur le point générique d'une composante irréductible de $\spec \cal A$. Il existe donc un unique $i$ tel que ${\cal Z}_j\to \spec \cal A$ se factorise par un morphisme ${\cal Z}_j\to\spec {\cal A}/\got p_i$ ; cette factorisation est unique et ${\cal Z}_j\to\spec {\cal A}/\got p_i$ est dominant. Comme ${\cal Z}_j$ est intègre et normal, ${\cal Z}_j\to\spec {\cal A}/\got p_i$ se factorise d'une unique manière par une flèche ${\cal Z}_j\to\spec {\cal B}_i$.

\medskip
Il découle immédiatement de ce qui précède que ${\cal A}\to {\cal C}$ se factorise d'une unique manière par un morphisme d'anneaux ${\cal B}\to {\cal C}$. Comme ${\cal A}\to {\cal B}$ est borné, et comme $\cal B$ est munie de sa structure de ${\cal A}$-algèbre de Banach finie, ${\cal B}\to {\cal C}$ est borné, ce qui achève la démonstration.~$\Box$  

\subsection*{Existence de la normalisation : le cas général}

\medskip
\deux{vocparacomp} {\bf Topologie générale.} Rappelons qu'un espace topologique $X$ est dit {\em paracompact} s'il est séparé et si chacun de ses recouvrements ouverts peut être raffiné en un recouvrement localement fini ; il revient au même, lorsque $X$ est séparé et localement compact, de demander qu'il soit réunion disjointe d'ouverts dénombrables à l'infini. Cette équivalence est établie par exemple dans \cite{topg}, I, \S 10, no. 12 ; la preuve qui y est donnée de l'implication réciproque montre même, sans le mentionner explicitement, que si $X$ est localement compact et paracompact et si $\cal B$ est une base de voisinages {\em non nécessairement ouverts} de $X$, alors tout recouvrement ouvert de $X$ peut être raffiné en un recouvrement localement fini constitué d'éléments de $\cal B$. 

\medskip
\deux{existnormx} {\bf Théorème.} {\em Tout espace analytique admet une normalisation.} 

\medskip
{\em Démonstration.} Soit $X$ un espace analytique.

\setcounter{cptbis}{0}

\medskip
\trois{xsepparcomp} {\bf Supposons que $X$ est paracompact}. Il possède alors un recouvrement localement fini $(X_i)$ par des domaines affinoïdes ; chacun d'eux est compact et donc fermé dans $X$. Pour tout $i$, l'espace affinoïde $X_i$ possède une normalisation $X'_i$ d'après la proposition \ref{NOR}.\ref{normaliaff} ci-dessus. Soient $i$ et $j$ deux indices. Pour tout domaine analytique $U$ de $X_i\cap X_j$, les espaces $X'_i\times_{X_i} U
$ et $X'_j\times_{X_j} U$ sont d'après le lemme \ref{NOR}.\ref{normdom} deux normalisations de $U$ ; ils existe donc un {\em unique} $U$-isomorphisme entre eux, ce qui permet de mettre en \oe uvre le procédé de recollement de \cite{brk2}, prop. 1.3.3, $b)$ ; on obtient un espace analytique $X'$ muni d'un morphisme vers $X$ et d'une collection d'isomorphismes $X'\times_X X_i\simeq X'_i$. 

\medskip
L'espace $X'$ est G-recouvert par les $X'_i$ qui sont normaux, et il est donc normal. Comme $X'_i\to X_i$ est quasi-dominant pour tout $i$, la flèche $X'\to X$ est quasi-dominante d'après le \ref{NOR}.\ref{qdomestloc}.\ref{qdomglocb}. 

\medskip
Soit $Z$ un espace analytique normal, et soit $Z\to X$ un morphisme quasi-dominant. Pour tout $i$ et pour tout domaine analytique $U$ de $X_i$, le domaine analytique $Z\times_XU$ de $Z$ est normal, et $Z\times_XU \to U$ est quasi-dominant par le \ref{NOR}.\ref{pseudom}.\ref{qdomprofib}. Dès lors, $Z\times_XU \to U$ se factorise d'une unique manière par la normalisation de $U$, laquelle s'identifie en vertu du lemme \ref{NOR}.\ref{normdom} au produit fibré $X'_i\times_{X_i} U$, autrement dit à $X'\times_XU$. L'existence et l'unicité de ces factorisations locales assurent l'existence et l'unicité d'une factorisation globale de $Z\to X$ par un morphisme $Z\to X'$ ; ainsi, $X'$ est une normalisation de $X$. 

\medskip
\trois{xnormacasgen} {\bf Le cas général.} On ne fait plus d'hypothèse sur $X$. Il possède un recouvrement $(X_i)$ par des ouverts paracompacts et topologiquement séparés ; en vertu  du \ref{NOR}.\ref{existnormx}.\ref{xsepparcomp}, chaque $X_i$ possède une normalisation $X'_i$. On construit une normalisation de $X$ en recollant les $X'_i$, par une méthode analogue à celle suivie  au \ref{NOR}.\ref{existnormx}.\ref{xsepparcomp}, mais en se fondant cette fois-ci sur sur \cite{brk2}, prop. 1.3.3, $a)$.~$\Box$ 

\medskip
{\bf Si $X$ est un espace analytique, il possède donc une normalisation ; celle-ci étant unique à unique isomorphisme près, on parlera désormais de {\em la} normalisation de $X$.} 

\medskip
\deux{propnorm} Soit $X$ un espace $k$-analytique et soit $X'$ sa normalisation.

\medskip
\setcounter{cptbis}{0}
\trois{normfinidim} Le morphisme $X'\to X$ est fini et surjectif ; si $X$ est purement de dimension $d$ pour un certain entier $d$, alors $X'$ est purement de dimension $d$ ; si $L$ est une extension complète de $k$, alors $X'_L\to X_L$ est quasi-dominant. En effet, ces assertions sont locales pour la G-topologie sur $X$ (pour la dernière, cela est dû au \ref{NOR}.\ref{qdomestloc}.\ref{qdomglocb}), ce qui autorise à supposer que $X$ est $k$-affinoïde, auquel cas elles résultent de la proposition \ref{NOR}.\ref{normaliaff} et du \ref{NOR}.\ref{qdomdim}.

\medskip
\trois{normred} Étant normal, $X'$ est réduit, et $X'\to X$ se factorise donc d'une unique manière par $X\duit$. Soit $Z$ un espace analytique normal et soit $Z\to X\duit$ un morphisme quasi-dominant. Le morphisme composé $Z\to X\duit\to X$ est encore quasi-dominant (\ref{NOR}.\ref{qdometfermzar}.\ref{fermzarqdom}), et se factorise donc d'une unique manière par une flèche $Z\to X'$. Comme $X\duit\to X$ est une immersion fermée, l'application naturelle $\mbox{Hom}_{X\duit}(Z,X')\to \mbox{Hom}_X(Z,X')$ est bijective. En conséquence, $Z\to X\duit$ se factorise d'une unique manière par $X'$ ; autrement dit, $X'\to X\duit$ identifie $X'$ à la normalisation de $X\duit$.

\subsection*{Normalisation et composantes irréductibles} 

\medskip
\deux{connirred} {\bf Proposition.} {\em Soit $X$ un espace analytique normal. L'ensemble des composantes irréductibles de $X$ coïncide avec celui de ses composantes connexes.} 

\medskip
{\em Démonstration.} Comme l'ensemble des composantes irréductibles de $X$ est G-localement fini, il suffit de montrer que celles-ci sont deux à deux disjointes. Soit $x\in X$ et soient $Y$ et $Z$ deux composantes irréductibles de $X$ contenant $x$. Soit $V$ un domaine affinoïde de $X$ contenant $x$. Soit $S$ (resp. $T$) une composante irréductible de $Y\cap V$ (resp. $Z\cap V$) contenant $x$. Les fermés de Zariski $S$ et $T$ de $V$ en sont deux composantes irréductibles. Comme $\spec {\cal A}_V$ est normal, ses composantes irréductibles sont deux à deux disjointes, et il en va donc de même de celles de $V$ ; en conséquence, $S=T$. Mais on a alors $Y=\overline{S}^{X\zar}=\overline{T}^{X\zar}=Z$.~$\Box$ 

\medskip
\deux{redpasnorm} {\bf Lemme.} {\em Soit $X$ un espace analytique réduit. Le lieu de normalité de $X$ est un ouvert de Zariski dense de $X$.} 

\medskip
{\em Démonstration.} On sait déjà qu'il s'agit d'un ouvert de Zariski. Concernant la densité, il suffit de la vérifier G-localement, et l'on peut donc supposer que $X$ est affinoïde. Soit $\cal A$ l'algèbre des fonctions analytiques sur $X$. Comme ${\cal A}$ est réduit, les anneau locaux des points génériques des composantes irréductibles de $\spec {\cal A}$ sont des corps ; en conséquence, le lieu de normalité de $\spec \cal A$ rencontre chacune de ses composantes irréductibles.~$\Box$ 

\medskip
\deux{irrnormconn} {\bf Proposition.} {\em Soit $X$ un espace analytique. Si $X$ est irréductible, sa normalisation est connexe.}

\medskip
{\em Démonstration.}  Soit $d$ la dimension de $X$ et soit $X'$ sa normalisation. Comme elle coïncide avec la normalisation de $X\duit$ (\ref{NOR}.\ref{propnorm}.\ref{normred}), on peut remplacer $X$ par $X\duit$ (ce qui ne modifie pas son caractère irréductible) et donc le supposer réduit. L'espace $X'$ est purement de dimension $d$ (\ref{NOR}.\ref{propnorm}.\ref{normfinidim}). Supposons que l'on ait $X'=X'_1\cup X'_2$  où les $X'_i$ sont deux ouverts disjoints et non vides de $X'$. Soit $i\in\{1,2\}$.  L'ouvert non vide $X'_i$ de $X'$ est de dimension $d$, et est par ailleurs un fermé de Zariski de $X'$. Son image sur $X$ est donc un fermé de Zariski de dimension $d$ de ce dernier ; puisque $X$ est lui-même irréductible et de dimension $d$, la flèche $X'_i\to X$ est surjective. 

\medskip
Ceci valant pour $i=1$ et pour $i=2$, tout point de $X$ a au moins deux antécédents sur $X'$. Le lieu de normalité de $X$, au-dessus duquel le morphisme de normalisation est un isomorphisme, est donc vide. Mais $X$ est non vide et réduit, et le lemme \ref{NOR}.\ref{redpasnorm} conduit alors à une contradiction.~$\Box$ 

\medskip
\deux{desccompnorm} {\bf Théorème.} {\em Soit $X$ un espace analytique et soit $f:X'\to X$ sa normalisation.

\medskip

$i)$ Si $\Omega$ est une composante connexe de $X'$ alors $f(\Omega)$ est une composante irréductible de $X$, et le morphisme naturel $\Omega\to f(\Omega)\duit$ (dont l'existence est assuré par le caractère réduit de $\Omega$) identifie $\Omega$ à la normalisation de $f(\Omega)\duit$  ;

\medskip
$ii)$  l'application $\Omega\mapsto f(\Omega)$ établit une bijection $\pi_0(X')\simeq \irr X$ .}

 \medskip
 {\em Démonstration.} Soit $(X_i)$ la famille des composantes irréductibles de $X$. Pour tout $i$, notons $X'_i$ la normalisation de $X_{i,red}$. En vertu du \ref{NOR}.\ref{propnorm}.\ref{normfinidim} et de la proposition \ref{NOR}.\ref{irrnormconn}, il suffit de montrer que le $X$-espace $\coprod X'_i$ s'identifie à la normalisation de $X$.
 
 \medskip
 L'espace $\coprod X'_i$ est normal. Pour tout indice $i$, le morphisme $X'_i\to X_{i,red}$ est quasi-dominant, et $X'_i\to X$ est dès lors quasi-dominant (\ref{NOR}.\ref{qdometfermzar}.\ref{fermzarqdom}) ; il s'ensuit que $\coprod X'_i\to X'$ est quasi-dominant  (\ref{NOR}.\ref{qdomestloc}.\ref{qdomglocs}.) 
 
 \medskip
 Soit $Z$ un espace analytique normal, et soit $Z\to X$ un morphisme quasi-dominant. Nous allons montrer que  $Z\to X$ se factorise d'une unique manière par $\coprod X'_i$, ce qui permettra de conclure. Il suffit de montrer que pour toute composante connexe $S$ de $Z$, la flèche $S\to X$ se factorise d'une unique manière par $\coprod X'_i$. Soit donc $S$ une composante connexe de $Z$ ; le lemme \ref{NOR}.\ref{connirred} assure que $S$ est irréductible. Par ailleurs $S$ est un ouvert de $Z$ ; en conséquence, l'espace analytique $S$ est normal et $S\to X$ est quasi-dominant. 
   
 \medskip
Puisque $S$ est irréductible et puisque $S\to X$ est quasi-dominant, l'adhérence de Zariski de l'image de $S$ dans $X$ est égale à $X_\ell$ pour un certain $\ell$ ; il en découle que si $j\neq \ell$, l'image de $S$ n'est pas incluse dans $X_j$. Par connexité de $S$, toute flèche  $S\to \coprod X'_i$ se factorise par $X'_j$ pour un certain $j$ ; de ce qui précède, on déduit alors que toute $X$-flèche  $S\to \coprod X'_i$ se factorise par $X'_\ell$.  Il suffit donc de démontrer que $S\to X$ se factorise de manière unique par $X'_\ell$.

\medskip
L'espace $S$ est réduit, et son image sur $X$ est incluse dans $X_\ell$. Il en découle que $S\to X$ se factorise d'une manière unique par $X_{\ell,red}$. Comme $X_{\ell,red}\to X$ est une immersion fermée, $\mbox{Hom}_{X}(S,X'_\ell)=\mbox{Hom}_{X_{\ell,red}}(S,X'_\ell)$ ; il suffit donc de démontrer que $S\to X_{\ell,red}$ se factorise d'une unique manière par $X'_\ell$.

\medskip
Or la flèche composée $S\to X_{\ell,red}\to X$ est quasi-dominante  ; il s'ensuit que $S\to X_{\ell,red}$ est  quasi-dominante (\ref{NOR}.\ref{qdometfermzar}.\ref{fermzarqdom}). L'espace $S$ étant normal, $S\to X_{\ell,red}$ se factorise de manière unique par la normalisation $X'_\ell$ de $X_{\ell,red}$.~$\Box$ 

\medskip
De ce qui précède et du lemme \ref{COMP}.\ref{compouvdense} découle immédiatement le corollaire qui suit.

\medskip
\deux{imrecouvdense} {\bf Corollaire.} {\em Soit $X$ un espace analytique et soit $X'$ sa normalisation. Soit $U$ un ouvert de Zariski dense de $X$. Le produit fibré $X'\times_XU$ est un ouvert de Zariski dense de $X'$.}~$\Box$ 

\medskip
{\bf Le théorème \ref{NOR}.\ref{desccompnorm} et le corollaire \ref{NOR}.\ref{imrecouvdense} seront utilisés dans la suite sans rappel ni justification.} 

\subsection*{Normalisation et anneaux locaux} 

\medskip
\deux{normannloc} {\bf Lemme.} {\em Soit $X$ un bon espace analytique et soit $X'$ sa normalisation. Soit $x\in X$ et soient $x'_1,\ldots,x'_r$ ses antécédents sur $X'$. L'anneau $\prod {\cal O }_{X',x'_i}$ s'identifie naturellement à la normalisation de ${\cal O}_{X,x}$.}

\medskip
{\em Démonstration.} On procède en plusieurs étapes. 

\medskip
\setcounter{cptbis}{0}
\trois{xaffqi} On peut supposer que $X$ est affinoïde ; notons $\cal A$ l'algèbre correspondante. Soient $\got p_1,\ldots, \got p_s$ les idéaux premiers minimaux de ${\cal O}_{X,x}$. Le lemme \ref{RAP}.\ref{idminoxx} permet de se ramener, quitte à restreindre $X$, au cas où il existe des idéaux $\got q_1,\ldots,\got q_s$ de $\cal A$ possédant les propriété suivantes, $F_i$ désignant pour tout $i$ le lieu des zéros de $\got q_i$ sur $X$ : 

\medskip
$i)$ $\got q_i{\cal O}_{X,x}=\got p_i$ pour tout $i$ ; 

\medskip
$ii)$ pour tout voisinage affinoïde $V$ de $x$ dans $X$ le point $x$ n'appartient qu'à une composante irréductible $G_i$ de $F_i\cap V$ ; les $G_i$ sont deux à deux distinctes et sont exactement les composantes irréductibles de $V$ contenant $x$.

\medskip
\trois{regals} Soit $V$ un voisinage affinoïde de $x$ dans $X$ ; posons $V'=V\times_X X'$. D'après $ii)$, il y a exactement $s$ composantes irréductibles de $V$ qui passent par $x$ ; il y a donc exactement $s$ composantes connexes de $V'$ dont l'image sur $V$ contient $x$. Ceci valant pour tout $V$, on a $r=s$. 

\medskip
\trois{caspartrun} Traitons maintenant un cas particulier, celui où ${\cal O}_{X,x}$ est intègre ; notons qu'alors $r=s=1$, et l'on désignera simplement par $x'$ l'unique antécédent de $x$ sur $X'$. L'anneau local ${\cal O}_{X',x'}$ est normal, et fini sur ${\cal O}_{X,x}$.

\medskip
Soit $f$ appartenant au noyau de ${\cal O}_{X,x}\to {\cal O}_{X',x'}$. Soit $\Omega$ un voisinage affinoïde de $x$ dans $X$ sur lequel $f$ est définie et soit $\Omega'$ l'image réciproque de $\Omega$ sur $X'$. Comme l'image $f'$ de $f$ dans ${\cal O}_{X',x'}$ est nulle le lieu des zéros de $f'$ contient un voisinage de $x'$ dans $\Omega'$. Par normalité de $\Omega'$, il contient la composante connexe $\Omega'_0$ de $x'$ dans $\Omega'$ (laquelle est irréductible). Le lieu des zéros de $f$ contient donc l'image de $\Omega'_0$, qui est une composante irréductible de $\Omega$ passant par $x$ ; l'assertion $ii)$ ci-dessus garantit l'unicité d'une telle composante (rappelons que $r=s=1$) qui est donc un voisinage de $x$ dans $X$. Il en découle que $f$ est nilpotente au voisinage de $x$, et donc nulle lorsqu'on la voit comme appartenant à l'anneau local {\em intègre} ${\cal O}_{X,x}$ ; ainsi,  ${\cal O}_{X,x}\hookrightarrow {\cal O}_{X',x'}$. 

\medskip
Soit $f$ appartenant à ${\cal O}_{X',x'}$. Il existe un voisinage affinoïde $\Omega$ de $x$ dans $X$ tel que $f$ soit définie sur l'image réciproque $\Omega'$ de $\Omega$. Par la construction de la normalisation dans le cas affinoïde, il existe $g\in{\cal A}_{\Omega}$ qui n'est pas un diviseur de zéro et telle que $gf$ provienne d'une fonction appartenant à ${\cal A}_{\Omega}$. Comme $g$ n'est pas diviseur de zéro dans ${\cal A}_{\Omega}$, le lieu des zéros de $g$ ne contient pas l'unique composante irréductible de $\Omega$ à laquelle appartient $x$ ; il ne contient donc aucun ouvert non vide de cette dernière, et est dès lors non nul lorsqu'on le voit comme élément de ${\cal O}_{X,x}$. L'image de $f$ dans le corps des fractions de ${\cal O}_{X',x'}$ provient donc du corps des fractions de ${\cal O}_{X,x}$. Ainsi ${\cal O}_{X,x}\hookrightarrow {\cal O}_{X',x'}$ apparaît comme une injection finie à but normal entre anneaux intègres qui induit un isomorphisme au niveau des corps de fractions ; par conséquent, ${\cal O}_{X',x'}$ s'identifie à la normalisation de ${\cal O}_{X,x}$.

\medskip
\trois{gigprimei} Revenons au cas général et aux notations correspondantes. Soit $V$ un voisinage affinoïde de $x$ dans $X$ ; posons $V'=V\times_X X'$ et considérons les fermés $G_i$ introduits lors de l'énoncé de l'assertion $ii)$ ; {\em on les munit de leur structure réduite}. Pour tout $i$, notons $G'_i$ la normalisation de $G_i$, identifiée à une composante connexe de $V'$. Puisque $s=r$, on peut renuméroter les $x'_i$ de sorte que $x'_i \in G'_i$ pour tout $i$. 

\medskip
Fixons $i$. L'anneau ${\cal O}_{X',x'_i}$ coïncide avec ${\cal O}_{G'_i,x'_i}$. Comme le fermé de Zariski réduit $G_i$ est défini au voisinage de $x$ par l'idéal $\got q_i$, l'anneau local ${\cal O}_{G_i,x}$ est canoniquement isomorphe à $({\cal O}_{X,x}/\got p_i)\duit$ (lemme \ref{COMP}.\ref{redannloc}), qui n'est autre que ${\cal O}_{X,x}/\got p_i$ puisque $\got p_i$ est premier. Le cas particulier traité au \ref{NOR}.\ref{normannloc}.\ref{caspartrun} assure que ${\cal O}_{G'_i,x'_i}$ s'identifie à la normalisation de l'anneau local intègre ${\cal O}_{G_i,x}$, c'est-à-dire que ${\cal O}_{X',x'_i}$ s'identifie à la normalisation de ${\cal O}_{X,x}/\got p_i$. Ceci valant pour tout $i$, le produit $\prod {\cal O}_{X',x'_i}$ s'identifie à la normalisation de ${\cal O}_{X,x}$.~$\Box$

\subsection*{Normalisation et extension des scalaires}

\medskip
\deux{normextscal} {\bf Proposition.} {\em Soit $X$ un espace $k$-analytique et soit $L$ une extension complète de $k$. Soit $X'$ la normalisation de $X$ et soit $Y$ celle de $X_L$.  Il existe un unique $X$-morphisme $Y\to X'_L$ ; il est fini, surjectif, et c'est un isomorphisme si  et seulement si $X'_L$ est normal.}

\medskip
{\em Démonstration.}  On procède en plusieurs étapes.

\medskip
\setcounter{cptbis}{0}
\trois{existyxl} {\em Existence et finitude du morphisme}. Par définition, $Y\to X_L$ est quasi-dominant, et $X_L\to X$ est quasi-dominant d'après le \ref{NOR}.\ref{qdometplat}.\ref{qdomext}. Par conséquent, $Y\to X$ est quasi-dominant. Comme $Y$ est normal, $Y\to X$ se factorise par un unique morphisme $Y\to X'_L$, qui est fini puisque $Y\to X_L$ et $X'_L\to X_L$ le sont. 

\medskip
\trois{suryxl} {\em Surjectivité du morphisme.} On peut supposer que $X$ est réduit. Soit $U$ le lieu de normalité de $X$ et soit $U'$ son image réciproque sur $X'$. Le lemme \ref{NOR}.\ref{redpasnorm} assure que $U$ est dense dans $X$, le corollaire \ref{NOR}.\ref{imrecouvdense} garantit que $U':=U\times_XX'$ est un ouvert de Zariski dense de $X'$. Par définition de $U$, la flèche $U'\to U$ est un isomorphisme ; dès lors $U'_L\simeq U_L$ ; notons que la densité de $U'$ dans $X'$ entraîne celle de $U'_L$ dans $X'_L$, et que $U'_L$ n'est autre que l'image réciproque de $U_L$ sur $X'_L$. Comme $Y\to X_L$ est surjective, l'image de $Y$ sur $X_L$ contient $U_L$ ; comme $X'_L\to X_L$ induit un isomorphisme $U'_L\simeq U_L$, l'image de $Y$ sur $X'_L$ contient $U'_L$. Cette image est par ailleurs, en vertu de la finitude de $Y\to X_L$, un fermé de Zariski de $X'_L$. L'ouvert $U'_L$ de $X'_L$ étant dense pour la topologie de Zariski, $Y\to X'_L$ est surjective. 

\medskip
\trois{bijyxl} {\em Condition nécessaire et suffisante pour avoir affaire à un isomorphisme.} Soit  $Z$ un espace analytique normal et soit $Z\to X_L$ une flèche quasi-dominante. La composée $Z\to X_L\to X$ est quasi-dominante et $Z$ est normal ; il existe dès lors un unique $k$-morphisme $Z\to X'$ par lequel $Z\to X$ se factorise, et donc un unique $L$-morphisme $Z\to X'_L$ par lequel $Z\to X_L$ se factorise. Il en découle, $X'_L\to X_L$ étant quasi-dominant en vertu du \ref{NOR}.\ref{propnorm}.\ref{normfinidim},  que $X'_L$ s'identifie à la normalisation de $X_L$ si et seulement si il est normal.~$\Box$

\section{Autour de la régularité géométrique}\label{UNIREG}

{\em Afin d'éviter des distinctions fastidieuses entre la caractéristique nulle et la caractéristique positive, on note désormais par $p$ {\em l'exposant caractéristique} de $k$. Si $p=1$ ({\em i.e.} si $k$ est de caractéristique nulle) alors $k^{1/p^n}$ désignera simplement $k$. La {\em clôture parfaite} de $k$ est la réunion des $k^{1/p^n}$.}

\subsection*{La validité sur $k^{1/p}$ entraîne la validité géométrique}
\setcounter{cpt}{0}

\deux{clotparfparf} {\bf Lemme.} {\em Le complété $L$ de la clôture parfaite de $k$ est parfait.}

\medskip
{\em Démonstration.} Soit $x\in L$. Il existe une suite $(x_i)$ d'éléments de $\bigcup k^{1/p^n}$ qui converge vers $x$ ; la suite $(x_i^{1/p})$ est une suite d'éléments de $\bigcup k^{1/p^n}$ qui est de Cauchy et dont la limite $y$ dans $L$ vérifie l'égalité $y^p=x$.~$\Box$

\bigskip
\deux{omegareg} {\bf  Lemme.} {\em Soit $X$ un espace $k$-affinoïde et soit $x\in X$. 

\medskip
\begin{itemize}
 \item [$i)$] $\dim{\hres(x)}\Omega^1_{X/k}\otimes\hres(x)\geq \dim x X$ ;
\item[$ii)$] si $\dim{\hres(x)}\Omega^1_{X/k}\otimes\hres(x)=\dim x X$, alors ${\cal O}_{X,x}$ est régulier ;
\item[$iii)$] si $|k^{*}|\neq\{1\}$, si ${\cal O}_{X,x}$ est régulier et si $\hres(x)=k$, on a l'égalité $\dim{\hres(x)}\Omega^1_{X/k}\otimes\hres(x)=\dim xX$.
\end{itemize}}

\medskip
{\em Démonstration.} Soit $\cal A$ l'algèbre des fonctions analytiques sur $X$. Si $L$ est une extension complète de $k$, la flèche $X_{L}\to X$ est un morphisme fidèlement plat d'espaces annelés. On peut donc supposer, quitte à étendre les scalaires, que $\hres(x)=k$ et que $|k^{*}|\neq\{1\}$. Il existe alors une immersion fermée de $X$ dans un polydisque compact {\em de centre $x$} ; on en déduit l'existence d'un voisinage strictement $k$-affinoïde de $x$ dans $X$. Cette remarque permet de restreindre $X$ de sorte qu'il soit lui-même strictement $k$-affinoïde. On a alors l'égalité $\dim{\tiny \mbox{Krull}}{\cal O}_{X,x}=\dim{x}X$ (\ref{RAP}.\ref{dimanal}.\ref{oxxkrull}). 

\medskip
Soit $\got{m}$ l'idéal maximal de ${\cal O}_{X,x}$. La dérivation $f\mapsto f-f(x)$ induit un isomorphisme $\Omega^1_{X/k}\otimes_{x}k\simeq \got{m}/\got{m}^{2}$, et le lemme s'ensuit immédiatement.~$\Box$ 

\medskip
\deux{proprgeom} {\bf Proposition.} {\em Soit $X$ un espace $k$-analytique et soit $x\in X$. Soit $\mathsf P$ une propriété appartenant à l'ensemble $\cal R$ défini au {\rm \ref{GAGA}.\ref{listepropri}}. Les propositions suivantes sont équivalentes : 

\medskip
$i)$ il existe une extension complète et parfaite $L$ de $k$ telle que pour tout point $y$ de $X_{L}$ situé au-dessus de $x$, l'espace $X_{L}$ satisfasse $\mathsf P$ en $y$ ;

\medskip
$ii)$ pour toute extension complète $F$ de $k$ et pour tout point $z$ de $X_{F}$ situé au-dessus de $X$, l'espace $X_{F}$ satisfait $\mathsf P$ en $y$. 

\medskip
\noindent
Si de plus $\mathsf P$ est la régularité, ces deux propositions équivalent  à la troisième ci-dessous :

\medskip
$iii)$ $\dim{\hres(x)}\Omega^1_{X/k}\otimes\hres(x)=\dim{x}X$.} 

\medskip
{\em Démonstration.} Montrons que $i)\Rightarrow ii)$. Supposons que $i)$ est vraie, soit $F$ une extension complète de $k$ et soit $z$ un point de $X_F$ situé au-dessus de $x$. Soit $\KK$ une extension complète composée de $F$ et $L$ et soit $t$ un point de $X_\KK$ situé au-dessus de $z$ ; soit $y$ l'image de $t$ sur $X_L$. Comme $L$ est parfait, $\KK$ est une extension anaytiquement séparable de $L$ ; l'espace $X_L$ satisfaisant $\mathsf P$ en $y$, l'espace $X_{\KK}$ satisfait $\mathsf P$ en $t$ ; on en déduit que $X_F$ satisfait $\mathsf P$ en $z$. 

\medskip
L'implication $ii)\Rightarrow i)$ est triviale. 

\medskip
Plaçons-nous maintenant dans le cas où $\mathsf P$ est la régularité. Supposons que $iii)$ est vraie. Soit $L$ le complété de la clôture parfaite de $k$. Soit $y$ un point de $X_{L}$ situé au-dessus de $x$ ; grâce à l'égalité $\dim{\hres(y)} \Omega^1_{X_{L}/L}\otimes\hres(y)=\dim{y}X_{L}$ et en vertu du lemme ~\ref{UNIREG}.\ref{omegareg}, l'espace $X_{L}$ est régulier en $y$ ; comme $L$ est parfait (lemme \ref{UNIREG}.\ref{clotparfparf}), l'assertion $i)$ est ainsi établie. Supposons maintenant que $ii)$ est vraie. Soit $F$ une extension complète de $\hres(x)$ dont la valeur absolue n'est pas triviale. L'espace $X_F$ possède un $F$-point $y$ au-dessus de $x$. Comme $X_F$ est régulier en $y$ en vertu de l'hypothèse $ii)$, on a d'après le lemme~\ref{UNIREG}.\ref{omegareg} l'égalité $\dim{\hres(y)}\Omega^1_{X_F/F}\otimes\hres(y)=\dim{y}X_F$ ; on en déduit que  $\dim{\hres(x)}\Omega^1_{X/k}\otimes\hres(x)=\dim{x}X\; ;$ ainsi, $iii)$ est vérifiée.~$\Box$

\medskip
\deux{defql} {\bf Définition.} On dit que l'espace $X$ satisfait {\em géométriquement} la propriété $\mathsf P$ en $x$ si les propriétés équivalentes $i)$ et $ii)$ sont vérifiées ; on dit que $X$ satisfait géométriquement $\mathsf P$ s'il la satisfait géométriquement en chacun de ses points. On parlera ainsi d'espace géométriquement régulier, géométriquement normal ({\em i.e.} $S_2$ et géométriquement $R_1$, d'après le critère de Serre) géométriquement réduit ({\em i.e.} $S_1$ et géométriquement $R_0$, d'après le critère de Serre), etc. En référence à la propriété $iii)$, on emploiera le plus souvent l'expression «quasi-lisse» en lieu et place de «géométriquement régulier». 

\medskip
\deux{remgeomparf} {\bf Remarque.} Si $k$ est parfait, alors par définition $X$ satisfait $\mathsf P$ en $x$ si et seulement si il la satisfait géométriquement en $x$. 

\medskip
\deux{qlouvzar} {\bf Proposition.} {\em Soit $X$ un espace $k$-analytique. Le lieu de quasi-lissité $U$ de $X$ est un ouvert de Zariski de $X$, la dimension locale $d$ est une fonction localement constante sur $U$, et $\Omega^1_{U/k}$ est {\rm G}-localement libre de rang $d$.} 

\medskip
{\em Démonstration.} Les propriétés énoncées étant de nature G-locale, on peut supposer que $X$ est $k$-affinoïde. On désigne par $\cal A$ l'algèbre des fonctions analytiques sur $X$ et par $\cal X$ le spectre de $\cal A$. Soit $(X_{i})$ la famille finie  des composantes irréductibles de $X$ ; pour chaque $i$, on note $d_{i}$ la dimension de $X_{i}$ et ${\cal X}_{i}$ le fermé de Zariski de $\cal X$ qui correspond à $X_{i}$. Si $x$ est un point quasi-lisse de $X$, alors ${\cal O}_{X,x}$ est régulier, et en particulier intègre ; il s'ensuit, d'après le lemme~\ref{RAP}.\ref{localintegre}, que $x$ n'est situé que sur une seule composante irréductible de $X$ ; si $X_i$ est la composante en question, on a $\dim x X=d_i$.

\medskip
Soit ${\cal V}$ l'ouvert de ${\cal X}$ formé des points qui ne sont situés que sur une seule composante irréductible ; on a ${\cal V}=\coprod{\cal V}\cap{\cal X}_{i}$. Pour tout $i$ et tout point $\bf x$ de ${\cal V}\cap{\cal X}_{i}$, le lemme~\ref{UNIREG}.\ref{omegareg} assure que $\dim{\kappa({\bf x})}\Omega^1_{{\cal A}/k}\otimes \kappa({\bf x})\geq d_{i}$ ; par semi-continuité supérieure du rang, l'ensemble ${\cal U}_{i}$ des points $\bf x$ de ${\cal V}\cap{\cal X}_{i}$ en lesquels on a égalité est un ouvert de Zariski de $\cal X$. 

\medskip
D'après ce qui précède, $U$ s'identifie à l'image réciproque de $\coprod {\cal U}_i$, et c'est donc bien un ouvert de Zariski de $X$, sur lequel la dimension locale est par construction localement constante. L'anneau ${\cal O}_{X,x}$ est régulier pour tout point $x$ de $U$ ; il en résulte que le schéma $\coprod {\cal U}_{i}$ est régulier. Le faisceau cohérent induit par $\Omega^1_{{\cal A}/k}$ sur le schéma noethérien localement intègre $\coprod {\cal U}_{i}$ ayant un rang localement constant, il est localement libre. On en conclut que $\Omega^1_{U/k}$ est localement libre ; par construction, son rang en tout point $x$ de $U$ coïncide avec $\dim xU$.~$\Box$

\medskip
Le lemme suivant sera utilisé à plusieurs reprises par la suite. 

\medskip
\deux{homeozar} {\bf Lemme.} {\em Soit $n\in \NN$ et soit $X$ un espace $k$-analytique. L'application naturelle $X_{k^{1/p^n}}\to X$ induit un homéomorphisme entre les espaces sous-jacents {\em munis de leurs topologies de Zariski.}}

\medskip
{\em Démonstration.} Cette application étant bijective ({\em cf.} \ref{RAP}.\ref{sensxl}.\ref{raddensehomeo}) et continue pour les topologies de Zariski, il reste à s'assurer qu'elle transforme tout fermé de Zariski de $X_{k^{1/p^n}}$ en un fermé de Zariski de $X$. La question est locale pour la G-topologie sur $X$, que l'on peut donc supposer $k$-affinoïde. Soit $F$ un fermé de Zariski de $X_{k^{1/p^n}}$. Il est défini comme le lieu d'annulation simultanée d'une famille $(f_i)$ de fonctions analytiques sur $X_{k^{1/p^n}}$ ; on peut également le décrire comme le lieu d'annulation simultanée des $f_i^{p^n}$. Mais  $f_i^{p^n}$ provient pour tout $i$ d'une fonction analytique $g_i$ sur $X$ ; l'image de $F$ sur $X$ est le lieu d'annulation simultanée des $g_i$, qui est un fermé de Zariski de $X$.~$\Box$ 

\medskip
L'énoncé de la proposition qui suit a été inspiré à l'auteur par le lemme 3.3.1 et de la proposition 3.3.6 de \cite{comprig}. 

\medskip
\deux{kunsurp} {\bf Proposition.} {\em Soit  $X$ un espace $k$-analytique, soit $x$ un point de $X$ et soit $y$ l'unique point de $X_{k^{1/p}}$ situé au-dessus de $x$. Soit $\mathsf P$ une propriété appartenant à $\cal R$. Les propositions suivantes sont équivalentes : 

\medskip
\begin{itemize}
\item[$i)$] $X$ satisfait géométriquement $\mathsf P$ en $x$ ; 

\item[$ii)$] $X_{k^{1/p}}$ satisfait $\mathsf P$ en $y$.

\end{itemize}

\medskip
En particulier, $X$ est quasi-lisse en $x$ si et seulement si $X_{k^{1/p}}$ est régulier en $y$.}

\medskip
{\em Démonstration.} L'implication $i)\Rightarrow ii)$ est évidente, et la dernière assertion est simplement une reformulation de l'équivalence  $i)\iff ii)$ lorsque $\mathsf P$ est la régularité. Il reste à établir $ii)\Rightarrow i)$ ; on suppose donc que $ii)$ est vérifiée. Soit $U$ le lieu de validité de $\mathsf P$ sur $X_{k^{1/p}}$. C'est un ouvert de Zariski de $X_{k^{1/p}}$ ; son image $V$ sur $X$ est un ouvert de Zariski de $X$  d'après le lemme \ref{UNIREG}.\ref{homeozar} ci-dessus; en remplaçant $X$ par $V$, on se ramène au cas où $X_{k^{1/p}}$ satisfait $\mathsf P$, et l'on cherche à montrer que $X_L$ satisfait $\mathsf P$ pour toute extension complète $L$ de $k$. On peut supposer que $X$ est $k$-affinoïde. 

\medskip
Soit $\bf r$ un polyrayon déployant $X$ (rappelons que par {\em définition} un tel $\bf r$ est en particulier $k$-libre). On a $(k_{\bf r})^{1/p}=(k^{1/p})_{\bf r^{1/p}}$, et $(k_{\bf r})^{1/p}$ est donc une extension analytiquement séparable de $k^{1/p}$. On en déduit que $X_{(k_{\bf r})^{1/p}}$ satisfait $\mathsf P$. Soit $L$ une extension complète de $k$ et soit $K$ une extension complète de $k$ composée de $L$ et $k_{\bf r}$ ; si $X_K$ satisfait $\mathsf P$, alors $X_L$ satisfait $\mathsf P$. Quitte à remplacer $k$ par $k_{\bf r}$ et $X$ par $X_{\bf r}$, on peut par conséquent faire l'hypothèse que $|k^{*}|\neq\{1\}$ et que $X$ est strictement $k$-affinoïde. 

\medskip
Soit $Y$ le lieu de non-régularité de $X_{k^{1/p}}$. C'est un fermé de Zariski de $X_{k^{1/p}}$ dont on note $d$ la codimension de Krull. Son image $Z$ sur $X$ est un fermé de Zariski de $X$ de codimension de Krull égale à $d$ (lemme \ref{UNIREG}.\ref{homeozar}) ; pour toute extension complète $L$ de $k$, la codimension de Krull de $Z_L$ dans $X_L$ vaut également $d$.

\medskip
La propriété $\mathsf P$ est ou bien la régularité, auquel cas $Y$ et $Z$ sont vides, ou bien la propriété $R_m$ pour un certain $m\in \NN$, auquel cas $d>m$. Dans les deux situations, il suffit de démontrer que $X_L-Z_L$ est régulier pour toute extension complète $L$ de $k$, autrement dit que $X-Z$ est quasi-lisse. On va montrer que tout domaine strictement $k$-affinoïde $W$ de $X-Z$ est quasi-lisse, ce qui permettra de conclure. Soit $W$ un tel domaine et soit $\cal A$ l'algèbre de fonctions correspondante ; comme $W_{k^{1/p}}$ est inclus dans $X_{k^{1/p}}-Y$, il est régulier.

\medskip
\setcounter{cptbis}{0}
\trois{finsepfinie} {\em Si $F$ est une extension finie purement inséparable de $k$, l'algèbre ${\cal A}_{F}$ est régulière}. C'est évident si $p=1$ ; sinon, soit $(\alpha_{1},\ldots,\alpha_{n})$ une $p$-base de $F$ sur $k$. Pour tout $i$, notons $n_{i}$ le plus petit entier non nul tel que $\alpha_{i}^{p^{n_{i}}}\in k$, et $\beta_{i}$ l'élément $\alpha_{i}^{p^{n_{i}-1}}$ de $F$. Soit $F_{0}$ le sous-corps de $F$ engendré par $k$ et les $\beta_{i}$. Il se plonge dans $k^{1/p}$ ; on en déduit que ${\cal A}_{F_{0}}$ est régulier. Par ailleurs $\Omega^1_{{\cal A}_{F}/{\cal A}}$ est un ${\cal A}_{F}$-module libre de rang $n$ (les $\mathsf{d}\alpha_{i}$ en constituent une base), et $\Omega^1_{{\cal A}_{F_{0}}/{\cal A}}$ est un ${\cal A}_{F_{0}}$-module libre de rang $n$ (les $\mathsf{d}\beta_{i}$ en constituent une base). Par le critère de régularité de Kiehl (\cite{kie}, Satz. 2.2 ; {\em cf.} \cite{comprig}, th. 1.1.1), ${\cal A}_{F}$ est régulière.

\medskip
\trois{fquelconque} {\em Si $F$ est une extension finie de $k$, l'algèbre ${\cal A}_{F}$ est régulière}. Pour le vérifier, on peut agrandir $F$ et la supposer normale, elle admet alors un dévissage $k\subset F_{1}\subset F$, où $F_{1}$ est purement inséparable sur $k$ et où $F$ est séparable sur $F_{1}$. L'algèbre ${\cal A}_{F_{1}}$ est régulière d'après ce qui précède, et la ${\cal A}_{F_{1}}$-algèbre étale ${\cal A}_{F}$ est donc régulière. 

\medskip
\trois{concluunivreg} {\em Conclusion.} Soit $w$ un point rigide de $W$. Le~\ref{UNIREG}.\ref{kunsurp}.\ref{fquelconque} ci-dessus assure que $W_{\hres(w)}$ est régulier ; comme $W_{\hres(w)}$ possède un $\hres(w)$-point au-dessus de $w$, le lemme~\ref{UNIREG}.\ref{omegareg} entraîne la quasi-lissité de $W$ en $w$. Le lieu de non quasi-lissité de $W$ est un fermé de Zariski, dont on vient de voir qu'il ne contient aucun point rigide. Par le {\em Nullstellensatz}, il est vide.~$\Box$ 

\medskip
\deux{lieugeomreg} {\bf Corollaire.} {\em Soit $X$ un espace $k$-analytique et soit $\mathsf P$ une propriété appartenant à $\cal R$. Le lieu des points de $X$ en lesquels $\mathsf P$ est géométriquement satisfaite est un ouvert de Zariski de $X$ ; sa formation commute à toute extension du corps de base.}

\medskip
{\em Démonstration.} La première assertion découle du fait que le lieu de validité de $\mathsf P$ sur $X_{k^{1/p}}$ est un ouvert de Zariski de  ce dernier, de la proposition ci-dessus, et du lemme \ref{UNIREG}.\ref{homeozar} ; la seconde est triviale.~$\Box$

\subsection*{Caractère réduit et extension des scalaires}

{\em Pour tout entier positif $n$, on notera $k_n$ l'extension $k^{1/p^n}$ de $k$. L'on désignera par $k_\infty$ le complété de la clôture parfaite de $k$, c'est-à-dire de la réunion des $k_n$ ; le corps $k_\infty$ est parfait d'après le lemme} \ref{UNIREG}.\ref{clotparfparf}.

\medskip
Le théorème qui suit figure essentiellement dans l'article  \cite{comprig} de Conrad (lemme 3.3.1) ; il y est énoncé dans le cadre strictement analytique, mais ses arguments s'appliquent ici ; nous en donnons la preuve pour la commodité du lecteur. 

\medskip
\deux{defredfin} {\bf Théorème (Conrad).} {\em Soit $X$ un espace $k$-analytique quasi-compact. Il existe $n$ tel que l'espace $k_n$-analytique $Y:=(X_{k_n})\duit$ soit géométriquement réduit ; pour toute extension complète $L$ de $k_n$, l'espace $(X_L)\duit$ s'identifie à $Y_ L$. }

\medskip
{\em Démonstration.} Par le lemme \ref{COMP}.\ref{xredl}, la seconde assertion résulte de la première ; il en découle que si celle-ci est vraie pour un certain $n$, elle l'est pour tout entier supérieur ou égal à $n$. 

\medskip
Il suffit donc, pour établir le théorème, de démontrer la première assertion G-localement sur $X$ ; l'on peut dès lors supposer que ce dernier est $k$-affinoïde, et l'on note $\cal A$ son algèbre de fonctions.

\medskip
Pour tout $n$, notons ${\cal I}_n$ l'idéal de ${\cal A}_{k_n}$ formé des éléments nilpotents. Comme $({\cal I}_n{\cal A}_{k_\infty})$ est une suite croissante d'idéaux de l'anneau noethérien ${\cal A}_{k_\infty}$, il existe $n$ tel que ${\cal I}_{n+1}{\cal A}_{k_\infty }={\cal I}_n{\cal A}_{k_\infty }$.

\medskip
Par platitude de ${\cal A}_{k_{n+1}}\to {\cal A}_{k_{\infty}}$, l'idéal  ${\cal I}_{n+1}{\cal A}_{k_\infty}$ (resp. ${\cal I}_n{\cal A}_{k_\infty}$) s'identifie à ${\cal I}_{n+1}\otimes_{{\cal A}_{k_{n+1}}} {\cal A}_{k_\infty}$ (resp. (${\cal I}_n{\cal A}_{k_{n+1}})\otimes_{{\cal A}_{k_{n+1}}} {\cal A}_{k_\infty}$). 

\medskip
En conséquence, $({\cal I}_{n+1}/{\cal I}_n{\cal A}_{k_{n+1}})\otimes_{{\cal A}_{k_{n+1}}} {\cal A}_{k_\infty}=0$. La {\em fidèle} platitude de ${\cal A}_{k_{n+1}}\to {\cal A}_{k_\infty}$ assure alors que ${\cal I}_{n+1}={\cal I}_n{\cal A}_{k_{n+1}}$. 

\medskip
De ce fait $({\cal A}_{k_n}/{\cal I}_n)\otimes_{{\cal A}_{k_n}}{\cal A}_{k_{n+1}}$ coïncide avec ${\cal A}_{k_{n+1}}/{\cal I}_{n+1}$, et est dès lors réduit. Autrement dit, $((X_{k_n})\duit)_{k_{n+1}}$ est réduit ;  compte-tenu de la proposition \ref{UNIREG}.\ref{kunsurp} et du fait que le caractère réduit est, par le critère de Serre, la conjonction des propriétés $R_0$ et $S_1$, cela signifie exactement que $(X_{k_n})\duit$  est géométriquement réduit.~$\Box$

\subsection*{Normalité et extension des scalaires}

\medskip
Le théorème qui suit figure essentiellement dans l'article  \cite{comprig} de Conrad (th. 3.3.6) ; il y est énoncé dans le cadre strictement analytique, mais ses arguments s'appliquent ici ; nous en donnons la preuve pour la commodité du lecteur.

\medskip
\deux{defnormfin} {\bf Théorème (Conrad).} {\em Soit $X$ un espace $k$-analytique quasi-compact. Il existe  $n$ tel que la normalisation $Y$ de $X_{k_n}$, vue comme un espace $k_n$-analytique, soit géométriquement normale ; pour toute extension complète $L$ de $k_n$, la normalisation de $X_L$ s'identifie à $Y_ L$. }

\medskip
{\em Démonstration.} Par la proposition \ref{NOR}.\ref{normextscal}, la seconde assertion résulte de la première ; il en découle que si celle-ci est vraie pour un certain $n$, elle l'est pour tout entier supérieur ou égal à $n$. 

\medskip
Il suffit donc, pour établir le théorème, de démontrer la première assertion G-localement sur $X$ ; l'on peut dès lors supposer que ce dernier est $k$-affinoïde. Le théorème \ref{UNIREG}.\ref{defredfin} permet, quitte à remplacer $k$ par $k_n$ pour un entier $n$ convenable, de faire l'hypothèse que $X\duit$ est géométriquement réduit ; en substituant $X\duit$ à $X$, on se ramène finalement au cas où $X$ lui-même est géométriquement réduit ; on note $\cal A$ son algèbre de fonctions. Pour tout $n\in \NN\cup\{\infty\}$, on désigne par $Y_n$ la normalisation de $X_{k_n}$ et par ${\cal B}_n$ l'algèbre des fonctions analytiques sur $Y_n$ ; la proposition \ref{NOR}.\ref{normextscal} fournit pour tout couple $(n,m)$ d'éléments de $\NN\cup\{\infty\}$ tel que $n\leq m$ un morphisme fini et surjectif $Y_m\to Y_n\hotimes_{k_n}k_m$. 

\medskip
Fixons $n\in \NN$. Soit $U$ le lieu de normalité de $X_{k_n}$ ; c'est un ouvert de Zariski de $X_{k_n}$ qui est dense d'après le lemme \ref{NOR}.\ref{redpasnorm} ; par le corollaire \ref{NOR}.\ref{imrecouvdense}, son image réciproque $V$ sur $Y_n$ est dense dans $Y_n$ ; par construction, $Y_n\to X_{k_n}$ induit un isomorphisme $V\simeq U$.

\medskip
Il résulte de tout ceci que $V_{k_\infty}$ est un ouvert de Zariski dense de $Y_{n,k_\infty}$ et que $Y_{n,k_\infty}\to X_{k_\infty}$ induit un isomorphisme $V_{k_\infty}\to U_{k_\infty}$. Puisque $X$ est géométriquement réduit, $U_{k_\infty}$ est réduit ; dès lors, $V_{k_\infty}$ est réduit aussi, et comme il s'agit d'un ouvert de Zariski dense de de $Y_{n,k_\infty}$, ce dernier possède la propriété $R_0$. Par ailleurs l'espace $Y_n$, étant normal, possède par le critère de Serre la propriété $S_2$, et $Y_{n,k_\infty}$ la possède donc également ; il possède {\em a fortiori} la propriété $S_1$ ; en conclusion, $Y_{n,k_\infty}$ satisfait les propriétés $R_0$ et $S_1$ ; par le critère de Serre, il est réduit. 

\medskip
La surjectivité de $Y_\infty\to Y_{n,k_\infty}$, combinée au caractère réduit de $Y_{n,k_\infty}$, entraîne l'injectivité de ${\cal B}_n\otimes_{{\cal A}_{k_n}}{\cal A}_{k_\infty}\to {\cal B}_{\infty}$. 

\medskip
On dispose en particulier d'une suite croissante $({\cal B}_n\otimes_{{\cal A}_{k_n}}{\cal A}_{k_\infty})$ de sous-${\cal A}_{k_\infty}$-modules de ${\cal B}_\infty$. Par noethérianité de ce dernier (qui est un ${\cal A}_{k_\infty}$-module fini), cette suite est stationnaire. Il existe donc un entier $n$ tel que l'on ait l'égalité ${\cal B}_n\otimes_{{\cal A}_{k_n}}{\cal A}_{k_\infty}={\cal B}_{n+1}\otimes_{{\cal A}_{k_{n+1}}}{\cal A}_{k_\infty}$. Autrement dit, la flèche naturelle $${\cal B}_n\otimes_{{\cal A}_{k_n}}{\cal A}_{k_{n+1}}\to {\cal B}_{n+1}$$ devient un isomorphisme après application du foncteur $\bullet\otimes_{{\cal A}_{k_{n+1}}}{\cal A}_\infty$. La $\cal{A}_{k_{n+1}}$-algèbre ${\cal A}_\infty$ étant fidèlement plate, ${\cal B}_n\otimes_{{\cal A}_{k_n}}{\cal A}_{k_{n+1}}\simeq {\cal B}_{n+1}$ ; en conséquence, ${\cal B}_n\otimes_{{\cal A}_{k_n}}{\cal A}_{k_{n+1}}$ est normal.

\medskip
On vient de démontrer que $Y_{n,k_{n+1}}$ est normal ; compte-tenu de la proposition \ref{UNIREG}.\ref{kunsurp} et du fait que la normalité  est, par le critère de Serre, équivalente à la conjonction des propriétés $R_1$ et $S_2$, cela signifie exactement que l'espace $k_n$-analytique $Y_n$ est géométriquement normal.~$\Box$ 

\subsection*{Préservation de l'irréductibilité par passage à la clôture parfaite} 

\medskip
\deux{perfclos} {\bf Proposition.} {\em Soit $X$ un espace $k$-analytique irréductible. L'espace $X_{k_\infty}$ est irréductible.} 

\medskip
{\em Démonstration.} On procède en deux temps. On utilisera implicitement tout au long de la preuve le fait que pour tout espace $k$-analytique $W$, l'application naturelle $W_{k_\infty}\to W$ induit un homéomorphisme entre les espaces topologiques sous-jacents ({\em cf.} \ref{RAP}.\ref{sensxl}.\ref{raddensehomeo}).

\setcounter{cptbis}{0}
\medskip
\trois{perfcloscomp} {\bf Le cas où $X$ est quasi-compact.} Le théorème \ref{UNIREG}.\ref{defnormfin} assure l'existence d'un entier $n$ tel que la normalisation $Y$ de $X_{k_n}$ soit un espace $k_n$-analytique géométriquement normal. Comme $X$ est irréductible, $X_{k_n}$ l'est aussi en vertu du lemme \ref{UNIREG}.\ref{homeozar} ; il en résulte que $Y$ est connexe ; dès lors, $Y_{k_\infty}$ est connexe. Comme $Y$ est géométriquement normal, $Y_{k_\infty}$ s'identifie par la proposition \ref{NOR}.\ref{normextscal} à la normalisation de $X_{k_\infty}$ ; la connexité de $Y_{k_\infty}$ entraîne donc l'irréductibilité de $X_{k_\infty}$. 

\medskip
\trois{perfclosgen} {\bf Le cas général.} Choisissons un G-recouvrement $(X_i)$ de $X$ par des domaines affinoïdes ; pour tout $i$, notons $(X_{i,j})$ la famille des composantes irréductibles de $X_i$. Il découle du cas compact traité au \ref{UNIREG}.\ref{perfclos}.\ref{perfcloscomp} ci-dessus que $X_{i,j,k_\infty}$ est irréductible pour tout $(i,j)$ ; en conséquence, les $X_{i,j,k_\infty}$ sont, à $i$ fixé, les composantes irréductibles de $X_{i,k_\infty}$. 

\medskip
Soit $Z$ une composante irréductible de $X_{k_\infty}$ et soit $T$ son image sur $X$. Pour tout $i$, l'intersection $Z\cap X_{i,k_\infty}$ est réunion de certaines des $X_{i,j,k_\infty}$ ; on en déduit que $T\cap X_i$ est pour tout $i$ réunion de certaines des $X_{i,j}$. De ce fait, $T$ est un fermé de Zariski de $X$, qui contient au moins une composante irréductible d'un domaine affinoïde de $X$ (puisque $Z$ est non vide) ; l'espace $X$ étant irréductible, $T=X$ et l'on a donc $Z=X_{k_\infty}$ ; l'espace $X_{k_\infty}$ est par conséquent irréductible, ce qu'il fallait démontrer.~$\Box$

\section{Connexité et irréductibilité géométriques} \label{GEO}
\setcounter{cpt}{0}

{\em On désigne toujours par $p$ l'exposant caractéristique de $k$. On fixe une clôture algébrique $k^a$ de $k$ et l'on note $\ka$ son complété.}

\subsection*{Quelques généralités} 

{\em Il arrivera fréquemment, dans ce qui suit, que l'on traite implicitement les composantes irréductibles d'un espace analytique comme des sous-ensemble analytique fermés de celui-ci, sans prendre la peine d'en préciser la structure ; bien entendu, nous ne nous le permettrons que dans des situations où cette ambiguïté n'a aucune incidence.}

\medskip
\deux{conntestaff} {\bf Lemme.}  {\em Soit $X$ un espace $k$-analytique et soit $(X_i)$ un G-recouvrement de $X$ par des domaines analytiques connexes. Soit $L$ une extension complète de $k$ telle que $X_{i,L}$ soit connexe pour tout $i$. L'application $\pi_0(X_L)\to \pi_0(X)$ est bijective.} 

\medskip
{\em Démonstration.} Soit ${\cal R}$ (resp. ${\cal R}_L$) la relation d'équivalence sur $X$ (resp. $X_L$) définie comme suit : pour tout couple $(x,y)$ de points de $X$ (res. $X_L$) on a $x\;{\cal R}\;y$ (resp. $x\;{\cal R}_L\;y$) si et seulement si il existe une famille finie $i_1,\ldots,i_r$ d'indices et une famille finie $x=x_1,x_2,\ldots, x_r=y$ de points de $X$ (resp. $X_L$) telle que $x_t\in X_{i_t}$ (resp. $x_t\in X_{i_t,L}$) pour tout $t$. Chacun des $X_i$ (resp. des $X_{i,L}$) étant connexe, l'ensemble $\pi_0(X)$ (resp. $\pi_0(X_L)$) s'identifie naturellement au quotient $X/{\cal R}$ (resp. $X_L/{\cal R}_L$). 

\medskip
Pour tout couple d'indices $(i,j)$ l'application $$(X_{i,L}\cap X_{j,L})=(X_i\cap X_j)_L\to X_i\cap X_j$$ est surjective ; on en déduit que $X_L\to X$ induit une bijection $X/{\cal R}\simeq X_L/{\cal R}_L$, et partant une bijection 
$\pi_0(X_L)\simeq \pi_0(X)$.~$\Box$

\medskip
\deux{introcx} Soit $X$ un espace $k$-analytique {\em réduit}. Nous noterons $\got s(X)$ l'algèbre des fonctions analytiques sur $X$ qui annulent G-localement un polynôme non nul et {\em séparable}\footnote{Dans \cite{brk4}, Berkovich a introduit l'anneau $\got c(X)$ qui est défini de façon analogue, mais sans condition de séparabilité du polynôme ; pour les questions qui nous intéressent, il sera suffisant de travailler avec $\got s(X)$, qui a de surcroît l'avantage d'être, contrairement à $\got c(X)$, utilisable même lorsque $X$ n'est pas supposé réduit.} à coefficients dans $k$ ; il est immédiat que $V\mapsto \got s(V)$ est un sous-faisceau en algèbres de ${\cal O}_{X\grot}$, et que $X$ est vide (resp. connexe) si et seulement si $\got s(X)=0$ (resp. n'a pas d'idempotent non trivial). Les deux lemmes ci-dessous sont essentiellement dus à Berkovich dans le cas strictement analytique (\cite{brk4}, lemmes 8.1.1 et 8.1.4) ; nous les démontrons néanmoins pour la commodité du lecteur. 

\medskip
\deux{cxcorps} {\bf Lemme (Berkovich).} {\em Soit $X$ un espace $k$-analytique connexe et non vide. L'anneau $\got s(X)$ est une extension finie séparable de $k$.}

\medskip
{\em Démonstration.} Puisque $X$ est non vide la $k$-algèbre $\got s(X)$ est non nulle. Nos allons tout d'abord établir que chaque élément de $\got s(X)$ annule un polynôme irréductible et séparable à coefficients dans $k$, ce qui prouvera que $\got s(X)$ est une extension algébrique séparable de $k$.

\medskip
\setcounter{cptbis}{0}
\trois{qannulef} Soit $f\in \got s(X)$ et soit $\psi: X\to \Aff^{1,an}_k$ le morphisme qu'elle induit. Soit $(X_i)$ un G-recouvrement de $X$ par des domaines analytiques tels que pour tout $i$ la fonction $f_{|X_i}$ annule un polynôme non nul et séparable $P_i$ à coefficients dans $k$. Pour tout $i$, la flèche $\psi_{|X_i}: X_i\to  \Aff^{1,an}_k$ se factorise par un morphisme $\psi'_i $ de $X_i$ vers $\coprod {\cal M}(k[\tau]/(P_{i,j}))$ où les $P_{i,j}$ sont les facteurs irréductibles unitaires de $P_i$, qui sont séparables et tous de multiplicité $1$ ; pour tout $j$, l'image réciproque $X_{i,j}$ de $ {\cal M}(k[\tau]/(P_{i,j}))$ par $\psi'_i$ est un ouvert fermé, et donc un domaine affinoïde, de $X_i$ ; par la définition même de $X_{i,j}$, l'on a $P_{i,j}(f_{|X_{i,j}})=0$. La famille $(X_{i,j})$ constitue un G-recouvrement de $X$. Si $(i,j)$ et $(i',j')$ sont tels que $P_{i,j}\neq P_{i',j'}$ l'existence d'une relation de Bezout entre $P_{i,j}$ et $P_{i',j'}$ force $X_{i,j}\cap X_{i',j'}$ à être vide. 

\medskip
Comme $X\neq \emptyset$ il existe $(i,j)$ tel que $X_{i,j}\neq\emptyset$ ; nous allons montrer que 
$P_{i,j}(f)=0$. Il suffit de vérifier que $P_{i,j}(f_{|X_{i',j'}})=0$ pour tout $(i',j')$. Soit donc $(i',j')$ un couple d'indices. Si $X_{(i',j')}=\emptyset$ il n'y a rien à démontrer. Sinon, il existe par connexité de $X$ une famille finie $(i,j)=(i_1,j_1),(i_2,j_2),\ldots,(i_r,j_r)=(i',j')$ telle que $X_{(i_l,j_l)}\cap X_{(i_{l+1},j_{l+1})}\neq \emptyset$ pour tout $l$ compris entre $1$ et $r-1$. Par ce qui précède, on a $P_{(i_l,j_l)}=P_{(i_{l+1},j_{l+1})}$ pour tout $l$ compris entre $1$ et $r-1$ ; par conséquent, $P_{i,j}=P_{i',j'}$ et $P_{i,j}(f_{|X_{i',j'}})=0$, ce qu'on souhaitait établir. 

\medskip
Ainsi, $f$ annule bien le polynôme irréductible et séparable $P_{i,j}$ et $\got s(X)$ est dès lors une extension algébrique séparable de $k$. 

\medskip
\trois{cxfinisurk} Montrons que $\got s(X)$ est de dimension finie sur $k$. Soit $V$ un domaine affinoïde non vide de $X$ et soit $\bf r$ un polyrayon déployant $V$. Par le {\em Nullstellensatz} l'espace $V_{\bf r}$ possède un point $y$ qui est $k_{\bf r}$-rigide. Si $E$ est une extension finie de $k$ contenue dans $\got s(X)$ alors $E_{\bf r}=E\otimes_kk_{\bf r}$ est un corps, et l'évaluation $\got s(X)\hookrightarrow \hres(y)$ induit un $k_{\bf r}$-morphisme $E_{\bf r}\hookrightarrow \hres(y)$. On a en conséquence $\dim k E=\dim {k_{\bf r}}E_{\bf r}\leq \dim {k_{\bf r}}\hres(y)$. Ceci valant pour toute extension finie de $k$ contenue dans l'extension algébrique $\got s(X)$ de $k$, cette dernière est elle-même finie.~$\Box$

\medskip
\deux{cxfonct} Si $X$ un espace $k$-analytique, si $L$ est une extension complète de $k$, si $Y$ est un espace $L$-analytique et si $\pi : Y\to X$ est un morphisme la flèche naturelle ${\cal O}_{X\grot}(X)\to {\cal O}_{Y\grot}(Y)$ envoie $\got s(X)$ dans $\got s(Y)$. 

\medskip
\deux{cxred} Si $X$ est un espace $k$-analytique le morphisme naturel de ${\cal O}_{X\grot}(X)$ vers ${\cal O}_{X_{\rm red, G}}(X\duit)$ induit un isomorphisme $\got s(X)\simeq \got s(X\duit)$ ; on le vérifie G-localement à l'aide de la propriété d'unique relèvement infinitésimal des morphismes étales. 

\medskip
\deux{remcxcorps} {\bf Remarque.} Soit $X$ un espace $k$-analytique connexe et non vide, et soit $F$ un sous-corps de $\got s(X)$ contenant $k$ ; on peut voir $X$ comme un espace $F$-analytique. Il est immédiat que le corps $\got s(X)$ est indépendant du fait que l'on considère $X$ comme un espace $k$-analytique ou comme un espace $F$-analytique.

\medskip
\deux{remcxr} {\bf Remarque.} Soit $X$ un espace $k$-analytique connexe et non vide et soit $\bf r$ un polyrayon $k$-libre. L'espace $k_{\bf r}$-analytique $X_{\bf r}$ est lui aussi non vide ; il est connexe comme on le voit en se ramenant au cas affinoïde (c'est possible grâce au lemme \ref{GEO}.\ref{conntestaff}) et en utilisant alors alors le lemme \ref{RAP}.\ref{tenseurkrintegre}. Dès lors $\got s(X_{\bf r})$ est bien défini, et est une extension finie séparable de $k_{\bf r}$ d'après le lemme \ref{GEO}.\ref{cxcorps}. 

\medskip
\deux{cxcompxf} {\bf Lemme (Berkovich).} {\em Soit $X$ un espace $k$-analytique connexe et non vide. Les propositions suivantes sont équivalentes : 

\medskip
$i)$ pour toute extension finie séparable $F$ de $k$, l'espace $X_F$ est connexe ;

\medskip
$ii)$ pour toute extension finie $F$ de $k$, l'espace $X_F$ est connexe ; 

\medskip
$iii)$ $\got s(X)=k$.} 

\medskip
{\em Démonstration.} Montrons que $i)\Rightarrow ii)$. On suppose donc que $i)$ est vraie. Soit $F$ une extension finie de $k$ et soit $F\sep$ la fermeture séparable de $k$ dans $F$. Par hypothèse, $X_{F\sep}$ est connexe ; comme le morphisme fini $X_F\to X_{F\sep}$ induit un homéomorphisme entre les espaces topologiques sous-jacents ({\em cf.} \ref{RAP}.\ref{sensxl}.\ref{raddensehomeo}) $X_F$ est connexe.

\medskip
Montrons que $ii)\Rightarrow iii)$. On suppose donc que $ii)$ est vraie. L'espace $X$ est muni d'une flèche surjective vers ${\cal M}(\got s(X))$ ; dès lors $X_{\got s(X)}$ est  muni d'une flèche surjective vers ${\cal M}(\got s(X)\otimes_k \got s(X))$. La connexité de $X_{\got s(X)}$ entraîne donc celle de ${\cal M}(\got s(X)\otimes_k \got s(X))$, c'est-à-dire encore celle du schéma $\spec (\got s(X)\otimes_k \got s(X))$ ; il en résulte que $\got s(X)=k$. 

\medskip
Montrons que $iii)\Rightarrow i)$. Supposons que $iii)$ est vraie. Soit $F$ une extension finie séparable de $k$ ; quitte à agrandir $F$,  on peut supposer que c'est une extension galoisienne de $k$. Posons $G=\mbox{Gal}(F/k)$. 

\setcounter{cptbis}{0}
\medskip
\trois{galtranscomp} Le groupe $G$ opère transitivement sur toute  fibre de $X_F\to X$ (la fibre de ce morphisme en un point $x$ de $X$ s'identifiant à ${\cal M}(F\otimes_k\hres(x))$ ) ; de plus, $X_F\to X$ est fini et plat, et donc ouvert et fermé ; chacune des composantes connexes de $X_F$ se surjecte par conséquent sur $X$ ; on déduit de tout ceci que $G$ opère transitivement sur $\pi_0(X_F)$. Notons que ce fait entraîne la finitude de $\pi_0(X_F)$.

\medskip
\trois{fixestab} Soient $Y_1,\ldots,Y_r$ les composantes connexes de $X_F$. Soit $H$ le stabilisateur de $Y_1$ dans $G$. Soit $f\in F^H$ ; l'action de $G$ sur $\pi_0(X_F)$ étant transitive, il existe pour tout $i>1$ un élément $g_i$ de $G$ tel que $g_i(Y_1)=Y_i$. La fonction $(f,g_2(f),\ldots,g_r(f))\in {\cal O}_{X_{F,\rm G}}(X_F)=\prod {\cal O}_{X_{F,\rm G}}(Y_i)$ ne dépend pas du choix des $g_i$ et est invariante sous l'action de $G$. Elle provient de ce fait d'une fonction $f_0\in {\cal O}_{X\grot}(X)$ : on le voit en se ramenant au cas affinoïde, auquel cas cela résulte de la théorie de Galois {\em schématique}. Si $P$ désigne le polynôme minimal sur $k$ de l'élément $f$ de $F$ alors $P$ annule $(f,g_2(f),\ldots,g_r(f))$ et donc $f_0$, par injectivité de ${\cal O}_{X\grot}(X)\to {\cal O}_{X_{F,\rm G}}(X_F)$ (laquelle se démontre par fidèle platitude après réduction triviale au cas affinoïde) ; ainsi, $f_0$ appartient à $\got s(X)$ et donc à $k$ en vertu de l'hypothèse $iii)$ purement inséparable sur $k$ ; par conséquent, $f\in k$. 

\medskip
\trois{conclugalcomp}{\em Conclusion.} Le sous-corps $F^H$ de $F$ est d'après ce qui précède égal à $k$, ce qui signifie que $H=G$. Ce dernier fixe par conséquent $Y_1$ ; puisqu'il agit par ailleurs transitivement sur $\pi_0(X_F)$, on a $\pi_0(X_F)=\{Y_1\}$ ; ainsi, $X_F$ est connexe, ce qu'il fallait démontrer.~$\Box$

\medskip
\deux{extconnexe} {\bf Lemme.} {\em Soit $X$ un espace $k$-analytique et soit $F$ une extension complète de $k$. Les propositions suivantes sont équivalentes : 

\medskip
$i)$ pour toute extension complète $L$ de $k$, l'espace $X_L$ est connexe (resp. irréductible) ;

\medskip
$ii)$ pour toute extension complète $L$ de $F$, l'espace $X_L$ est connexe (resp. irréductible).}

\medskip
{\em Démonstration.} Il est clair que $i)\Rightarrow ii)$. Supposons que $ii)$ est vraie, et soit $L$ une extension complète de $k$. Soit $K$ une extension complète de $k$ composée de $L$ et $F$. L'espace $X_K$ est connexe (resp. irréductible) puisque l'on suppose que $ii)$ est vraie. Comme $X_K\to X_L$ est surjective, $X_L$ est connexe (resp. irréductible).~$\Box$

\subsection*{Étape technique : preuve de certains résultats sous des hypothèses restrictives}

\medskip
Les résultats qui suivent seront établis plus tard en toute généralité, par une méthode consistant essentiellement à se ramener aux cas particuliers traités ci-dessous. 

\medskip
\deux{kafini} {\bf Lemme (Berkovich).} {\em Soit $X$ un espace $k$-analytique.

\medskip
$i)$ Si $X_{\ka}$ est connexe, alors $X_F$ est connexe pour toute extension finie $F$ de $k$ ;

\medskip
$ii)$ si $X_F$ est connexe pour toute extension finie $F$ de $k$ et si $X$ possède un point rigide, alors $X_{\ka}$ est connexe.}

\medskip
{\em Démonstration.} Montrons $i)$. Supposons que $X_{\ka}$ est connexe et soit $F$ une extension finie de $k$. Choisissons un plongement $F\hookrightarrow k^a$. L'application continue induite $X_{\ka}\to X_F$ est surjective et $X_{\ka}$ est connexe par hypothèse ; il en résulte que $X_F$ est connexe.

\medskip
Montrons $ii)$. On suppose que $X_F$ est connexe pour toute extension finie $F$ de $k$ et que $X$ possède un point rigide $x$. Si $Y$ est un domaine analytique {\em compact} de $X$, la bijection continue naturelle $Y_{\ka}\to \lim\limits_\leftarrow Y_F$, où $F$ parcourt l'ensemble $\cal S$ des extensions finies  de $k$ incluses dans $\ka$, est un homéomorphisme (\ref{RAP}.\ref{sensxl}.\ref{densehomeolimproj}).

Pour tout couple $(F,F')\in {\cal S}^2$ avec $F\subset F'$, l'application $Y_{F'}\to Y_F$ est finie et plate, et en particulier ouverte ; il en découle formellement que pour tout $K\in \cal S$ l'application $Y_{\ka}\simeq  \lim\limits_\leftarrow Y_F\to Y_K$ est ouverte ; elle est par ailleurs fermée compte-tenu du fait que $Y_{\ka}$ est compact et que $Y_F$ est topologiquement séparé ($Y$ lui-même l'étant par {\em définition} de la compacité).

\medskip
Ceci valant pour tout domaine analytique compact $Y$ de $X$, l'application $X_{\ka}\to X_F$ est ouverte et fermée pour tout $F\in \cal S$. Il existe $F\in{\cal S}$ et un point $x'$ sur $X_F$ situé au-dessus de $x$ tel que $\hres(x')=F$. Supposons que $X_{\ka}$ soit réunion disjointe de deux ouverts non vides $U_1$ et $U_2$. Par hypothèse, $X_F$ est connexe. L'application $X_{\ka}\to X_F$ étant ouverte et fermée, $U_1\to X_F$ et $U_2\to X_F$ sont surjectives. Il en résulte que $x'$ a au moins deux antécédents sur $X_{\ka}$. Or comme $\hres(x')=F$, le point  $x'$ n'a qu'un antécédent sur $X_{\ka}$, ce qui est absurde. On en déduit que $X_{\ka}$ est connexe. ~$\Box$ 

\medskip
Le théorème qui suit est dû à Conrad (\cite{comprig}, th. 3.2.1) ; la preuve que nous proposons est nouvelle. 

\medskip
\deux{xkaconnexe} {\bf Théorème (Conrad).} {\em Supposons que $|k^*|\neq\{1\}$. Soit $X$ un espace strictement $k$-affinoïde. Les propositions suivantes sont équivalentes :

\medskip

$i)$ $X_L$ est connexe pour toute extension complète $L$ de $k$ ;

\medskip
$ii)$ $X_{\ka}$ est connexe.} 

\medskip
{\em Démonstration.} Il est immédiat que $i)\Rightarrow ii)$. Montrons la réciproque. On suppose donc que $ii)$ est vraie. En vertu du lemme \ref{GEO}.\ref{extconnexe}, il suffit de montrer que $X_L$ est connexe pour toute extension complète $L$ de $\ka$. On peut donc, quitte à remplacer $k$ par $\ka$, supposer que $k$ est algébriquement clos. On fixe une extension complète $L$ de $k$. On note $\cal A$ l'algèbre des fonctions $k$-analytiques sur $X$.

\medskip
On raisonne par l'absurde. Supposons que $X_L$ n'est pas connexe. L'anneau ${\cal A}_L$ possède dès lors un idempotent non trivial $g$ ; l'on a donc $g^2-g=0, ||g||=1$ et $||g-1||=1$, où $||.||$ désigne la (semi)-norme spectrale de ${\cal A}_L$. Soit $\eta$ un réel appartenant à $]0;1[\cap |k^*|$. Le sous-anneau de ${\cal A}_L$ engendré par $\cal A$ et $L$ est dense relativement à $||.||$, et chacun de ses éléments vit dans un sous-anneau de ${\cal A}_L$ engendré par $\cal A$ et un sous-corps de $L$ de type fini sur $k$. Il existe en conséquence : 

\medskip
$\bullet$ un sous-corps $L_0$ de $L$ de type fini sur $k$, dont on note $\widehat{L_0}$ le complété ;

\medskip
$\bullet$ un élément $g_0$ de ${\cal A}_{\widehat{L_0}}$ tel que $||g_0^2-g_0||<\eta, ||g_0||=1$ et $||g_0-1||=1$. 

\medskip
Précisons, concernant cette dernière inégalité, qu'elle repose sur le fait que si $a\in {\cal A}_{L_0}$ alors la semi-norme spectrale de l'élément $a$ de ${\cal A}_{L_0}$  coïncide avec celle de son image dans ${\cal A}_L$, en raison par exemple de la surjectivité de la flèche ${\cal M}({\cal A}_L)\to {\cal M}({\cal A}_{L_0})$. 

\medskip
Soit $\cal Y$ un $k$-schéma intègre et de type fini dont le corps des fonctions est isomorphe à $L_0$. Le complété $\widehat{L_0}$ s'identifie à $\hres(y)$ pour un certain $y\in{\cal Y}\an$ situé au-dessus du point générique de $\cal Y$ ; l'algèbre ${\cal A}_{\widehat{L_0}}$ est naturellement isomorphe à l'algèbre des fonctions analytiques sur $q^{-1}(y)$ où $q: {\cal Y}\an\times_k X\to {\cal Y}\an$ est la première projection.

\medskip
En vertu de la densité de $\kappa(y)$ dans $\hres(y)$, il existe un voisinage strictement $k$-affinoïde $Y$ de $y$ dans ${\cal Y}\an$, et une fonction $\gamma$ sur $Y\times_k X$, tels que l'on ait $$||\gamma_{|q^{-1}(y)}^2-\gamma_{|q^{-1}(y)}||<\eta, ||\gamma_{|q^{-1}(y)}||=1\;\mbox{et}\;||\gamma_{|q^{-1}(y)}-1||=1.$$

\medskip
Soit $U_0$ (resp. $U_1$) l'ouvert de $Y\times_k X$ défini par la conjonction d'inégalités $|\gamma^2-\gamma|<\eta$ et $|\gamma|<\eta$ (resp. $|\gamma^2-\gamma|<\eta$ et $|\gamma-1|<\eta$). Soit $Z$ le domaine affinoïde de $Y\times_k X$ défini par l'inégalité $|\gamma^2-\gamma|\geq \eta$. Par choix de $\gamma$, les trois conditions suivantes sont satisfaites : 

\medskip
$\bullet$ pour tout $x\in q^{-1}(y)$ l'on a $|\gamma^2(x)-\gamma(x)|<\eta$ ; ainsi, $y\notin q(Z)$ ;

\medskip
$\bullet$ il existe $x\in q^{-1}(y)$ tel que $|\gamma(x)|=1$ ; comme $|\gamma^2(x)-\gamma(x)|<\eta$, l'on a $|\gamma(x)-1|<\eta$ ; ainsi, $y\in q(U_1)$ ;

\medskip
$\bullet$ il existe $x\in q^{-1}(y)$ tel que $|\gamma(x)-1|=1$ ; comme $|\gamma^2(x)-\gamma(x)|<\eta$, l'on a $|\gamma(x)|<\eta$ ; ainsi, $y\in q(U_0)$.

\medskip
Soit $T$ le complémentaire du compact $q(Z)$ de $Y$ et soit $T'$ un voisinage strictement $k$-affinoïde de $y$ dans $T$. Comme $y\in q(U_0)\cap q(U_1)\cap T'$, il existe un domaine strictement $k$-affinoïde $U'_0$ (resp. $U'_1$) de $U_0$ (resp. $U_1$) qui est inclus dans $q^{-1}(T')$ et qui rencontre $q^{-1}(y)$. L'intersection $q(U'_0)\cap q(U'_1)$ est donc une partie non vide de $T'$ ; l'algèbre strictement $k$-affinoïde ${\cal A}_{U'_0}\hotimes_{{\cal A}_{T'}}{\cal A}_{{U'_1}}$ est de ce fait non nulle. Par le {\em Nullstellensatz}, elle possède un quotient isomorphe à $k$. Ceci entraîne l'existence d'un $k$-point $t$ appartenant à $$T'\cap q(U'_0)\cap q(U'_1)\subset T\cap q(U_0)\cap q(U_1).$$ 

\medskip
Soit $x\in q^{-1}(t)$. Comme $t\in T$, l'on a $|\gamma^2(x)-\gamma(x)|<\eta$ ; on en déduit (compte-tenu du fait que $\eta<1$ et grâce à l'inégalité ultramétrique) que l'on a $|\gamma(x)|<\eta$ ou $|\gamma(x)-1|<\eta$, ces deux inégalités étant exclusives l'une de l'autre. La fibre $q^{-1}(t)$ est donc égale à la réunion disjointe de ses ouverts $q^{-1}(t)\cap U_0$ et $q^{-1}(t)\cap U_1$. Le point $t$ appartenant à $q(U_0)\cap q(U_1)$, lesdits ouverts sont tous deux non vides. En conséquence, $q^{-1}(t)$ n'est pas connexe. Or $t$ est un $k$-point ; dès lors, $q^{-1}(t)\simeq X$, lequel est connexe par hypothèse. On aboutit ainsi à une contradiction.~$\Box$ 

\subsection*{Extension des scalaires à $k_{\bf r}$}

\deux{cxrcxr} {\bf Proposition.} {\em Soit $X$ un espace $k$-affinoïde connexe et non vide ; soit $\DD$ un produit fini de $k$-disques et de $k$-couronnes compacts, et soient $\cal A$ et $\cal B$ les algèbres de fonctions analytiques respectives de $X$ et $\DD$ ; on note $\bf T$ la famille des fonctions coordonnées sur $\DD$.  Supposons qu'il existe une $\cal B$-algèbre finie étale $\cal C$ telle que $X\times_k \DD\to \DD$ se factorise par ${\cal M}({\cal C})$. Si $c\in \cal C$ et si l'on écrit son image dans ${\cal A}\hotimes_k {\cal B}$ sous la forme $\sum a_I{\bf T}^I$, où les $a_I$ sont des éléments de $\cal A$, alors chacun des $a_I$ appartient à $\got s(X)$.}

\medskip
{\em Démonstration.} On procède en plusieurs étapes.

\medskip
\setcounter{cptbis}{0}
\trois{cxrcxrstr} {\em Le cas où $|k^*|\neq \{1\}$ et où $X$ est strictement $k$-affinoïde.} Soit $\DD'$ l'espace ${\cal M}({\cal C})$. Comme $X\to {\cal M}(k)$ se factorise par ${\cal M}(\got s(X))$, l'on dispose d'un diagramme commutatif $\diagram & \DD'_{\got s(X)}\ddto \rto& \DD_{\got s(X)}\ddto \\ X\times_k \DD\urto\drto& &\\& \DD'\rto & \DD\enddiagram$, où le morphisme composé $X\to \DD_{\got s(X)}$ est la projection naturelle $X\times_k\DD\simeq X\times_{\got s(X)}\DD_{\got s(X)}\to \DD_{\got s(X)}$. On peut donc, en se fondant sur la remarque \ref{GEO}.\ref{remcxcorps}, remplacer $k$, $\DD$ et $\DD'$ par $\got s(X)$, $\DD_{\got s(X)}$ et $\DD'_{\got s(X)}$, et se ramener ainsi au cas où $\got s(X)=k$ ; il découle alors du lemme \ref{GEO}.\ref{cxcompxf}, du lemme \ref{GEO}.\ref{kafini} et du théorème \ref{GEO}.\ref{xkaconnexe} que  $X_L$ est connexe pour toute extension complète $L$ de $k$.

\medskip
Soit $\DD''$ l'image (ensembliste) de $X\times_k\DD$ sur $\DD'$. Soit $\omega\in \DD''$, soit $z$ un antécédent de $\omega$ sur $X\times_k \DD$ et soit $t$ l'image de $\omega$ sur $\DD$. La fibre $Z$ de $X\times_k \DD$ en $t$ est un espace $\hres(t)$-analytique qui s'identifie à $X_{\hres(t)}$ ; en conséquence $Z_L$ est connexe pour toute extension complète $L$ de $\hres(t)$. 

\medskip
Soient $\omega=\omega_1,\ldots,\omega_n$ les antécédents de $t$ sur $\DD'$ ; l'espace $\hres(t)$-analytique $Z$ se factorise par $\coprod {\cal M}(\hres(\omega_i))$ ; comme il est connexe et comme $z$ s'envoie sur $\omega$, il se factorise par ${\cal M}(\hres(\omega))$ ; ceci entraîne que $\omega$ est le seul antécédent de $t$ sur $\DD'$ qui appartienne à $\DD''$ ; en conséquence, $\DD''\to \DD$ est injective.

\medskip
Comme $Z_L$ est connexe pour toute extension complète $L$ de $\hres(t)$, le lemme  \ref{GEO}.\ref{cxcompxf} assure que l'extension finie séparable $\hres(\omega)$ de $\hres(t)$ est triviale. Il existe en conséquence un voisinage affinoïde $V$ de $\omega$ dans $\DD'
$ et un voisinage affinoïde connexe  $U$ de $t$ dans $\DD$ tel que $\DD'\to \DD$ induise un isomorphisme $V\simeq U$ ; l'espace $\DD'\times_\DD U$ étant fini et étale sur $U$, le domaine affinoïde $V$ est nécessairement une composante connexe de $\DD'\times_\DD U$. La projection  $X\times_k U\to U$ est compacte, à image connexe, et ses fibres sont connexes, chacune d'elle étant de la forme $X_L$ pour une certaine extension complète $L$ de $k$ ; dès lors  $X\times_k U$ est connexe. L'image réciproque de $V$ sur ce dernier en est un ouvert fermé qui est non vide, puisqu'il contient $z$. C'est en conséquence $X\times_k U$ tout entier. 

\medskip
Cela signifie que $X\times_k U\to U$ se factorise par $V$ ; la flèche $V\to U$ étant un isomorphisme, $X\times_kU\to V$ est surjectif ; de ce fait $V\subset \DD''$, ce qui montre à la fois que $\DD''$ est un ouvert de $\DD'$ et que $\DD''\to \DD$ est un isomorphisme local.

\medskip
Par ailleurs, on a vu plus haut que $\DD''\to \DD$ est injectif, et $\DD''\to \DD$ est surjectif puisque $X\times_k \DD\to \DD$ se factorise par $\DD''$ ; le morphisme $\DD''\to \DD$ est un isomorphisme local bijectif, c'est donc un isomorphisme. La flèche $X\times_k \DD\to \DD'$ se factorise en conséquence par $\DD$ ; autrement dit, ${\cal C}\to{\cal A}\hotimes_k\cal B$ se factorise par $\cal B$ ; l'image de la $\cal B$-algèbre $\cal C$ dans  ${\cal A}\hotimes_k {\cal B}$ est dès lors égale à $\cal B$, c'est-à-dire à l'ensemble des séries qui s'écrivent $\sum a_I{\bf T}^I$ avec les $a_I$ dans $k$ ; la démonstration dans le cas où $|k^*|\neq \{1\}$ et où $X$ est strictement $k$-affinoïde est donc terminée. 

\medskip
\trois{cxrcxrgen} {\em Le cas général.} Soit $c\in \cal C$, écrivons son image dans $ {\cal A}\hotimes_k {\cal B}$ sous la forme $\sum a_I {\bf T}^I$ où les $a_I$ appartiennent à $\cal A$. Soit $\bf r$ un polyrayon déployant $X$. Fixons $I$. Il résulte du cas particulier traité au \ref{GEO}.\ref{cxrcxr}.\ref{cxrcxrstr} ci-dessus que $a_I$, vu comme un élément constant de ${\cal A}_{\bf r}$, annule un polynôme unitaire séparable à coefficients dans $k_{\bf r}$. Soit $\mathsf A$ la sous-algèbre de $\cal A$ engendrée par $a_I$ ; le produit tensoriel $\mathsf A\otimes_k k_{\bf r}$ s'injecte dans ${\cal A}_{\bf r}$ (on peut le voir directement à l'aide de l'expression des éléments de ${\cal A}_{\bf r}$ comme des séries, ou bien utiliser le lemme \ref{RAP}.\ref{prodalgplat}) et s'identifie à la sous-algèbre de cette dernière engendrée par $a_I$ ; il en résulte que la $k_{\bf r}$-algèbre $\mathsf A\otimes_k k_{\bf r}$ est finie étale ; en conséquence, $\mathsf A$ est une $k$-algèbre finie étale.~$\Box$ 

\medskip
\deux{cxextkr}{\bf Corollaire.} {\em Soit $X$ un espace $k$-analytique connexe et non vide et soit $\bf r$ un polyrayon $k$-libre. Le corps $\got s(X_{\bf r})$ est naturellement isomorphe à $\got s(X)_{\bf r}$.}

\medskip
{\em Démonstration.} Rappelons que $X_{\bf r}$ est lui-même non vide et connexe, et que $\got s(X_{\bf r})$ est donc bien défini, et est une extension finie séparable de $k_{\bf r}$ (remarque \ref{GEO}.\ref{remcxr}). Soit $f\in \got s(X_{\bf r})$ ; écrivons-la $f=\sum a_I{\bf T}^I$, où les $a_I$ sont des fonctions analytiques sur $X$. Le but est de montrer que chacun des $a_I$ appartient à $\got s(X)$ ; c'est une propriété G-locale, et on peut donc supposer que $X$ est affinoïde. Le morphisme $X_{\bf r}\to {\cal M}(k_{\bf r})$ se factorise par ${\cal M}(k_{\bf r}[f])$ ; comme $k_{\bf r}[f]$ est une $k_{\bf r}$-algèbre finie étale, la proposition ci-dessus assure que $a_I\in\got s(X)$ pour tout $I$.~$\Box$

\subsection*{La connexité géométrique}
\medskip
\deux{conngeomfond} {\bf Théorème.} {\em Soit $X$ un espace $k$-analytique. Les propositions suivantes sont équivalentes.

\medskip
$i)$ $\got s(X)=0$ ou $\got s(X)=k$ ;

\medskip
$ii)$ il existe une extension complète algébriquement close $\KK$ de $k$ telle que $X_\KK$ soit connexe ; 

\medskip
$iii)$ $X_L$ est connexe pour toute extension complète $L$ de $k$ ;

\medskip
$iv)$ $X_{\ka}$ est connexe ; 

\medskip
$v)$ $X_F$ est connexe pour toute extension finie $F$ de $k$ ;

\medskip
$vi)$ $X_F$ est connexe pour toute extension finie séparable $F$ de $k$.}

\medskip
{\em Démonstration.} Montrons que $i)\Rightarrow ii)$. On suppose que $i)$ est vraie. Si $\got s(X)=0$ alors $X$ est vide et $ii)$ est vraie. Supposons que  $\got s(X)=k$. Comme $\got s(X)$ est un corps, $X$ est connexe et non vide.  Il existe de ce fait un polyrayon $k$-libre $\bf r$ tel que $X_{\bf r}$ possède un point $k_{\bf r}$-rigide. Le corollaire \ref{GEO}.\ref{cxextkr} permet d'identifier $\got s(X_{\bf r})$ à $\got s(X)_{\bf r}$ ; ainsi $\got s(X_{\bf r})=k_{\bf r}$. Si $\KK$ désigne le complété d'une clôture algébrique de $k_{\bf r}$, les lemmes \ref{GEO}.\ref{cxcompxf} et  \ref{GEO}.\ref{kafini} assurent que $X_{\KK}$ est connexe, ce qui prouve $ii)$. 

\medskip
Montrons que $ii)\Rightarrow iii)$. On suppose que $ii)$ est vraie. Grâce au lemme \ref{GEO}.\ref{extconnexe} il suffit de montrer que $X_L$ est connexe pour toute extension complète $L$ de $\KK$ ; on peut donc supposer $k$ algébriquement clos et $X$ connexe. Le lemme \ref{GEO}.\ref{conntestaff} permet de se ramener au cas où $X$ est $k$-affinoïde. Si $X$ est vide, il n'y a rien à démontrer. Supposons $X\neq \emptyset$ ; remarquons qu'alors $\got s(X)=k$ puisque $k$ est algébriquement clos. Soit $\bf r$ un polyrayon déployant $X$. Le corollaire \ref{GEO}.\ref{cxextkr} assure que $\got s(X_{\bf r})=\got s(X)_{\bf r}=k_{\bf r}$. Soit $\Bbb F$ le complété d'une clôture algébrique de $k_{\bf r}$. En vertu des lemmes \ref{GEO}.\ref{cxcompxf} et  \ref{GEO}.\ref{kafini}, l'espace $X_{\Bbb F}$ est connexe ; par le théorème \ref{GEO}.\ref{xkaconnexe}, l'espace $X_L$ est connexe pour tout extension complète $L$ de $k_{\bf r}$. Le lemme \ref{GEO}.\ref{extconnexe} garantit alors que $X_L$ est connexe pour tout extension complète $L$ de $k$ ; ainsi, $iii)$ est vérifiée.

\medskip
L'implication $iii)\Rightarrow iv)$ est triviale ; l'implication $iv)\Rightarrow v)$ est l'assertion $i)$ du lemme \ref{GEO}.\ref{kafini} ; l'implication $v)\Rightarrow vi)$ est triviale ; l'implication $vi)\Rightarrow i)$ est triviale si $X$ est vide et résulte sinon du lemme \ref{GEO}.\ref{cxcompxf}.~$\Box$ 

\medskip
\deux{defgeomconn} On dira qu'un espace $k$-analytique $X$ est {\em géométriquement connexe} s'il satisfait les conditions équivalentes du théorème ci-dessus. 

\subsection*{L'irréductibilité géométrique} 

\deux{irrgeomfond} {\bf Théorème.} {\em Soit $X$ un espace $k$-analytique et soit $X'$ sa normalisation. Les propositions suivantes sont équivalentes.

\medskip
$i)$ $\got s(X')=k$ ;

\medskip
$ii)$ il existe une extension complète algébriquement close $\KK$ de $k$ telle que $X_\KK$ soit irréductible ; 

\medskip
$iii)$ $X_L$ est irréductible pour toute extension complète $L$ de $k$ ;

\medskip
$iv)$ $X_{\ka}$ est irréductible ; 

\medskip
$v)$ $X_F$ est irréductible pour toute extension finie $F$ de $k$ ;

\medskip
$vi)$ $X_F$ est irréductible pour toute extension finie séparable $F$ de $k$.}

\medskip
{\em Démonstration.} Elle repose essentiellement sur le théorème précédent. On procède en plusieurs étapes. 

\medskip
\setcounter{cptbis}{0} 
\trois{unequivsix} {\em Montrons que $i)\iff vi)$}. Si $F$ est une extension finie séparable de $k$, elle est analytiquement séparable et $X'_F$ est donc normal ; d'après la proposition \ref{NOR}.\ref{normextscal}, $X'_F$ s'identifie alors à la normalisation de $X_F$, et $X_F$ est dès lors irréductible si et seulement si $X'_F$ est connexe et non vide. L'assertion $vi)$ équivaut par conséquent à dire que $X'_F$ est connexe et non vide pour toute extension finie séparable $F$ de $k$ ; et cette dernière affirmation est-elle même vraie, par l'équivalence entre les propriétés $i)$ et $vi)$ du théorème \ref{GEO}.\ref{conngeomfond}, si et seulement si $\got s(X')=k$.

\medskip
\trois{unimpliquedeux} {\em Montrons que $i)\Rightarrow ii)$}. On suppose donc que $i)$ est vraie. Soit $k_\infty$ le complété de la clôture parfaite de $k$ et soit $Y$ la normalisation de $X_{k_\infty}$. Soit $E$ une extension finie séparable de $k_\infty$. Le lemme de Krasner assure que $E$ est isomorphe à $k_\infty\otimes_k F$ pour une certaine extension finie séparable $F$ de $k$ ; puisque $E$ est parfait et contient un sous-corps dense et radiciel sur $F$, c'est le complété d'une clôture radicielle de $F$. 

\medskip
Puisque $i)$ est vraie, $vi)$ est vraie par le \ref{GEO}.\ref{irrgeomfond}.\ref{unequivsix} et $X_F$ est donc irréductible ; il découle alors de la proposition \ref{UNIREG}.\ref{perfclos} que $X_E$ est irréductible. L'espace $Y_E$ étant normal par séparabilité de $E$ sur $k_\infty$, il s'identifie à la normalisation de $X_E$ par la proposition  \ref{NOR}.\ref{normextscal} ;  on en déduit que $Y_E$ est connexe.

\medskip
Ceci vaut pour toute extension finie séparable $E$ de $k_\infty$. Si $\KK$ désigne le complété d'une clôture algébrique de ce dernier, il résulte de l'équivalence des propriétés $vi)$ et $iv)$ du théorème \ref{GEO}.\ref{conngeomfond} que $Y_\KK$ est connexe. Comme $k_\infty$ est parfait, $\KK$ en est une extension analytiquement séparable et $Y_\KK$ est par conséquent normal ; il s'identifie donc à la normalisation de $X_\KK$ par la proposition \ref{NOR}.\ref{normextscal}, et l'irréductibilité de $X_\KK$ en découle. 

\medskip
\trois{deuximpliquetrois} {\em Montrons que $ii)\Rightarrow iii)$.} On suppose donc que $ii)$ est vraie ; pour montrer $iii)$ il suffit, grâce au lemme \ref{GEO}.\ref{extconnexe}, de démontrer que $X_L$ est irréductible pour toute extension complète $L$ de $\KK$.

\medskip
Soit $Z$ la normalisation de $X_\KK$ ; comme $X_\KK$ est irréductible, $Z$ est connexe. Soit $L$ une extension complète de $\KK$ ; puisque $\KK$ est algébriquement clos, $Z_L$ est normal et s'identifie donc, en vertu de la proposition \ref{NOR}.\ref{normextscal}, à la normalisation de $X_L$. Par ailleurs, l'équivalence entre les propriétés $iii)$ et $iv)$ du théorème \ref{GEO}.\ref{conngeomfond} garantit que $Z_L$ est connexe. Il s'ensuit que $X_L$ est irréductible, ce qu'on souhaitait établir. 

\medskip
\medskip
\trois{suitetefin} {\em Fin de la démonstration.} L'implication $iii)\Rightarrow iv)$ est triviale. Supposons que $iv)$ est vraie, et soit $F$ une extension finie de $k$. Choisissons un plongement $F\hookrightarrow \ka$. La flèche induite $X_{\ka}\to X_F$ est surjective, et $X_F$ est donc irréductible, ce qui établit $v)$. L'implication $v)\Rightarrow vi)$ est triviale ; l'implication $vi)\Rightarrow i)$ a été vue au  \ref{GEO}.\ref{irrgeomfond}.\ref{unequivsix}.~$\Box$ 

\medskip
\deux{defgeomirr} On dira qu'un espace $k$-analytique $X$ est {\em géométriquement irréductible} s'il satisfait les conditions équivalentes du théorème ci-dessus.

\medskip
\deux{remmero} {\bf Remarque.} Conservons les notations du théorème ci-dessus. Si $X$ est réduit, alors $\got s(X')$ peut s'interpréter comme l'anneau des fonctions méromorphes sur $X$ qui annulent G-localement un polynôme unitaire séparable à coefficients dans $k$.

\subsection*{Corps de définition d'une composante géométrique} 

\medskip
Dans le cas strictement analytique, une partie des résultats qui suivent figurent déjà dans l'article \cite{comprig} de Conrad (cor. 3.2.3 et th. 3.4.2).

\medskip
\deux{provenance} Si $X$ est un espace $k$-analytique, si $L$ est une extension complète de $k$ et si $E$ est un sous-corps complet de $L$ on dira par abus qu'une composante connexe (resp. irréductible) $T$ de $X_L$ {\em provient de $E$} s'il existe une composante connexe (resp. irréductible) $S$ de $X_E$ telle que $T$ soit $X_E$-isomorphe à $S\hotimes_EL$ ; une telle $S$ est alors nécessairement égale à l'image de $T$ sur $X_E$, et est donc uniquement déterminée.

\medskip
\deux{corpsdefconn} {\bf Théorème.} {\em Soit $X$ un espace $k$-analytique et soit $L$ une extension complète de $k$.

\medskip
$i)$ Il existe une bijection naturelle $\pi_0(X_L)\simeq \coprod\limits_{S\in \pi_0(X)}\pi_0(S_L)$, modulo laquelle $\pi_0(S_L)$ correspond pour tout $S$ au sous-ensemble de $\pi_0(X_L)$ formé des composantes dont l'image sur $Z$ est {\em incluse} dans $S$. 

\medskip
$ii)$ Soit $S\in \pi_0(X)$ et écrivons la $L$-algèbre finie étale $\got s(S)\otimes_kL$ sous la forme $\prod L_i$, où les $L_i$ sont des corps ; l'espace $\got s(S)$-analytique $S$ est géométriquement connexe ; l'espace $S_L$ s'identifie canoniquement à $\coprod S\hotimes_{\got s(S)}L_i$, et les $S\hotimes_{\got s(S)}L_i$ sont précisément ses composantes connexes ; l'on définit ainsi une bijection fonctorielle en $L$ entre $\pi_0(S_L)$ et $\spec \left(\got s(S)\otimes_k L\right)$ ; en particulier $\pi_0(S_L)$ est fini de cardinal inférieur ou égal à $[\got s(S):k]$. 

\medskip
$iii)$ Soit $S\in \pi_0(X)$. Si $T$ est une composante connexe de $X_L$ qui appartient à $\pi_0(S_L)$ modulo la bijection de $i)$ son image sur $X$ est en fait {\em égale} à $S$.

\medskip
$iv)$ Soit $S\in \pi_0(X)$. Si $L$ déploie $\got s(S)$ alors toute composante connexe de $S_L$ est un espace $L$-analytique géométriquement connexe. 

\medskip
$v)$ Soit $T\in \pi_0(X_L)$. L'ensemble des sous-corps complets $E$ de $L$ contenant  $k$ tels que $T$ provienne de $E$ admet un plus petit élément, qui est fini et séparable sur $k$ ; on l'appelle le {\em corps de définition} de $T$. 

\medskip
$vi)$ Soit $\cal E$ l'ensemble des sous-corps complets $E$ de $L$ contenant $k$ tels que {\em toute} composante connexe $T$ de $X_L$ provienne de $E$ ; l'ensemble $\cal E$ admet un plus petit élément $F$ ; la fermeture algébrique de $k$ dans $F$ est dense dans $F$ ; si $\pi_0(X)$ est fini alors $\pi_0(X_L)$ est fini et $F$ est fini et séparable sur $k$.

\medskip
$vii)$ Si $\pi_0(X)$ est un ensemble fini il existe une extension finie et séparable $K$ de $k$ telle que toute composante  connexe de $X_K$ soit un espace $K$-analytique géométriquement connexe.}

\medskip
\deux{rembijconn} {\bf Commentaire.} Il résulte de l'énoncé $i)$ du théorème ci-dessus que l'ensemble $\cal E$ de son énoncé $vi)$ peut également être défini comme l'ensemble des sous-corps complets $E$ de $L$ contenant $k$ tels que $\pi_0(X_L)\to \pi_0(X_E)$ soit bijectif, ainsi que comme l'ensemble des sous-corps complets $E$ de $L$ contenant $k$ tels que $S\hotimes_EL$ soit connexe pour toute composante connexe $S$ de $X_E$.

\medskip
\deux{demtheoconn} {\em Démonstration du théorème} \ref{GEO}.\ref{corpsdefconn}. 

\setcounter{cptbis}{0}
\medskip
\trois{pixpis} {\em Preuve de $i)$.} C'est immédiat, compte-tenu du fait que si $S$ est une composante connexe de $X$ son image réciproque sur $X_L$ s'identifie à $S_L$.

\medskip
\trois{compsigmal} {\em Preuve de $ii)$.} L'espace $\got s(S)$-analytique $S$ est géométriquement connexe d'après le théorème \ref{GEO}.\ref{conngeomfond} et la remarque \ref{GEO}.\ref{remcxcorps}. L'on a $$S_L=S\hotimes_kL=S\hotimes_{\got s(S)}(\got s(S)\otimes_kL)\simeq \coprod S\hotimes_{\got s(S)}L_i.$$ Puisque l'espace $\got s(S)$-analytique $S$ est géométriquement connexe, chacun des $S\hotimes_{\got s(S)}L_i$ est connexe, et les $S\hotimes_{\got s(S)}L_i$ sont donc les composantes connexes de $S_L$, ce qui suffit à démontrer $ii)$. 

\medskip
\trois{imcomconn} {\em Preuve de $iii)$.} Soit $T$ une composante connexe de $X_L$ appartenant à $\pi_0(S_L)$. Il résulte de $ii)$ que $T$ est de la forme $S\hotimes_{\got s(S)}K$ pour une certaine extension complète $K$ de $\got s(S)$, et le morphisme naturel $T\to S$ est donc surjectif. 

\medskip
\trois{ldeploie} {\em Preuve de $iv)$.} Supposons que $L$ déploie $\got s(S)$, soit $T$ une composante connexe de $S_L$ et soit $K$ une extension complète de $L$. D'après $ii)$ la composante $T$ est de la forme $S\hotimes_{\got s(S)}L$ relativement à un certain $k$-plongement de $\got s(S)$ dans $L$ ; il en découle que $T\hotimes_LK$ s'identifie à $S\hotimes_{\got s(S)}K$, lequel est connexe en vertu de la connexité géométrique de l'espace $\got s(S)$-analytique $S$, qui découle de $ii)$ ; par conséquent, $T$ est un $L$-espace analytique géométriquement connexe.

\medskip
\trois{corpsdedef} {\em Preuve de $v)$.} Soit $S$ la composante connexe de $X$ telle que $T\in \pi_0(S_L)$ (notons qu'en vertu de $iii)$, $S$ n'est autre que l'image de $T$ sur $X$). Comme $\got s(S)$ est une extension finie séparable de $k$, il existe un polynôme irréductible unitaire $P$ à coefficients dans $k$ tel que $\got s(S)$ soit isomorphe à $k[\tau]/P(\tau)$. En vertu de $ii)$ la composante $T$ correspond à un facteur irréductible unitaire $Q$ du polynôme $P$ dans $L[X]$ ; en réutilisant $ii)$, on voit immédiatement que le sous-corps de $L$ engendré par $k$ et les coefficients de $Q$ satisfait les propriétés requises.

\medskip
\trois{pilpif} {\em Preuve de $vi)$.} Le corps $F$ égal au complété de la sous-extension de $L$ engendrée par les corps de définition ({\em cf.} $v)$ ) des différentes composantes irréductibles de $X_L$ est clairement le plus petit élément de $\cal E$. Supposons que $\pi_0(Z)$ est fini ; comme $\pi_{0}(Z_L)$ est égal à $\bigcup\limits_{S\in\pi_0(Z)}\pi_0(S_L)$ d'après $i)$ et comme chacun des $\pi_0(S_L)$ est fini d'après $ii)$, l'ensemble $\pi_0(Z_L)$ est lui aussi fini ; le corps $F$ est alors fini et séparable sur $k$ d'après sa construction. 

\medskip
\trois{kgeomconn} {\em Preuve de $vii)$.} Supposons $\pi_0(X)$ fini ; il existe une extension finie galoisienne $K$ de $k$ dans laquelle se plongent chacun des $\got s(S)$ où $S$ parcourt $\pi_0(X)$. Si $S$ est une composante connexe de $X$ il résulte de $iv)$ que les composantes connexes de $S_K$ sont toutes des espaces $K$-analytiques géométriquement connexe ; on déduit de $i)$ que toute composante connexe de $X_K$ est un espace $K$-analytique géométriquement connexe.~$\Box$ 

\medskip
\deux{explicanorm} {\bf Passage de l'étude des composantes connexes à celle des composantes irréductibles.} Nous nous proposons maintenant de décrire le comportement des composantes {\em irréductibles} d'un espace  $k$-analytique par extension des scalaires ; l'idée consiste essentiellement à se ramener, au moyen de la normalisation, au cas des composantes {\em connexes} que l'on vient de traiter. Un fait semble s'y opposer : la formation de la normalisation ne commute pas, en général, au changement de corps de base ; on peut néanmoins contourner l'obstacle en remarquant que les propriétés cruciales de la normalisation qui nous intéressent ici (prop. \ref{NOR}.\ref{connirred} et th. \ref{NOR}.\ref{desccompnorm} $ii)$ ) 
sont quant à elle préservées par extension des scalaires ; ce fait est l'objet de la proposition ci-dessous. 

\medskip
\deux{presquenorm} {\bf Proposition.} {\em Soit $X$ un espace $k$-analytique et soit $f:X'\to X$ sa normalisation. Soit $F$ une extension complète de $k$ et soit $K$ une extension complète de $F$ ; on désigne par $\pi$ et $\pi'$ les  morphismes naturels $X_K\to X_F$ et $X'_K\to X'_F$.

\medskip
$i)$ L'on a $\irr {X'_F}=\pi_0(X'_F)$. 

\medskip
$ii)$ Si $S$ est une composante irréductible de $X'_F$  alors $f_F(S)$ est une composante irréductible de $X_F$, et $S\mapsto f_F(S)$ établit une bijection $\irr{X'_F}\simeq \irr{X_F}$.

\medskip
$iii)$ Si $T$ est une composante irréductible de $X_K$ et si $T'$ est l'unique composante connexe de $X'_K$ telle que $T=f_K(T')$ alors $\pi(T)$ coïncide avec $f_F(\pi'(T'))$ et est une composante irréductible de $X_F$. 

\medskip
$iv)$ Si $S$ est une composante irréductible de $X_F$ et si $S'$ est l'unique composante connexe de $X'_F$ telle que $S=f_F(S')$ alors $S\hotimes_FK=\pi^{-1}(S)=\bigcup f_K(S'_i)$ où les $S'_i$ sont les composantes connexes de $S'\hotimes_FK={\pi'}^{-1}(S')$.}

\medskip
{\em Démonstration.} On peut supposer que $X$ est réduit. Son lieu de normalité $U$ en est ouvert de Zariski dense dont on notera $U'$ l'image réciproque sur $X'$. 

\medskip
\setcounter{cptbis}{0} 
\trois{compirrextnorm} Soit $S$ une composante connexe de $X'$ ; c'est un espace $k$-analytique normal. Soit $E$ une extension finie séparable de $\got s(S)$ ; l'espace $S\hotimes_{\got s(S)}E$ est connexe d'après le théorème \ref{GEO}.\ref{conngeomfond} et la remarque \ref{GEO}.\ref{remcxcorps} ; puisque $S$ est normal et puisque $E$ est séparable sur $\got s(S)$, l'espace $S\hotimes_{\got s(S)}E$ est normal ; comme il est connexe, il est irréductible. Ceci valant pour toute extension finie séparable $E$ de $\got s(S)$, on déduit du théorème \ref{GEO}.\ref{irrgeomfond} que pour toute extension complète $\KK$ de $\got s(S)$, l'espace $S\hotimes_{\got s(S)}\KK$ est irréductible. L'espace $S_F$ s'identifie à $S\hotimes_{\got s(S)} (\got s(S)\otimes_k F)$ ; il est donc somme disjointe d'espaces dont chacun est de la forme $S\hotimes_{\got s(S)}\KK$ pour une certaine extension complète $\KK$ de $\got s(S)$, et est dès lors irréductible d'après ce qui précède. L'assertion $i)$ s'ensuit aussitôt. 

\medskip
\trois{compxprimex} {\em Preuve de $ii)$.} Soit $\Sigma$ une composante connexe (ou irréductible, ce qui revient au même comme on vient de le voir) de $X'_F$ et soit $S$ la composante connexe de $X'$ au-dessus de laquelle $\Sigma$ est située. L'espace $X'$ étant la normalisation de $X$, le sous-ensemble $T:=f(S)$ de $X$ en est une composante irréductible ; soit $d$ sa dimension $k$-analytique. Les espaces $S$ et $T$ sont purement de dimension $d$ ; les espaces $S_F$ et $T_F$ sont donc purement de dimension $d$ ; comme $\Sigma$ est une composante connexe de $S_F$, on a $\dim{}\Sigma=d$ ; par finitude du morphisme $f_F$, le sous-ensemble $f_F(\Sigma)$ de $X_F$ en est un fermé de Zariski irréductible de dimension $d$, situé sur $T_F$ ; c'est donc une composante irréductible de ce dernier, et {\em a fortiori} une composante irréductible de $X_F$ (lemme \ref{COMP}.\ref{irrcompextscal}), ce qui achève de montrer la première assertion de $ii)$. 

\medskip
Établissons maintenant la seconde. Soit $T$ une composante irréductible de $X_F$. La surjectivité de $f_F$ et ce qui précède impliquent l'existence d'une composante irréductible $S$ de $X'_F$ telle que $T=f_F(S)$. Par ailleurs la densité de $U'$ dans $Z'$ entraîne celle de $U'_F$ dans $X'_F$ (cor. \ref{COMP}.\ref{ouvdenseextscal}) et $S$ rencontre de ce fait $U'_F$ . Or $U'_F\to U_F$ est un isomorphisme, ce qui implique qu'un point de $U_F$ a un {\em et un seul} antécédent par $f_F$. Soit $S'$ une composante irréductible de $X'_F$ telle que $f_F(S')=T$ et soit $t$ un point de $U'_F\cap S$ . Comme $t$ est l'unique antécédent de $f_F(t)$, qui a par hypothèse un antécédent sur $T'$, on a $t\in T'$. En conséquence $U'_F\cap T\subset T'$. Comme $U'_F\cap T$ est un ouvert de Zariski non vide, et donc dense, de l'espace irréductible $T$, on a $T\subset T'$ et donc $T'=T$, ce qui achève de montrer $ii)$. 

\medskip
\trois{irrcompfonct} {\em Preuve de $iii)$.} L'égalité $\pi(T)=f_F(\pi'(T'))$ est formelle ; l'assertion $iii)$ du théorème \ref{GEO}.\ref{corpsdefconn} garantit que $\pi'(T')$ est une composante connexe de $X_F$, et $f_F(\pi'(T'))$ est donc bien une composante irréductible de $X_F$ d'après $ii)$. 

\medskip
\trois{irrcompinvfonc} On sait déjà (lemme \ref{COMP}.\ref{irrcompextscal}) que $\pi^{-1}(S)$ est réunion finie de composantes irréductibles de $X_K$ ; en vertu de $ii)$, $\pi^{-1}(S)$ est donc de la forme $\bigcup f_K(\Sigma_j)$, où les $\Sigma_j$ sont des composantes connexes de $X'_K$ ; il est par conséquent égal à la réunion des $f_K(\Sigma)$ où $\Sigma$ parcourt l'ensemble des composantes connexes de $X'_K$ telles que $\pi(f_K(\Sigma))$ soit incluse dans $S$.

\medskip
Soit $\Sigma$ une composante connexe de $X'_K$ ; posons $\Theta=\pi'(\Sigma)$ ; comme on l'a vu au \ref{GEO}.\ref{presquenorm}.\ref{irrcompfonct}, $\Theta$ est une composante connexe de $X'_F$, et $f_F(\Theta)=\pi(f_K(\Sigma))$ est une composante irréductible de $X_F$ ; on en déduit que $\pi(f_K(\Sigma))$ est {\em incluse} dans $S$ si et seulement si elle est {\em égale} à $S$, et que cela se produit si et seulement si $\Theta$ est égale à $S'$, autrement dit si et seulement si $\Sigma$ est l'une des $S'_i$ ; il en découle que $\pi^{-1}(S)=\bigcup f_K(S'_i)$, ce qu'il fallait démontrer.~$\Box$ 

\medskip
\deux{rappirrfonct} {\bf Quelques mots sur le foncteur $L\mapsto \irr X_L$}. Soit $X$ un espace $k$-analytique, soit $L$ une extension complète de $k$ et soit $K$ une extension complète de $L$. Soit $T$ une composante irréductible de $X_K$. Comme $X_K\to X_L$ est quasi-dominant (\ref{NOR}.\ref{qdometplat}.\ref{qdomext}), l'image de $T$ sur $X_L$ est situé sur une et une seule composante irréductible de $X_L$ (on peut aussi le déduire, de manière plus élémentaire, du lemme \ref{COMP}.\ref{irrcompextscal}) ; on dispose donc d'une application $\irr {X_K}\to \irr {X_L}$ ; cette construction fait de $L\mapsto \irr {X_L}$ un foncteur en $L$ ; remarquons que nous venons simplement de rappeler dans un cas particuier ce qui a été mentionné en toute généralité au \ref{NOR}.\ref{pseudom}.\ref{irrfonct}. 

\medskip
\deux{corpsdefirr} {\bf Théorème.} {\em Soit $X$ un espace $k$-analytique et soit $L$ une extension complète de $k$.

\medskip
$i)$ Il existe une bijection naturelle $\irr {X_L}\simeq \coprod\limits_{S\in \irr X}\irr {S_L}$, modulo laquelle $\irr S_L$ correspond pour tout $S$ au sous-ensemble de $\irr {X_L}$ formé des composantes dont l'image sur $Z$ est {\em incluse} dans $S$. 

\medskip
$ii)$ Soit $S\in \irr {X}$ et soit $S'$ sa normalisation ; il existe une bijection fonctorielle en $L$ entre $\irr {S_L}$ et $\spec \left(\got s(S')\otimes_k L\right)$ ; en particulier $\irr {S_L}$ est fini de cardinal inférieur ou égal à $[\got s(S'):k]$. 

\medskip
$iii)$ Soit $S\in \irr X$. Si $T$ est une composante irréductible de $X_L$ qui appartient à $\irr {S_L}$ modulo la bijection de $i)$ son image sur $X$ est en fait {\em égale} à $S$.

\medskip
$iv)$ Soit $S\in \irr X$ et soit $S'$ sa normalisation. Si $L$ déploie $\got s(S')$ alors toute composante irréductible de $S_L$ est un espace $L$-analytique géométriquement irréductible. 

\medskip
$v)$ Soit $T\in \irr {X_L}$. L'ensemble des sous-corps complets $E$ de $L$ contenant  $k$ tels que $T$ provienne de $E$ admet un plus petit élément, qui est fini et séparable sur $k$ ; on l'appelle le {\em corps de définition} de $T$. 

\medskip
$vi)$ Soit $\cal E$ l'ensemble des sous-corps complets $E$ de $L$ contenant $k$ tels que {\em toute} composante irréductible $T$ de $X_L$ provienne de $E$ ; l'ensemble $\cal E$ admet un plus petit élément $F$ ; la fermeture algébrique de $k$ dans $F$ est dense dans $F$ ; si $\irr X$ est fini alors $\irr X_L$ est fini et $F$ est fini et séparable sur $k$.

\medskip
$vii)$ Si $\irr X$ est un ensemble fini il existe une extension finie et séparable $K$ de $k$ telle que toute composante  irréductible de $X_K$ soit un espace $K$-analytique géométriquement irréductible.}

\medskip
\deux{rembijirr} {\bf Commentaire.} Il résulte de l'énoncé $i)$ du théorème ci-dessus que l'ensemble $\cal E$ de son énoncé $vi)$ peut également être défini comme l'ensemble des sous-corps complets $E$ de $L$ contenant $k$ tels que $\irr{X_L}\to \irr {X_E}$ soit bijectif, ainsi que comme l'ensemble des sous-corps complets $E$ de $L$ contenant $k$ tels que $S\hotimes_EL$ soit irréductible pour toute composante irréductible $S$ de $X_E$.

\medskip
\deux{demtheoirr} {\em Démonstration du théorème} \ref{GEO}.\ref{corpsdefirr}. 

\setcounter{cptbis}{0}
\medskip

\trois{irrsirrsl} {\em Preuve de $i)$.} Il s'agit simplement d'une reformulation du lemme \ref{COMP}.\ref{irrcompextscal}. 

\medskip
\trois{imcompirr} {\em Preuve de $iii)$.} Cette assertion a déjà été établie : elle est contenue dans l'énoncé $iii)$ de la proposition \ref{GEO}.\ref{presquenorm} ci-dessus. 

\medskip
\trois{irrsprimel} {\em Preuve de $ii)$, $iv)$, $v)$ et $vi)$.} On les déduit formellement de la proposition \ref{GEO}.\ref{presquenorm} et des énoncés correspondants du théorème \ref{GEO}.\ref{corpsdefconn}.~$\Box$ 

\subsection*{À propos des espaces géométriquement réduits} 

Le but de ce qui suit est d'expliquer brièvement, par un contre-exemple simple ({\em cf.} \cite{comprig}, \S 3.3) pourquoi les équivalences entre les propriétés $iii)$ et $v)$ des théorèmes \ref{GEO}.\ref{conngeomfond} et \ref{GEO}.\ref{irrgeomfond}, respectivement relatifs à la connexité et à l'irréductibilité, ne subsistent pas lorsqu'on s'intéresse à la propriété d'être réduit. {\em On suppose que la caractéristique de $k$ est non nulle ; elle coïncide donc avec l'exposant caractéristique $p$ de $k$. } 

\medskip
\deux{fpegalsomme} {\bf Lemme.} {\em Soit $(a_i)_{i\in \NN}$ une suite d'éléments de $k$ telle que $|a_i|$ tende vers zéro quand $i$ tend vers l'infini. Posons ${\cal A}=k\{ T\} [\tau ]/(\tau^p-\sum a_i T^{pi}$). Les propositions suivantes sont équivalentes :

\medskip
$i)$ l'anneau $\cal A$ est réduit ;

\medskip
$ii)$ il existe $i$ tel que $a_i\notin k^p$.}

\medskip
{\em Démonstration.} Soit $F$ le corps des fractions de $k\{T\}$. Le morphisme d'anneaux $k\{T\}\to \cal A$ est fini et fidèlement plat. Comme $k\{T\}$ est régulier ({\em cf.} lemme \ref{EXC}.\ref{disquereg}), il est de Cohen-Macaulay et $\cal A$ satisfait la propriété $S_1$. En vertu du critère de Serre, il est donc réduit si et seulement si il satisfait $R_0$, c'est-à-dire si et seulement si l'anneau local artinien ${\cal A}\otimes_{k\{T\}}F$ est réduit, soit encore si et seulement si $\sum a_i T^{pi}\notin F^p$. D'autre part $$\sum a_i T^{pi}\in F^p\iff\sum a_i T^{pi}\in k\{T\}^p\iff a_i\in k^p \;\mbox{pour tout}\;i,$$ la première équivalence résultant de la normalité de $k\{T\}$ et la seconde d'un calcul immédiat ; le lemme s'ensuit aussitôt.~$\Box$

\medskip
\deux{exgeomreginfini} {\bf Un exemple.} Soit $(r_i)$ une suite de réels strictement positifs tendant vers zéro quand $i$ tend vers l'infini. Prenons pour $k$ le complété de $\FF_p(a_i)_{i\in \NN}$ pour la valeur absolue qui envoie un poynôme $\sum \lambda_I {\bf a}^I$ sur $\max |\lambda_I |{\bf r}^I$. La famille $(a_i)$ est par construction une $p$-base topologique de $k$ sur $k^p$, et est en particulier libre sur le corps $k^p$. Soit $X$ l'espace $k$-affinoïde ${\cal M}(k\{ T\} [\tau ]/(\tau^p-\sum a_i T^{pi})\;)$. Il résulte du lemme précédent que $X_{\ka}$ n'est pas réduit. 

\medskip
Soit $K$ une extension finie de $k$. Comme $[K^p:k^p]=[K:k]$, le $k^p$-espace vectoriel $K^p$ est de dimension finie, et il en va de même de son sous-$k^p$-espace vectoriel $K^p\cap k$. Par conséquent, il existe $i$ tel que $a_i$ n'appartienne pas à $K^p$. Le lemme ci-dessus assure alors que $X_K$ est réduit. 

\medskip
\deux{contrexgnor} {\bf Remarque.} Nous renvoyons le lecteur au paragraphe 3.3 de \cite{comprig} pour un contre-exemple dans le même esprit concernant la normalité : Conrad y construit, pour tout nombre premier impair $p$, un corps ultramétrique complet $F$ de caractéristique $p$ et un espace $F$-analytique $X$ qui est géométriquement réduit, qui est tel que $X_K$ soit normal pour toute extension finie $K$ de $F$, mais qui n'est pas géométriquement normal. 

\section{Quelques propriétés stables par produit} \label{PROD}

\setcounter{cpt}{0}

\deux{prodnorm} {\bf Théorème.}  {\em Soit $\mathsf P$ une propriété des anneaux locaux noethériens appartenant  à l'un des trois ensembles ${\cal S}$, ${\cal Q}$ et ${\cal R}$ définis au {\rm \ref{RAP}.\ref{listepropri}}. Soit $X$ un espace $k$-analytique, soit $L$ une extension complète de $k$ et soit $Y$ un espace $L$-analytique ; soit $Z$ l'espace $L$-analytique $X\times_kY$. Soit $z\in Z$ et soient $x$ et $y$ ses images respectives sur $X$ et $Y$. Si $X$ satisfait $\mathsf P$ géométriquement en $x$ et si $Y$ satisfait $\mathsf P$ géométriquement en $y$, alors $Z$ satisfait $\mathsf P$ géométriquement en $z$.}

\medskip
{\em Démonstration.} Comme $X\times_k Y\simeq X_L\times_LY$ on peut remplacer $k$ par $L$ et donc supposer que $L=k$. On peut vérifier l'assertion après extension des scalaires, et donc faire l'hypothèse que $k$ est algébriquement clos ; il suffit de traiter le cas où $X$ et $Y$ sont affinoïdes.

\setcounter{cptbis}{0}
\medskip
Soit ${\cal X}$ (resp. $\cal Y$, resp. $\cal Z$) le spectre de l'algèbre des fonctions analytiques sur $X$ (resp. $Y$, resp. $Z$) et soit $\bf x$ (resp. $\bf y$, resp. $\bf z$) l'image de $x$ (resp. $y$, resp. $z$) sur ${\cal X}$ (resp. ${\cal Y}$, resp. ${\cal Z}$).

\medskip
\trois{reducoxxcorps} {\em Réduction au cas où ${\cal O}_{{\cal X},\bf x}$ et ${\cal O}_{{\cal Y},\bf y}$ sont des corps.} Le schéma $\cal Z$ est plat sur $\cal X$ par le lemme \ref{RAP}.\ref{prodalgplat} ; comme ${\cal O}_{{\cal X},\bf x}$ satisfait $\mathsf P$, il suffit de prouver que pour tout idéal premier $\got p$ de ${\cal O}_{{\cal X},\bf x}$ l'anneau ${\cal O}_{{\cal Z},\bf z}\otimes_{{\cal O}_{{\cal X},\bf x}}\kappa(\got p)$ satisfait $\mathsf P$ ; fixons un tel $\got p$.  

\medskip
Soit $\cal T$ le sous-schéma fermé irréductible et {\em réduit} de $\cal X$ dont le point générique $\bf t$ correspond à $\got p$ ; désignons par ${\cal W}$ le produit fibré ${\cal Z}\times_{\cal X}\cal T$. L'anneau ${\cal O}_{{\cal Z},\bf z}\otimes_{{\cal O}_{{\cal X},\bf x}}\kappa(\got p)$ satisfait $\mathsf P$ si et seulement si pour toute générisation $\bf w$ de $\bf z$ située au-dessus de $\bf t$ l'anneau local ${\cal O}_{{\cal W},\bf w}$ satisfait $\mathsf P$ ; soit donc $\bf w$ une telle générisation.

\medskip
Soit $T$ le sous-espace analytique fermé de $X$ correspondant à $\cal T$ ; remarquons que le schéma $\cal W$ s'identifie au spectre de l'anneau des fonctions analytiques sur $W:=Z\times_XT=Y\times_k T$. Choisissons un antécédent $w$ de $\bf w$ sur $W$ ; son image $y'$ sur $Y$ est une générisation de Zariski de $y$, et $Y$ satisfait donc $\mathsf P$ en $y'$ ; l'image $t$ de $w$ sur $T$ est située au-dessus du point générique $\bf t$ du schéma intègre $\cal T$. Il nous suffit maintenant de montrer que $W$ satisfait $\mathsf P$ en $w$ ; on est ainsi revenu au problème initial, mais $(X,Y,x,y,z)$ a été remplacé par $(T,Y,t,y',w)$. Le bilan de cette opération est le suivant : {\em on peut, en modifiant $X$, en conservant $Y$ et en remplaçant $y$ par l'une de ses générisations de Zariski, se ramener au cas où ${\cal O}_{{\cal X},\bf x}$ est un corps.} On réapplique le procédé, mais cette fois-ci à $Y$ ; l'espace $X$ n'est plus modifié ; le point $x$ est éventuellement remplacé par l'une de ses générisations de Zariski, mais celle-ci est encore nécessairement située au-dessus de $\bf x$, puisque  ${\cal O}_{{\cal X},\bf x}$ est un corps. 

\medskip
On se retrouve finalement bien dans la situation où ${\cal O}_{{\cal X},\bf x}$ et ${\cal O}_{{\cal Y},\bf y}$ sont des corps.

\medskip
\trois{utiliserlisse} {\em Preuve du théorème lorsque ${\cal O}_{{\cal X},\bf x}$ et ${\cal O}_{{\cal Y},\bf y}$ sont des corps. } Le corps $k$ étant algébriquement clos, $X$ est quasi-lisse en $x$ et $Y$ est quasi-lisse en $y$ par la prop. \ref{UNIREG}.\ref{proprgeom}. Les propriétés de base du faisceau des formes différentielles (\cite{brk2}, \S 3.3) assurent que $$\dim{\hres(z)}\Omega^1_{Z/k}\otimes \hres(z)=\dim{\hres(x)}\Omega^1_{X/k}\otimes \hres(x)+\dim{\hres(y)}\Omega^1_{Y/k}\otimes \hres(y).$$ Ce dernier terme est égale à $\dim x X+\dim y Y$ en vertu de la quasi-lissité de $X$ en $x$ et de $Y$ en $y$, et donc à $\dim z Z$ d'après le lemme \ref{RAP}.\ref{dimproduit}; le lemme \ref{UNIREG}.\ref{omegareg} permet alors de conclure que $Z$ est régulier en $z$ ; il y satisfait {\em a fortiori} $\mathsf P$, et le théorème est démontré.~$\Box$ 

\medskip
\deux{remvasgeome} {\bf Commentaire.} Il est vraisemblable que le théorème ci-dessus reste vrai si l'on se contente de supposer que $Y$ satisfait simplement ({\em i.e.} non nécessairement géométriquement) la propriété $\mathsf P$ en $y$ ; nous espérons le déduire d'un travail en cours sur la platitude en géométrie de Berkovich. 

\medskip
\deux{connprodaff} {\bf Lemme.} {\em Soient $X$ et $Y$ deux espaces $k$-analytiques connexes. Supposons qu'il existe un G-recouvrement $(X_i)$ de $X$ et un G-recouvrement $(Y_j)$ de $Y$ tels que $X_i\times_k Y_j$ soit connexe pour tout couple $(i,j)$. Alors $X\times_k Y$ est connexe.}

\medskip
{\em Démonstration.} Soient $z$ et $z'$ deux points de $X\times_kY$, soient $x$ et $x'$ (resp. $y$ et $y'$) leurs images respectives sur $X$ (resp. $Y$). Soit $i$ (resp. $j$, resp. $i'$, resp. $j'$) un indice tel que $x\in X_i$ (resp. $y\in Y_j$, resp. $x'\in X_{i'}$, resp. $y'\in Y_{j'}$). Comme $X$ est connexe, il existe une famille finie $i_0=i,i_1,\ldots,i_r=i'$ d'indices telle que $X_{i_l}\cap X_{i_{l+1}}\neq \emptyset$ pour tout $l$ compris entre $0$ et $r-1$. De même, il existe une famille finie $j_0=j,j_1,\ldots,j_s=j'$ d'indices telle que $Y_{j_\lambda}\cap Y_{j_{\lambda+1}}\neq \emptyset$ pour tout $\lambda$ compris entre $0$ et $s-1$. Munissons $\{0,\ldots,r\}\times \{0,\ldots,s\}$ de  l'ordre lexicographique pour lequel (par exemple) la première coordonnée donne le ton ; dans la famille $(X_{i_l}\times_kY_{j_\lambda})_{(l,\lambda)}$, chacun des domaines rencontre le suivant, le premier contient $z$ et le dernier $z'$. Comme ils sont tous connexes, $z$ et $z'$ sont situés sur la même composante connexe de $X\times_k Y$.~$\Box$ 

\medskip
\deux{prodconnirr} {\bf Théorème.} {\em Soit $X$ un espace $k$-analytique, soit $L$ une extension complète de $k$ et soit $Y$ un espace $L$-analytique. Supposons $X$ géométriquement connexe (resp. géométriquement irréductible) et $Y$ connexe (resp. irréductible). Le produit $X\times_k Y$ est alors connexe (resp. irréductible).}

\medskip
{\em Démonstration.} Comme $X\times_k Y\simeq X_L\times_LY$ on peut remplacer $k$ par $L$ et donc se ramener au cas où $L=k$ ; on peut supposer, pour l'assertion relative à la connexité, que $X$ et $Y$ sont non vides. 

\medskip
\setcounter{cptbis}{0}
\trois{xygeomconnirr} {\em Réduction au cas où $X$ et $Y$ sont tous deux géométriquement connexes (resp. géométriquement irréductibles).} Soit $K$ une extension finie galoisienne de $k$ déployant $\got s(Y)$ (resp. $\got s(Y')$, où $Y'$ est la normalisation de $Y$), et soit $\Gamma$ le groupe de Galois de $K$ sur $k$. D'après l'énoncé $ii)$ du théorème \ref{GEO}.\ref{corpsdefconn} (resp. l'énoncé $ii)$ du théorème \ref{GEO}.\ref{corpsdefirr}), le $\Gamma$-ensemble des composantes connexes (resp. irréductibles) de $Y_K$ est en bijection avec le $\Gamma$-ensemble des idéaux maximaux de $K\otimes_k \got s(Y)$ (resp. $K\otimes_k \got s(Y')$) et est en particulier {\em transitif}. 

\medskip
Soit $\{Y_i\}_i$ l'ensemble des composantes connexes (resp. irréductibles) de $Y_K$. Supposons que l'on ait montré que $X_K\times_K Y_i$ est connexe (resp. irréductible) pour tout $i$. Sous cette hypothèse ${\cal E}:=\{X_K\times_K Y_i\}_i$ est un $\Gamma$-ensemble fini qui est constitué d'ouverts fermés non vides, connexes et deux à deux disjoints de $(X\times_k Y)_K$ (resp. de fermés de Zariski irréductibles de $(X\times_k Y)_K$ deux à deux non comparables pour l'inclusion) qui recouvrent $(X\times_kY)_K$ ; les éléments de $\cal E$ sont donc les composantes connexes (resp. irréductibles) de $(X\times_kY)_K$ ; comme $\Gamma$ agit transitivement sur $\{Y_i\}_i$, il agit transitivement sur  $\cal E$ et $X\times_k Y$ est de ce fait connexe (resp. irréductible).

\medskip
{\em Il suffit donc, pour établir notre théorème, de prouver que $X_K\times_K Y_i$ est connexe (resp. irréductible) pour tout $i$.} Comme $K$ déploie $\got s(Y)$ (resp. $\got s(Y')$), il résulte de l'énoncé $iv)$ du théorème \ref{GEO}.\ref{corpsdefconn} (resp. de l'énoncé $iv)$ du théorème \ref{GEO}.\ref{corpsdefirr}) que les $Y_i$ sont géométriquement connexes (resp. irréductibles).

\medskip
L'on peut donc se contenter de démontrer que le produit de deux espaces analytiques géométriquement connexes (resp. géométriquement irréductibles) sur un même corps de base est connexe (resp. irréductible) ; autrement dit, nous nous sommes ramené au cas où $X$ et $Y$ sont tous deux géométriquement connexes (resp. géométriquement irréductibles).

\medskip
\trois{prodconn} {\em Preuve de l'assertion relative à la connexité lorsque $X$ et $Y$ sont tous deux géométriquement connexes. } Quitte à étendre les scalaires, on peut supposer  que $k$ est algébriquement clos. Le lemme \ref{PROD}.\ref{connprodaff} ci-dessus permet de se ramener au cas où $X$ et $Y$ sont affinoïdes. Pour tout point $y$ de $Y$ la fibre $(X\times_kY)_y$ s'identifie à $X_{\hres(y)}$ et est donc connexe. La projection $X\times_k Y\to Y$ étant compacte, à image connexe et à fibres connexes, l'espace $X\times_k Y$ est connexe.

\medskip
\trois{prodirr} {\em Preuve de l'assertion relative à l'irréductibilité lorsque lorsque $X$ et $Y$ sont tous deux géométriquement irréductibles.} Quitte à étendre les scalaires, on peut supposer que $k$ est algébriquement clos. Comme $X$ et $Y$ sont irréductibles, leurs normalisations $X'$ et $Y'$ sont connexes ; le produit $X'\times_k Y'$ est donc connexe en vertu du \ref{PROD}.\ref{prodconnirr}.\ref{prodconn} ci-dessus ; le théorème \ref{PROD}.\ref{prodnorm} assure que $X'\times_k Y'$ est normal ; en conséquence, $X'\times_k Y'$ est irréductible. La flèche $X'\times_kY'\to X\times_k Y$ étant surjective, $X\times_k Y$ est irréductible.~$\Box$ 

\medskip
\deux{corcomprod} {\bf Corollaire.} {\em Soit $X$ un espace $k$-analytique géométriquement connexe et non vide (resp. géométriquement irréductible), soit $L$ une extension complète de $k$ et soit $Y$ un espace $L$-analytique. Soit $(Y_i)$ la famille des composantes connexes (resp. irréductibles) de $Y$. Les $X\times_k Y_i$ sont alors les composantes connexes (resp. irréductibles) de $X\times_k Y$.}

\medskip
{\em Démonstration.} Traitons tout d'abord l'assertion relative à la connexité. Les $X\times_k Y_i$ sont des ouverts fermés de $X\times_k Y$ qui sont non vides, connexes en vertu du théorème précédent, et deux à deux disjoints ; comme ils recouvrent $X\times_k Y$, ce sont ses composantes connexes. 

\medskip
Traitons maintenant l'assertion relative à l'irréductibilité. Les $X\times_k Y_i$ sont des fermés de Zariski de $X\times_k Y$ qui le recouvrent, sont deux à deux non comparables pour l'inclusion, et sont irréductibles en vertu du théorème précédent. Comme $\{X\times_k Y_i\}_i$ est un ensemble G-localement fini de fermés de Zariski de $X\times_k Y$ d'après le \ref{COMP}.\ref{testglocfin}, les $X\times_k Y_i$ sont les composantes irréductibles de $X\times_k Y$.~$\Box$

\end{document}